
  \magnification 1200

  \input pictex
  \input miniltx
  \input color
  \input amssym


  \font \bbfive = bbm5
  \font \bbeight = bbm8
  \font \bbten = bbm10
  \font \rs = rsfs10 \font \rssmall = rsfs10 scaled 833  
  \font \eightbf = cmbx8
  \font \eighti = cmmi8 \skewchar \eighti = '177
  \font \fouri = cmmi5 scaled 800 
  \font \eightit = cmti8
  \font \eightrm = cmr8
  \font \eightsl = cmsl8
  \font \eightsy = cmsy8 \skewchar \eightsy = '60
  \font \eighttt = cmtt8 \hyphenchar \eighttt = -1

  \font \sixi = cmmi6 \skewchar \sixi = '177
  \font \sixrm = cmr6
  \font \sixsy = cmsy6 \skewchar \sixsy = '60
  \font \tensc = cmcsc10

  \scriptfont \bffam = \bbeight
  \scriptscriptfont \bffam = \bbfive
  \textfont \bffam = \bbten

  \newskip \ttglue

  \def \eightpoint {\def \rm {\fam 0 \eightrm }\relax
  \textfont 0 = \eightrm \scriptfont 0 = \sixrm \scriptscriptfont 0 = \fiverm
  \textfont 1 = \eighti \scriptfont 1 = \sixi \scriptscriptfont 1 = \fouri
  \textfont 2 = \eightsy \scriptfont 2 = \sixsy \scriptscriptfont 2 = \fivesy
  \textfont 3 = \tenex \scriptfont 3 = \tenex \scriptscriptfont 3 = \tenex
  \def \it {\fam \itfam \eightit }\relax
  \textfont \itfam = \eightit
  \def \sl {\fam \slfam \eightsl }\relax
  \textfont \slfam = \eightsl
  \def \bf {\fam \bffam \eightbf }\relax
  \textfont \bffam = \bbeight \scriptfont \bffam = \bbfive \scriptscriptfont \bffam = \bbfive
  \def \tt {\fam \ttfam \eighttt }\relax
  \textfont \ttfam = \eighttt
  \tt \ttglue = .5em plus.25em minus.15em
  \normalbaselineskip = 9pt
  \def \MF {{\manual opqr}\-{\manual stuq}}\relax
  \let \sc = \sixrm
  \let \big = \eightbig
  \let \rs = \rssmall
  \setbox \strutbox = \hbox {\vrule height7pt depth2pt width0pt}\relax
  \normalbaselines \rm }

  \def \setfont #1{\font \auxfont =#1 \auxfont }
  \def \withfont #1#2{{\setfont {#1}#2}}


  \def \text #1{{\mathchoice {\hbox {\rm #1}} {\hbox {\rm #1}} {\hbox {\eightrm #1}} {\hbox {\sixrm #1}}}}
  \def \varbox #1{\setbox 0\hbox {$#1$}\setbox 1\hbox {$I$}{\ifdim \ht 0< \ht 1 \scriptstyle #1 \else \scriptscriptstyle #1 \fi }}

  \def \rsbox #1{{\mathchoice {\hbox {\rs #1}} {\hbox {\rs #1}} {\hbox {\rssmall #1}} {\hbox {\rssmall #1}}}}
  \def \mathscr #1{\rsbox {#1}}

  \font \rstwo = rsfs10 scaled 900
  \font \rssmalltwo = rsfs10 scaled 803
  \def \smallmathscr #1{{\mathchoice {\hbox {\rstwo #1}} {\hbox {\rstwo #1}} {\hbox {\rssmalltwo #1}} {\hbox {\rssmalltwo #1}}}}


  \def \TRUE {Y}

  \def \ifundef #1{\expandafter \ifx \csname #1\endcsname \relax }

  \def \undefrule {\kern 2pt \vrule width 2pt height 5pt depth 0pt \kern 2pt}
  \def \UndefLabels {}
  \def \possundef #1{\ifundef {#1}\undefrule {\eighttt #1}\undefrule
    \global \edef \UndefLabels {\UndefLabels #1\par }
  \else \csname #1\endcsname \fi }


  \newcount \secno \secno = 0
  \newcount \stno \stno = 0
  \newcount \eqcntr \eqcntr = 0

  \def \define #1#2{\global \expandafter \edef \csname #1\endcsname {#2}}
  \long \def \error #1{\medskip \noindent {\bf ******* #1}}
  \def \fatal #1{\error {#1\par Exiting...}\end }

  \def \advseqnumbering {\global \advance \stno by 1 \global \eqcntr =0}

  \def \current {\ifnum \secno = 0 \number \stno \else \number \secno \ifnum \stno = 0 \else .\number \stno \fi \fi }


  \def \deflabel #1#2{\relax
    \ifundef {#1PrimarilyDefined}\relax
      \define {#1}{#2}\relax
      \define {#1PrimarilyDefined}{#2}\relax
      \immediate \write 1 {\string \newlabe l {#1}{#2}}\relax
    \else
      \edef \old {\csname #1\endcsname }\relax
      \edef \new {#2}\relax
      \if \old \new \else \fatal {Duplicate definition for label ``{\tt #1}'', already defined as ``{\tt \old }''.}\fi
      \fi }

  \def \label #1 {\deflabel {#1}{\current }}

  \def \lbldeq #1 $$#2$${\ifundef {InsideBlock}\advseqnumbering
    \edef \lbl {\current }\else
    \global \advance \eqcntr by 1
    \edef \lbl {\current .\number \eqcntr }\fi
    \deflabel {#1}{\lbl }
    $$
    #2
    \eqno {(\lbl )}
    $$}

  \def \split #1.#2.#3.#4;{\global \def \parone {#1}\global \def \partwo {#2}\global \def \parthree {#3}\global \def
\parfour {#4}}
  \def \NA {NA}
  \def \ref #1{\split #1.NA.NA.NA;(\possundef {\parone }\ifx \partwo \NA \else .\partwo \fi )}


  \newcount \bibno \bibno = 0

  \def \Bibitem #1 #2; #3; #4 \par {\smallbreak
    \global \advance \bibno by 1
    \item {[\possundef {#1}]} #2, {``#3''}, #4.\par
    \ifundef {#1PrimarilyDefined}\else
      \fatal {Duplicate definition for bibliography item ``{\tt #1}'',
      already defined in ``{\tt [\csname #1\endcsname ]}''.}
      \fi \ifundef {#1}\else \edef \prevNum {\csname #1\endcsname } \ifnum \bibno =\prevNum \else \error {Mismatch
        bibliography item ``{\tt #1}'', defined earlier (in aux file ?) as ``{\tt \prevNum }'' but should be ``{\tt
        \number \bibno }''.  Running again should fix this.}  \fi \fi
    \define {#1PrimarilyDefined}{#2}\relax
    \immediate \write 1 {\string \newbi b {#1}{\number \bibno }} }

  \def \jrn #1, #2 (#3), #4-#5;{{\sl #1}, {\bf #2} (#3), #4--#5}
  \def \Article #1 #2; #3; #4 \par {\Bibitem #1 #2; #3; \jrn #4; \par }

  \def \references {\begingroup \bigbreak \eightpoint \centerline {\tensc References} \nobreak \medskip \frenchspacing }


  \catcode `\@ =11
  \def \citetrk #1{{\bf \possundef {#1}}}
  \def \c@ite #1{{\rm [\citetrk {#1}]}}
  \def \sc@ite [#1]#2{{\rm [\citetrk {#2}\hskip 0.7pt:\hskip 2pt #1]}}
  \def \du@lcite {\if \pe@k [\expandafter \sc@ite \else \expandafter \c@ite \fi }
  \def \cite {\futurelet \pe@k \du@lcite }
  \catcode `\@ =12


  \long \def \Quote #1\endQuote {\begingroup \leftskip 35pt \rightskip 35pt \parindent 17pt \eightpoint #1\par \endgroup }
  \long \def \Abstract #1\endAbstract {\vskip 1cm \Quote \noindent #1\endQuote }

  \def \part #1#2{\vfill \eject \null \vskip 0.3\vsize
    \withfont {cmbx10 scaled 1440}{\centerline {PART #1} \vskip 1.5cm \centerline {#2}} \vfill \eject }


  \def \fix {\smallskip \noindent $\blacktriangleright $\kern 12pt}
  \def \newpage {\vfill \eject }

  \def \emph #1{{\it #1}\/}

  \def \state #1 #2\par {\begingroup \def \InsideBlock {} \medbreak \noindent \advseqnumbering {\bf \current .\enspace
#1.\enspace \sl #2\par }\medbreak \endgroup }

  \def \definition #1\par {\state Definition \rm #1\par }


  \newcount \CloseProofFlag
  \def \closeProof {\eqno \endproofmarker \global \CloseProofFlag =1}
  
  \long \def \Proof #1\endProof {\begingroup \def \InsideBlock {} \global \CloseProofFlag =0
     \medbreak \noindent {\it Proof.\enspace }#1
     \ifnum \CloseProofFlag =0 \hfill $\endproofmarker $ \looseness = -1 \fi
     \medbreak \endgroup }

  \def \quebra #1{#1 $$$$ #1}
  \def \explica #1#2{\mathrel {\buildrel \hbox {\sixrm #1} \over #2}}
  \def \explain #1#2{\explica {\ref {#1}}{#2}}
  
  \def \=#1{\explain {#1}{=}}

  \def \pilar #1{\vrule height #1 width 0pt}
  \def \stake #1{\vrule depth #1 width 0pt}

  \newcount \fnctr \fnctr = 0
  \def \fn #1{\global \advance \fnctr by 1
    \edef \footnumb {$^{\number \fnctr }$}\relax
    \footnote {\footnumb }{\eightpoint #1\par \vskip -10pt}}


  \def \item #1{\par \noindent \kern 1.1truecm\hangindent 1.1truecm \llap {#1\enspace }\ignorespaces }
  
  \def \Item #1{\smallskip \item {{\rm #1}}}

  \newcount \zitemno \zitemno = 0
  \def \izitem {\global \zitemno = 0}

  \def \zitemplus {\global \advance \zitemno by 1 \relax }
  \def \rzitem {\romannumeral \zitemno }
  \def \rzitemplus {\zitemplus \rzitem }
  \def \zitem {\Item {{\rm (\rzitemplus )}}}

  \newcount \nitemno \nitemno = 0
  
  \def \nitem {\global \advance \nitemno by 1 \Item {{\rm (\number \nitemno )}}}

  \newcount \aitemno \aitemno = -1
  \def \boxlet #1{\hbox to 6.5pt{\hfill #1\hfill }}
  \def \iaitem {\aitemno = -1}
  \def \aitemconv {\ifcase \aitemno a\or b\or c\or d\or e\or f\or g\or h\or i\or j\or k\or l\or m\or n\or o\or p\or q\or
    r\or s\or t\or u\or v\or w\or x\or y\or z\else zzz\fi }
  \def \aitem {\global \advance \aitemno by 1\Item {(\boxlet \aitemconv )}}


  \def \deflabeloc #1#2{\deflabel {#1}{\current .#2}\deflabel {Local#1}{#2}}
  \def \lbldzitem #1 {\zitem \deflabeloc {#1}{\rzitem }}
  \def \lbldaitem #1 {\aitem \deflabeloc {#1}{\aitemconv }}
  \def \aitemmark #1 {\deflabel {#1}{\aitemconv }}


  \def \ds {\displaystyle }
  \def \and {\mathchoice {\hbox {\quad and \quad }} {\hbox { and }} {\hbox { and }} {\hbox { and }}}

  \def \calcat #1{\,{\vrule height8pt depth4pt}_{\,#1}}
  \def \IFF {\kern 7pt\Leftrightarrow \kern 7pt}
  \def \IMPLY {\kern 7pt \Rightarrow \kern 7pt}
  \def \for #1{\mathchoice {\quad \forall \,#1} {\hbox { for all } #1} {\forall #1}{\forall #1}}
  \def \endproofmarker {\square }
  \def \"#1{{\it #1}\/} \def \umlaut #1{{\accent "7F #1}}
  \def \inv {^{-1}}
  \def \*{\otimes }
  \def \bfdef #1{\global \expandafter \edef \csname #1\endcsname {{\bf #1}}}
  \bfdef N \bfdef Z \bfdef C \bfdef R
  \def \exists {\mathchar "0239\kern 1pt }


  \IfFileExists {\jobname .aux}{\input \jobname .aux}{\null }
  \immediate \openout 1 \jobname .aux

  \def \close {\if \empty \UndefLabels \else
    \message {*** There were undefined labels ***} \medskip \noindent
    ****************** \ Undefined Labels: \tt \par \UndefLabels \fi
    \closeout 1
    \par \vfill \supereject \end }


  \def \medcup {\mathop {\mathchoice {\raise 1pt \hbox {$\mathchar "1353$}}{\mathchar "1353}{\mathchar "1353}{\mathchar
"1353}}}
  \def \medcap {\mathop {\mathchoice {\raise 1pt \hbox {$\mathchar "1354$}}{\mathchar "1354}{\mathchar "1354}{\mathchar
"1354}}}
  \def \clauses #1{\def \crr {\vrule width 0pt height 10pt depth 5pt\cr }\left \{ \matrix {#1}\right .}
  \def \cl #1 #2 #3 {#1, & \hbox {#2 } #3\hfill \crr }

  %
  %

  \def \startsection #1 \par
    {\goodbreak \bigbreak
    \begingroup
    \global \edef \secname {#1}\relax
    \global \advance \secno by 1
    \stno = 0
    \AddToTableOfContents {\number \secno .}{#1}}

  \def \sectiontitle \par
    {\noindent {\bf \number \secno .\enspace \secname .}
    \nobreak \medskip
    \noindent }

  \def \endsection {\endgroup }

  %
  %

  \def \ToCf {\TRUE }
  \if \ToCf \TRUE \def \AddToTableOfContents #1#2{{\let \the =0\edef \a {\write 2{\string \inde x #1; #2; \the \count 0;}}\a }}
            \else \def \AddToTableOfContents #1#2{}\fi
  \def \condinput #1 {\IfFileExists {#1}{\input #1}{*** MISSING FILE #1.  RUNNING AGAIN MIGHT FIX THIS.}}
  \def \index #1; #2; #3;{\line {\hbox to 20pt{\hfill #1} #2 \ \dotfill \ \ #3}}
  \def \tableOfContents {\centerline {CONTENTS} \bigskip \IfFileExists {c}{\input c}{*** MISSING FILE c.  RUNNING AGAIN MIGHT FIX THIS.}ontents.aux \if \ToCf \TRUE \openout 2 contents.aux \fi }

  %
  %
  \def \dotedline #1 ... #2{\smallskip \line {\qquad #1\quad \dotfill \quad #2\qquad }}
  \def \PrintSymbol #1; #2; #3;{\dotedline \hbox to 1.5cm{\hfill $#1$\hfill }\ #2 ... {#3}}

  \def \SymbolBank {}

  \def \newsymbol #1#2{#1\global \edef \SymbolBank {\SymbolBank \PrintSymbol #1; #2; \number \secno .\number \stno ;}}

  \def \SymbolIndex {\begingroup \eightpoint \vskip 1.5cm \goodbreak
    \centerline {\bf SYMBOL INDEX} \bigskip
    \PrintSymbol \hbox {SYMBOL} ; DESCRIPTION; LOCATION; \medskip
    \AddToTableOfContents {}{Symbol index} \SymbolBank \endgroup }

  \def \PrintConcept #1; #2;{\dotedline #1 ... {#2}}

  \def \ConceptBank {}

  \def \newConcept #1#2{\emph {#1}\def \tempVar {#2}\ifx \empty \tempVar \def \tempVar {#1}\else \def \tempVar {#2}\fi
    \global \edef \ConceptBank {\ConceptBank \PrintConcept \tempVar ; \number \secno .\number \stno ;}}

  \def \ConceptIndex {\begingroup \eightpoint \vskip 1.5cm \goodbreak
    \centerline {\bf CONCEPT INDEX}
    \bigskip \PrintConcept CONCEPT; LOCATION;
    \medskip \AddToTableOfContents {}{Concept index}
    \ConceptBank
    \endgroup }

  \def \PART #1{\newpage \setbox 0\hbox {\bf PART #1}
    \dimen 0 \wd 0
    \advance \dimen 0 by 20pt
    \centerline {\hbox to \dimen 0{\hrulefill }}
    \medskip \centerline {\box 0}
    \centerline {\hbox to \dimen 0{\hrulefill }}
    \vskip 2cm
    \write 2{\string \medskip \centerline {\eightrm PART #1} \medskip}
    }

  \def \labelarrow #1{\setbox 0\hbox {\ \ $#1$\ \ }\ {\buildrel \textstyle #1 \over {\hbox to \wd 0 {\rightarrowfill }}}\ }
  \def \subProof #1{\medskip \noindent #1\enspace }
  \def \itmProof (#1) {\subProof {(#1)}}
  \def \itemImply #1#2{\subProof {#1$\Rightarrow $#2}}
  \def \itmImply (#1) > (#2) {\itemImply {(#1)}{(#2)}}
  \def \refl #1{\ref {Local#1}}
  \def \clspan {\overline {\hbox {span} \pilar {6pt}}\ }

  \def \Frac #1#2{{#1\over #2}}


  \def \scrCcL #1{{\mathscr C}_c(#1, \Lbpar )}
  \def \scrCzL #1{{\mathscr C}_0(#1, \Lbpar )}
  \def \Cg {C^*(\G , \Lbpar )}
  \def \Cgf {\EssAlg {\G } {\Lb \kern 1pt}}
  \def \red {{\hbox{\sixrm red}}}
  \def \Cgr {C^*_\red (\G , \Lbpar )}
  \def \Cw {C^*\big (\G _N,\WeylBund \big )}
  \def \Cwr {C^*_\red \big (\G _N,\WeylBund \big )}
  \def \G {{\cal G}}
  \def \Gz {\G ^{(0)}}
  \def \CzGz {C_0\big (\Gz \big )}
  \def \CcGz {C_c\big (\Gz \big )}
  \def \CurlyCcGz {{\mathscr C}_c\big (\Gz \big )}
  \def \ess {{\hbox {\sixrm ess}}}
  \def \EssAlg #1#2{C^*_\ess (#1, #2)}

  \def \Hilb {{\cal H}}
  \def \BH {{\cal B}(\Hilb )}
  \def \subsz #1{_{\varbox {#1}}}
  \def \supsz #1{^{\varbox {#1}}}
  \def \E #1;{E\subsz {#1}}
  \def \EV #1#2{\kern 0.5pt\big \langle #1, #2\big \rangle \kern 0.5pt}
  \def \Fix #1{\hbox {\sl Fix}\kern 1pt(#1)}
  \def \Frac #1#2{{#1\over #2}}
  \def \H #1;{H\subsz {#1}}
  \def \INCL #1#2{\left \langle #2\,,\, #1\right \rangle }
  \def \Incl #1#2{\big \langle #2, #1\big \rangle }
  \def \Isot #1;{B(#1)}
  \def \L #1;{L\subsz {#1}}
  \def \Npt #1;{N_{#1}}
  \def \O #1{\Omega (#1)}
  \def \Orb #1{\text {Orb}(#1)}
  \def \PIA #1;{C\subsz {#1}}
  \def \Weylbisg #1{C_0\big (#1,\WeylBund \big )}
  \def \bd #1{\partial \kern 1pt \Fix {#1}}
  \def \cclosint #1#2{C[#1, #2]}
  \def \clsr #1{\overline {#1\pilar {7.8pt}}}\relax
  \def \ev #1#2{\kern 0.5pt\langle #1, #2\rangle \kern 0.5pt}
  \def \evx #1{\ev {#1}x}
  \def \ft #1#2#3{\E #1#2;(#3)}
  \def \iFix #1{\hbox {\sl F\i x}\kern -7pt\raise 2pt\hbox {$^\circ $}\kern
3pt(#1)}
  \def \incl #1#2{ \langle #2, #1\kern 0.5pt\rangle }
  \def \int #1{#1^\circ }
  \def \li #1{\lim _{#1\to \infty }}
  \def \p #1;{p\subsz {#1}}
  \def \q #1;{q\subsz {#1}}

  \def \BIJ {\Isot I, J;}
  \def \BJ {\Isot J;}
  \def \BJI {\Isot J, I;}
  \def \Bun {\mathscr B}
  \def \Bx {\Isot x;}
  \def \Byx {\Isot y, x;}
  \def \Lbsys {\smallmathscr {L}}
  \def \Lb {\Lbsys \kern 1pt}
  \def \Lbpar {\Lbsys \kern 2pt}
  \def \Lbsubs {\Lbsys \kern 2pt}
  \def \EllTwo {\ell ^2(\G _x,\Lbpar )}
  \def \CIJ {\PIA IJ; }
  \def \CJ {\PIA J; }
  \def \Cx {\PIA x; }
  \def \Cyx {\PIA yx;}
  \def \EE {{\cal E}}
  \def \Ex {\E x; }
  \def \Eyx {\E yx;}
  \def \qess {q_\ess }
  \def \HIJ {\H IJ;}
  \def \ISG {\goth {S}_N}
  \def \LIJ {\L IJ;}
  \def \LJ {\L J;}
  \def \Lx {\L x;}
  \def \Lyx {\L yx;}
  \def \Mult {{\cal M}}
  \def \NoA {N\!A}
  \def \Nyx {\Npt yx; }
  \def \Pf {E_{\text {\sixrm free}}}
  \def \WeylBund {\Lbsys \kern 0.5pt_N }
  \def \Xt {X^{(2)}}
  \def \circle {{\bf T}}
  \def \free {\kern -1pt_{\scriptscriptstyle f}}
  \def \mnrm {_{\hbox {\sixrm min}}}
  \def \pij {\p IJ;}
  \def \px {\p x; }
  \def \pyx {\p yx;}
  \def \qij {\q IJ;}
  \def \Norm #1#2{N(#2,#1)}
  \def \Dom #1{\text {dom}(#1)}
  \def \dom #1{\text {dom}(#1)}
  \def \Ran #1{\text {ran}(#1)}
  \def \ran #1{\text {ran}(#1)}
  \def \Ker {\text {Ker}}

  \def \dstext #1{\quad \text {#1} \quad }
  \def \nbimod {\hbox {Bimod}_A(n)}
  \def \bimod #1{\hbox {Bimod}_A(#1)}
  
  \def \pd #1#2{\theta _{#2}(#1)}

  \centerline {\bf CHARACTERIZING GROUPOID C*-ALGEBRAS OF}
  \smallskip
  \centerline {\bf NON-HAUSDORFF \'ETALE GROUPOIDS}

  \vfill
  \centerline {\tensc R. Exel\footnote {$^{\ast }$}{\eightrm Universidade Federal de Santa Catarina
and University of Nebraska--Lincoln. Partially supported by CNPq.}  and David R. Pitts\footnote {$^{\ast \ast }$}{\eightrm
University of Nebraska--Lincoln.  Partially supported by the Simons Foundation, grant \#316952.}
  }

  \vfill
  \Abstract
  \def \Lb {\hbox {\rssmalltwo L}\kern 2pt}\relax
  \def \Lbpar {\Lb }\relax
  Given a not-necessarily Hausdorff, topologically free, twisted \'etale groupoid $(\G , \Lb )$, we consider its \emph
{essential groupoid C*-algebra}, denoted $\Cgf $, obtained by completing $\scrCcL \G $ with the smallest among all
C*-seminorms coinciding with the uniform norm on $C_c(\Gz )$.  The inclusion of C*-algebras
  $$
  \INCL {\Cgf }{C_0(\Gz ) }
  $$
  is then proven to satisfy a list of properties characterizing it as what we call a \emph {weak Cartan} inclusion.  We
then prove that every weak Cartan inclusion $\incl BA$, with $B$ separable, is modeled by a topologically free, twisted
\'etale groupoid, as above.  In our second main result we give a necessary and sufficient condition for an inclusion of
C*-algebras $\incl BA$ to be modeled by a twisted \'etale groupoid based on the notion of \emph {canonical states}.  A
simplicity criterion for $\Cgf $ is proven and many examples are provided.
  \footnote {}{\vskip 1pt \eightrm \hskip 3.5pt MSC (2020): 46L05, 46L45, 22A22, 46L55}
  \endAbstract

\vfill\vfill\vfill\vfill\vfill\vfill
\eject

\tableOfContents

\startsection Introduction

\sectiontitle

In their landmark paper \cite{FM}, Feldman and Moore gave methods for describing certain von Neumann algebras using
twisted measured equivalence relations arising from a Cartan MASA.  Their result may be thought of as a deep and
far-reaching generalization of the process of describing all linear maps on a finite dimensional Hilbert space using a
fixed orthonormal basis.  Nearly a decade later, Kumjian~\cite{Kumjian} defined the notion of a diagonal $A$ in a
C*-algebra $B$, and provided a theory of twisted topological equivalence relations which allowed him to classify the
inclusion $A\subseteq B$ in terms of such twisted topological equivalence relations.  However, because Kumjian's results
require that each pure state of $A$ uniquely extends to a pure state on $B$, there are natural settings for which his
results do not apply.  In a 2007 paper, Renault~\cite{RenaultCartan} introduced a notion of Cartan MASA in a C*-algebra,
and gave an elegant classification theorem for C*-algebras containing a Cartan MASA completely analagous to the
Feldman-Moore results for von Neumann algebras.  Renault's result replaces the measured equivalence relations found in
the von Neumann algebra setting with an appropriate class of twisted locally compact groupoids which are Hausdorff and
\'etale.

This is a satisfying state of affairs.  However, there still are desirable settings where one has an inclusion of
C*-algebras $A\subseteq B$ (which we write as $\incl BA$) which satisfies the same regularity condition as a Cartan
subalgebra, but which may not satisfy other conditions found in the definition of a Cartan MASA.  For example, the
crossed product constructions found in \cite[Section~6.1]{SRI} give a large class of examples where $A$ is a regular
MASA in $B$, but which do not have a conditional expectation $E:B\rightarrow A$.  In such settings, it is possible to
use the constructions of Kumjian and Renault to produce a twisted groupoid $(\G,\Lb)$, but the underlying groupoid $\G$
is necessarily non-Hausdorff (see \cite[Theorem~4.4]{SRITwo}).  It is not possible to describe such inclusions $\incl
BA$ using the methods of Renault and Kumjian.

The present work has several main goals.  The first is to give a groupoid model for inclusions of C*-algebras satisfying
certain properties inspired by, although much weaker than, the properties characterizing Cartan inclusions, as defined
by Renault in \cite {RenaultCartan}.  The second main goal is intimately related to the first: we define the \emph
{essential groupoid C*-algebra} associated to a not necessarily Hausdorff, topologically free, \'etale groupoid $\G $,
equipped with a Fell line bundle $\Lb $.  Taken together, these two results provide a classification for weak Cartan
inclusions in the spirit of Renault and Kumjian using twisted, but possibly not Hausdorff, \'etale groupoids.  Our
results apply in considerable generality, and include the regular MASA inclusions $\incl BA$ mentioned in the previous
paragraph, but they also apply to many regular inclusions where $A$ is an abelian subalgebra of $B$ which is not maximal
abelian.  Let us be more specific.

Let $\G$ be a topologically free, \'etale groupoid, equipped with a Fell line bundle $\Lb $.  While we do not assume
$\G$ is Hausdorff, we do assume its unit space is Hausdorff.  Well-known constructions associate two norms on the
$*$-algebra $\scrCcL \G $ formed by the continuous local cross-sections of $\Lb $: the \emph {maximum} and the \emph
{reduced} norms, respectively (see \ref {DefCCLG} for the precise definition).  Their completions yield the full and the
reduced twisted groupoid C*-algebras of $\G $, usually denoted by $\Cg $ and $\Cgr $.

When $\G $ is topologically free, that is, when the units of $\G $ with trivial isotropy form a dense set, we show in
Theorem \ref {MinNorm} that there is a very canonical third choice, namely the \emph {smallest} among all C*-seminorms
on $\scrCcL \G $ coinciding with the uniform norm on $\CcGz $.  The completion of $\scrCcL \G $ under this C*-seminorm
is the \emph {essential groupoid C*-algebra} referred to above, which we denote by
  $$
  \Cgf .
  $$

Inspired by an earlier version of the present work, Kwasniewski and Meyer \cite {KM} introduced a different notion of
essential C*-algebra which applies to groupoids which are not necessarily topologically free. Their construction is
described in detail in Section \ref{KMSection}.

  In case $\G $ is Hausdorff, $\Cgf$ turns out to be isomorphic to $\Cgr $, but otherwise this is not necessarily so.
In fact we show in \ref {LargestIdeal} that, in the general case, there exists a largest ideal
  $$
  \Gamma _\red \trianglelefteq \Cgr ,
  $$
  having a trivial intersection with $\CzGz $, which we call the \emph {{gray} ideal}, and such that
  $$
  \Cgf \simeq \Cgr /\Gamma _\red .
  $$

The pair of C*-algebras
  \lbldeq ModelInclu
  $$
  \Incl {\Cgf }{\CzGz }
  $$
  is then taken to be the groupoid model for the inclusions mentioned above.  Given an inclusion of C*-algebras
  \lbldeq AbsInclu
  $$
  \incl BA,
  $$
  that is, $A$ and $B$ are C*-algebras with $A\subseteq B$, satisfying certain natural properties, we then show that
there exists a topologically free, twisted \'etale groupoid $\G $, such that \ref {AbsInclu} and \ref {ModelInclu} are
isomorphic inclusions.

The reader familiar with Kumjian's work \cite {Kumjian} on C*-diagonals, or Renault's paper \cite {RenaultCartan} on
Cartan subalgebras, will be correct in guessing that regularity \cite [Definition 4.1]{RenaultCartan} is an important
condition in our characterization, but although also a natural guess, maximal commutativity is not included.  In fact,
our level of generality goes much beyond that, even allowing us to treat situations where both $A$ and $B$ are
commutative C*-algebras.  To mention but one of the commutative examples fitting our theory, the inclusion of $C(S^1)$
in $C([0,1])$ obtained by viewing functions on the circle as periodic functions on the interval, also falls under our
grasp.

Instead of maximal commutativity, one of the main properties appearing in our characterization is related to the unique
extension of pure states, a property which also plays a predominant role in \cite {Kumjian}.  To be precise, we assume
that the first component of our inclusion, namely the C*-algebra $A$, is commutative, and hence has the form $C_0(X)$,
where $X$ is a locally compact Hausdorff space. For each $x$ in $X$, we consider the pure state $\varphi _x$ on $A$
given by evaluation at $x$. It is well known that $\varphi _x$ admits an extension to a state on $B$, but such an
extension is often not unique.  When uniqueness holds, we say that $x$ is a \emph {free point}.  This is in analogy with
the groupoid case in which uniqueness holds if and only if the isotropy group of $x$ is trivial, as proven in \ref
{CharacFreePts}.  Still in analogy with groupoids, we say that a regular inclusion $\incl BA$ is \emph {topologically
free} if the set of free points is dense in $X$.  Indeed, in the above context of a twisted \'etale groupoid we have
that $\G $ is topologically free if and only if either one of the inclusions
  $$
  \Incl {\Cg }{\CzGz }, \and
  \Incl {\Cgr }{\CzGz }
  $$
  are topologically free, in the above sense.

However topological freeness, in itself, is not sufficient for our purposes, requiring a deeper analysis of the question
of the uniqueness of pure state extensions, which is what we do by introducing the notion of \emph {relatively free
points}.  To explain what we mean by this, let $b$ be a fixed element in $B$.  A point $x$ in $X$ is then called \emph
{free relative to} $b$, provided all extensions of $\varphi _x$ coincide on $b$.  It is obvious that a free point is
also free relative to every $b$ but, on the other hand, there may be many more points which are only free relative to a
given $b$.

An element of $B$ whose set of relatively free points has a dense interior is called \emph {smooth}.  This notion is
again inspired by groupoids since the better behaved normalizers of \ref {ModelInclu}, namely the cross-sections of $\Lb
$ which are continuous and compactly supported on some open bisection, are smooth whenever $\G $ is topologically free
\ref {Sexy}.

This brings us to the precise definition of one of the two crucial properties involved in our characterization, namely
we say that a regular inclusion $\incl BA$ is \emph {smooth} when there exists a $*$-semigroup $N$ of smooth normalizers
such that $\NoA $ spans a dense subspace of $B$.

The reader acquainted with \cite {Kumjian} and \cite {RenaultCartan} is probably anticipating some further condition
relating to the existence of a faithful conditional expectation, but unless we want to restrict our study to Hausdorff
groupoids, which we do not, conditional expectations must not be required.  Instead, we prove that, for a topologically
free inclusion $\incl BA$, the set
  $$
  \Gamma = \medcap _{x\in F} BJ_x,
  $$
  where $F$ is the set of free points, and $J_x$ is the kernel of $\varphi _x$, is a two sided ideal of $B$ which we
call the \emph {{gray} ideal}.

Incidentally, in case $\G $ is topologically free, we prove in \ref {FullRedGray} that the inverse image of the {gray}
ideal of $\Cgr $ under the \emph {left regular representation}
  $$
  \Lambda :\Cg \to \Cgr ,
  $$
  coincides with the {gray} ideal of $\Cg $.  This produces an isomorphism of quotient algebras,
  \medskip
  \begingroup \noindent \hfill \beginpicture
  \setcoordinatesystem units <0.035truecm, -0.02truecm>
  \setplotarea x from -100 to 200, y from -30 to 125
  \put {\null } at -100 -30
  \put {\null } at -100 125
  \put {\null } at 200 -30
  \put {\null } at 200 150
  \put {$\Cg $} at 0 0
  \put {$\Cgr $} at 100 0
  \arrow <0.15cm> [0.25, 0.75] from 32 0 to 65 0
  \put {$\Lambda $} at 50 -12
  \put {$\ds \Frac {\Cg }{\Gamma }$} at 0 100
  \put {$\ds \Frac {\Cgr }{\Gamma _\red }$} at 100 100
  \arrow <0.15cm> [0.25, 0.75] from 0 23 to 0 65
  \arrow <0.15cm> [0.25, 0.75] from 100 23 to 100 65
  \put {$\simeq $} at 50 100
  \put {$= $} at 145 100 \put {$\Cgf $} at 187 100
  \endpicture \hfill \null \endgroup

For Hausdorff groupoids the {gray} ideal of $\Cgr $ vanishes by \ref {EssGpdHausdorff.ii}, and this is perhaps the
reason why it has not been intensely studied before, but for a general topologically free inclusion it may well be
nonzero.  In any case, for topologically free inclusions, the {gray} ideal is always the largest ideal of $B$ having a
trivial intersection with $A$, as proved in \ref {LargestIdeal}.

This property of the {gray} ideal is of course related to the \emph {ideal intersection property}, meaning the property
describing the fact that every nonzero ideal of $B$ has a nonzero intersection with $A$.  For topologically free
inclusions we then conclude that the ideal intersection property is equivalent to the vanishing of the {gray} ideal \ref
{RelTFs}.

We therefore say that $\incl BA$ is a \emph {{light} inclusion} when the {gray} ideal of $B$ vanishes, and this turns
out to be the second fundamental property in our characterization.

In order to give a unified name to the objects under our radar we say that $\incl BA$ is a \emph {weak Cartan inclusion}
when it is both smooth and {light}.  Theorems \ref {MainOne} and \ref {MainTwo}, two of our main results, may then be
precisely described by saying that the class of weak Cartan inclusions $\incl BA$, with $B$ separable, coincides with
the class of inclusions \ref {ModelInclu}, with $\G $ topologically free and second countable.

In our third main result, namely Theorem \ref {MainThree}, we give a characterization of (possibly exotic) twisted
groupoid C*-algebras via a condition (existence of canonical states) which holds for any twisted \'etale groupoid by
\ref {CanoInGpd}.  This is therefore likely to be the first necessary and sufficient condition for an inclusion of
C*-algebras to be modeled by a twisted \'etale groupoid.

As a byproduct of our theory, in \ref {GpdSimplicityCriteria} we give a very clear cut criterion for the simplicity of
the essential groupoid C*-algebra, which states that if $\G $ is topologically free, then $\Cgf $ is simple if and only
$\G $ is minimal.

A purely algebraic version of \ref {GpdSimplicityCriteria} has been obtained by Nekrashevych in \cite [Proposition
4.1]{Nekr}, where it is also shown that the algebraic version of the {gray} ideal of a certain groupoid algebra
associated to the Grigorchuk group \cite [Example 4.5]{Nekr} is nontrivial.

Unfortunately this sheds little light on the questions listed near the end of the introduction in \cite {CEP}, except
maybe to reinterpret them in terms of the {gray} ideal.  From our point of view, the second question in that list asks
for conditions on a twisted, topologically free \'etale groupoid for the {gray} ideal of $\Cgr $ to vanish.

Our modest contributions to this question are two: one positive and one negative.
  On the bright side, in \ref {MainThree.ii} we are able to give sufficient conditions for a {light} inclusion to be
modeled by the reduced (as opposed to the essential) C*-algebra of its Weyl groupoid, in terms of the existence of
canonical states.  The reduced groupoid C*-algebra is then shown to have a vanishing {gray} ideal.  Our negative result,
on the other hand, consists in showing that the vanishing of the {gray} ideal is not just determined by the structure of
$\G $, as it is influenced by the line bundle as well!  Precisely, in section \ref {GrayIdealTwisted} we give an example
of a non-Hausdorff \'etale groupoid which admits two Fell line bundles, such that the {gray} ideal of $\Cgr $ vanishes
for one such line bundle, but not for the other.

This paper is the result of a collaboration which took place during an extended visit of the first named author to the
Mathematics Department of the University of Nebraska -- Lincoln, and he would like to express his sincere thanks to all
of the members of the Department, especially to his hosts D.~Pitts and A.~Donsig, for their warm hospitality.

\endsection

\PART {I -- INCLUSIONS}

\startsection Local modules

  \def \BIK {\Isot I, K;}
  \def \BJK {\Isot J, K;}
  \def \CIK {\PIA IK; }
  \def \CJI {\PIA JI; }
  \def \CJK {\PIA JK; }
  \def \dist {\hbox {dist}}
  \def \pik {\p IK;}
  \def \pji {\p JI;}
  \def \pjk {\p JK;}
  \def \ui {u \supsz I}
  \def \uj {u \supsz J}

\sectiontitle

One of the most basic aspects of our techniques, to be developed in this and the following section, is the study of
the algebraic relationship between ideals in a C*-algebra $A$, on the one hand, and a larger C*-algebra $B$
containing $A$, on the other.  During this section we will therefore assume the following:

\state {Standing Hypotheses I} \label StandingOne
  \rm
  \izitem
  \zitem $\newsymbol {\incl {B}{A}}{inclusion of C*-algebras}$ will be a fixed \newConcept {inclusion}{inclusion of
C*-algebras} of C*-algebras, that is, $B$ is a C*-algebra and $A$ is a closed $*$-subalgebra of $B$,
  \zitem $I$ and $ J$ will be fixed, closed 2-sided ideals of $A$,
  \zitem for $K\in \{I, J\}$, we will let $\{u\supsz K_i \}_i $ be an approximate unit for the corresponding ideal
and we will put $v\supsz K_i = 1-u\supsz K_i $.

We are making no assumption about the existence of units in our algebras.  Thus the occurrence of ``1'' above is to
be interpreted in $\Mult (B)$, the multiplier algebra of $B$.  In other words, $v\supsz K_i\in \Mult (B)$.

Also notice that the two approximate units above may be assumed to be indexed on the same directed set.

In this work, whenever we speak of the product of two or more linear subspaces of $B$, such as $IBJ$ in the next
lemma, we will always be referring to the \emph {closed linear span} of the set of individual products.  In many
cases the Cohen-Hewitt factorization Theorem \cite [32.22]{CH} implies that the closed linear span is actually no
bigger than the set of individual products, but this will not always be of great importance to us.  Also, when
speaking of the sum of two linear subspaces, we will always be referring to the closure of the set of all sums.

\state Lemma
  Given any $b$ in $B$, one has that
  \izitem
  \lbldzitem LeftIdeal $\ds b\in IB\IFF \li i \ui _i b= b\IFF \li i v\supsz I_i b= 0$,
  \lbldzitem RightIdeal $\ds b\in BJ\IFF \li i b\uj _i = b\IFF \li i bv\supsz J_i = 0$,
  \lbldzitem Hered $\ds b\in IBJ\IFF \li i \ui _i b\uj _i = b$,
  \lbldzitem LeftRightHered $IBJ= IB\cap BJ$.

\Proof The verification of \refl {LeftIdeal} and \refl {RightIdeal} are left to the reader.  Given $b$ in $IB\cap
BJ$, we have for all $i$ that
  $$
  \|b-\ui _i b\uj _i \| \leq
  \|b-\ui _i b\| + \|\ui _i b-\ui _i b\uj _i \| \leq
  \|b-\ui _i b\| + \|\ui _i \|\|b-b\uj _i \|,
  $$
  which converges to zero by \refl {LeftIdeal} and \refl {RightIdeal}, so
  \lbldeq DoisLados
  $$
  b= \li i \ui _i b\uj _i .
  $$

  This shows that $b\in IBJ$, and hence that $IB\cap BJ\subseteq IBJ$.  Since the reverse inclusion is obvious, we
have proven \refl {LeftRightHered}.

Given $b$ in $IBJ$, we then have that $b\in IB\cap BJ$, so the above argument yields \ref {DoisLados}, hence showing
($\Rightarrow $) of \refl {Hered}.  The converse is also clear, so the proof is complete.
  \endProof

\state Proposition \label JCB
  Setting
  $$
  \newsymbol {\CIJ }{defined as $\big \{c\in B: cJ\subseteq IB, \ Ic\subseteq BJ\big \}$}
  = \big \{c\in B: cJ\subseteq IB, \ Ic\subseteq BJ\big \},
  $$
  one has that $A\subseteq \PIA JJ;$, and
  $$
  IJ\subseteq I\CIJ = \CIJ J= IBJ\subseteq \CIJ .
  $$

\Proof
  The fact that $A\subseteq \PIA JJ; $ being trivial, let us next prove that $IBJ\subseteq \CIJ $.  For this, pick
$x\in I$, $b\in B$, and $y\in J$, and observe that
  $$
  xbyJ\subseteq xB\subseteq IB, \and Ixby\subseteq By\subseteq BJ,
  $$
  so $xby\in \CIJ $, proving that $IBJ\subseteq \CIJ $.  Given $x$ in $I$ and $c$ in $\CIJ $, observe that
  $$
  xc\in Ic\subseteq BJ,
  $$
  so we see that
  $I\CIJ \subseteq BJ$.  Since one obviously has $I\CIJ \subseteq IB$, it follows that
  $$
  I\CIJ \subseteq IB\cap BJ\={LeftRightHered} IBJ.
  $$

  Having already proven that
  $$
  I\CIJ \subseteq IBJ\subseteq \CIJ ,
  $$
  we may left-multiply everything by $I$ to get
  $$
  I\CIJ \subseteq IBJ\subseteq I\CIJ ,
  $$
  whence $I\CIJ = IBJ$.  The proof that $\CIJ J= IBJ$ is done along similar lines.

In order to prove the last remaining inclusion, notice that
  $$
  IJ\subseteq IBJ\subseteq \CIJ ,
  $$
  so
  $$
  IJ= IIJ\subseteq I\CIJ .
  \closeProof
  $$
  \endProof

\definition \label DefIsoMod
  The \newConcept {local module}{}, or the \newConcept {isotropy module}{} of the inclusion $\incl {B}{A}$ at the
pair of ideals $(I, J) $ is the quotient Banach space
  $$
  \newsymbol {\BIJ }{local module, defined as $\CIJ /\HIJ $} = \Frac {\CIJ }{\HIJ }.
  $$
  where
  $$
  \newsymbol {\HIJ }{defined as $I\CIJ = \CIJ J= IBJ$} := I\CIJ = \CIJ J= IBJ.
  $$
  We will denote the quotient map by
  $$
  \newsymbol {\pij }{projection from $\CIJ $ to $\BIJ $} : \CIJ \to \BIJ .
  $$
  Finally, in case $I=J$, we will write $\Isot J;$ as a short hand for $\Isot J, J;$, and we will refer to it as the
\newConcept {local algebra}{}, or the \newConcept {isotropy algebra}{} of the inclusion $\incl {B}{A}$ at the ideal
$J$.

\state Proposition
  If, in addition to $I$ and $J$, we are given a third ideal $K\trianglelefteq A$, then
  \izitem
  \zitem $\CIJ \CJK \subseteq \CIK $,
  \zitem $(\CIJ )^*=\CJI $,
  \lbldzitem Juxtap there exists a bi-linear map (denoted simply by juxtaposition)
  $$
  \BIJ \times \BJK \to \BIK ,
  $$
  such that
  $$
  \pij (a) \pjk (b) = \pik (ab), \for a\in \CIJ , \for b\in \CJK ,
  $$
  \lbldzitem StarBun there exists an anti-linear map \quad
  $$
  *: \BIJ \to \BJI ,
  $$
  such that
  $$
  \pij (a)^* = \pji (a^*), \for a\in \CIJ .
  $$

\Proof (i)\enspace
  Given $c_1\in \CIJ $ and $c_2 \in \CJK $, we have
  $$
  c_1c_2K\subseteq c_1JB\subseteq IBB\subseteq IB,
  $$
  and similarly
  $$
  Ic_1c_2 \subseteq BJc_2 \subseteq BBK\subseteq BK,
  $$
  so $c_1c_2$ lies in $\CIK $.

\itmProof (ii) Given $c$ in $\PIA IJ;$, we have that
  $$
  c^*I= (Ic)^* \subseteq (BJ)^* = JB,
  $$
  and similarly $Jc^* \subseteq BI$, so $c^*$ is in $\PIA JI;$.

\itmProof (iii) Multiplying both sides of (i) on the left by $I$, and then again on the right by $K$ gives,
  $$
  I\CIJ \CJK \subseteq I\CIK , \and \CIJ \CJK K\subseteq \CIK K= I\CIK .
  $$

  Denoting by ``$\mu $'' the corresponding restriction of the multiplication operation of $B$, we then see that the
composition of maps
  $$
  \CIJ \times \CJK \labelarrow \mu \CIK \labelarrow \pik \Frac {\CIK } {I\CIK } = \BIK ,
  $$
  vanishes on
  $$
  \HIJ \times \CJK , \quad \hbox {and on}\quad \CIJ \times \H JK;,
  $$
  so it factors through the Cartesian product of quotient spaces
  $$
  \Frac {\CIJ }{\HIJ } \times \Frac {\CJK }{\H JK; } = \BIJ \times \BJK ,
  $$
  producing the required bi-linear map.

\itmProof (iv)
  We have
  $$
  (I\CIJ )^* = \CIJ ^*I^* = \CJI I= J\CJI ,
  $$
  so the composition
  $$
  \CIJ \labelarrow * \CJI \labelarrow \pji \Frac {\CJI }{J\CJI } = \BJI ,
  $$
  vanishes on $\HIJ $, and hence it factors through the quotient producing the desired anti-linear map.
  \endProof

In the special case that $I=J=K$, the above result implies that $\PIA JJ;$ is a C*-algebra, and \ref {JCB} says that
that $\H JJ;$ is a two-sided ideal in $\PIA JJ;$.  Consequently the local algebra $\Isot J, J;$ is indeed a
C*-algebra.

\state Proposition \label viKvi
  Setting
  $$
  \newsymbol {\LIJ }{defined as $IB+ BJ$} = IB+ BJ,
  $$
  (closed linear span) we have for any given $ b$ in $ B$, that
  $$
  b\in \LIJ \IFF \li i v\supsz I_i bv\supsz J_i =
  0 \IFF b= \li i \ui _i b+ b\uj _i - \ui _i b\uj _i .
  $$

\Proof The second ``$\Leftrightarrow $'' is evident, so we prove only the first one, beginning with ``$\Rightarrow
$''.  Observing that
  $$
  \li i v\supsz I_i x= 0 = \li i yv\supsz J_i, \for x\in I, \for y\in J,
  $$
  we deduce that for every $x\in I$, $y\in J$, and every $ a, c\in B$, the element $ b= xa+ cy$, satisfies
  $$
  \li i v\supsz I_i bv\supsz J_i = 0.
  $$
  Since the set of all $ b$'s considered above spans a dense subspace of $ \LIJ $, and since the linear maps
  $$
  b\in \LIJ \mapsto v\supsz I_i bv\supsz J_i \in B
  $$
  are uniformly bounded, the conclusion follows.  In order to prove the first ``$\Leftarrow $'', suppose that a
given $ b$ in $B$ satisfies the condition in the hypothesis.  Then
  $$
  0 =
  \li i v\supsz I_i bv\supsz J_i =
  \li i (1- \ui _i) b(1- \uj _i) =
  \li i b- \ui _i b- b\uj _i + \ui _i b\uj _i,
  $$
  so that
  $$
  b= \li i \ui _i b+ b\uj _i - \ui _i b\uj _i \in IB+ BJ= \LIJ .
  \closeProof
  $$
  \endProof

\state Proposition \label IsometryGlob
  Besides the map $\pij $ of \ref {DefIsoMod}, let
  $$
  \newsymbol {\qij }{projection from $B$ to $\LIJ $} :B\to \Frac {B} {\LIJ },
  $$
  denote the quotient map.  Then:
  \izitem
  \lbldzitem IsometryOne For every $b$ in $B$, one has that $\Vert \qij (b)\Vert = \ds \li i \Vert v\supsz
I_ibv\supsz J_i\Vert $,
  \lbldzitem IsometryTwo There is an isometry
  $$
  \Psi :\BIJ \to \Frac {B}{\LIJ },
  $$
  such that $\Psi \big (\pij (c)\big ) = \qij (c)$, for all $c$ in $\CIJ $.

\Proof (i)\enspace Given $b$ in $B$ and $x$ in $\LIJ $, we have that
  $$
  \Vert b-x\Vert \geq \limsup _{i\to \infty }\Vert v\supsz I_ibv\supsz J_i-v\supsz I_ixv\supsz J_i\Vert \={viKvi}
\limsup _{i\to \infty }\Vert v\supsz I_ibv\supsz J_i\Vert .
  $$
  Taking the infimum as $x$ range in $\LIJ $, we deduce that
  $$
  \Vert \qij (b)\Vert =
  \dist (b, \LIJ ) \geq
  \limsup _{i\to \infty }\Vert v\supsz I_ibv\supsz J_i\Vert \geq
  \liminf _{i\to \infty }\Vert v\supsz I_ibv\supsz J_i\Vert \quebra =
  \liminf _{i\to \infty }\Vert b-\ui _ib-b\uj _i+\ui _ib\uj _i\Vert \geq \Vert \qij (b)\Vert .
  $$
  Equality therefore holds throughout, so $\ds \li i\Vert v\supsz I_ibv\supsz J_i\Vert $ exists and coincides with
$\Vert \qij (b)\Vert $.

\itmProof (ii)
  Since $\HIJ \subseteq \LIJ $, the composition
  $$
  \CIJ \hookrightarrow B\ {\buildrel \qij \over \longrightarrow }\ \Frac {B}{\LIJ }
  $$
  vanishes on $\HIJ $, and therefore it factors through the quotient providing a continuous linear map $\Psi $, as
in the statement, with $\Vert \Psi \Vert \leq 1$, such that
  $$
  \Psi \big (\pij (c)\big ) = \qij (c), \for c\in \CIJ .
  $$
  We will next prove that $\Psi $ is isometric, which is to say that
  $$
  \Vert \pij (c)\Vert =\Vert \qij (c)\Vert , \for c\in \CIJ .
  $$
  Since $\Vert \Psi \Vert \leq 1$, it is clear that $\Vert \qij (c)\Vert \leq \Vert \pij (c)\Vert $.  On the other
hand, the reverse inequality is proved as follows:
  $$
  \Vert \pij (c)\Vert = \Vert \pij (v\supsz I_i cv\supsz J_i)\Vert \leq \Vert v\supsz I_i cv\supsz J_i\Vert .
  $$
  Taking the limit as $i\to \infty $, we get by (i) that $\Vert \pij (c)\Vert \leq \Vert \qij (c)\Vert $, showing
that $\Psi $ is indeed an isometry.
  \endProof

\state Corollary \label Intersection
  One has that $ \CIJ \cap \LIJ = \HIJ $.

\Proof In order to prove $``\supseteq $'', notice that $I\CIJ \subseteq \CIJ $ by \ref {JCB}, while
  $$
  \HIJ = I\CIJ \subseteq IB\subseteq \LIJ .
  $$
  Regarding the reverse inclusion, choose $c$ in $\CIJ \cap \LIJ $.  Then
  $$
  0 = \qij (c) = \Psi \big (\pij (c)\big ),
  $$
  so $\pij (c)=0$, which says that $c\in \HIJ $, as required.
  \endProof

An elementary case of \ref {IsometryOne} which will be useful later is as follows:

\state Lemma \label DistIdeal
  Let $A$ be any C*-algebra and let $J$ be an ideal in $A$.  Also let $\{u_i\}_i$ be an approximate unit for $J$ and
put $v_i = 1-u_i\in \Mult (A)$.  Then
  $$
  \Vert \pi (a)\Vert = \li i \Vert v_iav_i\Vert , \for a\in A,
  $$
  where $\pi :A\to A/J$ denotes the quotient map.

\Proof Plugging in $(A,A,J,J)$ in place of $(B,A,I,J)$ in \ref {StandingOne} observe that
  $$
  \LIJ = JA+ AJ= J,
  $$
  so the map $\qij $ of \ref {IsometryGlob} becomes $\pi $ and the conclusion follows immediately from \ref
{IsometryOne}.
  \endProof

\endsection

\startsection Regular ideals and the localizing projection

  \def \EIJ {\E IJ;}
  \def \EJ {\E J;}
  \def \pj {\p J;}
  \def \qj {\q J;}

\sectiontitle

At the very beginning of this section we will keep the same standing hypothesis as before, namely \ref
{StandingOne}.  The purpose here is to introduce a property for pairs of ideals which will be of fundamental
importance later on.

\definition \label DefineRegIdeal
  We shall say that $(I, J)$ is a \newConcept {regular pair of ideals}{} if the map $\Psi $ introduced in \ref
{IsometryTwo} is surjective.  In case $(J, J)$ is a regular pair, we shall simply say that $J$ is a \newConcept
{regular ideal}{}.

\state Remark \rm
  Recall that for a subset $E\subseteq B$, $E^\perp =\{b\in B: bE=Eb=(0)\}$.  The term `regular ideal' is also
used to denote an ideal $J\subseteq B$ such that $J=(J^\perp )^\perp $.  These terms differ, as can be seen by
considering an inclusion of the form $(C(X), C(X))$.  In this setting, every ideal $J\subseteq C(X)$ is regular in
the sense of \ref {DefineRegIdeal}, yet need not satisfy $J=(J^\perp )^\perp $.  Confusion will not arise, for it
will be clear from the context which meaning is being used.

A useful characterization of regular pairs is as follows:

\state Proposition \label RegSum
  $(I, J) $ is regular if and only if
  $$
  B= \CIJ + \LIJ .
  $$

\Proof
  Given $b$ in $B$ notice that $\qij (b)$ lies in the range of $\Psi $ if and only if there exists $c$ in $\CIJ $
such that $\qij (b) = \qij (c)$, which in turn is equivalent to saying that
  $b\in \CIJ + \LIJ $.
  \endProof

From now on we will boost our standing hypotheses by assuming that:

\state {Standing Hypotheses II} \label StandingTwo
  \rm
  Besides the assumptions of \ref {StandingOne}, from now on we will suppose that $(I, J) $ is a regular pair.

\medskip Under the present assumptions we then have that $\Psi \inv $ is a well defined map so we may consider the
map
  $$
  \EIJ :B\to \BIJ ,
  $$
  given by
  \lbldeq DefineExpecE
  $$
  \EIJ = \Psi \inv \circ \qij .
  $$

\definition
  The map $\newsymbol {\EIJ }{localizing projection} $ described above will be called the
  \newConcept {localizing projection}{} of the inclusion $\incl {B}{A}$ at $(I, J)$.

Useful properties of $\EIJ $ are as follows:

\state Proposition
  For every $c$ in $\CIJ $, and $b$ in $B$, one has that
  \izitem
  \lbldzitem Idempa $\EIJ (c) = \pij (c)$,
  \lbldzitem IdempaDois $b-c\in \LIJ \ \Leftrightarrow \ \qij (b) = \qij (c)\ \Rightarrow \ \EIJ (b) = \pij (c)$.

\Proof
  Left for the reader.
  \endProof

A first useful application of the methods above allows for the description of certain functionals on $B$.

\state Proposition \label {DPEigenGen}
  The correspondence $\ f\to f|_{\CIJ }\ $ establishes a bijection from the set of all continuous linear functionals
$f$ on $B$ vanishing on $\LIJ $, and the set of all continuous linear functionals on $\CIJ $ vanishing on $\HIJ $,
which in turn is obviously in one-to-one correspondence with the set of all continuous linear functionals on the
local module $\BIJ $.

\Proof
  Given that $f$ vanishes on $\LIJ $, it is clear that $f|_{\CIJ }$ vanishes on $\HIJ $ by \ref {Intersection}.
Moreover, if $f|_{\CIJ }=0$, then $f$ itself vanishes by \ref {RegSum}, so we see that our correspondence is
one-to-one.

In order to show surjectivity, let $g$ be a continuous linear functional on $\CIJ $ vanishing on $\HIJ $.  Then $g$
factors through the quotient producing a continuous linear functional $\dot g$ on $\CIJ /\HIJ =B(I,J)$, such that
$g(c)=\dot g\big (\pij (c)\big )$, for every $c$ in $\CIJ $.  Setting $f=\dot g\circ \EIJ $, it is clear that $f$
vanishes on $\LIJ $, and for every $c$ in $\CIJ $ we have that
  $$
  f(c) = \dot g\big (\EIJ (c)\big ) \={Idempa} \dot g\big (\pij (c)\big ) = g(c),
  $$
  so we see that $f|_{\CIJ }=g$.
  \endProof

\state Remark \label {eigenrelate}
  \rm
  Suppose $(A,B)$ is an inclusion with $A$ abelian, let $X=\hat A$, and suppose $I$ and $J$ are the maximal
ideals of $A$ corresponding to $y,x \in X$.  In this setting, the linear functionals vanishing on $\LIJ $ are the
eigenfunctionals on $B$ with range $y$ and source $x$ (see \cite [Definitions 2.1 and 2.2]{DonsigPitts}) and \cite
[Theorem~2.7]{DonsigPitts} also gives a description of the continuous linear functionals on the local module
$B(I,J)$.

\state Lemma \label CKinK
  If $I$, $J$ and $K$ are (not necessarily regular) ideals in $A$, then $\CIJ \L JK; \subseteq \L IK;$ and $\L
IJ;\PIA JK; \subseteq \L IK;$.

\Proof We have
  $$
  \CIJ \L JK; = \CIJ (JB+ BK) \quebra \subseteq \CIJ JB+ \CIJ BK\={JCB}
  I\CIJ B+ \CIJ BK\subseteq IB+ BK= \L IK; .
  $$
  This proves the first assertion and the remaining one follows by taking adjoints.
  \endProof

We will later obtain conditional expectations based on the localizing projections and the following result will be
the technical basis for that.

\state Proposition \label EisExpectation
  For every $c\in \PIA II; $, $d\in \PIA JJ;$ and $b\in B$, one has that
  $$
  \EIJ (cb) = \p II; (c)\EIJ (b), \and \EIJ (bd) = \EIJ (b)\p JJ; (d).
  $$

\Proof
  Given $c$ in $\PIA II; $, and $b$ in $B$, choose $c_1$ in $\CIJ $, and $k$ in $\LIJ $, such that $b= c_1 + k$, by
\ref {RegSum}.  We then have that
  $$
  cb= cc_1 + ck,
  $$
  and since $ck \in \PIA II; \LIJ \subseteq \LIJ $ by \ref {CKinK}, we conclude that
  $$
  \EIJ (cb) \={IdempaDois}
  \pij (cc_1) \={Juxtap}
  \p II; (c)\pij (c_1) \={IdempaDois}
  \p II; (c)\EIJ (b).
  $$
  This proves the first part of the statement and the second one may be proved similarly.
  \endProof

\fix For the remainder of this section we will focus on the special case in which $I=J$, therefore implicitly
assuming that $J$ is a regular ideal.  Accordingly, whenever the pair $(I,J)$ occurs in the notations introduced
above, we will replace it by $J$.  For example $\CIJ $ will be denoted by $\CJ $, and $\LIJ $ will be written as
$\LJ $.

\state Proposition \label VanishJK
  Let $\psi $ be a state on $B$.  Then $\psi $ vanishes on $J$ if and only if $\psi $ vanishes on $\LJ $.

\Proof
  The ``only if'' part is obvious because $J\subseteq \LJ $.  Conversely, if $\psi $ vanishes on $J$, then for every
$x$ in $J$, and every $b$ in $B$, we have by the Cauchy-Schwarz inequality that
  $$
  |\psi (xb)|^2 \leq \psi (xx^*)\ \psi (b^*b) = 0,
  $$
  so $\psi (xb) = 0$, and similarly $\psi (bx) = 0$. This proves that $\psi $ vanishes on $JB+ BJ= \LJ $.
  \endProof

\state Proposition \label UniqueExtension
  Let $\varphi $ be a state on $\CJ $ vanishing on $J$.  Then
  \izitem
  \zitem $\varphi $ vanishes on $\H J; $,
  \zitem $\varphi $ extends uniquely to a state $\psi $ on $B$.  Moreover, the extended state is given by
  $$
  \psi = \dot \varphi \circ \EJ ,
  $$
  where $\dot \varphi $ is the state on $\BJ =\CJ /\H J;$ obtained by passing $\varphi $ to the quotient.

\Proof
  The first point follows by applying \ref {VanishJK} to the inclusion $\incl {\CJ }{A}$.  Regarding (ii), let $\psi
$ be a state on $B$ extending $\varphi $.  Then $\psi $ necessarily vanishes on $J$ and hence also on $\LJ $ by \ref
{VanishJK}.  Given any $b$ in $B$, write $b= c+ k$, with $c\in \CJ $, and $k\in \LJ $, by \ref {RegSum}.  So
  $$
  \psi (b) = \psi (c+ k) = \psi (c) = \varphi (c) = \dot \varphi \big (\pj (c)\big ) \={Idempa}
  \dot \varphi \big (\EJ (b)\big ).
  \closeProof
  $$
  \endProof

\state Proposition \label EPositive
  $\EJ $ is a positive map from $B$ to $\BJ $.

\Proof Given any positive element $b$ in $B$, we need to show that $\EJ (b)$ is positive.  For this, it is enough to
show that $\rho \big (\EJ (b)\big ) \geq 0$, for all states $\rho $ of $\BJ $.  Given any such $\rho $, put $\varphi
= \rho \circ \pj $, so that $\varphi $ is a state on $\CJ $ vanishing on $J$, and clearly $\dot \varphi = \rho $.

Using the fact that every state on a closed $*$-subalgebra admits an extension to a state on the ambient algebra,
choose a state $\psi $ on $B$ extending $\varphi $.  By \ref {UniqueExtension} we have that $\psi = \dot \varphi
\circ \EJ $.  Therefore
  $$
  \rho \big (\EJ (b)\big ) =
  \dot \varphi \big (\EJ (b)\big ) =
  \psi (b) \geq 0.
  \closeProof
  $$
  \endProof

\definition \label DefineCondExp
  Let $D$ be a C*-algebra equipped with a $*$-homomorphism $\iota :A\to \Mult (D)$.  A positive, linear map
  $\EE : B\to D$ is called a \newConcept {generalized conditional expectation}{} if $\Vert \EE \Vert \leq 1$, and
  $$
  \EE (ab) = \iota (a) \EE (b), \and \EE (b) = \EE (b)\iota (a),
  $$
  for all $a$ in $A$ and all $b$ in $B$.  Should we want to highlight the role of the map $\iota $, we will also say
that $\EE $ is a \emph {generalized conditional expectation relative to $\iota $}.

\state Remark
  \rm
  In \cite [Definition 4.12]{Rieffel}, Rieffel also gives a notion of generalized conditional expectation.
While related, our definition of generalized conditional expectation differs from Rieffel's.

Observe that $A\subseteq \CJ $, so we may consider the $*$-homomorphism
  \lbldeq IncludeAinBJ
  $$
  \iota :A\to \Frac {\CJ }{\H J; } = \BJ
  $$
  given as the restriction of the quotient map $\pj $ to $A$.  The result below follows directly from our findings
so far:

\state Proposition \label GetCondExp
  The localizing projection $\EJ $ is a generalized conditional expectation from $B$ to $\BJ $.

In the context of \ref {DefineCondExp} it is sometimes important to require that $\iota $ be non-degenerate in the
sense that $\iota (A)D=D$.  However we currently do not have the means to prove the non-degeneracy of the map $\iota
$ introduced in \ref {IncludeAinBJ}, although this would change once we require, as we eventually will, that $A$
contains an approximate unit for $B$.

Approximate units in fact allow for the following convenient workaround, which applies quite generally.

\state Proposition \label Workaround
  Assuming that $A$ contains an approximate unit for $B$, let $\EE :B\to D$ be a generalized conditional expectation
as in \ref {DefineCondExp}.  Consider the closed $*$-subalgebra $\bar D$ of $D$ generated by $\EE (B)$.  Then
  \izitem
  \zitem $\iota (a)$ is a multiplier of $\bar D$, for every $a$ in $A$,
  \zitem viewing $\iota $ as a map from $A$ to $\Mult (\bar D) $, one has that $\iota $ is non-degenerate,
  \zitem $\EE $ is a generalized conditional expectation from $B$ to $\bar D$ relative to $\iota $,
  \zitem if $B$ is separable, then so is $\bar D$.

\Proof
  (i)\enspace Given $a$ in $A$, it follows from the definition that both $\iota (a)\EE (B)$ and $\EE (B)\iota (a)$
are contained in $\EE (B)$.  Therefore $\iota (a)\bar D$ and $\bar D\iota (a)$ are both contained in $\bar D $,
which is to say that $\iota (a)$ is a multiplier of $\bar D$.

\itmProof (ii) Choosing an approximate unit $\{w_i\}_i$ for $B$ contained in $A$, we claim that
  \lbldeq IotaWi
  $$
  \li i \iota (w_i)x= x, \for x\in \bar D.
  $$
  Indeed, if $x=\EE (b)$, for some $b$ in $B$, then
  $$
  \li i \iota (w_i)x=
  \li i \iota (w_i)\EE (b) =
  \li i \EE (w_ib) =
  \EE (b) = x.
  $$
  Since $\EE (B)$ generates a dense subalgebra of $\bar D$, and since $\Vert \iota (w_i)\Vert \leq 1$, it follows
that \ref {IotaWi} indeed holds for all $x$ in $\bar D$, and hence that $\iota (A)\bar D=\bar D$.

\itmProof (iii) Obvious.

\itmProof (iv) If $\{b_i\}_{i\in {\bf N}}$ is a dense subset of $B$, then $\{\EE (b_i)\}_{i\in {\bf N}}$ is a
countable set of generators for $\bar D $.
  \endProof

The following consequence of \ref {EPositive} will be very useful in the sequel.  In a sense this is the reason for
the introduction of local algebras and localizing projections.

\state Corollary \label CorrespondenceExtendedStates
  Given any state $\rho $ on $\BJ $, the composition
  $$
  \psi : B\labelarrow \EJ \BJ \labelarrow \rho {\bf C}
  $$
  is a state on $B$ vanishing on $J$.  Moreover, the correspondence $\rho \to \psi $ establishes a bijection from
the state space of $\BJ $ and the set of all states on $B$ vanishing on $J$.

\Proof
  Given a state $\rho $, as above, we have that $\rho \circ \EJ $ is a positive linear functional on $B$ by \ref
{EPositive}, which clearly vanishes on $J$, since $\qj $ also does.

Let us compute the norm of $\psi $. Since $\rho $ is a state, we have that $\Vert \rho \Vert =1$, so for every
$\varepsilon >0$, there exists some element $w$ in $\BJ $ such that $\rho (w)=1$, and $\Vert w\Vert <1+\varepsilon
$.  By definition of the quotient norm on $\BJ = \CJ / \H J; $, there exists some $c$ in $\CJ $, such that $\pj
(c)=w$, and $\Vert c\Vert <1+\varepsilon $.
  We then have that
  $$
  \psi (c) =
  \rho \big (\EJ (c)\big ) \={Idempa}
  \rho \big (\pj (c)\big ) =
  \rho (w) =1.
  $$
  Therefore
  $$
  1 = \psi (c) \leq \Vert \psi \Vert \Vert c\Vert < \Vert \psi \Vert (1+\varepsilon ),
  $$
  from where we conclude that $\Vert \psi \Vert >(1+\varepsilon )\inv $.  Since $\varepsilon $ is arbitrary we get
$\Vert \psi \Vert \geq 1$.  Clearly also $\Vert \psi \Vert \leq 1$, so $\psi $ is indeed a state on $B$, as
required.

In order to prove that our correspondence is injective it is enough to observe that $\EJ $ is a surjective map.  On
the other hand, to prove that our correspondence is surjective, pick any state $\psi $ on $B$ vanishing on $J$, and
hence also vanishing on $\LJ $ by \ref {VanishJK}.  The restriction of $\psi $ to $\CJ $ then vanishes on $\H J; $,
so it factors through the quotient, providing a state $\rho $ on $\BJ $, such that
  $$
  \rho \big (\pj (c)\big )= \psi (c), \for c\in \CJ .
  $$

  Given any $b$ in $B$, use \ref {RegSum} to choose $c$ in $\CJ $, such that $b-c\in \LJ $.  Then, recalling that
$\psi $ vanishes on $\LJ $, we get
  $$
  \psi (b) = \psi (c) = \rho \big (\pj (c)\big ) \={IdempaDois} \rho \big (\EJ (b)\big ),
  $$
  so $\psi =\rho \circ \EJ $, as desired.
  \endProof

\state Proposition \label PropoFree
  In addition to \ref {StandingTwo}, assume that $A$ contains an approximate unit for $B$.  Then the following are
equivalent:
  \izitem
  \zitem any state on $A$ vanishing on $J$ admits a unique state extension to $B$,
  \zitem $\pj (A)=\BJ $,
  \zitem $B= A+ \LJ $.

\Proof (i)$\Rightarrow $(ii)\enspace
  Assuming by contradiction that $\pj (A)\neq \BJ $, we have that the former is a closed, $*$-subalgebra of $\BJ
$. So, using the Hahn-Banach extension Theorem we may find a nonzero continuous linear functional $\varphi $ on $\BJ
$ vanishing on $\pj (A)$.  Defining
  $$
  \psi (b) = \varphi (b)+ \overline {\varphi (b^*)}, \for b\in \BJ ,
  $$
  we have that $\psi $ is a self-adjoint functional, which also vanishes on $\pj (A)$.  Letting $\psi _+$ and $\psi
_-$ be the components of the Hahn decomposition of $\psi $ into positive linear functionals we have that
  $$
  \psi =\psi _+-\psi _-,
  $$
  so $\psi _+$ and $\psi _-$ agree on $\pj (A)$.  Given an approximate unit $\{w_i\}_i$ for $B$ contained in $A$, it
is easy to see that $\{\pj (w_i)\}_i$ is an approximate unit for $\BJ $, so
  $$
  \Vert \psi _+\Vert = \sup _i \psi _+(\pj (w_i)) = \sup _i \psi _-(\pj (w_i))= \Vert \psi _-\Vert ,
  $$
  so $\psi _+$ and $\psi _-$ have the same norm which, upon multiplying by a suitable positive scalar, we might
assume is equal to 1.  Consequently $\psi _+$ and $\psi _-$ are distinct states on $\BJ $, coinciding on $\pj (A)$.

Observe that the states $\psi _+\circ \EJ $ and $\psi _-\circ \EJ $ coincide on $A$ because for all $a$ in $A$, one
has
  $$
  \psi _+\big (\EJ (a)\big ) \={Idempa} \psi _+\big (\pj (a)\big ) = \psi _-\big (\pj (a)\big ) = \psi _-\big (\EJ
(a)\big ).
  $$
  However, since $\EJ $ is onto $\BJ $ we have that $\psi _+\circ \EJ \neq \psi _-\circ \EJ $, contradicting (i).

\itmImply (ii) > (iii)
  Given any $b$ in $B$, use \ref {RegSum} to choose $c$ in $\CJ $, such that $b-c\in \LJ $.  By (ii) we then have
that
  $\pj (c) = \pj (a)$, for some $a$ in $A$, whence $c-a\in \H J; $.  We then have that
  $$
  b= a+ (c-a) + (b-c) \in A+ \H J; + \LJ \subseteq A+ \LJ .
  $$
  This proves (ii).

\itmImply (iii) > (i)
  Follows immediately from \ref {VanishJK}.
  \endProof

\endsection

\startsection Regular inclusions

\sectiontitle

The work done in the last two sections should be considered as preparatory for our discussion of regular inclusions,
which is the fundamental objective of this article.  In this short section we will review the most important
properties of this class of inclusions, mostly based on Kumjian's pioneering work on C*-diagonals \cite {Kumjian}.

\definition \label RegAndStuff
  Let $B$ be a C*-algebra and let $A$ be a closed $*$-subalgebra of $B$.
  \izitem
  \zitem An element $n\in B$ is said to be a \newConcept {normalizer}{} of $A$ in $B$, if
  $nAn^*\subseteq A$, and $n^*An\subseteq A$.
  \zitem The set of all normalizers of $A$ in $B$ will be denoted by $\newsymbol {\Norm BA }{set of all normalizers}
$.
  \zitem We will say that $A$ is a \newConcept {non-degenerate subalgebra}{} of $B$, if $A$ contains an approximate
unit for $B$.
  \zitem We will say that $\incl {B}{A}$ is a \newConcept {regular inclusion}{}, if $A$ is a non-degenerate
subalgebra and the closed linear span of $\Norm BA $ is dense in $B$.

Throughout this section we will adopt the assumptions listed below.

\state {Standing Hypotheses III} \label StandingThree
  \rm
  \izitem
  \zitem $\incl {B}{A}$ will be a fixed regular inclusion,
  \zitem we will moreover assume that $A$ is abelian and we will let $X$ denote the Gelfand spectrum of $A$.

\medskip Notice that \ref {RegAndStuff.iii} is equivalent to saying that $AB=B$, which in turn is to say that the
inclusion of $A$ into $B$ is non-degenerate.  This is of course related to our discussion right after \ref
{GetCondExp}.

Let us now list some important facts proved by Kumjian in \cite {Kumjian}.  In doing so, we will adopt the notation
  \lbldeq EvaluationAtX
  $$
  \newsymbol {\evx a}{evaluation of $a$ at $x$}
  $$
  to represent the value of a given element $a\in A$, seen as a function on $X$ via the Gelfand transform, applied
to a point $x\in X$.  Alternatively, if $x$ is seen as a character on $A$, this is meant to represent the value of
this character applied to $a$.

\state Proposition \label KnownFacts
  For every normalizers $n,m\in \Norm BA $, one has that
  \izitem
  \zitem $n^*n$ and $nn^*$ lie in $A$,
  \lbldzitem Beta there exists a homeomorphism
  $$
  \newsymbol {\beta _n}{partial homeomorphism determined by the normalizer $n$} :\Dom {n}\to \Ran {n},
  $$
  where
  $$
  \newsymbol {\Dom {n}}{domain of $\beta _n$} = \{x\in X: \evx {n^*n}\neq 0\}, \and
  \newsymbol {\Ran {n}}{range of $\beta _n$} = \{x\in X: \evx {nn^*}\neq 0\},
  $$
  such that
  $$
  \evx {n^*an} = \ev {a}{\beta _n(x)}\evx {n^*n}, \for a\in A, \for x\in \Dom n.
  $$
  \lbldzitem BetaDomain $\beta _{nm}= \beta _n\circ \beta _m$, where the composition is defined on the biggest
domain where it makes sense, namely on $\Dom {nm} = \beta _m\inv \big (\Ran m\cap \Dom n\big )$,
  \zitem $\beta _n\inv = \beta _{n^*}$.

The following elementary result will be used many times throughout this work.  We leave its proof to the reader.

\state Lemma \label Polys
  Let $n$ be any element of any C*-algebra.  Then there exists a sequence $\{r_i\}_i$ of real polynomials such that
  $$
  n= \li i r_i(nn^*) nn^*n= \li i nn^*nr_i(n^* n).
  $$

The following is a variant of \cite [Lemma~3.1]{PittsSR} (see also \cite [Lemma~2.13]{SRITwo}).  We include a proof
for convenience.

\state Lemma \label VanishDurin
  Given a normalizer $n$, let $\newsymbol {\Fix n}{set of fixed points for $\beta _n$} $ denote the set of fixed
points for $\beta _n$, namely
  $$
  \Fix n= \{x\in \Dom {n}: \beta _n(x)=x\}.
  $$
 Then
  $na=an$, for every $a$ in $A$ which vanishes on $\big (\Dom n\cup \Ran n\big )\setminus \Fix n$.

\Proof
  Given that $a$ vanishes on $\big (\Dom n\cup \Ran n\big )\setminus \Fix n$, we claim that
  \lbldeq DeuIdentidade
  $$
  \evx {n^*an} = \evx {a}\evx {n^*n}, \for x\in X.
  $$
  We will prove this first for
  $x$ in $\Fix n$, then for $x$ in $X\setminus \Dom n$, and finally for $x$ in $\Dom n\setminus \Fix n$.

  The case $x\in \Fix n$ follows immediately from \ref {Beta}, so we next pick $x$ in $X\setminus \Dom n$.  Then
$\evx {n^*n}=0$, so the right-hand-side of \ref {DeuIdentidade} vanishes and we must prove that so does the
left-hand-side.

Letting $\{r_i\}_i$ be as in \ref {Polys}, we have that
  $$
  \evx {n^*an} =
  \li i \EV {n^*ar_i(nn^*) nn^* n}{x} =
  \li i \EV {n^*ar_i(nn^*) n}{x} \evx {n^* n} = 0,
  $$
  as desired.

Finally, picking $x$ in $\Dom n\setminus \Fix n$, we have that $\evx a=0$, by hypothesis, so the right-hand-side of
\ref {DeuIdentidade} again vanishes, and we must once more prove that the left-hand-side vanishes as well.  In order
to do this we claim that $\beta _n(x)$ is in $\Ran n\setminus \Fix n$. Otherwise $\beta _n(\beta _n(x))=\beta
_n(x)$, implying that $\beta _n(x)=x$, which is ruled out by our assumption.  Therefore $\ev {a}{\beta _n(x)}=0$, so
also
  $$
  \evx {n^*an} = \ev {a}{\beta _n(x)}\evx {n^*n} = 0.
  $$
  Since
  $$
  X= \big (X\setminus \Dom n\big ) \ \cup \ \big (\Dom n\setminus \Fix n\big )\ \cup \ \Fix n,
  $$
  we have verified \ref {DeuIdentidade}.  We then conclude that $n^*an= an^*n$, so
  $$
  \Vert an- na\Vert ^2 = \Vert (an- na)^*(an- na)\Vert = 0,
  $$
  from where we see that $na=an$.
  \endProof

\state Proposition \label BetaIdentity
  Given a normalizer $n$, the following are equivalent:
  \izitem
  \zitem $\beta _n$ is the identity on $\Dom n$,
  \zitem $n\in A'$.

\Proof Assuming (i), notice that $\Fix n=\Dom n=\Ran n$, so $\big (\Dom n\cup \Ran n\big )\setminus \Fix n$ is
empty.  Therefore every element $a$ in $A$ satisfies the condition in the last sentence of \ref {VanishDurin} and
hence commutes with $n$.  Therefore $n\in A'$.

\itmImply (ii) > (i) Obvious.
  \endProof

We will later be interested in the $A$-$A$-bimodule generated by a single normalizer, so it is worth noticing the
following simple result.

\state Proposition \label SingleAndBiModules
  Given a normalizer $n$, one has that
  \izitem
  \zitem $\clsr {An} = \clsr {AnA} = \clsr {nA}$,
  \zitem if $m=na$, where $a$ is an element of $A$ not vanishing on any point of $\Dom n$, then $\clsr {An}=\clsr
{Am}$.

\Proof (i)\enspace
  Referring to the polynomials $r_i$ of \ref {Polys}, we have for every $a$ and $b$ in $A$, that
  $$
  anb=
  \li i ann^*nr_i(n^* n)b=
  \li i n(n^*an) r_i(n^* n)b\in \clsr {nA}.
  $$
  This proves that $\clsr {AnA}\subseteq \clsr {nA}$.  Since $n$ belongs to $\clsr {An}$ by \ref {RegAndStuff.iii},
the reverse inclusion is clear, so
  $$
  \clsr {AnA}=\clsr {nA},
  $$
  and the last equality in (i) follows similarly.

\itmProof (ii) Given that $a$ does not vanish on $\Dom n$, one easily proves that $\{(aa^*)^{1/k}\}_{k\in {\bf N}}$
acts as an approximate unit for $C_0(\Dom n)$ (although it might not belong there).  Since $n^*n$ lies in $C_0\big
(\Dom n\big )$, we may employ the polynomials $p_j^k$ and $q_j^k$ used in the proof of \ref {Polys}, to compute
  $$
  n^*n=
  \li k (aa^*)^{1/k}n^*n=
  \li k \li j p_j^k(aa^*)n^*n=
  \li k \li jaa^* q_j^k(aa^*)n^*n\in \clsr {aA}.
  $$
  By \ref {Polys} we also have that $n\in \clsr {nn^*nA}$, so
  $$
  \clsr {nA} \subseteq \clsr {nn^*nA} \subseteq \clsr {naA} = \clsr {mA},
  $$
  Since the reverse inclusion is obvious, (ii) is proved.
  \endProof

\endsection

\startsection Invariant ideals

  \def \HJ {\H J;}

\sectiontitle

We begin this short section by returning to \ref {StandingOne}, with stricter hypotheses explicitly stated when
needed.  The goal here is not so much to subsidize our main work on regular inclusions, but to collect some useful
facts that will be later used to prove some secondary results about simplicity.

\definition \label DefineInvar
  Let $\incl {B}{A}$ be an inclusion.  A closed 2-sided ideal $J\trianglelefteq A$ is said to be \newConcept
{invariant}{invariant ideal} if $BJ=JB$.

If $J$ in an invariant ideal of $B$, notice that
  $$
  \CJ =B, \and \HJ =\LJ =BJ=JB,
  $$
  and that the latter is an ideal in $B$.  Moreover
  $$
  \BJ = \Frac {\CJ }{\HJ } = \Frac {B}{\LJ }.
  $$
  The isometry $\Psi $ of \ref {IsometryTwo} is therefore the identity map, whence $J$ is regular.

\state Proposition \label GeraIntersec
  If $J\trianglelefteq A$ is an invariant ideal, then
  $
  J=
  (JB)\cap A.
  $

\Proof Clearly
  $$
  J=JJ\subseteq JB,
  $$
  so $J\subseteq (JB)\cap A$.  Conversely, if $a\in (JB)\cap A$, then, picking an approximate unit $\{u_i\}_i$ for
$J$, we have that
  $$
  a\={LeftIdeal}
  \li i u_i a\in J.
  \closeProof
  $$
  \endProof

\definition
  We will say that $A$ is \newConcept {$B$-simple}{} if $A$ has no nontrivial closed invariant 2-sided ideals.

\state Proposition \label BsimpleSimple
  If $B$ is simple, then $A$ is $B$-simple.

\Proof
  Let $J$ be an invariant ideal of $A$.  One then has that $JB$ is an ideal of $B$, so either $JB=B$, or $JB=\{0\}$.
Using \ref {GeraIntersec} we deduce that either $J=A$, or $J=\{0\}$.  So $A$ is $B$-simple.
  \endProof

For regular inclusions there is a standard procedure to produce invariant ideals:

\state Proposition \label InterInvar
  Suppose that $\incl {B}{A}$ is a regular inclusion, with $A$ abelian.
  If $K\trianglelefteq B$ is a closed 2-sided ideal, then $K\cap A$ is an invariant ideal.

\Proof
  Letting $J=K\cap A$, let us prove that $BJ\subseteq JB$ or, equivalently, that
  $$
  ba\in JB, \for b\in B, \for a\in J.
  $$
  Since our inclusion is regular, we may assume that $b=n$, where $n$ is a normalizer.  Using the polynomials $r_i$
of \ref {Polys} we then have
  $$
  na=
  \li i nn^*nr_i(n^*n)a=
  \li i nan^*nr_i(n^*n).
  $$
  Observing that $nan^*\in J$, we deduce from the calculation above that $na$ lies in $JB$, proving our claim.  By
taking adjoints we see that also $JB\subseteq BJ$, so $J$ is invariant.
  \endProof

Recall that an inclusion $\incl {B}{A}$ is said to satisfy the \emph {ideal intersection property} when every
nontrivial closed 2-sided ideal of $B$ has a nontrivial intersection with $A$.

\state Proposition \label Simplicity
  Suppose that $\incl {B}{A}$ is a regular inclusion, with $A$ abelian.  Suppose moreover that $\incl {B}{A}$
satisfies the ideal intersection property.  Then $A$ is $B$-simple, if and only if $B$ is simple.

\Proof
  One direction has already been taken care of in \ref {BsimpleSimple}, so let us concentrate on the other one.

Assuming that $A$ is $B$-simple, let $K$ be a closed 2-sided ideal of $B$, and put $J=K\cap A$.  Then $J$ is an
invariant ideal by \ref {InterInvar}, so either $J=A$, or $J=\{0\}$.  In the first case we have
  $$
  B= BA= BJ\subseteq BK\subseteq K,
  $$
  so $B=K$.  In the second case, we have that $K=\{0\}$ by the ideal intersection property.  So $B$ is simple.
  \endProof

\endsection

\startsection Extended multiplication for normalizers

\sectiontitle

Let $n\in \Norm BA$.  Recall that $\Dom n=\{x\in X:\evx {n^*n}\neq 0\}$, so $\Dom n$ is an open subset of $X$.
Notice that $\clsr {n^*nA}$ is precisely $C_0(\Dom n)$ and the multiplier algebra $M(\clsr {n^*nA})$ is $C^b(\Dom
n)$.  Similarly $\clsr {nn^*A}$ and $M(\clsr {nn^*A})$ are $C_0(\Ran n)$ and $C^b(\Ran n)$ respectively.

By \cite [Lemma~2.1]{SRI}, the map $nn^*a\mapsto n^*an$ extends uniquely to a $*$-isomorphism
  $$
  \theta _n: C_0(\Ran n)\rightarrow C_0(\Dom n)
  $$
  satisfying,
  \lbldeq thetanproperty
  $$
  n\theta _n(a)= an \quad \forall a\in C_0(\Ran n).
  $$

Let $\nbimod \subseteq B$ be the norm-closed $A$-bimodule generated by $n$.  By \ref {SingleAndBiModules}, if
$(u_\lambda )$ is an approximate unit for $C_0(\ran n)$ (respectively $C_0(\dom n)$, then for any $b\in \nbimod $,
  \lbldeq aumod
  $$
  u_\lambda b\rightarrow b \quad \text {(respectively } bu_\lambda \rightarrow b).
  $$
  In addition, $(\nbimod )^*(\nbimod )=C_0(\dom n)$ and $(\nbimod )(\nbimod )^*=C_0(\ran n)$.

Clearly $\nbimod $ is an $C_0(\Ran n)-C_0(\Dom n)$ bimodule.  We now show that the left and right multiplications
extend to the multiplier algebras, that is, we show $\nbimod $ is an $C^b(\Ran n)-C^b(\Dom n)$ bimodule.

We may assume that $B\subseteq \BH $ for some Hilbert space $\Hilb $. By the polar decomposition, exists a partial
isometry $w\in \BH $ such that $n=|n^*|w=w|n|$.  Then
  $$
  C^b(\Ran n)=\{T\in \BH : Tww^*=ww^*T=T
  $$
  and
  $$
  TC_0(\Ran n)\cup C_0(\Ran n)T\subseteq C_0(\Ran n)\}.
  $$
  There is an analogous description for $C^b(\Dom n)$.

\state Proposition \label ExtendedMult
  Let $n$ be a normalizer and let $f\in C^b\big (\Dom n\big )$ (resp.~$g\in C^b\big (\Ran n\big )$).  Then there
exists an element in $B$, which we denote by $nf$ (resp.~$gn$), such that
  $$
  \li i n(u_i f) = nf, \qquad
  \big (\hbox {resp.} \li i (gu_i)n= gn\big ),
  $$
  for every approximate unit $\{u_i\}_i$ of\/ $C_0\big (\Dom n\big )$ (resp.~$C_0\big (\Ran n\big )$).

\Proof For any $a\in A$, the map $a|n^*|\mapsto an$ is an isometry, hence extends to an isometry of $C_0(\Ran n)$
onto $\nbimod $ with inverse given by $b\mapsto bw^*$.  It follows that for $g\in C^b(\Ran n)$, $gnw^*=g|n^*|\in
C_0(\Ran n)$, so that $gn\in \nbimod $.  Similar considerations yield $nf\in \nbimod $ for every $f\in C^b(\Dom n)$.
 Thus, for any $g\in C^b(\Ran n)$ and $f\in C^b(\Dom n)$, we define the extended multiplication by
  $$
  gn:=g|n^*|w\in \nbimod ;\dstext {and} nf:=w|n|f\in \nbimod .
  $$
  The statements regarding approximate units follow from \ref {aumod}.
  \endProof

We now extend the definition of $\theta _n$ from $C_0(\ran n)$ to larger classes of functions on subsets of $X$.

\definition \label PartialCompositionUnderlined
  Let $S\subseteq X$ be any subset and $g$ is a complex-valued function defined on $S\cap \ran n$.  We define a
function $\theta _n(g)$ on the set
  $$
  \beta _n^{-1}(S\cap \ran n)=\{x\in \dom n: \beta _n(x)\in S\}
  $$
  using the formula,
  $$
  \evx {\theta _n(g)} =\ev g {\beta _n(x)}.
  $$

The extended notion of multiplication allows us to obtain \ref {thetanproperty} for $a\in C^b(\ran n)$.

\state Proposition \label Intertwiner
  Let $n$ be a normalizer.  For each $g$ in $C^b\big (\Ran n\big )$,
  $$
  gn= n\theta _n(g).
  $$

\Proof For $S=\ran n$, a computation gives
  $$
  \theta _n(nn^*g)=n^*gn\dstext {hence} n\theta _n (nn^*g)=(nn^*g)n.
  $$
  Noting that $\theta _n(nn^*)=n^*n$ we obtain, $n\theta _n(nn^*g)=nn^*n\theta _n(g)$, so $nn^*(n\theta
_n(g)-gn)=0$.  In particular, whenever $p$ is a polynomial vanishing at $0$, we obtain $p(nn^*) (n\theta _n (g)
-gn)=0$.  But $\{p(nn^*): p(0)=0\}$ contains an approximate unit for $C_0(\ran n)$.  Therefore, for any $g\in
C^b(\ran n)$,
  $$
  n\theta _n(g)=gn.
  $$
  \endProof

There is an extra case to be considered when two normalizers are involved and which will be relevant later on.  If
$n, m\in \Norm BA$, observe that $(\bimod n )(\bimod m )=\bimod {nm}$.  If $a\in C^b(\dom n)$, then $na\in \nbimod
$, so $(na)m\in \bimod {nm}$.  While $am$ can be defined using the polar decomposition of $m$, the product need not
belong to $B$, but as just observed, the product $nam$ does belong to $B$.  We wish to identify $c\in C^b(\dom
{nm})$ such that $nam=nmc$.

\state Proposition \label IntertwinerTwo
  Let $n$ and $m$ be normalizers and let $a\in C^b\big (\Dom n\big )$.  Then,
  $$
  nam= nm\theta _m(a).
  $$

\Proof Since $\theta _n^{-1}(a)\in C^b(\ran n)$, the restriction of of $\theta _n^{-1}(a)$ to $\Ran n \cap \beta
_n(\Dom m)$ is an element of $C^b(\Ran {nm})$.  Therefore, $\theta _{nm}(\theta _n^{-1}(a))\in C^b(\dom {nm})$ and
  $$
  nam=\theta _n^{-1}(a) nm=nm \theta _{nm}(\theta _n^{-1}(a)).
  $$
  For $x\in \dom {nm}$,
  $$
  \ev {\theta _{nm}(\theta _n^{-1}(a))} x= \ev {\theta _n^{-1}(a)} {\beta _{nm}(x)} =\ev {\theta _n^{-1}(a)} {\beta
_n(\beta _m(x))}= \ev {\theta _n(\theta _n^{-1}(a))} {\beta _m(x)} \quebra = \ev a {\beta _m(x)}= \ev {\theta _m(a)}
x.
  $$
  Thus $\theta _{nm}(\theta _n^{-1}(a))=\theta _m(a)$, and the result follows.
  \endProof

\endsection

\startsection Regularity of maximal ideals in regular inclusions

  \def \Hx {\H x;}
  \def \Hyx {\H yx;}
  \def \evy #1{\ev {#1}y}
  \def \qx {\q x; }
  \def \qyx {\q yx;}

\sectiontitle

We continue working under \ref {StandingThree}.  As the title suggests, our purpose here is to prove that every
maximal ideal of $A$ is regular, as defined in \ref {DefineRegIdeal}.  Since maximal ideals are in a one to one
correspondence with pure states, these will also play an important role.

\definition \label DataBasedOnPoint
  For each $x$ in $X$ we will denote by $\newsymbol {J_x}{ideal of functions vanishing at $x$}$ the ideal of $A$
given by
  $$
  J_x= \{a\in A: \evx a= 0\}.
  $$
  If $y$ is another point in $X$, we will also let
  $$
  \def \crr {, \vrule height 10pt depth 1pt width 0pt \cr }
  \matrix {
  \newsymbol {\Cyx }{same as $\PIA J_yJ_x;$}
    & = & \PIA J_yJ_x; \crr
  \newsymbol {\Hyx }{same as $\H J_yJ_x;$}
    & = & \H J_yJ_x; \crr
  \newsymbol {\Byx }{same as $B(J_y, J_x)$}
    & = & B(J_y, J_x) \crr
  } \hskip 2cm
  \matrix {
  \newsymbol {\pyx }{same as $\p J_yJ_x;$}
    & = & \p J_yJ_x; \crr
  \newsymbol {\Lyx }{same as $\L J_yJ_x;$}
    & = & \L J_yJ_x; \crr
  \newsymbol {\qyx }{same as $\q J_yJ_x;$}
    & = & \q J_yJ_x; \crr
  }
  $$
  and, in case $x=y$, we will denote these simply as
  $\newsymbol {\Cx }{same as $\PIA xx;$} $,
  $\newsymbol {\Hx }{same as $\H xx;$} $,
  $\newsymbol {\Bx }{same as $B(x, x)$} $,
  $\newsymbol {\px }{same as $\p xx;$} $,
  $\newsymbol {\Lx }{same as $\L xx;$} $ and
  $\newsymbol {\qx }{same as $\q xx;$} $ respectively.

As noted earlier in \ref {eigenrelate}, in the language of~\cite [Defintions~2.1 and 2.2]{DonsigPitts}, the
following characterizes the eigenfunctionals on $B$ with source $x$ and range $y$ as the linear functionals
vanishing on $\Lyx $.

\state Lemma \label VanishLyx
  Given $x, y\in X$, and given any continuous linear map $T $ from $B$ to some other normed space, one has that $T $
vanishes on $\Lyx $, if and only if
  $$
  T (ab) = \evy aT (b), \and T (ba) = T (b)\evx a, \for a\in A, \for b\in B.
  $$

\Proof Throughout this proof we fix an approximate unit $\{w_i\}_i$ for $B$ contained in $A$.
  Assuming that $T $ vanishes on $\Lyx $, let $a\in A$, and $b\in B$.  Then
  $$
  T (ab) =
  \li i T (aw_ib) =
  \li i T \big (aw_ib- \evy aw_ib+ \evy aw_ib\big ) \quebra =
  \li i T \big ((aw_i - \evy aw_i)b) + \evy aT (w_ib\big ) =
  \li i T \big (c_ib) + \evy aT (w_ib\big ) = \cdots
  $$
  where
  $$
  c_i:=aw_i - \evy aw_i.
  $$
  Observing that $c_i\in J_y$, we have that $c_ib\in J_yB\subseteq \Lyx $, so $T \big (c_ib) =0$, and hence the
above equals
  $$
  \cdots =
  \li i \evy aT (w_ib) =
  \evy aT (b).
  $$
  This shows the first identity displayed in the statement, and the second one follows by a similar argument.  In
order to prove the converse, pick $a_1\in J_y$, $a_2\in J_x$ and $b_1,b_2\in B$. Then $a_1b_1+b_2a_2$ is a typical
element of
  $J_yB+BJ_x= \Lyx $, and
  $$
  T (a_1b_1+b_2a_2) =
  \evy {a_1} T (b_1) +T (b_2)\evx {a_2} = 0,
  $$
  so the result follows from the fact that the elements considered span a dense subspace of $\Lyx $.
  \endProof

The consequence for states is as follows.

\state Lemma \label CondexForState
  Any state $\psi $ on $B$ vanishing on $J_x$ satisfies
  $$
  \psi (ab) = \evx a\psi (b) = \psi (ba), \for a\in A, \for b\in B.
  $$

\Proof If $\psi $ is a state vanishing on $J_x$, then it also vanishes on $\Lx $ by \ref {VanishJK}, so the
conclusion follows from the above result.
  \endProof

The next technical lemma will be of crucial importance. Its purpose is to understand the relationship between the
position of a given point $x$ relative to $\Dom n$, on the one hand, and the interplay between the normalizer $n$
and the various subspaces of $B$ determined by $J_x$, on the other.

\state Lemma \label NormaJota
  Given $x$ and $y$ in $X$, and given a normalizer $n$ in $\Norm BA $, one has that:
  \izitem
  \zitem If $x\notin \Dom {n}$, then $n\in BJ_x$.
  \zitem If $x\in \Dom {n}$, then $J_{\beta _n(x)} n\subseteq BJ_x$.
  \zitem If $x\in \Dom {n}$, and $\beta _n(x)\neq y$, then
  $n\in J_yB+ BJ_x= \Lyx $.

\Proof (i) If $\ev {n^*n}{x} = 0$, then surely $n^*n$ belongs to $J_x$, so
  $$
  n= \li k n(n^*n)^{1/k} \in BJ_x.
  $$

\itmProof (ii) For every $a$ in $J_{\beta _n(x)}$, notice that
  $$
  \ev {n^*an}{x} = \ev {a}{\beta _n(x)}\ev {n^*n}{x} = 0,
  $$
  so we see that $n^* an\in J_x$.  Letting $\{r_i\}_i$ be the sequence of polynomials given by \ref {Polys}, we have
  $$
  an=
  \li i ar_i(nn^*) nn^* n=
  \li i r_i(nn^*) nn^* an\in BJ_x.
  $$

\itmProof (iii)
  Assuming that $y\neq \beta _n(x)$, we may choose $v_1$ and $v'_2$ in $A$, such that
  $$
  \ev {v_1}y= 1, \quad \EV {v'_2}{\beta _n(x)} = \ev {n^*n}{x}\inv , \and v_1v'_2 = 0.
  $$

  Setting
  $
  v_2 = n^*v'_2n,
  $
  we then have
  $$
  \ev {v_2}{x} = \ev {n^*v'_2 n}{x} = \ev {v'_2}{\beta _n(x)}\ev {n^*n}{x} = 1,
  $$
  and
  $$
  v_1nv_2 =
  v_1nn^*v'_2n=
  v_1v'_2nn^*n= 0.
  $$

  We next put $u_1 = 1-v_1$, and $u_2 = 1 - v_2$, viewing them in the multiplier algebra of $B$.
  We then claim that
  $$
  u_1n\in J_yB, \and nu_2 \in BJ_x.
  $$
  To see this, let $\{w_i\}_i$ be an approximate unit for $B$ contained in $A$, and notice that
  $$
  u_1w_i = (1-v_1)w_i = w_i - v_1w_i,
  $$
  which belongs to $A$, and clearly vanishes on $y$, whence $u_1w_i \in J_y$. Therefore
  $$
  u_1n= \li i u_1w_in\in J_yB.
  $$
  The proof that $nu_2 \in BJ_x$ is similar.
  We then have that
  $$
  0 = v_1nv_2 =
  (1 - u_1)n(1 - u_2) = n- u_1n- nu_2 + u_1nu_2,
  $$
  so
  $$
  n= u_1n+ nu_2 - u_1nu_2 \in J_yB+ BJ_x.
  \closeProof
  $$
  \endProof

We now have the tools to prove the main result of this section.

\state Proposition \label NorGenCx
  Given $x, y\in X$, put
  $$
  \newsymbol {\Nyx }{set of normalizers sending $x$ to $y$} =\big \{n\in \Norm BA : x\in \Dom n, \hbox { and } \beta
_n(x)=y\big \}.
  $$
  Then
  \izitem
  \lbldzitem NxyInCxy $\Nyx \subseteq \Cyx $,
  \lbldzitem Qzero for every normalizer $n$, one has that $\qyx (n)\neq 0$, if and only if $n\in \Nyx $,
  \lbldzitem Ezero for every normalizer $n$, one has that $\Eyx (n)\neq 0$, if and only if $n\in \Nyx $,
  \lbldzitem NormalizersGenCx if $N$ is a subset of $\Norm BA $ with dense linear span, then $\Byx $ is the closed
linear span of $\pyx (N\cap \Nyx )$,
  \lbldzitem JpairReg $(J_y, J_x) $ is a regular pair of ideals.

\Proof Given $n$ in $\Nyx $, notice that $J_yn\subseteq BJ_x$ by \ref {NormaJota.ii}.  On the other hand it is well
known that $\beta _{n^*}$ is the inverse of $\beta _n$, so we have that $\beta _{n^*}(y)=x$, and then it again
follows from \ref {NormaJota.ii} that $J_xn^*\subseteq BJ_y$.  Taking adjoints this leads to $nJ_x\subseteq J_yB$,
so we see that $n$ lies in $\Cyx $, proving \refl {NxyInCxy}.

We next claim that
  \lbldeq QNulo
  $$
  \qyx (n)=0, \for n\in \Norm BA \setminus \Nyx .
  $$
  To prove this, observe that, for any such $n$, at least one of the conditions defining $\Nyx $ must fail, so we
suppose first that $x$ is not in $\Dom n$, which is to say that $\evx {n^*n}=0$.  Then by \ref {NormaJota.i} we have
that $n\in BJ_x\subseteq \Lyx $, so $\qyx (n)=0$, as desired.  If, on the other hand, $n$ lies in $\Dom n$, then
necessarily $\beta _n(x)\neq y$, so \ref {NormaJota.iii} gives
  $
  n\in \Lyx ,
  $
  whence $\qyx (n) = 0$.  This proves \ref {QNulo} and hence also the ``only if'' part of \refl {Qzero}.

In order to prove ``if'' part, suppose that $n$ is in $\Nyx $.  Then $n$ is also in $\Cyx $ by \refl {NxyInCxy}, so
\ref {IsometryTwo} gives $\Vert \pyx (n)\Vert =\Vert \qyx (n)\Vert $, and we see that it is enough to prove that
$\pyx (n)$ is nonzero which is the same as saying that $n$ is not in $\Hyx = \Cyx J_x$.  Otherwise, given any
approximate unit $\{u_i\}_i$ for $J_x$, one would have
  $$
  \evx {n^*n} \={RightIdeal}
  \li i \evx {n^*nu_i} =
  \li i \evx {n^*n}\evx {u_i} = 0.
  $$
  This is in contradiction with the fact that $x\in \Dom n$, so we see that $n\notin \Hyx $, whence $\pyx (n)$ is
nonzero, completing the proof of \refl {Qzero}.  Since $\Vert \qyx (n)\Vert =\Vert \Eyx (n)\Vert $ by \ref
{IsometryTwo}, we see that \refl {Ezero} follows from \refl {Qzero}.

Focusing now on point \refl {JpairReg} of the statement, we need to show that the map
  $$
  \Psi _{yx}:\Byx = \Frac {\Cyx }{\Hyx } \ \longrightarrow \ \Frac {B}{\Lyx },
  $$
  introduced in \ref {IsometryTwo} is surjective.

Let $N\subseteq \Norm BA $ be any subset whose linear span is dense in $B$, as in \refl {NormalizersGenCx}.  We then
have that $\qyx (B)$ is generated by the set of all $\qyx (n)$, as $n$ run in $N$.  Since the range of $\Psi _{yx}$
is closed by \ref {IsometryTwo}, it is then enough to show that $\qyx (n)$ lies in the range of $\Psi _{yx}$, for
every $n\in N$.  But this follows easily since $\qyx (n) = \Psi _{yx}(\pyx (n))$, for $n$ in $N\cap \Nyx $, by \refl
{NxyInCxy}, whereas $\qyx (n)=0$, for $n$ in $N\setminus \Nyx $.

Observe that the argument above in fact shows that $B/\Lyx $ is generated by $\qyx (N\cap \Nyx )$.  Since we know
that $\Psi _{yx}$ is an isometry from $\Byx $ onto $B/\Lyx $, it follows that $\Byx $ is generated by
  $$
  \Psi _{yx}\inv \big (\qyx (N\cap \Nyx )\big ) = \pyx (N\cap \Nyx ).
  $$
  This proves \refl {NormalizersGenCx}.
  \endProof

\definition
  Given $x$ in $X$, we will let $\Orb x$ be the {orbit} of $x$, namely
  $$
  \Orb x= \{\beta _n(x): n\in \Norm BA, \ \Dom n\ni x\}.
  $$

It should be noticed that $y$ is in $\Orb x$ if and only if $\Nyx $ is nonempty.  Moreover, when $y\notin \Orb x$,
the above result says that $\Byx =\{0\}$, and hence $B=\Lyx $.

Once we know that the pair of ideals $(J_y, J_x) $ is regular, the localizing projection $\E J_yJ_x;$ of \ref
{DefineExpecE} becomes available and we will refer to it by the simplified notation $\Eyx $, or just $\Ex $ when
$x=y$.

Even though we are not assuming either $A$ or $B$ to be unital algebras, the local algebras are always unital as we
will now prove.

\state Lemma \label BXUnital
  Let $x$ be in $X$.  Then
  \izitem
  \zitem given $y\in \Orb x$, $a\in A$, $h\in \Isot y, x;$, and $k\in \Isot x, y;$, one has that
  $$
  h \px (a) = h \Ex (a) = \evx ah, \and \px (a)k = \Ex (a)k = \evx ak,
  $$
  \zitem $\Isot x; $ is a unital C*-algebra,
  \zitem for every $a$ in $A$, one has that $\Ex (a) = \px (a) = \evx a1$.

\Proof (i)\enspace
  Let $\{w_i\}_i$ be an approximate unit for $B$ contained in $A$. Then, for every $a$ in $A$, we have that
  $$
  aw_i - \evx aw_i\in J_x\subseteq J_x\Cx ,
  $$
  so we see that
  $$
  \px (aw_i) = \evx a\px (w_i).
  $$
  Given $h$ in $\Isot y, x;$, choose $c$ in $\PIA yx;$, such that $\pyx (c)=h$, and note that
  $$
  h\Ex (a) \={Idempa}
  h\px (a) =
  \pyx (c)\px (a) =
  \li i \pyx (c)\px (aw_i ) \quebra =
  \evx a\li i \pyx (c)\px (w_i ) \={Juxtap}
  \evx a\li i \pyx (cw_i ) =
  \evx a\pyx (c) =
  \evx ah,
  $$
  proving the first part of (i).  The second part follows either by a similar argument, or by taking adjoints.

Choosing $a_0$ in $A$ such that $\evx {a_0}=1$, and applying (i) to the case that $y=x$, we deduce that
  $$
  h\px (a_0 ) = \px (a_0 )h = h, \for h\in \Isot x, x;,
  $$
  so it follows that $\px (a_0)$ is the unit of $\Isot x; $, proving (ii).

Regarding (iii), for all $a$ in $A$ we have that
  $a-\evx aa_0\in J_x$, so
  $$
  0 = \px \big (a-\evx aa_0\big ) = \px (a)-\evx a1= \Ex (a)-\evx a1.
  \closeProof
  $$
  \endProof

\endsection

\startsection Extension of pure states, relative free points and smooth normalizers

\sectiontitle

We shall continue to work under \ref {StandingThree}.  Nevertheless the reader might notice that, up to \ref
{MandatoryValues}, we need only assume the milder standing hypothesis that $\incl {B}{A}$ is an inclusion of
C*-algebras, with $A$ abelian, its spectrum being denoted $X$.

The purpose of this section is to conduct an in-depth study of the uniqueness of extensions of pure states of $A$ to
$B$.  The notions of free points and of smooth normalizers will spring from this study, and they will be revealed as
central ingredients in our main results.

Given $x$ in $X$, consider the associated pure state on $A$, namely
  \lbldeq Phix
  $$
  \newsymbol {\varphi _x}{evaluation state at $x$} : a\in A\mapsto \evx a\in {\bf C}.
  $$
  It is well known that $\varphi _x$ admits at least one extension to a state on $B$, and we will now discuss to
what extent such an extension is unique.

\definition \label DefineRelFree
  Given $b$ in B, we shall say that a given point $x$ in $X$ is \newConcept {free relative to $b$}{free point
relative to $b$}, if for every two states $\psi _1$ and $\psi _2$ extending $\varphi _x$, one has that $\psi
_1(b)=\psi _2(b)$.  The set of all such points will be denoted by $\newsymbol {F_b}{set of free points relative to
$b$}$.  We will moreover denote by $\varepsilon _b$ the complex valued function defined for every $x$ in $F_b$ by
  $$
  \newsymbol {\varepsilon _b(x)}{value of unique extension} = \psi (b),
  $$
  where $\psi $ is any state extension of $\varphi _x$.

Examples show that the set $F_b$ may be very poorly behaved.  In \ref {BadNormalizer}, below, we exhibit an example
in which $X$ is a closed interval of the real line while $F_b$ consists of the irrational numbers in this interval.

\state Proposition \label UniqueForBA
  Given $b$ in $B$ and $a$ in $A$, one has that
  \izitem
  \zitem $ F_{ab}=F_{ba}$,
  \zitem $F_b\subseteq F_{ab}$,
  \zitem $\varepsilon _{ba}=\varepsilon _{ab}=a\, \varepsilon _b$, on $F_b$.

\Proof
  Given $x$ in $F_{ab}$, let $\psi _1$ and $\psi _2$ be two states on $B$ extending $\varphi _x$.  Then
  $$
  \psi _1(ba) \={CondexForState}
  \psi _1(ab) =
  \psi _2(ab) \={CondexForState}
  \psi _2(ba),
  $$
  so $x$ lies in $F_{ba}$.  This shows that $F_{ab}\subseteq F_{ba}$, and the reverse inclusion follows similarly.

  Given $x$ in $F_b$, let $\psi _1$ and $\psi _2$ be two states on $B$ extending $\varphi _x$.  Then
  $$
  \psi _1(ab) \={CondexForState}
  \evx a\psi _1(b) =
  \evx a\psi _2(b) \={CondexForState}
  \psi _2(ab).
  $$
  This shows that $x$ lies in $F_{ab}$ and also that $\varepsilon _{ab}(x)=\evx a\varepsilon _b(x)$.  Therefore
$\varepsilon _{ab}=a\, \varepsilon _b$, on $F_b$.

The proof that $\varepsilon _{ba}=a\, \varepsilon _b$, on $F_b$ may be done in a similar way.
  \endProof

Even though $F_b$ may be badly behaved, $\varepsilon _b$ does not follow suit.

\state Proposition \label EbContinuous
  For every $b$ in $B$, one has that $\varepsilon _b$ is a continuous function on $F_b$.

\Proof
  Suppose by way of contradiction that $\{x_\lambda \}_{\lambda \in \Lambda }$ is a net in $F_b$, converging to some
$x$ in $F_b$, but such that $\{\varepsilon _b(x_\lambda )\}_{\lambda \in \Lambda }$ does not converge to
$\varepsilon _b(x)$.  Therefore there exists some positive number $r$ such that
  \lbldeq NoCOnvergence
  $$
  |\varepsilon _b(x_\lambda )-\varepsilon _b(x)|\geq r,
  $$
  for every $\lambda $ in a co-final subset $\Lambda _1\subseteq \Lambda $.  Observing that the subnet $\{x_\lambda
\}_{\lambda \in \Lambda _1}$ also converges to $x$, and upon replacing $\Lambda $ by $\Lambda _1$, we may assume
without loss of generality that \ref {NoCOnvergence} in fact holds for all $\lambda $ in $\Lambda $.

For each $\lambda $ in $\Lambda $ let us choose a state $\psi _{x_\lambda }$ on $B$ extending $\varphi _{x_\lambda
}$.  Using Alaoglu's Theorem we know that the set $\cal P$ formed by all positive linear functionals on $B$, of norm
no greater that one, is compact in the weak* topology.  Therefore there exists a subnet $\{\psi _{x_{\lambda _\mu
}}\}_\mu $ converging in $\cal P$, whose limit we denote by $\psi $.
  Given any $a$ in $A$, we then have that,
  $$
  \psi (a) =
  \li \mu \psi _{x_{\lambda _\mu }}(a) =
  \li \mu \ev a{x_{\lambda _\mu }} =
  \ev a{\li \mu x_{\lambda _\mu }} =
  \evx a,
  $$
  so we see that $\psi $ extends $\varphi _x$, and also that $\Vert \psi \Vert =1$, so $\psi $ is a state\fn
{Observe that the set of states in a non-unital C*-algebra is not necessarily compact.}.  Since $x$ lies in $F_n$,
we may use any state extension of $\varphi _x$ to compute $\varepsilon _b(x)$ so, in particular,
  $$
  \varepsilon _b(x) = \psi (b) = \li \mu \psi _{x_{\lambda _\mu }}(b) = \li \mu \varepsilon _b(x_{\lambda _\mu }),
  $$
  which is in contradiction with \ref {NoCOnvergence}.  This concludes the proof.
  \endProof

\definition \label DefineContSmooth
  We shall say that a given element $b$ in $B$ is:
  \izitem
  \zitem \newConcept {continuous} {continuous element} if $F_b$ is dense in $X$,
  \zitem \newConcept {smooth} {smooth element} if the interior of $F_b$ is dense in $X$.

The next result will find no later use in this work but it is perhaps an interesting curiosity.

\state Proposition
  The closure in $B$ of the set of smooth elements consists of continuous elements.

\Proof
  Let $\{b_i\}_{i\in \N }$ be a sequence of smooth elements converging to some $b$ in $B$.  By Baire's theorem we
have that
  $$
  Y:= \medcap _{i\in \N }F_{b_i}
  $$
  is dense.  The result will therefore be proved once we show that $Y\subseteq F_b$. So pick $x$ in $Y$, and let
$\psi _1$ and $\psi _2$ be states on $B$ extending $\varphi _x$.  Then
  $$
  \psi _1(b) = \li i \psi _1(b_i) = \li i \psi _2(b_i) = \psi _2(b),
  $$
  the middle equality due to the fact that $x$ lies in $F_{b_i}$.  This shows that $x$ lies in $F_b$, concluding the
proof.
  \endProof

From now on, and until the end of this section, we shall make effective use of \ref {StandingThree}.

\medskip Given a normalizer $n\in \Norm BA $, consider its associated partial homeomorphism
  $
  \beta _n:\Dom {n}\to \Ran {n},
  $
  and recall that $\Fix n$ denotes the set of fixed points for $\beta _n$.  As usual we denote the interior of $\Fix
n$ by \pilar {10pt}$\iFix n$, and we will put
  $$
  \bd n=\Fix n\setminus \iFix n.
  $$
  Since $\Fix n$ is closed relative to $\Dom {n}$, we see that $\bd n$ is the boundary of $\Fix n$ relative to $\Dom
{n}$.  However notice that $\bd n$ is not necessarily the boundary of $\Fix n$ relative to $X$, as it might not even
be closed there.

\state Proposition \label MandatoryValues
  Let $n$ be a normalizer in $\Norm BA $.  Then
  \izitem
  \zitem $X=\Fix n \cup F_n$,
  \zitem $\varepsilon _n$ vanishes on $X\setminus \Fix n$.

\Proof
  Pick $x$ in $X\setminus \Fix n$, and let $\psi $ be a state on $B$ extending $\varphi _x$.
  Observe that, since $\psi $ vanishes on $J_x$, it follows from \ref {VanishJK} that $\psi $ also vanishes on $\Lx
=J_xB+BJ_x$.

If $x\notin \Dom {n}$, we have by \ref {NormaJota.i} that $n\in BJ_x\subseteq \Lx $, so $\psi (n)=0$.  On the other
hand, if $x\in \Dom n$, then $\beta _n(x)\neq x$ by hypothesis so, plugging $y=x$ in \ref {NormaJota.iii}, again
leads to $n\in \Lx $, so $\psi (n)=0$.
  This shows that $x$ is free relative to $n$, hence proving (i).  Finally, since all states $\psi $ above were
shown to vanish on $n$, we have also proved (ii).
  \endProof

As uniqueness of state extensions is among our biggest worries, it is worth pointing out the following situation in
which uniqueness is hopeless.

\state Proposition \label TwoExtensions
  Given a normalizer $n$ in $\Norm BA $, one has that
  \izitem
  \zitem for every $x$ in $\Fix n$, there exists a state $\psi $ on $B$ extending $\varphi _x$, such that $\psi
(n)\neq 0$,
  \zitem $\bd n\cap F_n = \emptyset $,
  \zitem in particular, for each $x$ in $\bd n$, there are at least two distinct states on $B$ extending $\varphi
_x$.

\Proof (i)\enspace Given $x$ in $\Fix n$, we have that $n\in \Npt xx;$, so \ref {Ezero} implies that $\Ex (n)\neq
0$.  We may therefore find a state $\rho $ on $B(x) $ with $\rho \big (\Ex (n)\big )\neq 0$, so \ref
{CorrespondenceExtendedStates} implies that the state
  $$
  \psi :=\rho \circ \Ex
  $$
  vanishes on $J_x$, while clearly $\psi (n)\neq 0$.  The restriction of $\psi $ to $A$ is then a positive linear
functional which clearly vanishes on $J_x$.  One also has that $\psi |_A$ is a state on $A$ by \ref
{RegAndStuff.iii}, so it is clear that $\psi |_A=\varphi _x$.  This proves (i).

\itmProof (ii)
  Suppose by contradiction that $x$ lies in $\bd n\cap F_n$.
  Referring to the state $\psi $ given by (i), we then have that
  $$
  \varepsilon _n(x) = \psi (x)\neq 0.
  $$
  On the other hand, choosing a net $\{x_i\}_i$ in $\Dom n\setminus \Fix n$ converging to $x$, we have by \ref
{MandatoryValues.i} that $x_i\in F_n$, and then
  $$
  \varepsilon _n(x)= \li i \varepsilon _n(x_i)\={MandatoryValues.ii} 0,
  $$
  a contradiction.

\itmProof (iii) Follows immediately from (ii).
  \endProof

Putting together \ref {TwoExtensions.ii} and \ref {MandatoryValues.i}, observe that
  $$
  \bd n\subseteq X\setminus F_n\subseteq \Fix n.
  $$
  Definition \ref {DefineContSmooth} evidently suggests that free points are desirable, and the more of them, the
better.  For the special case of normalizers, we then have from the above that any possible trouble maker (i.e. a
non-free point) must reside in $\Fix n$.  This does not mean that every point of $\Fix n$ is a trouble maker but, as
seen above, those in the boundary definitely are!

\endsection

\startsection Free points

\sectiontitle

We shall continue to work under \ref {StandingThree}.  In the previous section we studied free points relative to a
given element of $B$.  Now we will be interested in another, non-relative version of freeness.

\definition \label DefFreePts
  \izitem
  \zitem A point $x$ in $X$ will be called \newConcept {free}{free point} when the state $\varphi _x$ on $A$ given
by evaluation at $x$ admits a unique extension to a state on $B$.
  \zitem The set of all free points of $X$ will henceforth be denoted by $\newsymbol {F}{set of all free points}$.

\medskip As already observed, the states on $B$ extending $\varphi _x$ are exactly those vanishing on $J_x$.
Consequently we have that $x$ is free if and only $J_x$ satisfies the equivalent conditions of \ref {PropoFree}.

The following reconciles the two version of freeness so far discussed.

\state Proposition \label XfAndXb
  Let ${\cal S}$ be a subset of $B$ which generates $B$ as an $A$-$A$-bimodule, meaning that $A{\cal S}A$ spans a
dense subspace of $B$.
  Then
  $$
  F= \medcap _{b\in {\cal S}} F_b.
  $$

\Proof We leave it for the reader to ponder the obvious reasons behind ``$\subseteq $'', while we concentrate on the
reverse inclusion.  Given any $x$ in $\medcap _{b\in {\cal S}} F_b$, and given any two states $\psi _1$ and $\psi
_2$ on $B$ extending $\varphi _x$, we have for every $b$ in ${\cal S}$, and $a$ and $a'$ in $A$ that
  $$
  \psi _1(aba')\={CondexForState}
  \evx {a}\psi _1(b)\evx {a'} =
  \evx {a}\psi _2(b)\evx {a'} \={CondexForState}
  \psi _2(aba'),
  $$
  so we see that $\psi _1$ and $\psi _2$ coincide on $A{\cal S}A$, and hence everywhere since $A{\cal S}A$ spans a
dense subspace.  So $x\in F$.
  \endProof

The following is a useful characterization of free points.

\state Proposition \label IsotropyOneDim
  Given $x$ in $X$, one has that $x$ is free, if and only if the isotropy algebra $\Bx $ is isomorphic to ${\bf C}$.

\Proof Follows immediately from \ref {PropoFree} and \ref {BXUnital.iii}.
  \endProof

There are several complex numbers associated to each free point, and each element of $B$, which turn out to
coincide:

\state Proposition \label ThreeValues
  Given $x$ in $F$, let $\psi _x$ be the unique state on $B$ extending $\varphi _{x}$, and recall that $\Ex $
denotes the localizing projection at the ideal $J_x$.  Then, identifying $\Bx $ with ${\bf C}$ in the obvious way,
we have that
  $$
  \Ex (b) = \psi _x(b) =\varepsilon _b(x), \for b\in B.
  $$

\Proof Observe that $\psi _x$ vanishes on $J_x$, so we may use the last sentence of \ref
{CorrespondenceExtendedStates} to produce a state $\rho $ on $\Bx $, such that $\psi _x=\rho \circ \Ex $.  However,
under the present hypothesis we have that $\Bx ={\bf C}$ in view of \ref {IsotropyOneDim}, so the state $\rho $ must
necessarily be the identity map, whence $\Ex = \psi _x$.

Plugging ${\cal S}=B$ in \ref {XfAndXb}, we deduce that $x$ is necessarily in $F_b$, for every $b$ in $B$, so the
second equality displayed in the statement follows from the definition of $\varepsilon _b$.
  \endProof

\state Proposition \label IntroFreeCondExp
  For $b$ in $B$, let $\Pf (b)$ be the bounded function on $F$ given by
  $$
  \Pf (b)\calcat x= \psi _x(b)
 , \for x\in F,
  $$
  where $\psi _x$ is the unique state on $B$ extending $\varphi _{x}$.  Then $\Pf (b)$ is continuous on $F$.

\Proof
  Follows immediately from \ref {EbContinuous} and the fact that $\Pf (b)$ is the restriction of $\varepsilon _b$ to
$F$, according to \ref {ThreeValues}.
  \endProof

We may therefore view $\Pf $ as a map
  $$
  \Pf : B\to C^b(F).
  $$

\state Proposition \label DefineFreeCondExp
  The map $\Pf $ defined above is a generalized conditional expectation from $B$ to $C^b(F)$, relative to the \emph
{restriction} map
  $$
  \iota :a\in A\mapsto a|_{F} \in C^b(F).
  $$
  We will therefore refer to $\newsymbol {\Pf }{free expectation}$ as the \newConcept {free expectation}{} for the
inclusion $\incl {B}{A}$.

\Proof
  It is evident that $\Pf $ is a positive, continuous map, with $\Vert \Pf \Vert \leq 1$.  Moreover, for all $a$ in
$A$, $b$ in $B$, and $x$ in $F$, one has that
  $$
  \Pf (ab)\calcat x= \psi _x(ab) \={CondexForState} \evx a\psi _x(b) = \evx a\Pf (b)\calcat x,
  $$
  so $\Pf (ab)=a\Pf (b)$ on $F$, and it is also clear that $\Pf (ba)= \Pf (b)a$.  This concludes the proof.
  \endProof

\endsection

\startsection Fourier coefficients

\sectiontitle

\label FourierSection
  We continue under \ref {StandingThree}.  As seen in \ref {JpairReg}, every pair of ideals of the form $(J_y, J_x)
$ is regular, hence yielding the localizing projection $\Eyx $.  In this section we would like to show that, for
every fixed $b$ in $B$, the correspondence
  $$
  (y,x)\in X\times X\mapsto \Eyx (b)\in \Byx
  $$
  behaves somewhat like the classical Fourier transform.  See also \cite [II.4.2]{Renault}.

As already observed after \ref {NorGenCx}, if $y$ is not in the orbit of $x$ under the action of the $\beta _n$,
then $\Byx =\{0\}$, so $\Eyx (b)$ vanishes for every $b$.  It is therefore sensible to discard this superfluous
information.

\state Proposition \label OrbEquiv
  Given $x$ and $y$ in $X$, the following are equivalent:
  \izitem
  \zitem $y\in \Orb x$,
  \zitem $\Nyx \neq \emptyset $,
  \zitem $\Byx \neq \{0\}$,
  \zitem there exists some $b$ in $B$, such that $\Eyx (b)\neq 0$.
  \medskip \noindent In addition, writing
  $$
  x\sim y,
  $$
  to express that the above equivalent conditions hold, we have that ``$\sim $'' is an equivalence relation on $X$.

\Proof
  Evidently (iv)$\Rightarrow $(iii), while (iii)$\Rightarrow $(ii) by \ref {NormalizersGenCx}.  The equivalence
(i)$\Leftrightarrow $(ii) follows immediately from the definitions, and we finally have that (ii)$\Rightarrow $(iv)
by \ref {Ezero}.  The last sentence in the statement is an obvious consequence of \ref {KnownFacts}.
  \endProof

We may therefore consider the set
  $$
  \Xt = \big \{(y, x)\in X\times X: x\sim y\big \},
  $$
  which, in strict technical terms, \emph {is} the equivalence relation introduced above.  As usual we may view $\Xt
$ as a groupoid under the multiplication operation such that
  $(w,z)(y,x)$ is defined if and only if $z=y$, in which case the product is set to be $(w,x)$.  We shall not be
concerned with topological aspects of $\Xt $, so we will see it simply as a discrete groupoid (we nevertheless call
the reader's attention to question \ref {FellBundleOverEqulvRel}).

\state Proposition \label BxyFellBun
  Let
  $$
  \Bun \kern 5pt= \kern -5pt \bigsqcup _{(y, x)\in \Xt } \kern -5pt \Byx ,
  $$
  (disjoint union) and, for each $x,y,z\in X$, such that both $(z,y)$ and $(y,x)$ lie in $\Xt $, consider the
bi-linear operation
  $$
  \Isot z, y; \times \Isot y, x; \to \Isot z, x;
  $$
  introduced in \ref {Juxtap} as well as the anti-linear map
  $$
  *:\Isot y, x; \to \Isot x, y;
  $$
  of \ref {StarBun}.  Then $\Bun $ is a Fell bundle \cite {KumjianFell} over the discrete groupoid $\Xt $.

\Proof
  Among the ten axioms listed in \cite [\S 2]{KumjianFell}, we check only the last two, leaving the remaining ones
for the reader.

Given $(y,x)\in \Xt $, and given any element $e\in \Isot y, x; $, we may write $e=\pyx (c)$, for some $c\in \PIA yx;
$.  So
  $$
  e^*e =
  \pyx (c)^* \pyx (c) \={StarBun}
  \p xy; (c^*) \pyx (c) \={Juxtap}
  \p xx; (c^*c) \quebra =
  \def \modc {|c\kern 0.5pt|}
  \p xx; (\modc ^2) =
  \p xx; (\modc )^* \p xx; (\modc ) \geq 0.
  $$
  This proves \cite [2.1.x]{KumjianFell}.  In addition, from the above it follows that
  \lbldeq OlhaOBreque
  $$
  \Vert e^*e\Vert =
  \Vert \p xx; (c^*c) \Vert = \cdots
  $$

  Observe that by \ref {JCB}, we have that
  $$
  J_x\subseteq J_x\PIA xx; = \PIA xx; J_x\trianglelefteq \PIA xx; .
  $$
  Choosing an approximate unit $\{u_i\}_i$ for $J_x$, we then have that $\{u_i\}_i$ is also an approximate unit for
$J_x\PIA xx; $, so if we set $ v_i=1-u_i, $ it follows from \ref {DistIdeal} that \ref {OlhaOBreque} equals
  $$
  \cdots = \li i \Vert v_ic^*cv_i \Vert =
  \li i \Vert cv_i \Vert ^2 \geq
  \li i \Vert \pyx (cv_i)\Vert ^2 =
  \Vert \pyx (c)\Vert ^2 =
  \Vert e\Vert ^2.
  $$
  This proves \cite [2.1.ix]{KumjianFell}.
  \endProof

We will shortly show that $\Bun $ is a saturated Fell bundle \cite [2.4]{KumjianFell} but, in order to do so, we
need the following technical fact:

\state Lemma \label BasisForSat
  Given $x$ and $y$ in $X$, and given $n$ in $\Npt yx;$, one has
  $$
  \Isot z,y; \pyx (n) = \Isot z,x;
  $$

\Proof
  The inclusion from left to right is proved as follows:
  $$
  \Isot z,y; \pyx (n) \subseteq
  \Isot z,y; \Isot y, x; \explain {Juxtap}\subseteq \Isot z, x;.
  $$
  On the other hand, given any $h\in \Isot z, x;$, write $h=\p zx;(c)$ for some $c$ in $\PIA zx;$.  Since $n^*$ is
in $\Npt xy;$, we have
  $$
  cn^*\in \PIA zx;\PIA xy;\subseteq \PIA zy;.
  $$
  Therefore
  $$
  \Isot z, y;\pyx (n)\ni
  \p zy;(cn^*)\pyx (n) \={Juxtap}
  \p zx;(cn^*n) =
  \p zx;(c)\p xx;(n^*n) \={BXUnital}
  \evx {n^*n} h.
  $$
  Since $\Isot z, y;\pyx (n)$ is a linear space and since $\evx {n^*n}\neq 0$, we conclude that $h\in \Isot z,
y;\pyx (n)$.
  \endProof

\state Proposition
  $\Bun $ is a saturated Fell bundle.

\Proof
  Given $x, y, z$ in $X$, such that $(z,y)$ and $(y,x)$ belong to $\Xt $, we need to prove that
  $$
  \Isot z,x; = \Isot z, y;\Isot y, x;
  $$
  (closed linear span). Choosing a normalizer $n$ in $\Npt yx;$, we get by \ref {BasisForSat} that
  $$
  \Isot z,x; =
  \Isot z,y; \pyx (n) \subseteq \Isot z, y;\Isot y, x;.
  $$
  Since the reverse inclusion is evident, the proof is complete.
  \endProof

The purpose of the next definition is to call attention to the relationship between the present context and
classical Fourier analysis.

\definition \label FTofb
  Given $b$ in $B$, the \newConcept {Fourier transform}{} of $b$ is the section of the bundle $\Bun $ given by
  $$
  \hat b:(y, x)\in \Xt \to \ft yxb\in \Byx .
  $$

In order to avoid duplication of notation we will continue to use $\ft yxb$, rather than $\hat b(y, x)$.

The next result is reminiscent of the fact that, for a periodic function $f$ on the real line, one has that
$\widehat {f^*}(n)= \overline {f(-n)}$.

\state Proposition \label FourAdjoin
  Given $b$ in $B$, and $(y,x)$ in $\Xt $, one has that
  $$
  \ft yx{b}^* = \ft xy{b^*}.
  $$

\Proof
  Pick $c$ in $\Cyx $ such that $b-c\in \Lyx $, and observe that
  $$
  b^*-c^*\in (\Lyx )^* = (J_yB+BJ_x)^* = BJ_y+ J_xB= \L xy;,
  $$
  so
  $$
  \ft xy{b^*} \={IdempaDois}
  \p xy; (c^*) \={StarBun}
  \pyx (c)^* \={IdempaDois}
  \ft yx{b}^*.
  \closeProof
  $$
  \endProof

Our next medium term goal is to prove a result related to the well known fact that the Fourier transform of the
product of two functions is the convolution product of the corresponding Fourier transforms.  This will require some
preparations so we begin by identifying what might be loosely thought of as the trigonometrical polynomials in our
context.

\state Proposition \label RowColFinite
  Let $B_0$ be the linear span (no closure) of $\Norm BA $ in $B$.  Then $B_0$ is a dense $*$-subalgebra of $B$.
Moreover, for every $b$ in $B_0$, and for every $x$ in $X$, the sets
  $$
  \{y\in X: \ft yxb\neq 0\}, \and \{z\in X: \ft xzb\neq 0\}
  $$
  are finite.

\Proof
  That $B_0$ is a $*$-subalgebra follows from the fact that $\Norm BA $ is closed under multiplication and adjoint.
It is dense in $B$ as a consequence of our assumption that $\incl {B}{A}$ is regular.  Writing any given $b$ in
$B_0$ as
  $$
  \ds b=\sum _{n\in F} n,
  $$
  where $F$ is a finite subset of $\Norm BA $, and given $y$ in $X$, suppose that
  $$
  0\neq \ft yxb= \sum _{n\in F} \ft yxn.
  $$
  It obviously follows that there exists some $n$ in $F$ such that $\ft yxn\neq 0$, so \ref {Ezero} implies that
$n\in \Nyx $, which is to say that $y=\beta _n(x)$, besides of course that $x\in \Dom n$.  The upshot is that $\ft
yxb$ can only be nonzero for $y$ in the finite set
  $$
  \{\beta _n(x): n\in F, \ \Dom n\ni x\},
  $$
  from where we see that the first set displayed in the statement is finite.  The proof that the same holds for the
second set is similar.
  \endProof

The announced version of the convolution formula in Fourier analysis may now be given, at least in a special case.

\state Proposition \label BabyConvol
  Given $a,b\in B_0$, and given $(y,x)\in \Xt $, one has that
  \lbldeq expectedconvform
  $$
  \ft yx{ab} = \sum _{z\in \Orb x} \ft yza\ft zxb.
  $$

\Proof
  By \ref {RowColFinite} we see that the above sum has at most finitely many nonzero terms, so convergence is not
yet an issue.

Since both sides of \ref {expectedconvform} represent linear maps in each variable $a$ and $b$, we may assume
without loss of generality that $a$ and $b$ are themselves normalizers, so we will refer to them as $m$ and $n$,
respectively.

Let us assume first that $\ft yx{mn} \neq 0$.  We then deduce from \ref {Ezero} that $mn\in \Nyx $ and that $\Eyx
(mn)=\pyx (mn)$.  Therefore
  $$
  y=\beta _{mn}(x)=\beta _m\big (\beta _n(x)\big ).
  $$
  Setting $z=\beta _n(x)$, we then have that $n\in \Npt zx; $ and $m\in \Npt yz; $, and moreover
  $$
  \ft yzm\ft zxn= \p yz; (m)\p zx; (n) \={Juxtap} \pyx (mn) = \Eyx (mn),
  $$
  and the conclusion will be reached once we show that the sum in the statement has no nonzero terms other than the
one discussed above.  To see this, assume that
  $$
  \ft y{z'}m\ft {z'}xn\neq 0,
  $$
  for some $z'$.  This implies in particular that $\ft {z'}xn\neq 0$ so, as argued above, $z'=\beta _n(x)=z$,
proving the desired uniqueness.

Assuming now that $\ft yx{mn} =0$, let us prove that the right-hand-side also vanishes.  Otherwise, there exists
some $z$ in $X$ such that $\ft yzm\ft zxn\neq 0$, so we again deduce from \ref {Ezero} that $\beta _m(z)=y$, and
$\beta _n(x)=z$.  Consequently $\beta _{mn}(x)=y$, and using \ref {Ezero} once more, we obtain $\ft yx{mn} \neq 0$,
a contradiction.
  \endProof

We will now extend a few other technical facts of classical Fourier analysis to be used in extending \ref
{expectedconvform} from elements of $B_0$ to arbitrary pairs of elements in $B$.

\state Lemma \label ClassTech
  Given $a$ and $b$ in $B$, $x$ in $X$, and given a finite subset $Y\subseteq \Orb x$ one has that
  \izitem
  \zitem (Bessel's inequality)
  $$
  \sum _{y\in Y} \ft xy{b^*}\ft yxb\leq \ft xx{b^*b},
  $$
  \zitem (Cauchy-Schwarz inequality)
  $$
  \Big \Vert \sum _{z\in Y} \ft yza\ft zxb\Big \Vert \leq
  \Big \Vert \sum _{z\in Y} \ft yza\ft yza^*\Big \Vert ^{1/2}
  \Big \Vert \sum _{z\in Y} \ft zxb^* \ft zxb\Big \Vert ^{1/2},
  $$
  for every $y$ in $\Orb x$,

\Proof
  (i)\enspace Assume first that $b$ lies in $B_0$.  We may then use \ref {BabyConvol} to get
  $$
  \ft xx{b^*b} =
  \sum _{y\in \Orb x} \ft xy{b^*}\ft yxb.
  $$
  Summing only for $y$ in $Y$ we get (i).  Since both sides of (i) represent continuous functions, and since $B_0$
is dense in $B$, the general case follows by taking limits.

\itmProof (ii)
  Considering $\Bun $ as a Fell bundle over the (discrete) groupoid $\Xt $, let $C_c(\Bun )$ be the $*$-algebra
formed by all finitely supported sections of $\Bun $, equipped with the operations described in \cite
[3.1]{KumjianFell}.  This admits a faithful representation on a Hilbert space $\Hilb $ by \cite [3.4]{KumjianFell}.
Identifying each $\Byx $ with its canonical image in $C_c(\Bun )$ we may then suppose that the $\Byx $ are
faithfully represented in $\BH $.

Since the $\Isot x, x; $ are C*-algebras, the above representations must necessarily be isometric on each $\Isot x,
x; $, and the same must be the case for all the $\Byx $ due to \cite [2.1.ix]{KumjianFell}.

The result then follows from the familiar inequality
  $$
  \Big \Vert \sum _{i=1}^k a_i b_i\Big \Vert \leq
  \Big \Vert \sum _{i=1}^k a_i a_i^*\Big \Vert ^{1/2}
  \Big \Vert \sum _{i=1}^k b_i^* b_i\Big \Vert ^{1/2},
  $$
  which holds whenever $a_1,\ldots ,a_k,b_1,\ldots b_k$ are operators on $\Hilb $.
  \endProof

The main result of this section is the following:

\state Theorem \label MatrixMult
  Given $a,b\in B$, and given $(y,x)\in \Xt $, one has that
  $$
  \ft yx{ab} = \sum _{z\in \Orb x} \ft yza\ft zxb,
  $$
  the series being unconditionally summable in $\Cyx $.

\Proof Let ${\cal Z}$ be the collection of all finite subsets of $\Orb x$, seen as a partially ordered set under
inclusion.  Evidently $\cal Z$ is a directed set so, given $a$ and $b$ in $B$, we may consider the net
  $\{S_Z(a, b)\}_{Z\in \cal Z}$, defined by
  $$
  S_Z(a, b)= \sum _{z\in Z} \ft yza\ft zxb, \for Z\in \cal Z.
  $$
  To say that the series in the statement is unconditionally summable to $\ft yx{ab}$ is thus to say that
  \lbldeq AdultConv
  $$
  \li Z S_Z(a,b) =
  \ft yx{ab}.
  $$

Observe that \ref {BabyConvol} already gives \ref {AdultConv}, provided $a$ and $b$ lie in $B_0$.  In order to
extend this to an arbitrary pair $(a, b)$, it therefore suffices to observe that $B_0$ is dense in $B$, but first we
need to show that the $S_Z$ form an equicontinuous family of functions.  In order to do so, note that
  $$
  \Vert S_Z(a,b)\Vert =
  \Big \Vert \sum _{z\in Z} \ft yza\ft zxb\Big \Vert \explain {ClassTech} \leq
  $$
  $$
  \leq
  \Big \Vert \sum _{z\in Z} \ft yza\ft yza^*\Big \Vert ^{1/2} \Big \Vert \sum _{z\in Z} \ft zxb^* \ft zxb\Big \Vert
^{1/2}
  \ \={FourAdjoin}
  $$
  $$
  =
  \Big \Vert \sum _{z\in Z} \ft yza\ft zy{a^*}\Big \Vert ^{1/2} \Big \Vert \sum _{z\in Z} \ft xz{b^*} \ft zxb\Big
\Vert ^{1/2}
  \ \explain {ClassTech.i}\leq
  $$
  $$
  \leq
  \Vert \ft yy{aa^*}\Vert ^{1/2} \Vert \ft zz{b^*b}\Vert ^{1/2} \leq
  \Vert a\Vert \Vert b\Vert .
  $$
  The equicontinuity of the $S_Z$ now follows easily.  This allows us to extend the conclusion of \ref {BabyConvol}
for all $a$ and $b$ in $B$, concluding the proof.
  \endProof

\state Corollary \label Parseval
  (Parseval's identity) For every $b$ in $B$ and for every $x\in X$, one has that
  $$
  \ft xx{b^*b} = \sum _{y\in \Orb x} \ft yxb^* \ft yxb
  $$
  the series being unconditionally summable in $\PIA xx; $.

\Proof
  Follows immediately from \ref {MatrixMult} upon letting $y=x$ and $a=b^*$.
  \endProof

\endsection

\startsection {Opaque} and {gray} ideals

\sectiontitle

We shall continue working under \ref {StandingThree}, namely assuming that $\incl {B}{A}$ is a regular inclusion and
that $A$ is abelian.

As already mentioned in the paragraph right after \ref {GetCondExp}, observe that the map $\iota :A\to \BJ $,
introduced in \ref {IncludeAinBJ}, requires extra hypotheses, not yet available there, to be a non-degenerate map.
However, in view of the fact that we are now assuming that $A$ is regular in $B$, and hence that $A$ contains an
approximate unit for $B$, it is clear that $\iota $ is a non-degenerate map.  Consequently $\E xx;$, which we have
already agreed to call $\Ex $, is a generalized conditional expectation in the stronger sense suggested after \ref
{GetCondExp}.

In a sense, our next result shows that, collectively, the $\Ex $ dominate all other generalized conditional
expectations.

\state Proposition \label MotherCondexp
  Let $D$ be a commutative C*-algebra equipped with a $*$-homo\-mor\-phism $\iota :A\to \Mult (D)$.  Also let $\EE
:B\to D$ be a generalized conditional expectation in the sense of \ref {DefineCondExp}.  If\/ $b\in B$ satisfies
$\Ex (b) = 0$ for all $x$ in $X$, then $\EE (b) = 0$.

\Proof Observe that, even though $\iota $ is not assumed to be non-degenerate, we may use \ref {Workaround} to
replace $D$ with a subalgebra $\bar D$, relative to which $\iota $ is non-degenerate.  In other words, we may
suppose without loss of generality that $\iota $ is non-degenerate.

Writing $\hat D$ for the spectrum of $D$, it is well known that $\Mult (D )$ may be naturally identified with the
C*-algebra $C^b(\hat D)$ formed by all bounded continuous functions on $\hat D$.  Moreover the non-degeneracy of
$\iota $ implies that there exists a proper continuous map $\theta :\hat D\to X$, such that
  $$
  \iota (a) = a\circ \theta , \for a\in A.
  $$
  This allows us to extend $\iota $ to the corresponding multiplier algebras, obtaining the map
  $$
  \tilde \iota : a\in \Mult (A)=C^b(X) \ \longmapsto \ a\circ \theta \in \Mult (D)=C^b(\hat D).
  $$

  Arguing by contradiction, suppose that $\EE (b)\neq 0$, so there is some $y$ in $\hat D$ such that $\ev {\EE
(b)}y\neq 0$.  Setting
  $$
  x= \theta (y), \and m = \big |\ev {\EE (b)}y\big |,
  $$
  observe that for all $v\in \Mult (A)$, one has that
  \lbldeq Bigvbv
  $$
  \|vbv\| \geq
  \|\EE (vbv)\| =
  \|\iota (v)\EE (b)\iota (v) \| \geq
  |\ev {\iota (v)\EE (b)\iota (v)}y| \quebra \geq
  |\ev {\iota (v)}y|\ |\ev {\EE (b)}y|\ |\ev {\iota (v)}y| =
  m|\ev {v}x|^2.
  $$
  Choosing an approximate unit $\{u_i \}_i $ for $J_x$, and putting $v_i = 1 - u_i $, we have that $\ev {v_i}x= 1$,
for every $i$, so \ref {Bigvbv} yields
  $$
  \|v_i bv_i\| \geq m.
  $$

On the other hand, since $\Ex (b) = 0$, one has that $b$ lies in $\Lx $, so \ref {viKvi} implies that
  $$
  \li i v_i bv_i = 0,
  $$
  a contradiction.
  \endProof

Let us now study the joint left-kernel of the $\Ex $.

\state Proposition \label IntroBlack
  Let
  $$
  \Delta = \big \{b\in B: \Ex (b^*b) = 0, \hbox { for all } x\in X\big \}.
  $$
  Then
  \izitem
  \zitem
  $
  \Delta = \medcap _{x\in X} BJ_x=
  \medcap _{x\in X} J_xB,
  $
  \zitem $\Delta $ is a closed two-sided ideal of $B$,
  \zitem if there exists a commutative C*-algebra $D$, and a faithful generalized conditional expectation $\EE :B\to
D$, in the sense of \ref {DefineCondExp}, then $\Delta =\{0\}$.

\Proof
  Given $b$ in $B$, and $x$ in $X$, let $\{u_i \}_i $ be an approximate unit for $J_x$ and set $v_i=1-u_i$. Then
  $$
  \Ex (b^*b) = 0 \iff
  b^*b\in \Lx \explain {viKvi}\iff
  \li i v_i b^*bv_i = 0 \quebra \iff
  \li i bv_i = 0 \explain {RightIdeal}\iff
  b\in BJ_x.
  $$
  This proves the first equality in (i) and hence it is clear that $\Delta $ is a left ideal.  In order to prove
that $\Delta $ is also a right ideal, it is enough to show that
  $$
  \Delta n\subseteq BJ_x,
  $$
  for every normalizer $n$, and for every $x$ in $X$.  Assuming that $\ev {n^*n}x= 0$, we have by \ref {NormaJota.i}
that $n\in BJ_x$, so
  $$
  \Delta n\subseteq \Delta BJ_x\subseteq BJ_x.
  $$

Otherwise, if $\ev {n^*n}{x}\neq 0$, let $y=\beta _n(x)$.  Then $J_yn\subseteq BJ_x$, by \ref {NormaJota.ii}, so
  $$
  \Delta n\subseteq BJ_yn\subseteq BBJ_x\subseteq BJ_x.
  $$
  This shows that $\Delta $ is a two-sided ideal.  It moreover follows that $\Delta $ is self-adjoint, from which
one deduces the second equality in (i).

In order to prove (iii), let $\EE $ be a faithful generalized conditional expectation and pick $b$ in $\Delta $.
  Then $\Ex (b^*b)=0$, for every $x$, whence $\EE (x^*x)=0$, by \ref {MotherCondexp}.  Since $\EE $ is faithful, we
deduce that $b=0$, as required.
  \endProof

\state Proposition \label BlackNoFourrier
  Given $b$ in $B$, a necessary and sufficient condition for $b$ to lie in $\Delta $ is that
  $$
  \ft yxb=0, \for (x,y)\in \Xt .
  $$

\Proof Given $b$ in $B$, we have for every $x$ in $X$ that
  $$
  \ft xx{b^*b} \={Parseval} \sum _{y\in \Orb x} \ft yxb^* \ft yxb,
  $$
  from where the conclusion follows immediately.
  \endProof

The above result may be interpreted as saying that $\Delta $ is invisible from the point of view of Fourier
coefficients.

\definition
  We shall say that $\newsymbol {\Delta }{{opaque} ideal} $ is the \newConcept {{opaque} ideal}{} relative to the
inclusion $\incl {B}{A}$.

\state Remark
  \rm
  V. Zarikian has shown us that it is possible for a regular inclusion to have a non-trivial {opaque} ideal, however the
verification is somewhat involved, so we do not supply an example here.  Here is a straightforward example of an
inclusion with a non-trivial {opaque} ideal; note however that this example is not a regular inclusion.  Let $\cal K$ be
the compact operators on $L^2[0,1]$ (Lebesgue measure), $A=C[0,1]$ represented as multiplication operators on
$L^2[0,1]$, and $B=A+\cal K$.  Then the {opaque} ideal for the inclusion $(A,B)$ is $\cal K$ (also see~\ref
{eigenrelate} and \cite [Example~2.8]{DonsigPitts}).

The following result resembles the fact that the diagonal subalgebra in a matrix algebra is maximal abelian.

\state Proposition \label DiagoComuta
  Given $b$ in $B$, the following are equivalent
  \izitem
  \zitem $\ft yxb=0$, whenever $x\neq y$,
  \zitem for all $a$ in $A$, the commutator $[a,b]=ab-ba$ lies in $\Delta $.

\Proof (i)$\Rightarrow $(ii)\enspace
  Given any $x$ and $y$ in $X$, we have by \ref {MatrixMult} that
  $$
  \ft yx{ab-ba} =
  \sum _{z\in \Orb x} \ft yza\ft zxb- \ft yzb\ft zxa\quebra =
  \ft yya\ft yxb- \ft yxb\ft xxa\={BXUnital.i}
  \ev ay\ft yxb- \ev ax\ft yxb.
  $$
  If $x=y$, this vanishes for obvious reasons, and otherwise it does so by hypothesis.  Consequently $ab-ba$ lies in
$\Delta $.

\itmImply (ii) > (i)
  Given that $x\neq y$, choose $a$ in $A$ such that $\ev ay=1$, while $\ev ax=0$.  Then the above computation gives
  $$
  0= \ft yx{ab-ba} =
  \ev ay\ft yxb- \ev ax\ft yxb= \ft yxb,
  $$
  as desired.
  \endProof

There is another relevant ideal of $B$ which we would like to study, but first we need to prove a technical result
regarding the set of $F$ of free points introduced in \ref {DefFreePts}.

\state Lemma \label InvarRel
  $F$ is invariant under the equivalence relation ``$\sim $'' introduced in \ref {OrbEquiv}.

\Proof
  Given that $x\sim y\in F$, choose a normalizer $n$ in $\Nyx $, so that $\beta _n(x)=y$.  Supposing that $\psi _1$
and $\psi _2$ are two state extensions of $\varphi _x$, consider the positive linear functionals $\rho _i$ on $B$
given by
  $$
  \rho _i(b) = \Frac {\psi _i(n^*bn)}{\evx {n^*n}}.
  $$
  For $a$ in $A$ we then have that
  $$
  \rho _i(a) =
  \Frac {\psi _i(n^*an)}{\evx {n^*n}} =
  \Frac {\evx {n^*an}}{\evx {n^*n}} \={Beta}
  \ev a{\beta _n(x)} =
  \ev ay= \varphi _y(a).
  $$
  so $\rho _i$ is a state, and $\rho _i$ extends $\varphi _y$.  Since $y$ is a free point we deduce that $\rho
_1=\rho _2$, which implies that
  $$
  \psi _1(n^*bn) = \psi _2(n^*bn), \for b\in B.
  $$
  Plugging in $nbn^*$ in place of $b$ in either side above, we get for every $i=1,2$, that
  $$
  \psi _i(n^*nbn^*n) \={CondexForState}
  \evx {n^*n}^2\psi _i(b).
  $$
  Taking into account that $\evx {n^*n}\neq 0$, we conclude that $\psi _1=\psi _2$, and hence that $\varphi _x$
admits a unique state extension.  Consequently $x\in F$, and the proof is complete.
  \endProof

In the following definition one should pay attention to the fact that $F$ will take the place previously occupied by
$X$.

\state Proposition \label IntroGray
  Let
  $$
  \Gamma = \big \{b\in B: \Ex (b^*b) = 0, \hbox { for all } x\in F\big \}.
  $$
  Then
  \izitem
  \zitem
  $
  \Gamma = \medcap _{x\in F} BJ_x=
  \medcap _{x\in F} J_xB,
  $
  \zitem $\Gamma $ is a closed two-sided ideal of $B$,
  \zitem a necessary and sufficient condition for a given $b$ in $B$ to lie in $\Gamma $ is that
  $$
  \ft yxb=0, \for x,y\in F.
  $$

\Proof
  The first equality in (i) is proved in the exact same way as the corresponding statement in \ref {IntroBlack.i}.
As before this implies that $\Gamma $ is a left ideal.  To prove that $\Gamma $ is also a right ideal we again
follow the argument given in the proof of \ref {IntroBlack}, observing that if $x$ is in $F$, then the element
$y=\beta _n(x)$ referred to there is also in $F$, thanks to \ref {InvarRel}.  With this, (i) and (ii) are proved as
before.

Finally, the proof of (iii) may be done as in \ref {BlackNoFourrier} by considering only elements of $F$, and
observing that if $x$ is in $F$, then $\Orb x\subseteq F$, by \ref {InvarRel}.
  \endProof

While the {opaque} ideal is invisible from the point of view of \emph {every} Fourier coefficient, the above result
says that $\Gamma $ has that same property, except that only \emph {free} Fourier coefficients do not see it.

\definition
  We shall say that $\newsymbol {\Gamma }{{gray} ideal} $ is the \newConcept {{gray} ideal}{} relative to the
inclusion $\incl {B}{A}$.

  Recalling that $\Pf $ denotes the free expectation introduced in \ref {IntroFreeCondExp}, the following is an
obvious consequence of the definitions:

\state Proposition \label FaithGray
  One has that
  $$
  \Gamma = \{b\in B:\Pf (b^*b)=0\},
  $$
  so $\Pf $ is a faithful expectation if and only if the {gray} ideal of $B$ vanishes.

Much has been said regarding inclusions $\incl {B}{A}$ satisfying the \emph {ideal intersection property}, namely
when every nontrivial ideal (always assumed closed and two-sided) of $B$ has a nontrivial intersection with $A$.
The following connects the {opaque} and {gray} ideals with their intersections with $A$.  When $\incl {B}{A}$ has
the unique pseudo-expectation property (see \cite [Section~1.3]{PZ}), \ref {LargestIdeal} shows that the {gray}
ideal is the left kernel of the pseudo-expectation, see \cite [Theorem~3.15]{SRI} and \cite [Theorem~6.5]{SRITwo}.

\state Proposition
  \izitem
  \lbldzitem BlackNoA The {opaque} ideal\/ $\Delta $ has trivial intersection with $A$.
  \lbldzitem GreiInterA $\Gamma \cap A$ consists of all elements of $A$ vanishing on $F$ (and hence also on $\clsr
{F}$).
  \lbldzitem LargestIdeal If $F$ is dense in $X$, then the {gray} ideal is the largest ideal of $B$ whose
intersection with $A$ is trivial.

\Proof \refl {BlackNoA}\enspace
  If $a$ is in $\Delta \cap A$, and if $\{u_i\}_i$ is an approximate unit for $J_x$, for a given $x$ in $X$, the
fact that $a$ belongs to $BJ_x$ implies that
  $$
  a\={RightIdeal} \li i au_i,
  $$
  so $\evx a=0$, and since $x$ is arbitrary, we deduce that $a=0$.

\itmProof (\LocalGreiInterA )
  Observe that for every $a$ in $A$, and for every $x$ in $F$, one has that
  $$
  \Ex (a^*a) \={ThreeValues} \psi _x(a^*a) = \evx {a^*a} = |\evx {a}|^2.
  $$
  The proof of \ref {GreiInterA} then follows.

\itmProof (\LocalLargestIdeal )
  Under the present hypothesis, \refl {GreiInterA} implies that $\Gamma \cap A=\{0\}$.
  On the other hand, suppose that $K$ is an ideal in $B$ with $K\cap A=\{0\}$.  We may then consider the quotient
inclusion
  $$
  \INCL {\Frac BK}{\Frac A{K\cap A}},
  $$
  which should actually be written as
  $$
  \INCL {\Frac BK}{A},
  $$
  given that the above denominator under $A$ is trivial.  Given any $x$ in $F$, consider the pure state $\varphi _x$
on $A$ given by evaluation at $x$.  Seeing $A$ as a subalgebra of $B/K$, we may then find a state $\rho $ on $B/A$
extending $\varphi _x$.  Writing
  $$
  \pi :B\to B/K,
  $$
  for the quotient map, observe that $\rho \circ \pi $ is a state on $B$ which clearly also extends $\varphi _x$.
Consequently $\rho \circ \pi $ coincides with the unique state extension of $\varphi _x$, which we have been calling
$\psi _x$.

Given $b$ in $K$ we then have
  $$
  \Ex (b^*b) \={ThreeValues} \psi _x(b^*b) = \rho \big (\pi (b^*b)\big ) = 0,
  $$
  since this holds for every $x$ in $F$, we conclude that $b\in \Gamma $, hence that $K\subseteq \Gamma $.
  \endProof

\endsection

\startsection Topologically free inclusions

\sectiontitle

We shall continue working under \ref {StandingThree}.

It is easy to see that the collection $\Norm BA $ of all normalizers of $A$ in $B$ forms a $*$-semigroup under
multiplication.  From now on we will be interested in its subsemigroups.

\definition \label DefineAdmissTotal
  We shall say that a $*$-subsemigroup $N\subseteq \Norm BA $ is
  \izitem
  \zitem \pilar {12pt} \newConcept {admissible}{admissible $*$-semigroup} when \
  $
  \clspan \{n^*na: n\in N,\ a\in A\} = A,
  $
  \zitem \pilar {12pt} \newConcept {generating} {generating $*$-semigroup} when \
  $
  \clspan \{na: n\in N,\ a\in A\} = B.
  $

Notice that the set referred to in the left-hand-side of the equality in (i), above, is precisely the closed ideal
of $A$ generated by $\{n^*n:n\in N\}$.  Regardless of whether or not $N$ is admissible, this ideal is clearly the
set of elements in $A$ which vanish outside the union of the open supports of the $n^*n$, as $n$ range in $N$.
  Thus, we have:

\state Proposition \label AdmissFullDom
  A $*$-subsemigroup $N\subseteq \Norm BA $ is admissible if and only if
  $$
  \medcup _{n\in N}\Dom n= X.
  $$

Regarding generating semigroups, recall that if $n$ is any normalizer, then the closed linear span of $nA$, $An$ and
$AnA$ coincide with each other by \ref {SingleAndBiModules}.  Therefore, for any subset $N\subseteq \Norm BA $, one
has that
  $$
  \clspan \NoA = \clspan AN= \clspan A\NoA .
  $$

For future reference we highlight the following:

\state Proposition \label BimodGen
  A $*$-subsemigroup $N\subseteq \Norm BA $ is generating if and only if
  $$
  B= \clspan A\NoA .
  $$

\state Proposition
  Every generating semigroup is admissible.

\Proof
  Let $N\subseteq \Norm BA $ be a generating semigroup and suppose, by way of contradiction, that
  $$
  J:=\clspan \{n^*na: n\in N,a\in A\}\neq A.
  $$

  Observing that $J$ is then a proper ideal in $A$, it must be contained in some maximal ideal, meaning that
$J\subseteq J_x$, for some $x$ in $X$.  So
  $$
  \evx {n^*n}=0, \for n\in N.
  $$

  Letting $\psi _x$ be any state on $B$ extending $\varphi _x$, Cauchy-Schwarz gives for all $a$ in $A$ that
  $$
  |\psi _x(an)|^2\leq |\psi _x(aa^*)|\ |\psi _x(n^*n)| =
  |\evx {aa^*}|\ |\evx {n^*n}| =0,
  $$
  This says that $AN$ is contained in the null space of $\psi _x$, contradicting the fact that $N$ is generating.
  \endProof

As the title of this section suggests, we will now introduce several variations of the notion of \emph {topological
freeness} for regular inclusions relevant to this work.

\definition \label ManyTopFree
  Let $\incl {B}{A}$ be a regular inclusion with $A$ abelian (as in \ref {StandingThree}). We shall say that $\incl
{B}{A}$ is a:
  \izitem
  \zitem \newConcept {topologically free}{topologically free inclusion} inclusion if the set $F$ of free points is
dense in the spectrum of $A$,
  \zitem \newConcept {smooth}{smooth inclusion} inclusion if there exists a generating $*$-semigroup consisting of
smooth normalizers,
  \zitem \newConcept {{light}}{{light} inclusion} inclusion if the {gray} ideal of $B$ vanishes,
  \zitem \newConcept {weak Cartan}{weak Cartan inclusion} inclusion if it is both smooth and {light}.

Notice that \ref {FaithGray} says that the {gray} ideal is precisely the left kernel of $\Pf $, so, being {light} is
equivalent to the free expectation $\Pf $ introduced in \ref {DefineFreeCondExp} being faithful.

\medskip

\state Remark \label {nonsepCartan}
  \rm
  In \cite {RenaultCartan}, Renault defines a Cartan inclusion as a regular inclusion $(A,B)$ with $B$ separable, $A$ a
MASA in $B$ which contains an approximate unit for $B$, and for which there exists a faithful conditional expectation of
$B$ onto $A$.  As in \cite {SRI}, the separability condition on $B$ is sometimes dropped from Renault's definition; such
inclusions are also called Cartan inclusions.  We caution the reader that in the absence of separability, a Cartan
inclusion need not be a weak Cartan inclusion in the sense of~\ref {ManyTopFree}.  For example, let $X$ be the closed
unit disk in the plane, let $\Gamma $ be the (discrete) group of all linear fractional transformations which carry $X$
to itself.  Then $(A,B) =\incl {C({X})\rtimes _\red \Gamma }{C({X})}$ satisfies the conditions of Renault's definition
of a Cartan inclusion except for the separability hypothesis.  However, because its set of free points is empty, it is
far from being topologically free, so it cannot be a {light} inclusion (by~\ref {RelTFs} below).  Thus this inclusion is
not a weak Cartan inclusion.  Despite this, it is an inductive limit of weak Cartan inclusions.

Besides the built-in implications ``\emph {weak Cartan}'' $\Rightarrow $ ``\emph {smooth}'', and ``\emph {weak
Cartan}'' $\Rightarrow $ ``\emph {{light}}'', the following explores the relationship between the above notions:

\state Proposition \label RelTFs
  Let $\incl {B}{A}$ be a regular inclusion with $A$ abelian. Then
  \izitem
  \zitem $\incl {B}{A}$ is a {light} inclusion if and only if it is topologically free and satisfies the ideal
intersection property,
  \zitem if $\incl {B}{A}$ is smooth and $B$ is separable, then $\incl {B}{A}$ is topologically free.

\Proof (i)\enspace
  Assume that $\incl {B}{A}$ is a {light} inclusion.  In order to prove that $F$ is dense in $X$ it is enough to
show that every $a$ in $A$ vanishing on $F$, must itself vanish.  Given that $a$ vanishes on $F$, we have that $a$
lies in $\Gamma $ by \ref {GreiInterA}, so $a=0$ by hypothesis.  This proves that $F$ is dense in $X$, so $\incl
{B}{A}$ is topologically free.

Once we know that $F$ is dense, the ideal intersection property follows from \ref {LargestIdeal}, while the converse
implication is also a consequence of \ref {LargestIdeal}.

\itmProof (ii)
  Let $N$ be a generating semigroup of smooth normalizers.  Since $B$ is separable, then so is $N$, hence there
exists a countable dense subset $N_0\subseteq N$.  By \ref {BimodGen} we have that $A\NoA $ spans a dense subspace
of $B$, whence so does $AN_0A$.  By \ref {XfAndXb} we then have that
  $$
  F= \medcap _{n\in N_0} F_n,
  $$
  and since the interior of each $F_n$ is dense in $X$, we obtain the conclusion by Baire's Theorem.
  \endProof

Regarding \ref {RelTFs.i} it is perhaps worth observing that the ideal intersection property alone does not imply
that an inclusion is {light}, a counterexample being given by $\Incl {B}{{\bf C}}$, where $B$ is simple and unital.
However we refer the reader to \ref {EssentialCommutTwo}, below, for a case in which these two properties are
equivalent.

For future reference let us summarize some of of our previous results using the terminology introduced above.

\state Proposition \label SummarySimplicity
  Let $\incl {B}{A}$ be a {light} inclusion. Then
  \izitem
  \zitem every nonzero ideal of $B$ has a nonzero intersection with $A$,
  \zitem $B$ is simple if and only if $A$ is $B$-simple.

\Proof
  The first point is just a rephrasing of the ideal intersection property which holds for {light} inclusions by \ref
{RelTFs}.  The second point then follows from
 \ref {Simplicity}.
  \endProof

Let us now study morphisms between inclusions.

\state Proposition \label TwoInclusions
  Let $\incl {B_1}{A}$ and $\incl {B_2}{A}$ be two inclusions of the same abelian C*-algebra $A$, with spectrum
denoted by $X$, and suppose that $\mu :B_1\to B_2$ is a surjective $*$-homomorphism whose restriction to $A$ is the
identity.  Then
  \izitem
  \zitem $\mu \big (\Norm {B_1}A \big )\subseteq \Norm {B_2}A $,
  \zitem if $A$ is a non-degenerate subalgebra of $B_1$, then $A$ is also a non-degenerate subalgebra of $B_2$,
  \lbldzitem RegToReg if $\incl {B_1}{A}$ is a regular inclusion, then so is $\incl {B_2}{A}$.
  \medskip \noindent Moreover, under the conditions of (iii), one also has that
  \iaitem
  \lbldaitem UnitoUni if $b\in B_1 $, then $F_b\subseteq F_{\mu (b)}$,
  \lbldaitem ConttoCont if $b\in B_1 $ is continuous, then so is $\mu (b)$,
  \lbldaitem SmothToSmooth if $b\in B_1 $ is smooth, then so is $\mu (b)$,
  \lbldaitem FreeToFree if $x$ is a $B_1$-free point of $X$, meaning free relative to $B_1$, then $x$ is also
$B_2$-free,
  \lbldaitem ProjNaMosca if $x$ is $B_1$-free, then
  $$
  \Ex ^1=\Ex ^2\circ \mu ,
  $$
  where
  $\Ex ^i$
  denotes the localizing projection of
  $\incl {B_i}{A}$ at the ideal $J_x$, for every $i=1,2$,
  \lbldaitem IdealQuase denoting the {gray} ideal of $B_i$ by $\Gamma _i$, for every $i=1,2$, one has that $\mu \inv
(\Gamma _2)\subseteq \Gamma _1$,
  \lbldaitem GreyKer $\Ker (\mu )\subseteq \Gamma _1$,
  \lbldaitem TopoFreetoTopoFree if $\incl {B_1}{A}$ is topologically free, then so is $\incl {B_2}{A}$,
    \global \advance \aitemno by 1
  \lbldaitem IItoII if $\incl {B_1}{A}$ is smooth, then so is $\incl {B_2}{A}$,
  \lbldaitem NoGreyIso if $\incl {B_1}{A}$ is {light}, then $\mu $ is an isomorphism,
  \lbldaitem IdealNaMosca if $\incl {B_1}{A}$ is topologically free, then $\mu \inv (\Gamma _2)=\Gamma _1$.

\Proof We leave the first three easy points to the reader.

  \itmProof (\LocalUnitoUni )
  Assuming that $x$ is not in $F_{\mu (b)}$, we have that $\varphi _x$ admits two state extensions to $B_2$, say
$\psi $ and $\psi '$, such that $\psi (\mu (b))\neq \psi '(\mu (b))$.  It then obvious that $\psi \circ \mu $ and
$\psi '\circ \mu $ are states on $B_1$ extending $\varphi _x$, taking on distinct values on $b$, so $x$ is not in
$F_b$.  This proves \refl {UnitoUni}.

\itmProof (\LocalConttoCont \&\LocalSmothToSmooth )
  Follow immediately from \refl {UnitoUni}.

\itmProof (\LocalFreeToFree )
  Follows immediately from \refl {UnitoUni} and \ref {XfAndXb}.

\itmProof (\LocalProjNaMosca )
  Assuming that $x$ is $B_1$-free, and hence also $B_2$-free by \refl {FreeToFree}, let $\psi _x$ be the unique
state extension of $\varphi _x$ to $B_2$.  Since $\psi _x\circ \mu $ is a state on $B_1$ extending $\varphi _x$, we
deduce that $\psi _x\circ \mu $ is the unique state on $B_1$ extending $\varphi _x$.  Consequently, for every $b$ in
$B_1$, we have
  $$
  \Ex ^2\big (\mu (b)\big ) \={ThreeValues}
  \psi _x\big (\mu (b)\big ) \={ThreeValues}
  \Ex ^1(b) .
  $$

\itmProof (\LocalIdealQuase )
  Given $b$ in $\mu \inv (\Gamma _2)$, we have for every $B_1$-free point $x$, that
  $$
  \Ex ^1(b^*b) \={LocalProjNaMosca}
  \Ex ^2\big (\mu (b^*b)\big ) =
  \Ex ^2\big (\mu (b)^*\mu (b)\big ) =
  0,
  $$
  so $b$ is in $\Gamma _1$.

\itmProof (\LocalGreyKer ) We have
  $$
  \Ker (\mu ) = \mu \inv (\{0\}) \subseteq \mu \inv (\Gamma _2) \explain {LocalIdealQuase}\subseteq \Gamma _1.
  $$

\itmProof (\LocalTopoFreetoTopoFree )
  Follows immediately from \refl {FreeToFree}.

\itmProof (\LocalIItoII )
  Follows immediately from \refl {SmothToSmooth}.

\itmProof (\LocalNoGreyIso )
  Follows immediately from \refl {GreyKer} and the fact that $\Gamma _1$ vanishes.

\itmProof (\LocalIdealNaMosca )
  We first claim that $\mu (\Gamma _1)\cap A=\{0\}$.  To see this let us pick $b$ in $\Gamma _1$ such that $a:=\mu
(b)\in A$.  It follows that $\mu (a-b)=0$, so
  $$
  a-b\in \Ker (\mu ) \explain {LocalGreyKer}\subseteq \Gamma _1,
  $$
  from where we see that $a$ lies in $\Gamma _1$, whence $a=0$ by \ref {LargestIdeal}.  This shows that indeed $\mu
(\Gamma _1)\cap A=\{0\}$ and, using \ref {LargestIdeal} once more, we conclude that $\mu (\Gamma _1)\subseteq \Gamma
_2$.  So
  $$
  \Gamma _1 \={LocalGreyKer} \Gamma _1 + \Ker (\mu ) = \mu \inv \big (\mu (\Gamma _1)\big ) \subseteq \mu \inv
(\Gamma _2).
  $$
  Combined with \refl {IdealQuase}, this concludes the proof.
  \endProof

Let us quickly discuss invariant ideals in the above context.

\state Proposition \label FunctorBSimple
  Let $\incl {B_1}{A}$ and $\incl {B_2}{A}$ be two inclusions of the same abelian C*-algebra $A$, and suppose that
$\mu :B_1\to B_2$ is a surjective $*$-homomorphism whose restriction to $A$ is the identity.  Let moreover $J$ be a
closed ideal in $A$.  Then:
  \iaitem
  \aitem if $J$ is $B_1$-invariant (that is, invariant relative to $B_1$, according to \ref {DefineInvar}), then $J$
is also $B_2$-invariant,
  \aitem if $\incl {B_1}{A}$ is regular, and $J$ is $B_2$-invariant, then $J$ is also $B_1$-invariant,
  \aitem if $A$ is $B_2$-simple, then it is also $B_1$-simple,
  \aitem if $\incl {B_1}{A}$ is regular, and $A$ is $B_1$-simple, then it is also $B_2$-simple.

\Proof (a)\enspace
  Applying $\mu $ to the identity $B_1J=JB_1$ gives the conclusion.

\itmProof (b) As already observed in the proof of \ref {InterInvar}, in order to prove that $B_1J\subseteq JB_1$, if
suffices to check that
  $$
  na\in JB_1,
  $$
  for every normalizer $n\in \Norm {B_1}A $, and for every $a$ in $J$.  Setting $m=\mu (n)$, we have that $m$ lies
in $\Norm {B_2}A $ by \ref {TwoInclusions.i}, and since $J$ is assumed to be $B_2$-invariant, we have for every $a$
in $J$, that
  $$
  ma\in B_2J= JB_2,
  $$
  so
  $$
  mam^* \in (JB_2)\cap A\={GeraIntersec}J.
  $$

  Since $\mu $ is the identity on $A$, and since $\mu (nan^*)=mam^*$, we conclude that $nan^*\in J$, whence,
employing the polynomials $r_i$ of \ref {Polys}, we have
  $$
  na= \li i nn^*nr_i(n^*n)a= \li i nan^*nr_i(n^*n) \in JB_1.
  $$
  This proves that $B_1J\subseteq JB_1$, and the reverse inclusion follows similarly, so $J$ is $B_1$-invariant.

\itmProof (c) Follows immediately from (a).

\itmProof (d) Follows immediately from (b).
  \endProof

Returning to \ref {StandingThree}, and assuming that $\incl {B}{A}$ is a topologically free inclusion, we have by
\ref {LargestIdeal} that $\Gamma \cap A$ is trivial, so the quotient map
  $$
  \pi :B\to B/\Gamma
  $$
  restricts to a one-to-one map on $A$, and then $A$ may be identified with the subalgebra $\pi (A)\subseteq
B/\Gamma $.  Setting
  $$
  \bar B= B/\Gamma ,
  $$
  we may then consider the inclusion $\incl {\bar B}{A}$, and therefore \ref {TwoInclusions} applies to the
inclusions $\incl {B}{A}$ and $\incl {\bar B}{A}$, together with the quotient map $\pi $.  Spelling out some of the
consequences, we arrive at an important result:

\state Theorem \label GetUltraFree
  Suppose that $\incl {B}{A}$ is a topologically free inclusion and let $\bar B$ be the quotient of $B$ by its
{gray} ideal.  Then $\incl {\bar B}{A}$ is {light}.

\Proof
  Stressing that the definition of topologically free inclusions requires regularity, we have that $\incl {\bar
B}{A}$ is regular by \ref {RegToReg}.
  That $\incl {\bar B}{A}$ is topologically free then follows from
  \ref {TopoFreetoTopoFree}.

Denoting by $\Gamma $ and $\bar \Gamma $ the {gray} ideals of $B$ and $\bar B$, respectively, we have that
  $$
  \bar \Gamma =
  \pi \big (\pi \inv (\bar \Gamma )\big ) \explain {IdealQuase}\subseteq \pi (\Gamma ) = \{0\},
  $$
  so the {gray} ideal of $\bar B$ vanishes, concluding the proof.
  \endProof

The next result describes what happens when we mod out the {gray} ideals of the two inclusions at the same time.

\state Proposition \label TwoQuotients
  Let $\incl {B_1}{A}$ be a topologically free inclusion, with $A$ abelian.  Also let $\incl {B_2}{A}$ be another
inclusion of the same C*-algebra $A$, equipped with a surjective $*$-homomorphism $\mu :B_1\to B_2$ whose
restriction to $A$ is the identity. Then
  \izitem
  \zitem $\incl {B_2}{A}$ is also a topologically free inclusion,
  \zitem letting $\bar B_1$ and $\bar B_2$ be the quotients of $B_1$ and $B_2$ by their respective {gray} ideals
$\Gamma _1$ and $\Gamma _2$, with quotient maps denoted by $\pi _1$ and $\pi _2$, one has that $\bar B_1$ and $\bar
B_2$ are isomorphic, via an isomorphism $\bar \mu $, such that the diagram
  \medskip
  \begingroup \noindent \hfill \beginpicture
  \setcoordinatesystem units <0.025truecm, -0.02truecm>
  \setplotarea x from -100 to 200, y from -30 to 125
  \put {\null } at -100 -30
  \put {\null } at -100 125
  \put {\null } at 200 -30
  \put {\null } at 200 125
  \put {$B_1$} at 0 0
  \put {$B_2$} at 100 0
  \arrow <0.15cm> [0.25, 0.75] from 20 0 to 80 0
  \put {$\mu $} at 50 -12
  \put {$\bar B_1$} at 0 100
  \put {$\bar B_2$} at 100 100
  \arrow <0.15cm> [0.25, 0.75] from 20 100 to 80 100
  \put {$\bar \mu $} at 50 118
  \arrow <0.15cm> [0.25, 0.75] from 0 20 to 0 80
  \put {$\pi _1$} at -16 50
  \arrow <0.15cm> [0.25, 0.75] from 100 20 to 100 80
  \put {$\pi _2$} at 116 50
  \put {$\simeq $} at 50 90
  \endpicture \hfill \null \endgroup
  \smallskip commutes.

\Proof
  The first point follows from \refl {RegToReg} and \refl {TopoFreetoTopoFree} in \ref {TwoInclusions}.  Regarding
(ii), consider the composition
  $$
  \tau : B_1 \labelarrow {\mu } B_2 \labelarrow {\pi _2} \bar B_2,
  $$
  and observe that
  $$
  \Ker (\tau ) =
  \Ker (\pi _2\circ \mu ) =
  \mu \inv \big (\Ker (\pi _2)\big ) =
  \mu \inv (\Gamma _2) \={IdealNaMosca}
  \Gamma _1,
  $$
  so $\tau $ factors through the quotient $\bar B_1=B_1/\Gamma _1$, providing the required map.
  \endProof

Let us conclude this section by briefly discussing inclusions in the commutative setting.

\fix For the remainder of this section we will therefore assume that both algebras in the inclusion referred to in
\ref {StandingThree} are commutative, and also unital for simplicity.  We may therefore assume that $A=C(X)$ and
$B=C(Y)$, where $X$ and $Y$ are compact Hausdorff spaces, and $A$ is seen as a subalgebra of $B$ by means of an
injective $*$-homomorphism of the form
  $$
  f\in C(X) \mapsto f\circ \pi \in C(Y),
  $$
  where $\pi :Y\to X$ is a continuous surjection.

\definition \label DefineEssSurj
  A continuous surjection
  $$
  \pi :Y\to X
  $$
  is said to be \newConcept {essential}{essential surjection} if there is no proper closed subset $Y'\subseteq Y$,
such that $\pi (Y')=Y$.

The following simple observation, whose proof we leave to the reader, connects this concept to the ideal
intersection property.

\state Proposition \label IIPCommut
  The inclusion
  $\incl {C(Y)}{C(X)}$ satisfies the ideal intersection property if and only if $\pi $ is an essential surjection.

The next result clarifies the meaning of free points in the commutative case.

\state Proposition \label FreeCommut
  A point $x$ in $X$ is free if and only if $\pi \inv (\{x\})$ is a singleton.  In this case,
  $$
  \Ex (b) = b(y), \for b\in C(Y),
  $$
  where $y$ is the unique element of $\pi \inv (\{x\})$.

\Proof If $y$ is any element in $\pi \inv (\{x\})$, it is clear that the evaluation state $\varphi _x$ on $C(X)$ is
extended by the state $\psi _y$ on $C(Y)$ given by evaluation at $y$.  Therefore, if $x$ is a free point of $X$, one
necessarily has that $\pi \inv (\{x\})$ is a singleton.

Conversely, due to the fact that $C(Y)$ is commutative, one has that
  $$
  \Lx = J_xC(Y) + C(Y)J_x= J_xC(Y),
  $$
  which is clearly a two-sided ideal in $C(Y)$, and one may prove that it consists precisely of the continuous
functions on $Y$ vanishing on $\pi \inv (\{x\})$.  Under the hypothesis that the latter is a singleton, say
  $$
  \pi \inv (\{x\})=\{y\},
  $$
  we deduce that $\Lx $ is a maximal ideal, namely the ideal of $C(Y)$ formed by the functions vanishing on $y$.

If $\psi $ is any state on $C(Y)$ extending $\varphi _x$, then $\psi $ vanishes on $\Lx $ by \ref {VanishJK}, so
$\psi $ must be the character of $C(Y)$ associated to $y$.  Therefore $\psi $ is the unique state extension of
$\varphi _x$.  This proves that $x$ is free, and the last sentence of the statement follows easily.
  \endProof

As observed in the paragraph immediately after the proof of \ref {RelTFs}, the ideal intersection property alone
does not imply topological freeness.  When $A=C(\bf T)$ and $B$ is the injective envelope of $A$ (see \ref
{injenvdes}), \cite [Theorem~3.2]{Gleason} together with \ref {IIPCommut} show $\incl B A$ has the ideal
intersection property.  Furthermore, for each $z\in \bf T$, the algebra
  $$
  B_z=\{f\in C({\bf T}\setminus \{z\}): \lim _{\delta \rightarrow 0^+} f(z e^{i\delta }) \text { and } \lim _{\delta
\rightarrow 0^-} f(z e^{i\delta }) \text { both exist}\}
  $$
  satisfies $A\subseteq B_z\subseteq B$.  As $z$ is not a free point for $\incl {B_z} A$, it is not a free point for
$\incl B A$ either; it follows that the set of free points for $\incl {B}{A}$ is empty.  This example shows
topological freeness does not follow from the ideal intersection property even for commutative inclusions.  However,
the following result shows that for most settings arising in practice, topological freeness follows from the ideal
intersection property.

\def \csub {\frak C}

\state Theorem \label EssMetrTopFree
  Suppose $\incl {C(Y)}{C(X)}$ satisfies the ideal intersection property.  The following statements hold.
  \iaitem
  \aitem For every $v\in C(Y)$, the free points relative to $v$, $F_v$, is a dense $G_\delta $ subset of $X$.
  \aitem Suppose {\it $C(Y)$ is countably generated over $C(X)$}, that is, there exists a countable set $\csub
\subseteq C(Y)$ such that $C(Y)=C^*(\csub \cup C(X))$.  Then the set of free points for $\incl {C(Y)}{C(X)}$
contains a dense $G_\delta $ set; in particular, $\incl {C(Y)}{C(X)}$ is a topologically free inclusion.
  \aitem When $Y$ is metrizable, $\incl {C(Y)}{C(X)}$ is a topologically free inclusion.

\Proof For each $x\in X$, let $\Delta (x)$ be the diameter of the set $\{v(y): \pi (y)=x\}\subseteq \bf C$ and
observe that $F_v=\{x\in X: \Delta (x)=0\}$.  For $\varepsilon >0$, let $W_\varepsilon =\{x\in X: \Delta (x)
<\varepsilon \}$.  We shall show that $W_\varepsilon $ is open and dense in $X$.

Suppose by contradiction that $W_\varepsilon $ contains a point $x$ which is not interior.  Then there is a net
$\{x_i\}_i$ converging to $x$, such that $\Delta (x_i)\geq \varepsilon $, for all $i$.  Choosing any $t$ such that
$\Delta (x)<t<\varepsilon $, we may find $y_i$ and $z_i$ in $Y$, such that $\pi (y_i)=\pi (z_i)=x_i$, and
$|v(y_i)-v(z_i)| > t$, for each $i$.

Since $Y$ is compact, and upon passing to subnets, we may assume that $y_i\to y$, and $z_i\to z$, for some $y$ and
$z$ in $Y$.  We then have that $\pi (y) = \lim _i\pi (y_i) = \lim _i x_i = x$, and likewise $\pi (z) = x$, so both
$y$ and $z$ lie in $\pi \inv (x)$.  Now this is a contradiction because
  $$
  \Delta (x) \geq |v(y)-v(z)| = \lim _i |v(y_i)-v(z_i)| \geq t > \Delta (x).
  $$

This shows that $W_\varepsilon $ is open.  To see that $W_\varepsilon $ is dense, suppose again by contradiction
that there is a nontrivial open set $U\subseteq X$, disjoint from $W_\varepsilon $.  Given any $x_0$ in $U$, pick
$y_0$ such that $\pi (y_0)=x_0$, and consider
  $$
  V = \big \{y\in Y: \pi (y)\in U, \ |v(y)-v(y_0)|<\varepsilon /3\big \}.
  $$
  Clearly $V$ is an open neighborhood of $y_0$ in $Y$, so $Z:= Y\setminus V$ is closed and proper.  We then claim
that $\pi (Z)=X$.  To see this, pick $x$ in $X$, and let us prove that $x$ lies in $\pi (Z)$.  Excluding the trivial
case in which $x\notin U$, let us assume that $x$ lies in $U$.  Therefore $\Delta (x)\geq \varepsilon $, so there
are $y$ and $z$ in $Y$, such that $\pi (y)=\pi (z)=x$, and $|v(y)-v(z)|>2\varepsilon /3$.  Notice that $y$ and $z$
cannot simultaneously lie in $V$, or else
  $$
  |v(y)-v(z)|< |v(y)-v(y_0)| + |v(y_0)-v(z)| < 2\varepsilon /3.
  $$
  Assuming without loss of generality that $y$ is not in $V$, we then have that $x = \pi (y) \in \pi (Z)$.
  This contradicts the fact that $\pi $ is essential, hence proving that $W_\varepsilon $ is dense.

Since $F_v=\bigcap _{n=1}^\infty W_{1/n}$ is dense by Baire's theorem, $F_v$ is a dense $G_\delta $ set.

\itmProof (b) Let $\csub \subseteq C(Y)$ be a countable set such that the C*-algebra generated by $C(X)\cup \csub $
is $C(Y)$.  Let $H$ be the set of all finite products of elements from $\csub \cup \csub ^*$.  Then $H$ is
countable, so (a) implies that the set $X_H=\bigcap _{v\in H} F_v$ is a dense $G_\delta $ set in $X$.  If $g\in
C(Y)$ and $\varepsilon >0$, we may find $n\in \bf N$, $\{f_k\}_{k=1}^n\subseteq C(X)$ and $\{v_k\}_{k=1}^n\subseteq
H$, so that
  $$
  \| g-g'\| <\varepsilon , \dstext {where} g'=\sum _{k=1}^n v_k f_k.
  $$
  Then $X_H\subseteq F_{g'}$, from which it follows whenever $x\in X_H$ and $y_1, y_2\in \pi ^{-1}(x)$,
$|g(y_1)-g(y_2)|< 2\varepsilon $.  As the choice of $\varepsilon >0$ is arbitrary, we conclude that $g(y_1)=g(y_2)$.
Therefore $X_H\subseteq F_g$ for each $g\in C(Y)$, so that $X_H\subseteq F$.  Thus, (a) holds.

\itmProof (c) Since $Y$ is metrizable, $C(Y)$ is separable, so (c) follows from (b).
  \endProof

As a consequence we have:

\state Corollary \label EssentialCommut
  If\/ $\incl {C(Y)}{C(X)}$ is a {light} inclusion, then $\pi $ is an essential surjection.  In case $C(Y)$ is
countably generated over $C(X)$, the converse also holds.

\Proof
  If $\incl {C(Y)}{C(X)}$ is a {light} inclusion, then it satisfies the ideal intersection property by \ref
{RelTFs.i}, so the first statement follows at once from \ref {IIPCommut}.  Regarding the second sentence in the
statement, the hypothesis implies that $\incl {C(Y)}{C(X)}$ satisfies the ideal intersection property by \ref
{IIPCommut}, and is topologically free by \ref {EssMetrTopFree}.  So the conclusion follows from \ref {RelTFs.i}.
  \endProof

A slightly different way to put the above is:

\state Corollary \label EssentialCommutTwo
  If\/ $C(Y)$ is countably generated over $C(X)$, then
  $\incl {C(Y)}{C(X)}$ is a {light} inclusion if and only if it satisfies the ideal intersection property.

\Proof Follows from \ref {EssentialCommut} and \ref {IIPCommut}.
  \endProof

\def \calM {{\cal M}}

We now consider smoothness of an inclusion $\incl {C(Y)}{C(X)}$ using the slight change of perspective from
considering a single smooth element to the set of all elements of $C(Y)$ smooth relative to a given open dense set.
Let $\Lambda $ be the collection of all open dense subsets of $X$ and for $\lambda \in \Lambda $, let
  $$
  \calM _\lambda =\{b\in C(Y): \lambda \subseteq F_b\}.
  $$
  Routine arguments show:
  \izitem
  \zitem $\calM _\lambda $ is a $C^*$-subalgebra of $C(Y)$ with $C(X)\subseteq \calM _\lambda \subseteq C(Y)$;
  \zitem every element of $\calM _\lambda $ is smooth for $\incl {C(Y)}{C(X)}$, hence
 $\incl {\calM _\lambda }{C(X)}$ is a smooth inclusion;
  \zitem if $\lambda _1, \lambda _2\in \Lambda $ and $\lambda _2\subseteq \lambda _1$, then $\calM _{\lambda
_1}\subseteq \calM _{\lambda _2}$; and
  \zitem $\bigcup _{\lambda \in \Lambda } \calM _\lambda $ is a self-adjoint subalgebra of $C(Y)$ which contains every
smooth element
for $\incl {C(Y)}{C(X)}$.

\state Proposition \label {Mlambda}
  Let $\calM =\clsr {\bigcup _{\lambda \in \Lambda } \calM _\lambda }$.  Then $\calM $ is the largest
$C^*$-subalgebra of $C(Y)$ satisfying $C(X)\subseteq \calM \subseteq C(Y)$ such that $\incl {\calM }{C(X)}$ is
smooth.  In particular, $\incl {C(Y)}{C(X)}$ is smooth if and only if $\calM =C(Y)$.

\Proof For $i\in \{1,2\}$, let $\lambda _i\in \Lambda $ and suppose $n_i\in \Norm {\calM _{\lambda _i}}{C(X)}$. Then
$n_i$ are smooth, and as $n_1n_2\in \calM _{\lambda _1\cap \lambda _2}$, we see $n_1n_2$ is also smooth.  It follows
that $\ds \bigcup _{\lambda } \Norm {\calM _\lambda }{C(X)}$ is a $*$-semigroup of smooth normalizers.  Since $\Norm
{\calM _{\lambda }}{C(X)}$ contains the unitary group of $\calM _\lambda $, it follows that the linear span of $\ds
\bigcup _{\lambda \in \Lambda } \Norm {M_\lambda }{C(X)}$ is dense in $\calM $.  Therefore, $\incl \calM {C(X)}$ is
smooth.

Suppose $B$ is a $C^*$-subalgebra of $C(Y)$ containing $C(X)$ and $\incl {B}{C(X)}$ is smooth.  Then for a smooth
normalizer $n\in \Norm {B}{C(X)}$, there exists $\lambda \in \Lambda $ so that $\lambda \subseteq F_n$.  Then $n\in
\calM _\lambda \subseteq \calM $.  Therefore, $B\subseteq \calM $.
  \endProof

In the presence of the ideal intersection property, more can be said as the next two results show.

\state Proposition \label {proj=smooth}
  Suppose $\incl {C(Y)}{C(X)}$ has the ideal intersection property.  Then the following statements hold.
  \izitem
  \zitem If $p\in C(Y)$ is a projection, then $p$ is smooth.
  \zitem If the linear span of the projections in $C(Y)$ is dense in $C(Y)$ (equivalently $Y$ is totally
disconnected), then $\incl {C(Y)}{C(X)}$ is smooth.

\Proof The projection $p$ is the characteristic function of a clopen set $G\subseteq Y$.  Setting $H=Y\setminus G$,
we then clearly have that $F_p= X\setminus \big (\pi (G)\cap \pi (H)\big )$, so it suffices to prove that $\pi
(G)\cap \pi (H)$ has empty interior.  Assuming by contradiction that there is a nonempty open set $U\subseteq \pi
(G)\cap \pi (H)$, observe that $\pi $ maps $G\cap \pi \inv (U)$ onto $U$.  Consequently $\pi $ maps the proper
closed subset $G\cup \big (H\setminus \pi \inv (U)\big ) $ onto $\pi (G)\cup (\pi (H)\setminus U)=X$, contradicting
the fact that $\pi $ is essential.

For the second statement, $\{I-2p: p\in \text {proj}(C(Y))\}$ is a group of smooth unitaries whose span contains
$\text {proj}(C(Y))$.  It follows that $\incl {C(Y)}{C(X)}$ is a smooth inclusion.
  \endProof

Finally we show that any inclusion with the ideal intersection property is an intermediate inclusion to a smooth
inclusion with the ideal intersection property.

\state Theorem \label {smoothext}
  Suppose $\incl {C(Y)}{C(X)}$ has the ideal intersection property.  Then there is a compact Hausdorff space $Z$
such that $C(X)\subseteq C(Y)\subseteq C(Z)$ and $\incl {C(Z)}{C(X)}$ is both smooth and has the ideal intersection
property.  If in addition, $Y$ is metrizable, then $Z$ can be chosen metrizable as well; when such a choice is made,
$\incl {C(Z)}{C(X)}$ is a weak Cartan inclusion.

\Proof By \cite [Theorem~3.2]{Gleason}, there is a extremally disconnected compact Hausdorff space $Z$ and a
continuous, essential, surjection $\alpha : Z\rightarrow Y$.  The map $C(Y)\ni b\mapsto b\circ \alpha $ allows us to
regard $C(Y)\subseteq C(Z)$.  As the composition $\pi \circ \alpha $ is an essential surjection of $Z$ onto $X$,
\ref {IIPCommut} shows that $\incl {C(Z)}{C(X)}$ has the ideal intersection property.  Since $Z$ is extremally
disconnected, it is totally disconnected, so \ref {proj=smooth} shows $\incl {C(Z)}{C(X)}$ is smooth.

Suppose now that $Y$ is metrizable and let $P$ be the Cantor set.  By the Hausdorff-Alexandroff Theorem, there
exists a continuous surjection $\phi : P\rightarrow Y$.  A Zorn's Lemma argument shows there exists a minimal closed
set $Z\subseteq P$ such that $\phi (Z)=Y$ (there may be many such $Z$).  Then $\phi |_Z$ is an essential surjection
of $Z$ onto $Y$.  Since $P$ is totally disconnected, $Z$ is a totally disconnected metric space.  As in the previous
paragraph, we obtain $C(Y)\subseteq C(Z)$ and the inclusion $\incl {C(Z)}{C(X)}$ is both smooth and has the ideal
intersection property.  That $C(Z)$ is separable is clear.  By \ref {RelTFs}, $\incl {C(Z)}{C(X)}$ is topologically
free, so it is a weak Cartan inclusion.
  \endProof

\endsection

\startsection Pseudo-expectations

\label PseudoExpSection

\sectiontitle

The set $F$ of free points of a regular inclusion $\incl BA$ is often rather badly behaved from a topological point
of view.  Illustrating this fact, in \ref {BadNormalizer} below, we will see an example in which $X$ is a closed
interval in the real line while $F$ consists of the set of irrational numbers in that interval.

However, if one adopts a more lenient point of view, the irrational numbers are perhaps not so badly behaved, since
they at least form a dense $G_\delta $ set.  In fact, in most instances in this work where we prove that an
inclusion is topologically free, namely in \ref {RelTFs.ii} and \ref {EssMetrTopFree} above, as well as \ref
{Sexy.iii$\Rightarrow $i} and \ref {WeylTopFree} below, the punch line comes as an application of Baire's theorem,
so the set of free points is actually proven to contain a dense $G_\delta $ set.

One might therefore suspect that topologically free inclusions whose set of free points satisfy the above property
might be worth considering, and this is in fact the case, as we would now like to show.

In order to describe what we have in mind, we need to use the concept of \emph {injective envelopes} and we refer
the reader to \cite [Section 2.1]{SRI} for a brief survey of the basic details.  Given a compact topological space
$X$, let us denote the injective envelope of $C(X)$ by $I\big (C(X)\big )$.  By a result of Dixmier \cite {Dixmier},
one has that
  \lbldeq injenvdes
  $$
  I\big (C(X)\big ) \simeq {\cal B}(X)/{\cal M}(X),
  $$
  where ${\cal B}(X)$ is the algebra of all bounded Borel functions on $X$, and ${\cal M}(X)$ is the ideal in ${\cal
B}(X)$ formed by the functions that vanish off a meager set.  The canonical inclusion
  \lbldeq CanonInclu
  $$
  \iota :C(X)\to I\big (C(X)\big ),
  $$
  an important ingredient of the injective envelope, may be described as the restriction to $C(X)$ of the quotient
map
  $$
  q:{\cal B}(X)\to {\cal B}(X)/{\cal M}(X).
  $$

Given any Borel subset $F\subseteq X$ containing a dense $G_\delta $ subset, consider the composition
  $$
  j: C(X) \to C^b(F) \to {\cal B}(X),
  $$
  where the leftmost arrow is the restriction map, and the rightmost arrow sends any $f$ in $C^b(F)$ to the
extension of $f$ to $X$ which vanishes off $F$.  Since the complement of $F$ in $X$ is meager, for every $f$ in
$C(X)$ one has that
  $f$ coincides with $j(f)$ modulo ${\cal M}(X)$, and this in turn implies that the composition
  $$
  C(X) \to C^b(F) \to {\cal B}(X) \labelarrow q {\cal B}(X)/ {\cal M}(X) \simeq I\big (C(X)\big )
  $$
  coincides with the canonical inclusion $\iota $ of \ref {CanonInclu}.  In particular we see that $\iota $ factors
through $C^b(F)$ and since this factorization is crucial for what we are about to do, let us name the relevant maps
in the following commutative diagram

  \begingroup \noindent \hfill \beginpicture
  \setcoordinatesystem units <0.015truecm, -0.015truecm>
  \setplotarea x from -60 to 260, y from -50 to 150
  \put {\null } at -60 -50
  \put {\null } at -60 150
  \put {\null } at 260 -50
  \put {\null } at 260 150
  \put {$C(X)$} at 0 0
  \put {$C^b(F)$} at 100 100
  \arrow <0.15cm> [0.25, 0.75] from 25 25 to 70 70
  \put {$\rho $} at 35.479 59.521
  \put {$I\big (C(X)\big )$} at 200 0
  \arrow <0.15cm> [0.25, 0.75] from 125 75 to 170 30
  \put {$\chi $} at 159.521 64.521
  \arrow <0.15cm> [0.25, 0.75] from 50 0 to 140 0
  \put {$\iota $} at 95 -12
  \endpicture \hfill \null \endgroup

  \bigskip

\state Proposition \label UniquePseudo
  Let $\incl BA$ be a regular inclusion with $A$ abelian and unital, and suppose that the set $F$ of free points in
the spectrum of $A$ contains a dense $G_\delta $ subset.  Denoting by $\Pf $ the free expectation introduced in \ref
{DefineFreeCondExp}, the composition
  $$
  B \labelarrow {\Pf } C^b(F) \labelarrow \chi I(A),
  $$
  is the unique generalized conditional expectation (as defined in \ref{DefineCondExp}) from $B$ to $I(A)$.

\Proof Writing $E=\chi \circ \Pf$ and using the commutative diagram above, it is easy to see that $E$ is a
generalized conditional expectation.

Before turning to the uniqueness question, let us discuss the map $\iota $ of \ref {CanonInclu} from the point of
view of spectra.  Letting $Y$ be the spectrum of $I\big (C(X)\big )$, and considering the map
  $$
  \kappa :Y\to X,
  $$
  dual to $\iota $, then $\kappa $ is surjective because $\iota $ is one-to-one.  We claim that $\kappa \inv (F)$ is
dense in $Y$.  In order to prove this let us denote by $Y'$ the closure of $\kappa \inv (F)$ in $Y$. We then have
that
  $$
  F \subseteq \kappa \big (\kappa \inv (F)\big ) \subseteq \kappa (Y'),
  $$
  so $\kappa (Y')$ is seen to be a compact set containing a dense subset, hence $\kappa (Y')=X$.  By \cite
[Corollary 2.18]{HadPaul} we have that $\kappa $ is essential, as defined in \ref {DefineEssSurj}, so $Y'=Y$,
proving that $\kappa \inv (F)$ is dense, as claimed.

Finally focusing on the uniqueness question, assume that $E':B\to I(A)$ is another generalized conditional expectation.  For each $y$
in $\kappa \inv (F)$, consider the positive linear functional $\psi _y$ on $B$ given by
  $$
  \psi _y(b) = \ev {E'(b)}y, \for b\in B.
  $$
  Given $a$ in $A$, we then have that
  $$
  \psi _y(a) =
  \ev {E'(a)}y =
  \ev {\iota (a)}y =
  \ev a{\kappa (y)} =
  \varphi _{\kappa (y)}(a),
  $$
  so we see that $\psi _y$ is a state on $B$ extending $\varphi _{\kappa (y)}$, and hence it is the unique such
state because $\kappa (y)$ is free, thanks to the choice of $y$ in $\kappa \inv (F)$.

Since the paragraph above applies just as well for $E$ in place of $E'$, we deduce that
  $$
  \ev {E(b)}y = \psi _y(b) = \ev {E'(b)}y, \for b\in B.
  $$
  Having seen that $\kappa \inv (F)$ is dense in $Y$, we conclude that $E(b)=E'(b)$, hence the required uniqueness.
  \endProof

Recall from
\cite [Definition 1.3]{SRI} and \cite [Definition 3.1]{KM}
that a generalized conditional expectation taking values in the injective envelope $I(A)$ is
called a \emph{pseudo-expectation}.   Adopting this terminology, the statement of \ref{UniquePseudo} then claims that there
exists a unique pseudo-expectation.

\endsection

\PART {II -- GROUPOIDS}

  \def \Nx {\Npt x;}
  \def \scrM {\smallmathscr {M}\, }
  \def \CcGx {C_c\big (\G (x), {\smallmathscr {M}\, } \big )}
  \def \CstarGx {C^*\big (\G (x), \scrM \big )}

\startsection \'Etale groupoids

\sectiontitle

So far our standing assumptions have always involved a fixed inclusion of C*-algebras.  From this point on we
will instead concentrate our attention on groupoids which will in due time  lead to a fundamental example of
inclusions of C*-algebras.  We therefore begin with a brief introduction to the theory of groupoids and their
C*-algebras, referring the reader to \cite{Renault}, \cite{Paterson},
\cite[Section 3]{actions} and \cite{Sims} for alternative
treatments and proofs of the statements left unproven in this section.

\definition \label succinctGpoid
A \emph{groupoid} is a set $\G$ equipped with
  an
``inverse map''
  $$
  \gamma \in \G \mapsto  \gamma^{-1}\in \G,
  $$
  as well as a ``multiplication operation''
  $$
  (\alpha , \beta ) \in   \G^{(2)} \mapsto  \alpha \beta \in \G,
  $$
  where $\G^{(2)}$ is a given subset of $\G \times \G$,
such that
  \izitem
  \zitem $(\gamma^{-1})^{-1} = \gamma$, for all $\gamma \in \G$,
  \zitem if $(\alpha,\beta)$ and $(\beta,\gamma)$ belong to
    $\G^{(2)}$, then $(\alpha\beta,\gamma)$ and $(\alpha,\beta\gamma)$
    belong to $\G^{(2)}$, and $(\alpha\beta)\gamma = \alpha(\beta\gamma)$,
  \lbldzitem GpdInvs $(\alpha ,\alpha ^{-1}) \in \G^{(2)}$ for all $\alpha  \in
    \G$, and for all $(\alpha ,\beta ) \in \G^{(2)}$, we have
  $$
  \alpha ^{-1} \alpha  \beta  = \beta ,  \and \alpha \beta \beta ^{-1} = \alpha .
  $$

This definition, taken from \cite{Hahn}, emphasizes the viewpoint
that groupoids should be seen as generalized groups,
where the multiplication operation is no longer everywhere defined.
These axioms are perhaps as succinct as possible and are
very useful when checking a specific example is a
groupoid.  However \ref{succinctGpoid} is perhaps not the ideal definition to give a
clear and intuitive picture of the notion of groupoid.  For the
reader acquainted with categories, a more substantive definition
asserts that a groupoid is a small category (that is, a category whose
objects form a set, rather than a class) and such that all morphisms
are invertible.

A \emph{unit} in a groupoid $\G$ is any element that can be written as $\gamma \gamma \inv$, for some $\gamma $ in $\G$.  Every unit  is
idempotent by \ref{GpdInvs}, and it is easy to see that, conversely, every idempotent element is necessarily a unit.
The set of all units in a groupoid $\G$, called the \emph{unit space},
is usually denoted by $\Gz$.  In symbols
  $$
  \Gz  =\big\{\gamma \gamma \inv: \gamma \in \G\big\} =\big\{u\in \G: (u, u)\in  \G^{(2)}, \ u= uu\big\}.
  $$

If $\gamma $ is any element of $\G$, the \emph{source} and \emph{range}  of $\gamma $ are respectively defined by
  $$
  s(\gamma ) = \gamma \inv\gamma , \and   r(\gamma ) = \gamma \gamma \inv.
  $$
  Observe that the source of $\gamma $ coincides with the range of $\gamma \inv$, and that both $s(\gamma )$ and $r(\gamma )$ are units.

Here is the relationship between the categorical definition of
groupoid with the definition found in~\ref{succinctGpoid}.  Given a
groupoid $\G$ as in~\ref{succinctGpoid}, one
considers the small category whose set of objects is $\Gz$, and for
any $x$ and $y$ in $\Gz$, one lets
  $$
  \text{Hom}(x,y) = \{\gamma \in \G: s(\gamma )=x, \ r(\gamma )=y\}.
  $$
  Taking the composition of morphisms to be the groupoid operation one easily proves that this is indeed a category in
which all morphisms are invertible.

\definition A \emph{topological groupoid} is a groupoid equipped with a topology relative to which both the inversion
and multiplication operations are continuous maps (where $\G^{(2)}$ is given the relative topology as a subset of $\G \times \G$).

Contrary to the point of view often adopted in certain branches
of topology, the topological groupoids considered here are not always assumed to be
Hausdorff, mainly in order not to rule out examples arising from many
important constructions.  Nevertheless, virtually all such examples
satisfy a weak separation property in that the unit space $\Gz$ is
Hausdorff for the relative topology\fn{This weak separation property
is sometimes replaced by the assumption that $\G$ is locally
Hausdorff, meaning that every point of $\G$ has a neighborhood which
is Hausdorff in the relative topology.}.

\definition An \emph{\'etale groupoid} is a topological groupoid $\G$
such that $\Gz$ is Hausdorff and locally compact, and such that the
range map
  $$
  r: \G \to  \Gz
  $$
  is a local homeomorphism.

Since $s(\gamma ) = r(\gamma \inv)$ for every $\gamma $ in $\G$, the source map  $s$ is also a local homeomorphism in any \'etale groupoid.

\state Proposition In any \'etale groupoid, the unit space
  $\Gz$ is an open subset of $\G$.

\fix From now on we shall fix an \'etale groupoid $\G$.

\definition A \emph{bisection} (also called a \emph{slice}) in $\G$ is any subset $U\subseteq G$ such that the source and range
maps are one-to-one when restricted to $U$.

Routine arguments show that the topology of  an \'etale groupoid admits a basis consisting of open bisections.

If $U$ is an open bisection in $\G$, notice that $r$ restricts to a homeomorphism
  $$
  r|_U: U \to  r(U)
  $$
  from $U$ to the open subset $r(U)\subseteq \Gz$.  Since the latter
is evidently locally compact and Hausdorff, then so is $U$.  In other
words, every open bisection is locally compact and Hausdorff.  It then
follows that \'etale groupoids are locally compact and locally
Hausdorff.

For $x$ in $\Gz$ we let
  $$
  \G _x= \{\gamma \in \G :s(\gamma )=x\}, \and
  $$
  $$
  \G (x) = \{\gamma \in \G :s(\gamma )=r(\gamma )=x\},
  $$
  the latter being called the \emph {isotropy group of\/ $\G $ at
$x$}.  The isotropy group is a group
under the induced multiplication operation, with $x$ playing the role
of the unit group element.

The union of all
isotropy groups
  $$
  \G '= \bigcup_{x\in \Gz} \G(x)  = \{\gamma \in \G :s(\gamma )=r(\gamma )\}
  $$
  is called the \emph {isotropy group bundle}.

A simple, yet enlightening example of a non-Hausdorff \'etale groupoid
is the \emph {two-headed snake}.  To construct it, let $1'$ be a
symbol, and let $X=[0,1]\cup\{1'\}$ with a base of open sets given by
the collection of sets,
$$\{(a,b)\cap [0,1]:  a < b \}\cup \{ (a,1)\cup \{1'\}:
0\leq a<1\}.$$
This set may be given the structure of a groupoid:  every
element of $[0,1]$ is a unit, and  the subset $\{1, 1'\}$ is
isomorphic to the cyclic group with two elements.

The reader interested in developing their intuition on \'etale
groupoids is urged to explore this example further, showing e.g.~that
it is indeed an \'etale groupoid, and that $[0,1]$ and $[0,1)\cup
\{1'\}$ are open bisections.

\endsection

\startsection Twists and line bundles

\sectiontitle

Associated to every \'etale groupoid $\G$ is a certain C*-algebra,
denoted $C^*(\G)$, which contains a canonical abelian subalgebra
isomorphic to $C_0(\Gz)$.   While inclusions of the form $(C_0(\Gz),
C^*(\G))$ encompass a wide variety of
examples,  we shall require a more general class of inclusions
arising from {\it twisted} groupoids.  The purpose of this short section is to
describe the notion of a twist over an \'etale groupoid.

The most widespread notion of twist (\cite[Section 4]{RenaultCartan},
see \cite[Definition~11.1.1]{Sims} for a less terse definition)
involves a groupoid extension
\lbldeq GrpdExt
  $$
  {\bf  T} \times  \Gz \to  \Sigma  \to  G,
  $$
  where $\bf T$ is the circle group.  In particular $\Sigma $ is a
principal $\bf T$-bundle, from which one may construct an associated
complex line bundle, which in turn may be shown to carry the
additional structure of a Fell bundle (see below for the full
definition).  Conversely, given such a Fell line bundle, the
associated circle bundle $\Sigma$ formed by vectors of norm-one
possesses the structure of a groupoid, yielding a groupoid
extension as in \ref{GrpdExt}.  In short, twists may be put in a
one-to-one correspondence with Fell line bundles, hence allowing us to
choose between these two formalisms when discussing twisted groupoid
C*-algebras.

In our opinion,  Fell line bundles are the  appropriate
formalism for the issues encountered in this work, and  we therefore
 emphasize them instead of twists.
Nevertheless the picture based on groupoid extensions is always
readily available through the associated circle bundle.

The notion of Fell bundles over groupoids was introduced by Kumjian in
\cite {KumjianFell}, motivated by \cite{Yg}.  Kumjian's original
definition was stated only for Hausdorff groupoids, but it only
requires only a slight modification to deal with the present
non-Hausdorff situation.  This modification is based on Buss and
Meyer's treatment of fields of Banach spaces over locally Hausdorff
spaces \cite [Appendix B]{BussMeyer}).

In what follows we give the precise definition for the special case of
interest here, namely when the base is an \'etale groupoid and the
fibers are one-dimensional.

\definition A \emph{Fell line bundle} over an \'etale groupoid $\G$  is a topological space
  $$
  \Lb = \bigsqcup_{\gamma \in \G} \Lb_\gamma
  $$
  formed by the disjoint union of one-dimensional complex Banach spaces $\Lb_\gamma $, such that, for every open bisection
$U\subseteq \G$, the restriction
  $$
  \Lbpar|_U = \bigsqcup_{\gamma \in U} \Lb_\gamma ,
  $$
  is a continuous Banach bundle \cite[II.13.4]{FD}.   Additionally,  defining $\Lbpar^{(2)}$ to be the subset of
$\Lb\times \Lb$ formed by the pairs $(\xi _1,\xi _2)\in  \Lb_{\gamma _1}\times \Lb_{\gamma _2}$, with  $(\gamma _1,\gamma _2)\in  \G^{(2)}$,
we are given a continuous ``multiplication'' operation
  $$
  \cdot : \Lbpar^{(2)} \to  \Lb,
  $$
  and a continuous ``adjoint'' operation
  $$
  *: \Lb \to  \Lb,
  $$
  such that, for every $\gamma _1, \gamma _2, \gamma _3\in \G$, with both $(\gamma _1,\gamma _2)$  and $(\gamma _2,\gamma _3)$ lying in $\G^{(2)}$, and for
every $\xi _i\in \Lb_{\gamma _i}$, where $i\in  \{1, 2, 3\}$, the
following hold:
  \izitem
  \zitem the multiplication operation maps $\Lb_{\gamma _1}\times \Lb_{\gamma _2}$ into $\Lb_{\gamma _1\gamma _2}$,
  \zitem the multiplication operation is bilinear on  $\Lb_{\gamma _1}\times \Lb_{\gamma _2}$,
  \zitem $(\xi _1\xi _2)\xi _3 = \xi _1(\xi _2\xi _3)$,
  \zitem $\|\xi _1\xi _2\| \leq  \|\xi _1\| \|\xi _2 \|$,
  \zitem the adjoint operation maps $\Lb_{\gamma _1}$ into $\Lb_{\gamma _1^{-1}}$,
  \zitem the adjoint operation is conjugate-linear on $\Lb_{\gamma _1}$,
  \zitem $\xi _1^{**} = \xi _1$,
  \zitem $(\xi _1 \xi _2)^* = \xi _2^* \xi _1^*$,
  \lbldzitem CstarProp $\|\xi _1^*\xi _1\| = \|\xi _1\|^2$,
  \lbldzitem Positiv $\xi _1^*\xi _1 \geq  0$.

\medskip Observe that the axioms implicitly state that $\Lb_u$ is a C*-algebra for every unit $u$ in $\Gz$.  Moreover,
for any $\xi _1$ as in the above definition, one has that $\xi _1^*\xi _1$ lies in the C*-algebra $\Lb_{\gamma _1\inv
\gamma _1}$. This said, the positivity referred to in \ref{Positiv} is to be interpreted with respect to that
C*-algebra.

Still referring to the C*-algebra $\Lb_u$, where $u$ is any given unit in $\Gz$, notice that it is necessarily
isomorphic to the only one-dimensional C*-algebra there is, namely $\C$.  This said, we see that $\Lb_u$ contains a
distinguished element, namely its unit.  From this it easily follows that the restriction of $\Lb$ to $\Gz$ is
necessarily the trivial one-dimensional vector bundle.  Nevertheless $\Lb$ is not necessarily trivial over
the entire groupoid $\G$.

Based on \cite [Proposition 3.4]{Hofmann} (see also \cite[Remark C18, vol 1]{FD}), some authors have generalized the
theory of Banach bundles to the upper-semicontinuous case, namely relaxing axiom \cite[II.13.4.(i)]{FD} by requiring
only that the norm be upper-semicontinuous.  For example, Buss and Meyer's treatment of fields of Banach spaces in
\cite{BussMeyer} actually applies to \emph{upper semi-continuous} fields.  However, for the special case of a Fell
bundle $\smallmathscr {B}$, this extra generality is unnecessary, because the continuity of the multiplication and
adjoint operations from \ref{CstarProp} imply that the norm on $\smallmathscr {B}$ must actually be continuous: indeed,
  $$
  x\in \smallmathscr {B}\mapsto x^*x\in \smallmathscr {B}_1\mapsto \|x^*x\|=\|x\|^2\in {\bf R}.
  $$

\endsection

\startsection The C*-algebra of a twisted groupoid

\label CstarALgGpd

\sectiontitle

We begin by describing the context within which we will conduct our work.

\state {Standing Hypotheses IV} \label StandingGPDS
  \rm
  We shall fix a (not necessarily Hausdorff), \'etale grou\-poid $\G $ whose unit space, denoted $\Gz $ or more often
$X$, is recalled to be Hausdorff and locally compact.  In addition we assume that $\G $ is equipped with a fixed Fell
line bundle $\Lb $.

  Since Fell line bundles are in one-to-one correspondence with twists, by a slight abuse of language we shall say that
the pair $(\G ,\Lbpar )$ is a \emph {twisted groupoid}.

If $U$ is an open bisection in $\G $, we will denote by $C_0(U, \Lbpar )$ the set of all continuous cross-sections
of the restricted line bundle $\Lbsys \kern 2pt|_U$ vanishing at $\infty $, while $C_c(U, \Lbpar )$ will denote the
subset of $C_0(U, \Lbpar )$ consisting of the compactly supported cross-sections.
  With the slightly different notations
  \lbldeq CZLCCL
  $$
  \newsymbol {\scrCzL U}{local cross-sections of line bundle vanishing at infinity}, \and
  \newsymbol {\scrCcL U}{compactly supported local cross-sections}
  $$
  we will denote the set of all cross-sections of $\Lb $ defined on the whole of $\G $ obtained by extending the
members of $C_0(U, \Lbpar )$ and $C_c(U, \Lbpar )$, respectively, to $\G $ by setting them to be zero off $U$.

One then defines
  \lbldeq DefCCLG
  $$
  \scrCcL \G
  $$
  to be the linear span of the union of the $\scrCcL U$, for all open bisections $U$.

As already mentioned, every bisection $U\subseteq \G$ is a locally compact Hausdorff  space and hence the usual tools of
analysis are conveniently available to study $C_0(U, \Lbpar )$,  and hence also its twin $\scrCzL U$.
The problems arising from the possible lack of Hausdorfness of $\G$ are then
tackled by focusing as much attention as possible on  $\scrCzL U$,  avoiding the consideration of more general
continuous  functions defined on the whole of $\G$.

We then make $\scrCcL \G $
into a $*$-algebra with the operations
  \lbldeq GpdOper
  $$
  (f*g)(\gamma ) = \sum _{\gamma _1\gamma _2=\gamma } f(\gamma _1)g(\gamma _2),
  \and
  (f^*)(\gamma ) = f(\gamma \inv )^*,
  $$
  for all $f,g$ in $\scrCcL \G $, and all $\gamma $ in $\G $.  The \emph {full twisted groupoid C*-algebra}, denoted
$\Cg $, is then defined to be the completion of $\scrCcL \G $ with the largest C*-norm.

In order to define the reduced algebra one adopts a different strategy, based on a collection of \emph {regular
representations}, defined as follows: for each $x$ in $X$, let
  $
  \EllTwo
  $
  denote the Hilbert space of square integrable cross-sections of the restriction of $\Lb $ to $\G _x$.

For $f$ in $\scrCcL \G $, and $\xi $ in $\EllTwo $, the formula
  \lbldeq RegRepForx
  $$
  \big (\pi _x(f) \xi \big )|_\gamma = \sum _{\gamma _1\gamma _2=\gamma } f(\gamma _1)\xi (\gamma _2), \for \gamma
\in \G (x),
  $$
  gives a well-defined element $\pi _x(f) \xi $ of $\EllTwo $.  The resulting operator $\pi _x(f)$ may then be shown
to be bounded, and the correspondence $f\mapsto \pi _x(f)$ becomes a $*$-representation of $\scrCcL \G $ on $\EllTwo
$.

The \emph {reduced twisted groupoid C*-algebra}, denoted $\Cgr $, is then defined to be the completion of $\scrCcL
\G $ with the norm
  \lbldeq ReducedNorm
  $$
  \Vert f\Vert _\red = \sup _{x\in X} \Vert \pi _x(f)\Vert , \for f\in \scrCcL \G .
  $$

It is well-known that both $\Cg $ and $\Cgr $ induce the uniform norm on $\scrCcL U$, for every bisection $U$, so
the relative closure of the latter in either one of these algebras is isomorphic to $\scrCzL U$, which will
therefore be identified as a subspace of both $\Cg $ and $\Cgr $.

For each $\gamma $ in $\G _x$ choose any element $\xi _\gamma $ in the fiber $\Lb _\gamma $ with $|\xi _\gamma |=1$,
and let $e_\gamma $ denote the cross-section of $\Lb $ over $\G _x$ mapping $\gamma $ to $\xi _\gamma $, and
vanishing everywhere else.  In addition, we stipulate that when $x\in \Gz $, $\xi _x$ is the unit of the fiber $\Lb
_x$.  One then has that $\{e_\gamma \}_{\gamma \in \G _x}$ is an orthonormal basis for $\EllTwo $.

Given any $f$ in $\scrCcL \G $, and given $\gamma $ in $\G _x$, observe that
  $$
  \big \langle \pi _x(f)e_x,e_\gamma \big \rangle =
  \big (\pi _x(f)e_x\big )\calcat \gamma e_\gamma (\gamma )^* \quebra =
  \sum _{\gamma _1\gamma _2=\gamma } f(\gamma _1)e_x(\gamma _2)e_\gamma (\gamma )^* =
  f(\gamma )e_x(x)e_\gamma (\gamma )^* =
  f(\gamma )\xi _x\xi _\gamma ^* =
  f(\gamma )\xi _\gamma ^*.
  $$
  Consequently,
  $$
  f(\gamma )=\big \langle \pi _x(f)e_x,e_\gamma \big \rangle \xi _\gamma , \for f\in \scrCcL \G ,
  $$
  from which it follows that
  \lbldeq BoundedValues
  $$
  |f(\gamma )|\leq \Vert f\Vert _\red ,
  $$
  and hence the correspondence
  $$
  f\in \scrCcL \G \mapsto f(\gamma )\in \Lb _\gamma ,
  $$
  is seen to extend to a continuous linear map $\varphi _\gamma
  :\Cgr \rightarrow \Lb _\gamma $. By abuse of language, we write
  \lbldeq ItIsAllFunctions
  $$
  f(\gamma ):= \varphi _\gamma (f), \for f\in \Cgr , \for \gamma \in \G ,
  $$
  so $f$ defines a section of $\Lb $, which is uniformly bounded but often discontinuous.  It is easy to see that if
$f$ is an element of $\Cgr $ lying in the kernel of every $\varphi _\gamma $, then $f=0$, so we may identify $f$
with the section of $\Lb $ it defines.  This practice will in fact be adopted throughout this work.

It may also be proved that the formulas \ref {GpdOper} hold for every $f$ and $g$ in $\Cgr $, so we may actually
think of $\Cgr $ as an algebra of cross-sections of $\Lb $ with the above operations.  The only tricky aspect of
doing so is that no simple criteria appears to exist in order to determine if a given cross-sections of $\Lb $
belongs to $\Cgr $.  This is akin to the difficulty in determining which bounded sequences $\{a_n\}_{n\in {\bf Z}}$
are given by the Fourier coefficients of a continuous function on the circle.

Since the reduced norm $\Vert \cdot \Vert _\red $ is of course smaller than the maximum norm, there exists a
$*$-homo\-mor\-phism
  \lbldeq IntForlLeftReg
  $$
  \Lambda :\Cg \to \Cgr ,
  $$
  extending the identity map on $\scrCcL \G $, often called the \newConcept {left regular representation}{} of $\Cg
$.  Given any $f$ in $\Cg $, we will occasionally also speak of $f(\gamma )$, actually meaning $\Lambda (f)(\gamma
)$, as above, although one should notice that if $f(\gamma )$ vanishes for every $\gamma $, all we can say is that
$\Lambda (f)=0$, which does not necessarily imply that $f=0$.

\endsection

\startsection Topologically free groupoids

\sectiontitle

Throughout this section we keep \ref{StandingGPDS} in force.
We shall then be concerned with the inclusion $\incl {B}{A}$, where
  $$
  A=C_0(X),
  $$
  and
  \lbldeq TwoAlgebras
  $$
  B= \left \{\matrix { \Cg \hfill \cr \pilar {15pt} \Cgr \hfill }\right .
  $$
  meaning that, unless indicated, $B$ could be either $\Cg $ or $\Cgr $.  It is well-known that
  $$
  \scrCzL {U}\subseteq \Norm BA,
  $$
  for any open bisection $U$ (whichever version of $B$ one happens to choose).

\definition \label DefineNG
  The union (rather than the linear span) of all $\scrCzL {U}$, as $U$ ranges in the collection of all open bisections of $\G $, will henceforth
be denoted by $\newsymbol {N_\G }{set of all normalizers supported on bisections} $.

It may be easily proved that $N_\G $ is $*$-subsemigroup of $\Norm BA $, and that
  $$
  A\subseteq N_\G ,\and \clspan N_\G =B,
  $$
  so
  $\incl {B}{A}$ is a regular inclusion, and $N_\G $ is a generating $*$-semigroup.

For Hausdorff groupoids, and under certain favorable hypotheses, $N_\G $ includes all normalizers \cite [Proposition
4.2]{RenaultCartan} but this is no longer true in the non-Hausdorff case, even though the other \emph {favorable
hypotheses} are kept in place \cite [3.2]{ExelNHaus}.

Given $n$ in $N_\G $, let $U_n$ be the open support of $n$, namely
  \lbldeq OpenSuport
  $$
  U_n=\{\gamma \in \G :n(\gamma )\neq 0\}.
  $$
  Choosing any open bisection $U$ such that $n\in \scrCzL {U}$, one sees that $U_n$ is an open subset of $U$, hence
$U_n$ is an open bisection as well, and $n\in \scrCzL {U_n}$.  The domain $\Dom n$ of $\beta _n$, as defined in \ref
{Beta}, is then easily seen to coincide with $s(U_n)$, and for each $x$ in $s(U_n)$, one has that
  \lbldeq BetaDois
  $$
  \beta _n(x)=r(\gamma ),
  $$
  where $\gamma $ is the unique element in $U_n$ such that $s(\gamma )=x$.

The following is the justification for the adoption of the term \emph {isotropy algebra} back in \ref {DefIsoMod}.

\state Theorem \label IdentifyIsotropyAlgebras
  Given a twisted groupoid $(\G ,\Lbpar )$, as in \ref {StandingGPDS}, pick any $x$ in $\Gz $, and let $\scrM $ be
the restriction of $\Lb $ to the isotropy group $\G (x)$.  Considering the inclusion
  $
  \INCL BA,
  $
  described in \ref{TwoAlgebras},
  denote by
  $$
  C=
  \left \{\matrix {
    \kern 5pt C^*\big (\G (x),\scrM \big ) \hfill \cr \pilar {15pt}
    C^*_\red \big (\G (x) , \scrM \big ) \hfill
    }\right .
  $$
  enforcing the first alternative for $C$ if and only if $B$ is chosen to be the first alternative in \ref
{TwoAlgebras}.
  Recalling that $B(x)$ denotes the isotropy algebra (see \ref {DefIsoMod} and \ref {DataBasedOnPoint}) of the above
inclusion relative to $x$, there is a *-homomorphism $\theta $ from $B(x)$ onto $C$, such that the following diagram
commutes:

  \begingroup \noindent \hfill \beginpicture
  \setcoordinatesystem units <0.018truecm, -0.018truecm>
  \setplotarea x from -150 to 170, y from -30 to 140
  \put {\null } at -150 -30
  \put {\null } at -150 140
  \put {\null } at 170 -30
  \put {\null } at 170 140
  \put {$B$} at 0 0
  \put {$B(x)$} at -100 100 \arrow <0.11cm> [0.5, 1.8] from -20 20 to -80 80
  \put {$E_x$} at -64.142 35.858 \arrow <0.11cm> [0.5, 1.8] from 20 20 to 80 80
  \put {$\rho $} at 58.485 41.515
  \put {$C$} at 100 100 \arrow <0.11cm> [0.5, 1.8] from -60 100 to 70 100
  \put {$\theta $} at 5 112
  \endpicture \hfill \null \endgroup

\noindent where $\rho $ is the isotropy projection introduced in \ref {IsotropyProjection}.  In addition,
  \izitem
  \zitem there are dense *-subalgebras of $B(x)$ and $C$, respectively, which are *-isomorphic under $\theta $,
  \lbldzitem IdIsoAlgsFull in case $B$ and $C$ are chosen to be the full C*-algebras (first choices above), then
$\theta $ itself is a *-isomorphism.

\Proof
  We shall initially work simultaneously with the different alternatives, observing that $\scrCcL \G $ and $\CcGx $
are dense *-subalgebras of $B$ and $C$, respectively, regardless of our choice.

  As usual we shall denote by $J_x$ the ideal of $\CzGz $ formed by the functions vanishing at $x$.
  Given $a$ in $J_x$, and $b$ in $\scrCcL \G $, notice that for every $g$ in $\G _x$, one has that
  $$
  (ab)(g) = \langle a, r(g)\rangle b(g) = \langle a, x\rangle b(g) = 0,
  $$
  and similarly $(ba)(g)=0$.  This shows that $ab$ and $ba$ vanish on $\G (x)$, whence $\rho $ vanishes on
  $$
  J_x\scrCcL \G +\scrCcL \G J_x.
  $$
  Consequently $\rho $ also vanishes on $L_x$, so it factors through $\Ex $, providing a bounded linear map $\theta
$ from $B(x)$ to $C$, such that the above diagram commutes.  Turning now to the proof of (i), we first claim that,
for all $b\in \scrCcL \G $, one has that
  \lbldeq DenseInjective
  $$
  \Ex (b)=0 \IFF \rho (b)=0.
  $$

  Since $\rho $ factors through $\Ex $, the implication ``$\Rightarrow $" holds trivially, so let us focus on the
the reverse implication by choosing $b$ in $\scrCcL \G $ such that $\rho (b)=0$.  In order to prove that $\Ex
(b)=0$, let us initially tackle the case in which $b\in \scrCcL {U}$, for some open bisection $U$.  Keeping in mind
that such a $b$ is a normalizer of $A$ in $B$, observe that the present case in turn breaks down into the following
three subcases:

\medskip \noindent {\tensc Subcase 1}: If $x\notin \Dom b$, then by \ref {NormaJota} we have that $b\in
BJ_x\subseteq L_x$, whence $\Ex (b)=0$, as desired.

\medskip \noindent {\tensc Subcase 2}: If $x\in \Dom {b}$, but $\beta _b(x)\neq x$, then again by \ref {NormaJota}
we have that $b\in L_x$, so $\Ex (b)=0$.

\medskip \noindent {\tensc Subcase 3}: If $x\in \Dom {b}$, and $\beta _b(x)= x$, then by \ref {BetaDois} there is
some $\gamma $ in $\G $ such that $b(\gamma )\neq 0$, and $s(\gamma ) = x = r(\gamma )$.  So we see that $b$ does
not vanish on the element $\gamma \in \G (x)$, hence contradicting the assumption that $\rho (b)=0$.  Consequently
this subcase cannot occur and hence may be discarded.

\medskip
  Returning to the general case, let us write
  $$
  b=\sum _{i\in I}n_i,
  $$
  where $I$ is a finite index set, and each $n_i$ lies in $\scrCcL {U_i}$, for some open bisection $U_i$.

Suppose that, for a given $i\in I$, one has that $\rho (n_i)=0$.  Then, by the case already proven above, one has
that $\Ex (n_i)=0$, and therefore $n_i$ contributes to neither side of \ref {DenseInjective}, so we may simply
ignore it.  In other words, we may assume, without loss of generality, that $\rho (n_i)\neq 0$, for every $i\in I$.
The eventual vanishing of either side of \ref {DenseInjective} should therefore be due to cancelation, rather than
the vanishing of each individual summand corresponding to each $n_i$.

Recalling that each $n_i$ is supported on the bisection $U_i$, there can be at most one element $g$ in $\G _x$ for
which $n_i(g)\neq 0$.  Moreover, since we are assuming that $n_i$ does not vanish identically on $\G (x)$, there
must in fact be exactly one $g$ in $\G (x)$ such that $n_i(g)\neq 0$.  Consequently, defining
  $$
  I_g = \{i\in I: n_i(g)\neq 0\},
  $$
  we see that $I$ splits as the disjoint union
  $$
  I = \bigsqcup _{g\in \G (x)}I_g.
  $$
  Accordingly, for each $g$ in $\G (x)$, let us define
  $$
  b_g = \sum _{i\in I_g}n_i,
  $$
  so evidently
  $b=\sum _{g\in \G (x)}b_g$.  Because $b_h(g)=0$, when $h\neq g$, it follows that
  $$
  0 = b(g) = \sum _{h\in \G (x)}b_h(g) = b_g(g),
  $$
  so in fact $b_g$ also vanishes at $g$, and hence $\rho (b_g)=0$.

We will next prove that $\Ex (b_g)=0$, for every $g$ in $\G (x)$, from where claim \ref {DenseInjective} will
follow.
  Fixing $g$ in $\G (x)$, define
  $$
  U_g=\bigcap _{i\in I_g}U_i.
  $$
  Since $n_i(g)\neq 0$, for every $i$ in $I_g$, we have that $g$ lies in every $U_i$, hence also in $U_g$.
Consequently
  $$
  x=r(g)\in r(U_g),
  $$
  so we may choose $a$ in $\CurlyCcGz $, such that $\langle a, x\rangle =1$, and such that $a$ vanishes outside some
compact set $K\subseteq r(U_g)$.
  Notice that, for every $i$ in $I_g$, and every $\gamma \in \G $, we have that
  $$
  (an_i)(\gamma ) = \langle a, r(\gamma )\rangle \, n_i(\gamma ),
  $$
  so we see that $an_i$ vanishes off
  $$
  L:=\{\gamma \in U_i:r(\gamma )\in K\}.
  $$
  Recalling that $r$ is a bijection from $U_i$ to $r(U_i)$, one has that $L$ does not depend on $i$, and also that
$L$ is a compact subset of $U_g$.

  \begingroup \noindent \hfill \beginpicture
  \setcoordinatesystem units <0.015truecm, 0.015truecm>
  \setplotarea x from -400 to 400, y from -180 to 160
  \put {\null } at -400 -180
  \put {\null } at -400 160
  \put {\null } at 400 -180
  \put {\null } at 400 160 \plot -260 -25 -160 -25 -160 25 -260 25 -260 -25 /
  \put {$L$} at -210 40
  \put {$g ^{\,\bullet } $} at -210 0 \plot 160 -25 260 -25 260 25 160 25 160 -25 /
  \put {$K$} at 210 40
  \put {$x ^{\,\bullet } $} at 210 0 \ellipticalarc axes ratio 3:2 63 degrees from 130 140 center at 0 0
  \arrow <0.15cm> [0.25, 0.75] from 125 142 to 130 140
  \put {$r$} at 0 180 \setdashes <1.5pt> \ellipticalarc axes ratio 3:4 360 degrees from -90 0 center at -210 0
  \put {$U_i$} at -210 186.667 \ellipticalarc axes ratio 4:3 360 degrees from -50 0 center at -210 0
  \put {$U_j$} at -390 0 \circulararc 360 degrees from -130 0 center at -210 0
  \put {$U_g$} at -280.711 70.711 \ellipticalarc axes ratio 3:4 360 degrees from 330 0 center at 210 0
  \put {$r(U_i)$} at 210 186.667 \ellipticalarc axes ratio 4:3 360 degrees from 370 0 center at 210 0
  \put {$r(U_j)$} at 405 0 \circulararc 360 degrees from 290 0 center at 210 0
  \put {$r(U_g)$} at 148.283 84.947
  \endpicture \hfill \null \endgroup

\noindent One then has that $an_i\in \scrCcL {U_g}$, for every $i$ in $I_g$, and so also
  $$
  ab_g = \sum _{i\in I_g}an_i \in \scrCcL {U_g}.
  $$
  We have already seen that $\rho (b_g)=0$, so also $\rho (ab_g)=0$.  Thanks to the first case already proved above,
we then have that $\Ex (ab_g)=0$, and hence
  $$
  0 = \Ex (ab_g) \={VanishLyx} \langle a, x\rangle \Ex (b_g) = \Ex (b_g),
  $$
  hence concluding the proof of claim \ref {DenseInjective} in the general case.

Defining
  $$
  B_0(x):=\Ex \big (\scrCcL \G \big ),
  $$
  observe that $B_0(x)$
  is a dense subset of $B(x)$, linearly spanned by the union of the sets $\Ex \big (\scrCcL U\big )$, for all open
bisections $U\subseteq \G $.  Observe also that, if $U$ is a bisection and $n\in \scrCcL U$, then by \ref
{NormaJota}, one has that $\Ex (n)=0$, unless $x\in \Dom n$ and $\beta _n(x)=x$, which is to say that
  $$
  n\in \Nx := \Npt xx;,
  $$
  according to the notation introduced in \ref {NorGenCx}.  Consequently $B_0(x)$ is linearly spanned by the set
  \lbldeq GenBZero
  $$
  \Ex \big (\Nx \cap \scrCcL \G \big ).
  $$

Since $\Nx $ is closed under multiplication and adjoints, and since $\Ex $ is multiplicative on $\Cx $ (which
contains $\Nx $ by \ref {NxyInCxy}), we conclude that $B_0(x)$ is a dense *-subalgebra of $B(x)$.

On the other hand, the image of $\scrCcL \G $ under $\rho $ clearly coincides with the dense *-subalgebra $\CcGx
\subseteq C$.  We may then consider the following restricted version of the diagram in the statement

  \begingroup \noindent \hfill \beginpicture
  \setcoordinatesystem units <0.02truecm, -0.02truecm>
  \setplotarea x from -150 to 170, y from -30 to 140
  \put {\null } at -150 -30
  \put {\null } at -150 140
  \put {\null } at 170 -30
  \put {\null } at 170 140
  \put {$\scrCcL \G $} at 0 0
  \put {$B_0(x)$} at -100 100 \arrow <0.11cm> [0.5, 1.8] from -20 20 to -75 75
  \put {$\Ex ^0$} at -61.642 33.358 \arrow <0.11cm> [0.5, 1.8] from 20 20 to 75 75
  \put {$\rho ^0$} at 61.642 33.358
  \put {$\CcGx $} at 100 100 \arrow <0.11cm> [0.5, 1.8] from -60 100 to 30 100
  \put {$\theta ^0$} at -15 112
  \endpicture \hfill \null \endgroup

\noindent the zero superscripts reminding us that these maps are the corresponding restrictions of the ones in the
main diagram.

Both $\Ex ^0$ and $\rho ^0$ are surjective by construction, hence so is $\theta ^0$.  Moreover $\theta ^0$ is
injective by claim \ref {DenseInjective}, so we see that $\theta ^0$ is a linear isomorphism.  We will next prove
that $\theta ^0$ is a *-isomorphism.

Recalling that \ref {GenBZero} spans $B_0(x)$, in order to prove that $\theta ^0$ is multiplicative it is enough to
check that
  $$
  \theta ^0\big (\Ex (n_1)\Ex (n_2)\big ) =
  \theta ^0\big (\Ex (n_1)\big ) \theta ^0\big (\Ex (n_2)\big ),
  $$
  whenever $n_i\in \Nx \cap \scrCcL {U_i}$, for $i=1, 2$, where the $U_i$ are open bisections.  Since $\Ex $ is
multiplicative on $\Nx $, the left hand side above coincides with $\theta ^0\big (\Ex (n_1n_2)\big )$, so our task
above translates to proving that
  \lbldeq RhoMult
  $$
  \rho (n_1n_2) = \rho (n_1) \rho (n_2).
  $$

Since $n_i\in \Nx $, we have that $\beta _n(x)=x$, which, according to \ref {BetaDois}, means that there is some
$g_i$ in $\G $ such that $n_i(g_i)\neq 0$, and $s(g_i)=x=r(g_i)$.  Such a $g_i$ is evidently in $\G (x)\cap U_i$,
and since $U_i$ is a bisection, we actually have that
  $$
  \G (x)\cap U_i = \{g_i\}.
  $$
  It follows that $g_1g_2 \in U_1U_2$, so
  $$
  \G (x)\cap (U_1U_2) = \{g_1g_2\},
  $$
  and, observing that $n_1n_2\in \scrCcL {U_1U_2}$, we also have
  $$
  (n_1n_2)(g_1g_2) = n_1(g_1)n_2(g_2),
  $$
  from where
  \ref {RhoMult} follows at once.
    We leave it for the reader to prove that $\theta ^0$ also preserves adjoints, so $\theta ^0$ is a *-isomorphism,
as claimed.  Summarizing we have the following diagram

  \begingroup \noindent \hfill \beginpicture
  \setcoordinatesystem units <0.015truecm, -0.015truecm>
  \setplotarea x from -150 to 180, y from -50 to 150
  \put {\null } at -150 -50
  \put {\null } at -150 150
  \put {\null } at 180 -50
  \put {\null } at 180 150
  \put {$B(x)$} at -100 0
  \put {$C$} at 100 0 \arrow <0.11cm> [0.5, 1.8] from -56 0 to 60 0
  \put {$\theta $} at 2 -18
  \put {$B_0(x)$} at -100 100
  \put {$\qquad \CcGx $} at 100 100 \arrow <0.11cm> [0.5, 1.8] from -56 100 to 36 100
  \put {$\theta ^0$} at -10 118
  \put {$\cup\kern1pt\vrule height 6pt width 0.4pt$} at -100 50
  \put {$\cup\kern1pt\vrule height 6pt width 0.4pt$} at 100 50
  \endpicture \hfill \null \endgroup

\noindent where the bottom map is a *-isomorphism of dense *-subalgebras, so it follows that $\theta $ must be, at
least, a *-homomorphism of C*-algebras.

This proves the statement up to and including (i), so we now focus on (ii), hence choosing $B$ and $C$ to be the
full C*-algebras.

\def \tnorm #1{\def \tri {\hbox {$|\kern -1.5pt|\kern -1.5pt|$}}\tri #1 \tri }

Being the restriction of the continuous map $\theta $, it is clear that $\theta ^0$ is continuous, and we will next
prove that $(\theta ^0)\inv $ is also continuous.  To see this, observe that the norm $\tnorm {{\kern 1pt \cdot
\kern 1pt}}$ defined on $\CcGx $ by
  $$
  \tnorm y = \|(\theta ^0)\inv (y)\|, \for y\in \CcGx ,
  $$
  is a C*-norm, hence dominated by the largest of all C*-norms on that algebra, namely the norm $\|{\kern 1pt \cdot
\kern 1pt}\|$ inherited by $\CstarGx $.  This is to say that
  $$
  \|(\theta ^0)\inv (y)\|\leq \|y\|,
  $$
  so $(\theta ^0)\inv $ is indeed continuous, as claimed.  From this it is now clear that $\theta $ is a
*-isomorphism of C*-algebras, concluding the proof.
  \endProof

  Unfortunately the conclusion of \ref{IdIsoAlgsFull} does not extend to reduced C*-algebras.    That is,  under the choice
$B=\Cgr$ in \ref{TwoAlgebras}, and $C=C^*_\red \big (\G (x) , \scrM \big )$ in \ref{IdentifyIsotropyAlgebras}, the map
  $$
  \theta  :  B(x) \to  C
  $$
  cannot be proven to be an isomorphism.  For a counter example, recall that Willett \cite{Willett} has found an example
of a Hausdorff,   \'etale,  non-amenable groupoid $\G$ whose maximal and reduced
C*-algebras coincide.  Within Willett's groupoid there is a unit $x$ whose isotropy group $\G(x)$ coincides with the
free group on
two generators.  Taking the trivial twist on $\G$, we then have that
  $$
  C=   C^*_\red(\G(x)) =   C^*_\red({\bf F}_2),
  $$
  while
  $$
  B(x) =
  C^*_\red(\G)(x)  =
  C^*(\G)(x)  \= {IdIsoAlgsFull}
  C^*(\G(x) ) =
  C^*({\bf F}_2),
  $$
  so $B(x)\not\simeq C$.

  The notions of \emph {freeness} for groupoids have evolved along the years and have been known by various names
which are not always consistent with one another in the literature.  Below we list some of these, although we will
only be concerned with \ref {DefTopoFreeGpd}, below.

\definition \label DefVariousTopoFree
  An \'etale groupoid $\G $ is said to be
  \izitem
  \zitem
    \newConcept {principal}{principal groupoid} when $\G (x)=\{x\}$, for every $x$ in $X$,
  \lbldzitem DefTopoFreeGpd
    \newConcept {topologically free}{topologically free groupoid}, or
    \newConcept {topologically principal}{topologically principal groupoid}, if the subset of $X$ formed by the
points $x$ for which $\G (x)=\{x\}$, is dense in $X$,
  \zitem
    \newConcept {essentially principal}{essentially principal groupoid} \cite [Definition II.4.3]{Renault}, when for
every invariant closed subset $Y\subseteq X$, the set of points $x$ in $Y$ for which $\G (x)=\{x\}$, is dense in
$Y$,
  \lbldzitem EffectiveGpd
    \newConcept {effective}{effective groupoid} \cite [page 137]{Mrcun}, when the interior of $\G '$ coincides with
$X$.

We observe that the expression \emph {essentially principal} has also been used in \cite [Definition
3.1]{RenaultCartan}, to mean \ref {DefVariousTopoFree.iv}, rather than \ref {DefVariousTopoFree.iii}.

An easy consequence of \ref {IdentifyIsotropyAlgebras} is the following characterization of free points of $X$
relative to the inclusion under analysis.

\state Corollary \label CharacFreePts
  \izitem
  \zitem A point $x$ in $X$ is free relative to $\incl BA$ if and only if the isotropy group $\G (x)$ coincides with
$\{x\}$.
  \lbldzitem TopoFreeFull $\Incl {\Cg }{C_0(X)}$ is a topologically free inclusion (see Definition \ref
{ManyTopFree}), if and only if\/ $\G $ is a topologically free groupoid.
  \lbldzitem TopoFreeRed $\Incl {\Cgr }{C_0(X)} $ is a topologically free inclusion if and only if\/ $\G $ is a
topologically free groupoid.
  \lbldzitem FreeEvalIsExp If $x$ is a free point in $X$, then the localizing projection $\Ex $ associated to $x$ is
given by
  $$
  \Ex (f) = f(x), \for f\in B.
  $$

\Proof (i)\enspace
  Follows immediately from \ref {IdentifyIsotropyAlgebras} and \ref {IsotropyOneDim}.

\itmProof (\LocalTopoFreeFull \ \& \LocalTopoFreeRed )
  Follow from (i).

\itmProof (\LocalFreeEvalIsExp )
  Observe that for every $x$ in $X$, the correspondence
  $$
  f\mapsto f(x)
  $$
  defines a state on $B$ which evidently extends the state $\varphi _x$ on $C_0(X)$ given by evaluation at $x$.  In
case $x$ is a free point, this must necessarily be the state $\psi _x$ of which \ref {ThreeValues} speaks, whence
(iv).
  \endProof

This therefore reconciles the classical notion of topological freeness with
  \ref {ManyTopFree.i}.
  However, for inclusions determined by groupoids and line bundles, as the ones we are currently studying, there is
more than meets the eye.

\state Proposition \label Sexy
  Consider the following statements:
  \izitem
  \zitem $\G $ is topologically free,
  \zitem the interior of\/ $\G '\setminus X$ is empty,
  \zitem every normalizer in $N_\G $ is smooth.
  \medskip \noindent Then (i)$\Rightarrow $(ii)$\Rightarrow $(iii).  If $\G $ is second countable then also
(iii)$\Rightarrow $(i).

\Proof
  (i)$\Rightarrow $(ii)\enspace
  Suppose by way of contradiction that $V$ is a nonempty open set contained in $\G '\setminus X$.  Then $s(V)$ is a
nonempty open subset of $X$ containing no free points, conflicting with the fact that the free points are dense,
hence proving (ii).

\itmImply (ii) > (iii) Pick $n$ in $N_\G $, and recall that $n$ lies in $\scrCzL {U_n}$, where $U_n$ is as in \ref
{OpenSuport}.
  It then follows from (ii) that the interior of
  $$
  W:= U_n\cap (\G '\setminus X)
  $$
  is empty.
  Observing that $\G '$ is closed and $X$ is open, it follows that $W$ is a closed subset of $U_n$, and hence that
$W$ is rare in $U_n$.  We leave it as an easy exercise for the reader to prove that
  $$
  \clsr {s(W)} \subseteq s(W) \cup \partial s(U_n),
  $$
  where the closure and boundary above are taken relative to $X$.

Still working within $X$, we next claim that
  $\clsr {s(W)}$ has empty interior.  To prove it, assume that $V$ is a nonempty open set contained in $\clsr
{s(W)}$, whence also
  $$
  V \subseteq s(W) \cup \partial s(U).
  $$

  Since $\partial s(U)$ is the boundary of an open set, it clearly has no interior, so $V$ cannot be entirely
contained in $\partial s(U)$.  Thus $V\setminus \partial s(U)$ is a nonempty open set contained in $s(W)$,
contradicting the fact that $W$ has empty interior.

The grand conclusion so far is that
  $$
  X\setminus \clsr {s(W)}
  $$
  is a dense open set, and we will conclude the proof by showing that all of the points in this set are free
relative to $n$.  We will actually prove a bit more, namely that every point not in $s(W)$ is free relative to $n$.

  For this let $x$ in $X\setminus s(W)$.  If $x$ is not in $\Fix n$, then \ref {MandatoryValues.i} takes care of our
needs, so assume that $x$ is in $\Fix n$, which is the same as saying that $\beta _n(x)=x$, so \ref {BetaDois}
implies that there exists some $\gamma $ in $U_n$ such that $s(\gamma )=x$ and $r(\gamma )=x$.  This $\gamma $ is
then in $U_n\cap \G '$, but it cannot be in $W$, since
  $$
  s(\gamma )= x\notin s(W).
  $$

  A quick glance at the definition of $W$, to be found near the beginning of this proof, will convince the reader
that $\gamma $ belongs to $X$.  Therefore $\gamma $ also lies in $U_n\cap X$.  Choosing some $v$ in $A$, supported
in $U_n\cap X$, such that $\evx v=1$, one may now show that $nv$ lies in $A$ so, if $\psi _1$ and $\psi _2$ are
states on $B$ extending $\varphi _x$, one has that
  $$
  \psi _1(n) =
  \psi _1(n) \evx v\={CondexForState}
  \psi _1(nv) = \varphi _x(nv) =
  \psi _2(nv) \={CondexForState}
  \psi _2(n) \evx v=
  \psi _2(n),
  $$
  proving that $x$ is indeed free relative to $n$.  This proves our claim that $X\setminus \clsr {s(W)} \subseteq
F_n$, so the smoothness of $n$ follows, verifying (iii).

\itmImply (iii) > (i) Assuming that $\G $ is second countable, we have that $\Cgr $ is separable, hence so is $N_\G
$.  One may therefore find a countable dense $*$-subsemigroup $N\subseteq N_\G $, which is therefore also
generating.  Using \ref {XfAndXb} we then have that
  $$
  F= \medcap _{n\in N} F_n,
  $$
  and the conclusion follows from (iii) and Baire's Theorem.
  \endProof

It is interesting to contrast \ref {DefVariousTopoFree.iv} with
  \ref {Sexy.ii}: it is obvious that the former implies the latter, that is
  $$
  \int {(\G ')}=X\quad \Rightarrow \quad \int {(\G '\setminus X)}=\emptyset ,
  $$
  but the converse fails in general (see \ref {CrazyGroupoids} and the discussion following it), although it holds
when $\G $ is Hausdorff, as the reader may easily prove based on the fact that $X$ is closed in $\G $.

We believe condition \ref {Sexy.ii}, above, has first been considered in \cite {BCFS}, even though only in the
context of Hausdorff groupoids, where it is equivalent to \ref {DefVariousTopoFree.iv}, as already mentioned.

\medskip As a consequence of \ref {Sexy}, we have the following improvement of \ref {CharacFreePts.ii\&iii}.

\state Corollary \label CharacFreePtsTwo
  Let $\incl {B}{A}$ denote any one of the inclusions
  $$
  \Incl {\Cg }{C_0(X)},\and \Incl {\Cgr }{C_0(X)}.
  $$
  If\/ $\G $ is a topologically free groupoid, then $\incl {B}{A}$ is a smooth inclusion.  If, in addition $\G $ is
second countable, then the converse also holds.

\Proof
  Assuming that $\G $ is topologically free we have by \ref {Sexy} that $N_\G $ consists of smooth normalizers,
meeting all of the conditions for a smooth inclusion.  If $\G $ is second countable, then $B$ is separable, so the
converse follows from \ref {RelTFs.ii}.
  \endProof

\endsection

\startsection The essential groupoid C*-algebra

\sectiontitle

In this section we will introduce a new version of the twisted C*-algebra for a topologically free groupoid, which
should be viewed as an alternative to the reduced algebra in the non-Hausdorff case.  Supporting this point of view
is the fact that these algebras may be abstractly characterized, much in the same way Renault has characterized
Cartan algebras in \cite {Renault}, as will be proved in our main result, Theorem \ref {MainTwo}, below.

Under \ref {StandingGPDS}, which we assume throughout this section, it is easy to see that the set $F$ of free
points is an invariant subset of $X$, in the sense that if $\gamma \in \G $, then
  \lbldeq InvarFree
  $$
  s(\gamma )\in F\IFF r(\gamma )\in F.
  $$

  Regardless of the perhaps badly behaved topology of $F$ (it is often not locally compact) we shall next consider
the reduction of $\G $ to $F$.

\definition
  We shall denote by $\G \free $ the reduction of $\G $ to $F$, namely
  $$
  \newsymbol {\G \free }{free component} = \{\gamma \in \G : s(\gamma )\in F\}.
  $$
  We shall refer to $\G \free $ as the \newConcept {free component} {free component of groupoid} of $\G $.

\state Proposition \label DescribeGray
  \izitem
  \zitem The {opaque} ideal of\/ $\Cgr $ is trivial.
  \zitem The {gray} ideal of $B$ (recall from \ref {TwoAlgebras} that $B$ denotes either the full or the reduced
groupoid C*-algebra) is given by
  $$
  \def \quad { }
  \matrix {
  \Gamma & = & \Big \{f\in B:& (f^*f)(x) = 0,& \hbox { for all } x\in F\Big \}\ \cr \pilar {16pt}
        & = & \Big \{f\in B:& f(\gamma ) = 0, & \hbox { for all } \gamma \in \G \free \Big \}.\cr \pilar {10pt}}
  $$

\Proof
  The map
  $$
  P_\red : f\in \Cgr \mapsto f|_X
  $$
  is a faithful generalized conditional expectation from $\Cgr $ to the algebra of all bounded functions on $X$,
hence the first point follows from \ref {IntroBlack.iii}.

The first equality in (ii) follows immediately from \ref {FreeEvalIsExp} and the definition of $\Gamma $ given in
\ref {IntroGray}.  Given $f$ in $B$, and given any $x$ in $F$, observe that
  $$
  (f^*f)(x) =
  \sum _{\gamma \inv \gamma =x} f^*(\gamma \inv )f(\gamma ) =
  \sum _{\gamma \in \G _x} |f(\gamma )|^2,
  $$
  so the second equality in (ii) follows.
  \endProof

The generalized conditional expectation referred to in the proof of \ref {DescribeGray.i} may fail to be faithful on
$\Cg $ if $\G $ is not amenable, hence the omission of the full groupoid C*-algebra there.

Since we are interested in {light} inclusions we will consider the following:

\definition \label DefineFreeGpd
  The \newConcept {essential groupoid C*-algebra}{} of $(\G ,\Lbpar )$ is the quotient C*-algebra
  $$
  \newsymbol {\Cgf }{essential groupoid C*-algebra} = \Cg /\Gamma ,
  $$
  where $\Gamma $ is the {gray} ideal.  The quotient map
  $$
  \newsymbol {\qess }{essential projection} :\Cg \to \Cgf
  $$
  will moreover be referred to as the \newConcept {essential projection}{}.

If the free component of $\G $ is too small, then clearly $\Cgf $ might lose a lot of information.  In the extreme
non-free case, namely when $\G $ is a group, the {gray} ideal is clearly the whole of $\Cg $, and hence $\Cgf $
reduces to zero.  One would therefore prefer to consider this construction only for groupoids whose free component
is not too small, such as for topologically free groupoids,  which is precisely the case we are interested in this work.

See section \ref{KMSection} below for an alternative version of the essential algebra in case $\G$ fails to be topologically free.

The occurrence of the full groupoid C*-algebra in the above definition might suggest that there is a reduced version
of $\Cgf $, built by moding out the {gray} ideal of $\Cgr $, instead.  However the following result shows that there
isn't.

\state Proposition \label FullRedGray
  If $\G $ is topologically free, the left regular representation
  $$
  \Lambda :\Cg \to \Cgr ,
  $$
  described in \ref {IntForlLeftReg}, factors through the respective {gray} ideals $\Gamma \trianglelefteq \Cg $,
and\/ $ \Gamma _\red \trianglelefteq \Cgr $, leading up to an isomorphism
  $$
  \Cgf = \Frac {\Cg }{\Gamma } \ \simeq \ \Frac {\Cgr }{\Gamma _\red }.
  $$

\Proof
  Since $\Incl {\Cg }{C_0(X)}$ is a topologically free inclusion by \ref {CharacFreePts.ii}, and since the
restriction of $\Lambda $ to $C_0(X)$ is the identity map, the result follows from \ref {TwoQuotients}.
  \endProof

In other words, even for non-amenable groupoids, there is no distinction between $\Cgf $ and its reduced version
${\Cgr }/{\Gamma _\red }$.

Recalling that both $\Cg $ and $\Cgr $ are defined to be the completion of $\scrCcL \G $ with appropriate C*-norms,
our next result presents another compelling reason to consider the essential groupoid C*-algebra; it is related to
the results on minimal norms found in~\cite [Section~7]{SRI}.

\state Theorem \label MinNorm
  Under \ref {StandingGPDS}, assume that $\G $ is topologically free.  Then, among the collection of all
C*-seminorms on $\scrCcL \G $ coinciding with the uniform norm on $C_c(X)$, there exists a smallest member, namely
the seminorm given by
  $$
  \Vert b\Vert \mnrm = \Vert \qess (b)\Vert , \for b\in \scrCcL \G .
  $$
  The completion of\/ $\scrCcL \G $ under this seminorm is consequently the essential groupoid algebra $\Cgf $.

\Proof
  Let $\Vert {\cdot }\Vert _1$ be any C*-seminorm on $\scrCcL \G $ coinciding with the uniform norm on $C_c(X)$, and
let $C^*_1(\G ,\Lbpar )$ be the corresponding completion.  Evidently the maximum C*-seminorm is bigger than $\Vert
{\cdot }\Vert _1$, so there is a $*$-homomorphism
  $$
  \Phi :\Cg \to C^*_1(\G ,\Lbpar )
  $$
  obtained by extending the identity map on $\scrCcL \G $.  By assumption $\Phi $ is isometric on $C_c(X)$, and
hence it is injective on $C_0(X)$.  Identifying the corresponding copies of $C_0(X)$, we may therefore assume that
$\Phi $ is actually the identity map on $C_0(X)$.

Since $\G $ is a topologically free groupoid, we have that $\Incl {\Cg }{C_0(X)}$ is a topologically free inclusion
by \ref {TopoFreeFull}, so we may apply \ref {TwoQuotients} to conclude that $\Phi $ factors through the quotients
providing an isomorphism
  $$
  \Frac {\Cg }{\Gamma } \ \simeq \ \Frac {C^*_1(\G ,\Lbpar )}{\Gamma _1},
  $$
  where $\Gamma _1$ is the {gray} ideal of $C^*_1(\G ,\Lbpar )$.  Denoting by $q_1$ the quotient map modulo $\Gamma
_1$, we then have for every $b$ in $\scrCcL \G $ that
  $$
  \Vert \qess (b)\Vert = \Vert q_1(b)\Vert \leq \Vert b\Vert _1,
  $$
  concluding the proof.
  \endProof

One concrete way to express the above seminorm $\|{\cdot }\|\mnrm $ is the following variation of the expression for
the reduced norm in \ref {ReducedNorm}, in which the supremum is taken over the set $F$ of free points, rather than
over $X$.

\state Proposition
  For every $f$ in $\scrCcL \G $, one has that
  $$
  \Vert f\Vert \mnrm = \sup _{x\in F} \Vert \pi _x(f)\Vert ,
  $$
  where the $\pi _x$ are as in \ref {RegRepForx}.  The completion of\/ $\scrCcL \G $ under the seminorm given by the
right-hand-side above is consequently the essential groupoid algebra $\Cgf $.

\Proof
  For each $x$ in $X$, it is clear that $\pi _x$ is bounded relative to the maximum norm on $\scrCcL \G $, and
therefore it extends to a $*$-representation of $\Cg $, which we will still denote by $\pi _x$, by abuse of
language.

  Given $\gamma $ in $\G _x$, a simple computation shows that
  $$
  \big \langle \pi _x(f)e_\gamma , e_\gamma \big \rangle = f\big (r(\gamma )\big ),
  $$
  for all $f$ in $\scrCcL \G $, and hence also for all $f$ in $\Cg $, by continuity.  In particular
  $$
  (f^*f)\big (r(\gamma )\big ) = \big \langle \pi _x(f^*f)e_\gamma , e_\gamma \big \rangle = \Vert \pi _x(f)e_\gamma
\Vert ^2.
  $$
  It then follows from \ref {DescribeGray.ii} and the invariance of $F$, mentioned in \ref {InvarFree}, that $f$
lies in the {gray} ideal if and only if $\pi _x(f)=0$, for every $x$ in $F$.  Put another way, the {gray} ideal is
the kernel of the representation of $\Cg $ given by
  $$
  \pi = \bigoplus _{x\in F}\pi _x,
  $$
  and hence we see that $\Vert \qess (f)\Vert = \Vert \pi (f)\Vert $.  Therefore
  $$
  \Vert f\Vert \mnrm \={MinNorm}
  \Vert \qess (f)\Vert =
  \Vert \pi (f)\Vert =
  \sup _{x\in F} \Vert \pi _x(f)\Vert ,
  $$
  concluding the proof.
  \endProof

 See section \ref {AltDefs} for further discussions regarding alternate constructions of the essential groupoid
C*-algebra.

As mentioned in the preamble for this section, the essential groupoid C*-algebra will later be characterized
abstractly, so here are the crucial algebraic properties of our model inclusion:

\state Theorem \label MainOne
  Under \ref {StandingGPDS}, and assuming that $\G $ is topologically free, the essential projection $\qess $
introduced in \ref {DefineFreeGpd} is injective on $C_0(X)$, and hence $C_0(X)$ may be identified with a subalgebra
of $\Cgf $.  Once this identification is made one has that
  $$
  \Incl {\Cgf }{C_0(X)}
  $$
  is a weak Cartan inclusion, as defined in \ref {ManyTopFree.iv}.

\Proof
  Since $\G $ is topologically free, we have that
  $$
  \incl {B}{A}:= \Incl {\Cg }{C_0(X)}
  $$
  is a topologically free inclusion by \ref {TopoFreeFull}.  Therefore the inclusion in the statement is {light} by
\ref {GetUltraFree}.  As seen in \ref {CharacFreePtsTwo}, one has that $\incl {B}{A}$ is smooth, so the conclusion
follows from \ref {IItoII}.
  \endProof

As the reader might have noticed, our main interest is to study groupoids which are not necessarily Hausdorff.
However one might of course wonder what happens if a Hausdorff groupoid is plugged into our construction.

\state Proposition \label EssGpdHausdorff
  Under \ref {StandingGPDS}, assume that $\G $ is Hausdorff and topologically free.  Then, considering the
inclusions described in \ref {TwoAlgebras}, one has that
  \izitem
  \zitem the {gray} ideal of\/ $\Cg $ coincides with the kernel of the left regular representation $\Lambda $
mentioned in \ref {IntForlLeftReg},
  \zitem the {gray} ideal of\/ $\Cgr $ vanishes,
  \zitem $\Cgf \simeq \Cgr $.

\Proof
  It is well known that the restriction of every $f$ in $\Cg $ to $X$ is a continuous map.  Points (i) and (ii) may
then be deduced from \ref {DescribeGray.ii}, and the density of the set $F$ of free points.  The last point then
follows from (ii).
  \endProof

\endsection

\startsection Kwasniewski and Meyer's version of the essential groupoid C*-algebra

\label KMSection

\def \supp {\text {supp}}

\sectiontitle

Inspired by an earlier version of the present work,  Kwasniewski and Meyer \cite {KM} introduced a different notion of the  essential
C*-algebra of a groupoid which is perhaps better suited when considering groupoids with small sets of free points.    Their
original definition was given in the context of {\it actions of inverse
semigroups on C*-algebras by Hilbert bimodules}, but it may be translated in  a more pedestrian  language when we only have
\'etale groupoids in mind.  We therefore dedicate this section to compare Kwasniewski and Meyer's version of
the essential C*-algebra with ours.

\fix The standing assumption for this section is   \ref {StandingGPDS}, but we shall add that the unit space of $\G $,
which we recall is also denoted by $X$, will be assumed to be second countable.

As described after \ref {ItIsAllFunctions}, any $f$ in $\Cgr $ may be viewed as a section of $\Lb $.  However,  unless $\G $ is
Hausdorff, there is no reason why such sections are continuous.  Our first task will be to show that they are at least Borel measurable.

The concept of Borel measurability for sections of Fell line bundles over non-Hausdorff groupoids perhaps deserves a few introductory
words:
initially we
recall from \cite [II.17.37]{FD} that a Banach bundle over a locally compact Hausdorff space, whose
fibers are all of the same finite dimension, is necessarily locally trivial.

Therefore our Fell line bundle $\Lb $ falls
under the scope of the above paragraph,
and is therefore locally trivial, as long as it is restricted to an open bisection (which is a locally compact
Hausdorff space).

In the special case in which $\Lbpar |_U$ is actually trivial,  where  $U$ is a given bisection,  one may choose a nowhere vanishing continuous
section $e$, leading to the isomorphism of bundles
  $$
  (\gamma , \lambda )\in U\times {\bf C}\mapsto \lambda e(\gamma )\in \Lbpar |_U,
  $$
  sometimes called a \emph {trivialization}.
  This said, any section $f$ of $\Lbpar |_U$ is necessarily given as
  \lbldeq TrivializeF
  $$
  f(\gamma ) = g(\gamma )e(\gamma ), \for \gamma \in U,
  $$
  where $g$ is a complex valued function on $U$, which is continuous if and only if $f$ is a continuous section.

\definition \label DefineBorelMeasurable
  Let $Y$ be an open subset of $\G $ (e.g.~an open bisection or perhaps $\G $ itself) and let $f$ be a section of $\Lbpar
|_Y$.  We shall say that $f$ is
  \emph {Borel measurable}\fn {Should the reader be worried about a possible problem regarding local Borel measurability,
we recall that $X$ is assumed to be second countable, and hence so are all open bisections. As a consequence,  any
\emph {locally Borel} subset is necessarily Borel.}  if, for every $\gamma \in Y$, there
is an open bisection $U\subseteq Y$, with $\gamma \in U$, and a continuous section $e$ of $\Lbpar |_U$, such that
  $$
  f(\gamma ) = g(\gamma )e(\gamma ), \for \gamma \in U,
  $$
  where $g$ is some Borel measurable, complex valued function on $U$.

Observe that we have not required $e$ to be nowhere vanishing in \ref {DefineBorelMeasurable}.  Nevertheless, upon
suitably reducing $U$, we may assume that $\Lbpar |_U$ is trivial and hence one could pick a nowhere vanishing
continuous section $e'$, in which case
  $$
  e(\gamma ) = h(\gamma )e'(\gamma ), \for \gamma \in U,
  $$
  for some complex valued continuous function $h$.  It follows that $f(\gamma )=g(\gamma )h(\gamma )e'(\gamma )$, for
all $\gamma $ in $U$, and clearly $gh$ is a Borel measurable function.  Summarizing, upon replacing $e$ by $e'$, and $g$
by $gh$, we may assume that $e$ is nowhere vanishing in \ref {DefineBorelMeasurable}.

\state Proposition Any $f$ in $\Cgr $, when viewed as a section of $\Lb $ as described above, is bounded and Borel
measurable.

\Proof The boundedness part follows immediately from \ref {BoundedValues}.  Given any $\gamma $ in $\G $, pick an open
bisection $U$ containing $\gamma $, and such that $\Lbpar |_U$ is trivial.  Next we choose a nowhere vanishing section
$e$ of
  $\Lbpar |_U$, which we may assume satisfies $\|e(\gamma )\| = 1$, for all $\gamma $ in $U$, since otherwise we may
replace $e(\gamma )$ by $e(\gamma )/\|e(\gamma \|$.  Letting ${\cal F}(U)$ be the Banach space of all bounded, complex
valued functions on $U$, consider the map
  $$
  R: \Cgr \to {\cal F}(U),
  $$
  given by
  $$
  R(f)(\gamma ) = f(\gamma )e(\gamma )^*, \for f\in \Cgr , \for \gamma \in U.
  $$
  Notice that $R(f)(\gamma )$ is actually an element of $\Lb _{\gamma \gamma \inv } = \Lb _{r(\gamma )}$, which we
tacitly, and canonically, identify with the field of complex numbers.  Clearly $R$ is a bounded linear map.

One therefore has that
  $$
  f(\gamma ) = R(f)(\gamma )\ e(\gamma ), \for \gamma \in U,
  $$
  so our task consists in proving that $R(f)$ is Borel measurable  for all $f$.

Observing that $\Cgr $ is the closed linear span of the
union of all $\scrCcL V$,  where $V$ ranges over all open bisections,  it suffices to prove that
  $$
  R\big (\scrCcL V\big )\subseteq {\cal B}(U),
  $$
  for all such $V$'s, where ${\cal B}(U)$ denotes the closed subspace of ${\cal F}(U)$ formed by all Borel measurable
functions.

Let us therefore pick  $f\in \scrCcL V$,  for some  $V$.  In this case
notice that $R(f)$ is continuous on the open subset of $U$ given by
  $$
  W:= V\cap U,
  $$
  while it vanishes on the relative complement $U\setminus W$.  The inverse image of any open subset of ${\bf C}$ under
$R(f)$ may then be described as the union of an open subset of $W$ with or without the closed set $U\setminus W$, which
is obviously a Borel subset of $U$.  This proves that $R(f)$ is Borel measurable.  \endProof

As a special case of the above result, namely taking $Y$ to be the unit space of $\G $,  which in turn we agreed to denote by $X$
in \ref {StandingGPDS},  we get a map
  $$
  \EE : f\in \Cgr \mapsto f|_X \in {\cal B}(X).
  $$
  We leave it for the reader to prove that $\EE $ is a generalized conditional expectation, as defined in   \ref {DefineCondExp}.

Recalling from section \ref {PseudoExpSection} that ${\cal M}(X)$ denotes the ideal of ${\cal B}(X)$ formed by the functions
vanishing off a meager, let
  $$
  I\big (C_0(X)\big ) = {\cal B}(X)/{\cal M}(X).
  $$
  Incidentally we remark that   $I\big (C_0(X)\big )$ is the injective envelope of $C_0(X)$.  The composition
  \lbldeq TheEssPseudo
  $$
  E_\ess :
  \Cgr \labelarrow {\EE } {\cal B}(X) \longrightarrow I\big (C_0(X)\big ),
  $$
  where the rightmost arrow indicates the quotient map, is therefore a pseudo-expectation.

Should we find ourselves under the
special conditions outlined in \ref {UniquePseudo}, then the above is necessarily the unique pseudo-expectation but, in
the general case, uniqueness is not guaranteed.  This pseudo-expectation is nevertheless of significance to us, so we shall give it a
special name.

\definition Given any twisted \'etale groupoid $(\G ,\Lbpar )$, as in \ref {StandingGPDS}, with second countable unit
space,  the pseudo-expectati\-on $E_\ess $ described in
\ref {TheEssPseudo} will be called the \emph {essential pseudo-expectation}.

\state Proposition \label LeftKerChar
  The left kernel of the essential pseudo-expectation $E_\ess $, namely
  $$
  \Omega = \big \{f\in \Cgr : E_\ess (f^*f)=0\big \},
  $$
  consists of all elements $f$ in $\Cgr $ such that, for
every open bisection $U\subseteq \G $, the restriction $f|_U$ vanishes off a meager subset of $U$.   Moreover $\Omega $ is a
two-sided ideal of $\Cgr $.

\Proof
  Fixing   $f$  in  $\Cgr $,  we have  by the Cauchy-Schwarz inequality,  that
  $$
  E_\ess (f^*f)=0 \IFF E_\ess (g^*f)=0, \for g\in \Cgr .
  $$
  Because the union of the $\scrCcL U$ span a dense subset of $\Cgr $, we then have
that
  \lbldeq CharactOmega
  $$
  f\in \Omega \IFF E_\ess (g^*f)=0, \for g\in \textstyle \bigcup \ \scrCcL U,
  $$
  where the union ranges over all open bisections $U$.
Fixing  such a  $U$, and choosing $g$ in $\scrCcL U$,
observe that, for every $x$ in $X$, one has
  $$
  (g^*f)(x) = \sum _{\gamma \inv \gamma =x} g^*(\gamma ^{-1}) f(\gamma )   = \sum _{\gamma \inv \gamma =x} {g(\gamma )}^* f(\gamma ).
  $$
  This clearly vanishes in case $x\notin s(U)$, and in case $x\in s(U)$, we have that
  $$
  (g^*f)(x) = {g(\gamma _x)}^* f(\gamma _x),
  $$
  where $\gamma _x$ is the (necessarily unique) element in $U$ whose source coincides with $x$.
  It follows that
  \lbldeq SupInBis
  $$
  \big \{x\in X: (g^*f)(x) \neq 0\big \} = \big \{x\in s(U):  g(\gamma _x) \neq 0,\text { and } f(\gamma _x)\neq 0\big \}.
  $$
  Denoting the support (no closure) of $f$ by
  $$
  \supp (f) = \big \{\gamma \in \G : f(\gamma )\neq 0\big \},
  $$
  we see that the set referred to in \ref {SupInBis} is contained in
  $$
  s\big (U\cap \supp (f)\big ).
  $$
Now suppose $f\in \Cgr $ and that for each open bisection $U\subseteq
\G$, $f|_U$ vanishes off a meager subset of $U$.  Then $U\cap \supp
(f)$ is a meager subset of $U$, and because $s$ establishes a
homeomorphism from $U$ to $s(U)$, $s\big (U\cap \supp (f)\big )$ is a
meager subset of $s(U)$.  Consequently $E_\ess (g^*f) = 0$, by \ref
{SupInBis}.  As $g$ is arbitrary, we conclude from \ref
{CharactOmega} that $f$ lies in $\Omega $.

Conversely, pick any $f$ in $\Omega $, and let $U$ be an open bisection of $\G $.  We must therefore prove that
$U\cap \supp (f)$ is meager in $U$.

For each $\gamma $ in $U$, choose $g_\gamma $ in $\scrCcL U$, such that $g_\gamma (\gamma )\neq 0$, as well as an open neighborhood $V_\gamma $  of $\gamma $
within $U$, where $g_\gamma $ is nowhere vanishing.  Clearly $\{V_\gamma \}_{\gamma \in U}$ is an open cover for $U$.

Observing that
$X$ is second countable  (and
hence Lindel\umlaut of) by assumption, so is $U$, and then we may extract a countable subcover $\{V_{\gamma _n}\}_{n\in {\bf N}}$ of the
above cover.

In order to prove that $U\cap \supp (f)$ is meager in $U$, it therefore suffices to prove that the same is true for
$V_{\gamma _n}\cap \supp (f)$, for every $n$.  Fixing $n$, and
in order to do away with the no longer useful subscript $\gamma _n$, we set
  $$
  g := g_{\gamma _n}, \and V := V_{\gamma _n}.
  $$
  Recall that $g\in \scrCcL U$,  $V\subseteq U$,  $g$ vanishes
nowhere on $V$, and that our task is to prove that $V\cap \supp (f)$ is
meager in $U$.  By hypothesis, we have that $E_\ess (g^*f)=0$, so
  $$
  \big \{x\in X: (g^*f)(x) \neq 0\big \}
  $$
  is meager in $X$.
Using \ref {SupInBis}  we deduce that
  $$
  \big \{x\in s(U):  g(\gamma _x) \neq 0,\text { and } f(\gamma _x)\neq 0\big \}
  $$
  is also meager in $X$.  Since the correspondence  $x\mapsto \gamma _x$ establishes a homeomorphism from $s(U)$ onto $U$, we further
conclude that
  $$
  \{\gamma \in U:  g(\gamma ) \neq 0,\text { and } f(\gamma )\neq 0\big \}
  $$
  is meager in $U$.  Since the latter clearly contains $V\cap \supp (f)$,  our task is accomplished.

That $\Omega$ is a left-ideal follows from the fact that $E_\ess$ is a
positive map and that for any elements $x$ and $y$ in a C*-algebra,
$y^*x^*xy\leq \|x\|^2y^*y$.   The alternative description of $\Omega $ provided by the first part of the statement may now be used to give a
routine proof
that $\Omega $ is a self-adjoint set,  hence also a right-ideal.
\endProof

We observe that the last statement of Proposition~\ref{LeftKerChar}
also follows
from \cite [4.11]{KM}.

\definition \label KMEss
  (\cite [Definition 4.4]{KM}) Given a twisted \'etale groupoid $(\G ,\Lbpar )$, as in \ref {StandingGPDS}, with second countable unit
space,  the Kwasniewski-Meyer
essential groupoid C*-algebra of $(\G ,\Lbpar )$ is defined to be the quotient of  $\Cgr $ by the left kernel of the
essential pseudo-expectation $E_\ess $.

We remark that, for any generalized conditional expectation whose left kernel is a two sided ideal, the induced
expectation on the   quotient is faithful.  Therefore, the Kwasniewski-Meyer algebra defined above  admits a faithful
pseudo-expectation.  This fact is also a special case of  \cite [4.11]{KM}.

As discussed in the last paragraph of section \ref {CstarALgGpd}, every element of the full C*-algebra $\Cg $ may also be
viewed as a section of $\Lb $.  Regardless of the draw back that this correspondence is no longer one-to-one,  we could have
developed the whole analysis above after  replacing  $\Cgr $ by $\Cg $.  In particular, the new essential
pseudo-expectation would be none other than the composition
  $$
  \Cg \labelarrow \Lambda \Cgr \labelarrow {E_\ess } I\big (C_0(X)\big ),
  $$
  where $\Lambda $ is the regular representation.   We would then define the Kwasniewski-Meyer algebra as the quotient of
$\Cg $ by the left kernel of $E_\ess \circ \Lambda $.  However, as everything factors through the reduced algebra,  it is clear that the latter
quotient is isomorphic to the one in \ref {KMEss}.  We should also remark that the original Definition,  namely \cite [4.4]{KM}, is closer in
spirit to this second point of view.

In what follows we explore the relationship between $\Cgf $ and \ref {KMEss}.

\state Theorem \label IsoEssAlgs
  Let $(\G ,\Lbpar )$ be a twisted \'etale groupoid as in \ref {StandingGPDS}, with second countable unit
space.  Suppose that $\G $ is not only
topologically free, but that there exists a dense $G_\delta $ subset $Y$ of $X$, contained in the set $F$ of free points for $\G $ (see \ref {CharacFreePts}).
Then there is a natural isomorphism between $\Cgf $ and the
Kwasniewski-Meyer essential groupoid C*-algebra of $(\G ,\Lbpar )$.

\Proof We should notice that, as discussed in the preamble of section
\ref {PseudoExpSection}, the present hypothesis, nominally stronger than topological freeness, often comes for free in
concrete examples of topologically free groupoids, so it is perhaps not too strong a restriction.

Addressing the proof, observe that
   both algebras under consideration are obtained as quotients of $\Cgr $  (see  \ref {FullRedGray} regarding $\Cgf $), so
it suffices to prove that the ideals relative to which these quotients are taken coincide.  In case of $\Cgf $ the ideal
is the {gray} ideal $\Gamma _\red $, so we need to prove that $\Gamma _\red $  coincides with the left kernel $\Omega $ of
$E_\ess $.

On the other hand,
from \ref {FaithGray} we have that the {gray} ideal is the left kernel of the free expectation $\Pf $ defined
in \ref {DefineFreeCondExp}. So,
in order to prove the statement, it is enough to verify that
  $$
  \Pf (f) = 0 \IFF E_\ess (f)=0, \for f\in \Cgr .
  $$

  Observe further  that the free expectation is precisely the restriction operator to the set of
free points $F$ by \ref {FreeEvalIsExp}.  This said, our task may be translated by saying that, for every $f$ in $\Cgr $,
  $$
  f|_F = 0 \IFF f|_X \text { vanishes on the complement of a meager set}.
  $$
  The implication ``$\Rightarrow $'' is then evident, since the hypothesis that $F$ contains a dense $G_\delta $ set clearly implies
that the complement of $F$ in $X$ is meager.

Conversely, suppose that $M$ is a meager subset of $X$, such that $f$ vanishes on $X\setminus M$,  and set
  $$
  U= \{x\in F: f(x)\neq 0\}.
  $$
  Recalling from \ref {EbContinuous} that $f$ is necessarily continuous on $F$, we conclude that $U$ is open relative to
$F$.  Evidently $U\subseteq M$, so $U$ is also meager in $X$.
The conclusion then follows immediately from the Lemma below.
\endProof

\state Lemma Let $X$ be a Baire space and let $Y\subseteq F\subseteq
  X$, where $Y$ is a dense $G_\delta $ set in $X$.  If\/ $U$ is a
  subset of $F$ which is open relative to $F$, and meager relative to
  $X$, then $U$ is empty.

\Proof We argue by contradiction.  Assume that $U$ is not empty.
  Write $U=V\cap F$, where $V$ is an open subset of $X$, and
  let $\{C_n\}_{n\in {\bf N}}$ be a sequence of closed subsets of $X$, each having empty interior, such that $U\subseteq \bigcup _{n\in {\bf N}}C_n$.

Routine arguments show $F$ is a Baire space, and since $U$ is open
in $F$, necessarily $U$ is also Baire.  Therefore there is some
$n$ such that $C_n\cap U$ has nontrivial interior relative to $U$,
that is, there is an open subset $W$ of $X$ such that
  $$
  \emptyset \neq W\cap U  \subseteq C_n.
  $$
As subsets of $X$,  $F$ is dense, $W\cap V$ is open, and $C_n$ is closed, so
$$W\cap V\subseteq\overline{W\cap V\cap F}=\overline{W\cap U}\subseteq
C_n.$$
This contradicts the fact that $C_n$ has empty interior.
\endProof

\endsection

\startsection The relative Weyl groupoid

  \def \preG {{\cal S}}

\sectiontitle

Our description of the Weyl groupoid and twist for a given inclusion of C*-algebras will be based on \cite {Kumjian}
and \cite {RenaultCartan}, but for various reasons we need to make some small adjustments.  This section should
therefore be seen as a survey of the well-known construction, dotted with what we hope are some small
improvements. As such, we will often omit proofs of the results which are either explicitly in \cite {Kumjian} and
\cite {RenaultCartan}, or may be easily adapted from there.

The main adjustment we need to make is to base the construction of the twisted Weyl groupoid on a chosen
$*$-subsemigroup
  $
  N\subseteq \Norm BA .
  $
  The resulting Weyl groupoid and twist will therefore depend on $N$.  In fact it will become apparent that the
dependence on $N$ is an important feature of the non-Hausdorff/expectation-less situation.  One could obviously take
$N=\Norm BA $, itself, in which case one gets the well known Weyl groupoid and twist.

Another significant difference with respect the usual Weyl groupoid is that we use a different notion of germs,
based on \cite [page 140]{Paterson} (see also \cite [Section 4]{actions} and \cite [Section 3.1]{BussExel}), rather
than the more standard notion of germs used by Renault in his definition of the Weyl groupoid.  In the case of
Cartan subalgebras, or for essentially principal groupoids, there is not much difference between the two notions
because the pseudo groups involved have few fixed points.  However, since we plan to treat much more general
situations, such as inclusions in which both algebras are commutative, the classical notion of germs seems
inappropriate.

We should also say that our notion of germs appears implicitly in Kumjian's proof of \cite [Theorem 3.1]{Kumjian}.

Another slight technical difference is that, rather than building the twist, and then looking at the associated line
bundle, we feel it is simpler, and perhaps a bit more elegant, to directly construct the line bundle.  It is well
known that twists and Fell line bundles are interchangeable \cite {DKR} but, after all, it is from the line bundle,
rather than the twist, that one constructs the twisted groupoid C*-algebra.
  Should the twist ever be necessary, one may easily recover it as the corresponding circle bundle.

We should also note that the groupoid coming out of our construction is not expected to be Hausdorff.

Our assumptions throughout this section will be as follows:

\state {Standing Hypotheses V} \label StandingChooseN
  \rm
  In addition to \ref {StandingThree}, namely that $\incl {B}{A}$ is a regular inclusion with $A$ abelian, we shall fix
a generating (see \ref {DefineAdmissTotal}) $*$-subsemigroup
  $
  N\subseteq \Norm BA .
  $

The reader might want to note that until \ref {MapFromGpg}, all of our results hold under the weaker assumption that
$N$ is admissible.

In order to describe the Weyl groupoid of $\incl {B}{A}$, we start by considering the set
  $$
 \preG = \big \{(n,x)\in N\times X: \evx {n^*n}\neq 0\big \},
  $$
  where we will introduce the following equivalence relation:

\definition \label DefineEquivGerm
  Given $(n_1,x_1)$ and $(n_2,x_2)$ in $\preG $, we shall say that $(n_1,x_1) \sim (n_2,x_2)$ provided $x_1=x_2$,
and there are $a_1$ and $a_2$ in $A$, such that
  $$
  n_1a_1=n_2a_2, \and \ev {a_1}{x_1} \neq 0 \neq \ev {a_2}{x_2}.
  $$

The underlying set for the Weyl groupoid is then defined to be the quotient
  $$
  \G =\preG /{\sim }
  $$
  and we shall denote by $[n,x]$ the equivalence class of any given element $(n,x)$ in $\preG $.

Given two elements, say $[m,y]$ and $[n,x]$ in $\G $, we define their product if and only if $\beta _n(x)=y$, in
which case we set
  $$
  [m,y] [n,x] = [mn,x],
  $$
  observing that
  $$
  \evx {n^*m^*mn} \={Beta}
  \ev {m^*m}{\beta _n(x)} \, \evx {n^*n} =
  \ev {m^*m}{y} \, \evx {n^*n} \neq 0,
  $$
  by hypothesis, so that $(mn, x)$ is indeed a member of $\preG $.
  The inverse operation may then be deduced from the above, but it may also be explicitly defined by
  $$
  [n,x]\inv = [n^*,\beta _n(x)].
  $$

Given any $n$ in $N$, and an open subset $U\subseteq \Dom n$, set
  \lbldeq subbase4top
  $$
  \Omega (n,U) = \big \{[n,x]: x\in U\big \}.
  $$
  The topology on $\G $ is then defined to be the one generated by the collection of all $\Omega (n,U)$.

As it turns out, the $\Omega (n,U)$ indeed form a basis for that topology.  The crucial check supporting this
assertion is done as follows: given
  $$
  \gamma \in \Omega (n_1,U_1)\cap \Omega (n_2,U_2),
  $$
  we may write
  $$
  \gamma = [n_1,x] = [n_2,x],
  $$
  so we may pick $a_1$ and $a_2$ as in \ref {DefineEquivGerm}.  Setting $m:= n_1a_1 = n_2a_2$, we then have that
  $$
  \gamma \in \Omega (m,V)\subseteq \Omega (n_1,U_1) \cap \Omega (n_2,U_2),
  $$
  where $V$ is any open neighborhood of $x$ such that $V\subseteq U_1\cap U_2$, and neither $a_1$ nor $a_2$ vanish
on $V$.

As in \cite {Kumjian} and \cite {RenaultCartan} one may now prove that $\G $ is an \'etale groupoid such that each
$\Omega (n,U)$ is an open bisection.  However, $\G $ need not be Hausdorff, as examples show.

Given any $\gamma =[n,x]$ in $\G $, observe that
  \lbldeq UnitInG
  $$
  \gamma \inv \gamma = [n^*,\beta _n(x)]\ [n,x] = [n^*n,x].
  $$
  In the particular case that $n\in A$, that is, $\gamma =[a,x]$ with
  $a\in A$ and $\evx a\neq 0$, one checks that
  $$
  \gamma \inv \gamma = [a^*a,x]=[a,x],
  $$
  so we see that $[a,x] $ is a unit of $\G $.

For a general normalizer $n$ in $N$, since $n^*n$ lies in $A$, our calculation in \ref {UnitInG} implies that every
unit of $\G $ is of the above form, hence
  $$
  \Gz =\big \{[a,x]: a\in A, \ \evx a\neq 0\big \}.
  $$
  The correspondence
  \lbldeq XAsUnitSpace
  $$
  [a,x]\in \Gz \mapsto x\in X,
  $$
  may then be shown to be a homeomorphism, so we may identify the unit space of $\G $ with $X$.

When checking that the above map is onto $X$, we need to show that for every $x$ in $X$, there exists some $a$ in
$N\cap A$, such that $\evx a\neq 0$, and since $N$ is admissible, we may indeed find some $a$ of the form $n^*n$,
with $n$ in $N$, by \ref {AdmissFullDom}.

\medskip In order to define the relevant line bundle we again bring the set $\preG $ into play, introducing one more
equivalence relation, this time on ${\bf C}\times \preG $.

\definition
  Given $(\lambda _1, n_1,x_1)$ and $(\lambda _2, n_2,x_2)$ in ${\bf C}\times \preG $, we shall say that
  $$
  (\lambda _1, n_1,x_1) \approx (\lambda _2, n_2,x_2)
  $$
  provided $x_1=x_2$, and there are $a_1$ and $a_2$ in $A$, such that
  $$
  n_1a_1=n_2a_2, \and {\lambda _1\over \ev {a_1}{x_1}} = {\lambda _2\over \ev {a_2}{x_2}},
  $$
  where the terms in the denominator are implicitly required not to vanish.  We also put
  $$
  \Lb ={{\bf C}\times \preG \over \approx }
  $$
  and we shall denote by $[\lambda , n,x]$ the equivalence class of any given element $(\lambda , n,x)$ in ${\bf
C}\times \preG $.

Using brackets to denote equivalence classes for both ``$\sim $'' and ``$\approx $'' will cause no confusion, since
the number of coordinates inside those brackets determines which equivalence relation is under consideration.

Observing that
  $$
  (\lambda _1, n_1,x_1) \approx (\lambda _2, n_2,x_2)\ \Rightarrow \ (n_1,x_1) \sim (n_2,x_2),
  $$
  we see that the formula
  $$
  \pi \big ([\lambda , n,x]\big )= [n,x], \for (\lambda , n,x)\in {\bf C}\times \preG ,
  $$
  gives a well defined function from $\Lb $ to $\G $.  Given any $\gamma $ in $\G $, one then defines the \emph
{fiber of\/ $\Lb $ over $\gamma $}, by
  $$
  \Lb _\gamma =\pi \inv (\{\gamma \}).
  $$

In order to give a sensible description of $\Lb _\lambda $, write $\gamma =[m,x]$ and notice that $\Lb _\gamma $
consists precisely of all elements $[\lambda ,n,x]$ such that $[n,x]=[m,x]$.  Given any such element $[\lambda ,
n,x]$, we may find $a$ and $b$ in $A$, with $\ev {a}{x} \neq 0 \neq \ev {b}{x}$, and $na=mb$.  Setting
  $$
  \mu = {\lambda \evx b\over \evx a},
  $$
  it is easy to see that
  $
  [\lambda , n,x] = [\mu , m,x],
  $
  so we see that
  $$
  \Lbsubs _{[m, x]} = \big \{[\mu ,m,x]: \mu \in {\bf C}\big \}.
  $$
  In addition, the correspondence
  $$
  \mu \in {\bf C}\mapsto [\mu , m,x]\in \Lbsubs _{[m, x]}
  $$
  is easily seen to be a bijection from the ``complex line'' to $\Lbsubs _{[m, x]}$, through which we give $\Lbsubs
_{[m, x]}$ the structure of a one-dimensional complex vector space.

Let us now describe the Fell bundle structure \cite {KumjianFell} of $\Lb $.  Besides the above linear structure on
each fiber, we need multiplication operations
  $$
  \Lb _{\gamma _1}\times \Lb _{\gamma _2} \to \Lb _{\gamma _1\gamma _2},
  $$
  for every pair $(\gamma _1,\gamma _2)$ in the set $\G ^{(2)}$ of \emph {composable pairs}, as well as adjoint
operations
  $$
  *:\Lb _\gamma \to \Lb _{\gamma \inv },
  $$
  for every $\gamma $ in $\G $.  Regarding the former, write $\gamma _i=[n_i,x_i]$, so that $x_1=\beta _{n_2}(x_2)$.
Any pair of elements in $\Lb _{\gamma _1}\times \Lb _{\gamma _2}$ may then be written uniquely as
  $$
  \Big ([\lambda _1, n_1,x_1],[\lambda _2, n_2,x_2]\Big ),
  $$
  for $\lambda _1,\lambda _2\in {\bf C}$, and we then set their product to be
  $$
  [\lambda _1, n_1,x_1] \cdot [\lambda _2, n_2,x_2]= [\lambda _1\lambda _2, n_1n_2,x_2].
  $$
  Likewise, given $[\lambda , n,x] \in \Lbsubs _{[n,x]}$, we define
  $$
  [\lambda , n,x]^* = \big [\, \bar \lambda , n^*,\beta _n(x)\big ].
  $$

The uniqueness of the representation $[\lambda , n,x]$ for any given element of $\Lb _\gamma $ only stands after one
chooses to represent $\gamma $ as $[n,x]$, so all of the above operations must still be shown not to depend on the
choice of representatives.  Let us therefore verify this, at least in the stickier case, which is
  \lbldeq StickyPoint
  $$
  [\lambda _1, n_1,x_1] = [\lambda '_1, n'_1,x_1] \IMPLY [\lambda _1\lambda _2, n_1n_2,x_2] = [\lambda _1'\lambda
_2, n_1'n_2,x_2].
  $$
  In the statement of our task it is, of course, implicit that
  $$
  x_1=\beta _{n_2}(x_2) \in \Dom {n_1}\cap \Dom {n_1'}\cap \Ran {n_2},
  $$
  and the hypothesis implies that there are $a$ and $a'$ in $A$, such that
  $$
  n_1a=n'_1a', \and {\lambda _1\over \ev {a}{x_1}} = {\lambda '_1\over \ev {a'}{x_1}}.
  $$

Choosing any $v$ in $A$, vanishing off $\Ran {n_2}$, and such that $\ev v{x_1}\neq 0$, we may replace $a$ and $a'$
respectively by $av$ and $a'v$; thus we may assume that both $a$ and $a'$ lie in $C_0\big (\Ran {n_2}\big )$.
Consequently
  $$
  \pd a{n_2}, \pd {a'}{n_2} \in C_0\big (\Dom {n_2}\big )\subseteq A,
  $$
  and then
  $$
  n_1n_2\pd a{n_2}\={Intertwiner}
  n_1an_2=
  n'_1a'n_2 \={Intertwiner}
  n'_1n_2\pd {a'}{n_2},
  $$
  while
  $$
  {\lambda _1\lambda _2\over \ev {\pd a{n_2}}{x_2}} =
  {\lambda _1\lambda _2\over \ev {a}{\beta _{n_2}(x_2)}} =
  {\lambda _1\lambda _2\over \ev {a}{x_1}} =
  {\lambda '_1\lambda _2\over \ev {a'}{x_1}} =
  {\lambda _1'\lambda _2\over \ev {a'}{\beta _{n_2}(x_2)}} =
  {\lambda _1'\lambda _2\over \ev {\pd {a'}{n_2}}{x_2}},
  $$
  proving \ref {StickyPoint}.

These operations may now be shown to satisfy all of the 10 axioms listed in \cite [2.1]{KumjianFell}, except for
their topological aspects, as we have not yet equipped $\Lb $ with any topology.

Even though each fiber $\Lb _\gamma $ is a one-dimensional vector space, it does not come equipped with a canonical
basis, meaning a canonical generator.  This is in contrast to the fibers over unit elements of $\G $, which are
invariant under the multiplication and adjoint operations, and hence possess the structure of a one-dimensional
C*-algebra.  In this case $\Lb _\gamma $ does indeed possess a very canonical linear generator, namely its unit.

Since the unit in a one-dimensional C*-algebra is its only nonzero idempotent element, we may search for it among
the solutions of the equation $\xi ^2=\xi $.

Given $\gamma $ in $\Gz $, say $\gamma =[a,x]$, with $a$ in $N\cap A$ satisfying $\evx a\neq 0$, write $\xi
=[\lambda ,a,x]$, so the above equation takes the form
  $$
  [\lambda ,a,x] = [\lambda ,a,x]^2 = [\lambda ^2,a^2,x],
  $$
  This means that $(\lambda ,a,x) \approx ( \lambda ^2,a^2,x)$, so, for suitable $b$ and $c$ in $A$, we have
  $$
  ab=a^2c, \and {\lambda \over \evx b} = {\lambda ^2\over \evx c}.
  $$
  The solution may then be found by observing that $\lambda \neq 0$, so
  $$
  \lambda =
  {\evx c\over \evx b} =
  {\evx {a^2}\evx c\over \evx {a^2}\evx b} =
  {\evx {a^2c}\over \evx a\evx {ab}} =
  {1\over \evx a}.
  $$
  The conclusion is then that the unit of the one-dimensional C*-algebra $\Lbsubs _{[a,n]}$ is given by
  \lbldeq FoundUnit
  $$
  {\bf 1} =[\evx a\inv ,a,x].
  $$

  This also allows us to compute the absolute value of any given $\xi =[\lambda ,a,x]$, as follows
  \lbldeq AbsVal
  $$
  \big |[\lambda ,a,x]\big | = \Big |\lambda \evx a[\evx a\inv ,a,x]\Big | = \big |\lambda \evx a\big |\ |{\bf 1}| =
\big |\lambda \evx a\big |.
  $$
  More generally, given an arbitrary element $\gamma =[n,x]$ in $\G $, and any
  $$
  \xi =[\lambda , n,x] \in \Lbsubs _{[n, x]},
  $$
  we have
  $$
  \vert \xi \vert ^2 = \big |\xi ^*\xi \big | =
  \big |\big [\bar \lambda , n^*,\beta _n(x)\big ] \big [\lambda , n,x\big ]\big | =
  \big |\big [\bar \lambda \lambda , n^*n,x\big ]\big | \={AbsVal}
  |\lambda |^2\evx {n^*n}.
  $$
  Consequently
  \lbldeq NormInFiber
  $$
  \big \vert [\lambda , n,x]\big \vert = |\lambda |\ \evx {n^*n}^{1/2}.
  $$

To provide a topology for $\Lb $, we shall introduce one based on \cite [II.13.18]{FD} and \cite [Appendix
B]{BussMeyer}, a construction that requires a collection of local cross-sections of $\Lb $.  In order to produce
such a collection, choose $n$ in $N$, and let
  $$
  a\in C^b\big (\Dom n\big ).
  $$
  We will often choose such $a$'s by restricting some member of $A$ to $\Dom n$, but we will see that it is useful
to introduce the following for any bounded continuous function on $\Dom n$, regardless of whether or not $a$ extends
to a continuous function on $X$. A hint about the relevance of this choice is that $C^b\big (\Dom n\big )$ is the
multiplier algebra of $C_0\big (\Dom n\big )$, as in \ref {ExtendedMult}.  With the notation of~\ref {subbase4top},
consider the bisection
  \lbldeq DefineOn
  $$
  \newsymbol {\O n}{bisection determined by a normalizer} := \Omega \big (n,\Dom n\big ),
  $$
  and let
  \lbldeq OneSection
  $$
  \newsymbol {\xi _{n,a}}{canonical cross-section of Weyl line bundle} \ :\ [n,x]\in \O n\ \mapsto \ \big [a(x),
n,x\big ]\in \Lbsubs _{[n, x]}.
  $$
  (The reader is cautioned that if $\gamma \in \Omega (n)$, this formula relies on representing $\gamma $ as $\gamma
=[n,x]$ for appropriate $x\in \Dom n$.)
  We evidently have that $\pi \circ \xi _{n,a}$ is the identity map on $\O n$, so $\xi _{n,a}$ is a local
cross-section of $\Lb $.  Observe that
  \lbldeq ComputeNorm
  $$
  \big \vert \xi _{n,a}([n,x])\big \vert =
  \big \vert [a(x), n,x]\big \vert \={NormInFiber}
  |a(x)|\ \evx {n^*n}^{1/2},
  $$
  which is a continuous function of the variable $x$, and hence also of the variable $[n,x]$, due to the fact that
$\O n$ is a bisection and hence $[n, x] \leftrightarrow x$ is a homeomorphism.  We may then use \cite [II.13.18]{FD}
to give $\Lb $ a unique topology, making each $\xi _{n,a}$ a continuous local cross-section.

The careful reader will have noticed that \cite [II.13.18]{FD} assumes that the base topological space is Hausdorff,
but the proof given there makes no fundamental use of that hypothesis and hence may be adapted to our present
situation in which $\G $ is not necessarily Hausdorff.  See also \cite [Appendix B]{BussMeyer}.

One may then prove that both the multiplication and the adjoint operations are continuous with respect to this
topology, so that $\Lb $ is a bona fide Fell line bundle over $\G $.

\definition
  Let $\incl {B}{A}$ be a regular inclusion with $A$ abelian, and let $N\subseteq \Norm BA $ be a generating
$*$-semigroup.  Then the groupoid $\G $ and the line bundle $\Lb $ introduced above are called respectively the
  \newConcept {Weyl groupoid}{} and the \newConcept {Weyl line bundle}{} relative to $N$.  Should we want to
emphasize the important fact that our construction depends on the choice of $N$, we will denote them respectively by
  $$
  \newsymbol {\G _N}{relative Weyl groupoid}, \and \newsymbol {\WeylBund }{relative Weyl line bundle} .
  $$

\endsection

\startsection Fell bundles over inverse semigroups

\sectiontitle

The standing assumptions for this section will be \ref {StandingChooseN}.  Our main purpose is to find a natural
$*$-homo\-mor\-phism from $\Cw $ to $B$.  The primary tool to be used here is the theory of Fell bundles over
inverse semigroups introduced by Sieben in unpublished work and later developed in \cite {cpisg}, to which the
reader is referred to for the definition and basic concepts.  We thus begin by discussing the relevant inverse
semigroup.

Speaking of the bisections defined in \ref {DefineOn}, given $n$ and $m$ in $N$, observe that \ref {KnownFacts.iii}
implies that
  $
  \O n\O m= \O {nm},
  $
  while \ref {KnownFacts.iv} gives
  $
  \O n\inv = \O {n^*}.
  $
  Therefore the collection of all bisections
  \lbldeq introISG
  $$
  \ISG = \big \{\O n:n\in N\big \}
  $$
  forms a $*$-subsemigroup of the inverse semigroup of all open bisections of $\G $.  Consequently $\ISG $ is an
inverse semigroup in its own right.

Recall from \cite [5.4.ii]{actions} and \cite [Definition 2.14]{BussExel} that, in order to verify that $\ISG $ is
\emph {wide}, one needs to prove that $\ISG $ is a cover of $\G $, and for every $U,V$ in $\ISG $, and every $\gamma
\in U\cap V$, there must exist $W$ in $\ISG $, such that $\gamma \in W\subseteq U\cap V$.

\state Proposition \label Wide
  The inverse semigroup $\ISG $ introduced in \ref {introISG} is wide.

\Proof
  By construction, it is clear that $\ISG $ covers $\G $.  Given $n$ and $m$ in $N$, let $\gamma \in \O n\cap \O m$,
so we may write
  $$
  \gamma = [n,x] = [m,x],
  $$
  for some $x$ in $\Dom n\cap \Dom m$.  The two representatives for $\gamma $ implicit above lead to a pair of
elements $a$ and $b$ in $A$, with $na=mb$, and $\ev {a}{x} \neq 0 \neq \ev {b}{x}$.  It follows that
  $$
  \gamma = [na,x] \in \O {na} \subseteq \O n\cap \O m,
  $$
  completing the proof.
  \endProof

Given the important role played by inverse semigroups here, many ingredients of our construction will have to be
parametrized by $\ISG $, although they are currently parametrized by $N$.  We will therefore need to gain a good
understanding of the correspondence $n\to \O n$.  We thus start with the following auxiliary device.

\state Lemma \label Smallerbisection
  Given $n$ and $m$ in $N$, the following are equivalent:
  \izitem
  \zitem $\O n\subseteq \O m$,
  \zitem $\Dom n\subseteq \Dom m$, and there are $a$ and $b$ in $C_0\big (\Dom n\big ) $, both of which are nonzero
on every point of\/ $\Dom n$, and such that $na=mb$.

\Proof (ii)$\Rightarrow $(i) Obvious.

\itmImply (i) > (ii)
  Since the source of $\O n$ is $\Dom n$, and similarly for $\O m$, we have that $\Dom n\subseteq \Dom m$.
  We next claim that, for every $x$ in $\Dom n$, there are $a_x$ and $b_x$ in $C_0\big (\Dom n\big ) $, such that
  $$
  0\leq a_x\leq 1, \quad \evx {a_x}=1,
  \and na_x=mb_x.
  $$
  To prove the claim, pick any $x$ in $\Dom n$, and notice that
  $[n, x]$ lies in $\O n$, hence also in $\O m$.  The only element of $\O m$ whose source coincides with that of
$[n, x]$ is $[m,x]$, so we must necessarily have
  $$
  [n, x] = [m, x],
  $$
  and hence we may choose $a$ and $b$ in $A$, such that $\evx {a}\neq 0\neq \evx {b}$, and $na=mb$.  Upon
multiplying both $a$ and $b$ by some $v$ in $C_0(\Dom n)$ satisfying $\evx v=\evx a\inv $, we may assume that both
$a$ and $b$ lie in $C_0\big (\Dom n\big )$, and $\evx a=1$.

We next consider the subsets of $X$ given by
  $$
  V = \big \{y\in X: |\ev ay|>2/3\big \},
  \and
  K = \big \{y\in X: |\ev ay|\geq 1/3\big \},
  $$
  and we notice that $V$ is open, $K$ is compact, and
  $$
  x\in V\subseteq K\subseteq \Dom n.
  $$
  Using Tietze's extension Theorem, let $f$ be a bounded function on $X$, such that
  $$
  f(y)={1\over \ev ay}, \for y\in K,
  $$
  and choose $v$ in $A$, supported on $V$, such that $\evx v=1$, and $0\leq v\leq 1$.  Then $af$ coincides with $1$
on the support of $v$, so
  $$
  a_x:= afv=v,
  $$
  and clearly
  $$
  na_x= nafv = mbfv = mb_x,
  $$
  where $b_x=bfv$, so $a_x$ and $b_x$ form the desired pair of functions, and the claim is thus verified.
  Let us introduce the open subset of $\Dom n$ given by
  $$
  U_x=\{y\in X:\ev {a_x}y> 0\},
  $$
  which clearly contains $x$.  Evidently $\{U_x\}_{x\in \Dom n}$ is an open cover for $\Dom n$, and we claim that it
admits a countable subcover.  To see this, write
  $$
  \Dom n = \medcup _{k\in {\bf N}}\{x\in X:\evx {n^*n}\geq 1/k\},
  $$
  from where we deduce that $\Dom n$ is a countable union of compact sets.  If, for each such compact set, we select
a finite subcover of the above cover, altogether we will get the claimed countable subcover for $\Dom n$, say
$\{U_{x_k} \}_{k\in {\bf N}}$.  Setting
  $$
  r_k=\max \{\Vert a_{x_k}\Vert , \Vert b_{x_k}\Vert \},
  $$
  we finally define
  $$
  a= \sum _{k=1}^\infty {a_{x_k}\over r_k2^k},
  \and
  b= \sum _{k=1}^\infty {b_{x_k}\over r_k2^k}.
  $$
  Again we have that
  $
  na=mb,
  $
  while $a$ is strictly positive on $\Dom n$.  Regarding $b$, notice that for every $x$ in $\Dom n$ we have that
  $$
  |\evx b|^2\evx {m^*m} = \evx {b^*m^*mb} = \evx {a^*n^*na} = |\evx a|^2\evx {n^*n} >0
  $$
  so $\evx b\neq 0$, and the proof is concluded.
  \endProof

As an immediate consequence of \ref {Smallerbisection} we have:

\state Corollary \label Samebisection
  Given $n$ and $m$ in $N$, the following are equivalent:
  \izitem
  \zitem $\O n=\O m$,
  \zitem $\Dom n=\Dom m$, and there are $a$ and $b$ in $A$, neither of which vanish on $\Dom n$, such that $na=mb$.

Based on the above one could introduce an equivalence relation on $N$ by saying that $n$ and $m$ are equivalent if
and only if the equivalent conditions in \ref {Samebisection} hold.  The quotient space will then evidently be in a
one-to-one correspondence with the inverse semigroup $\ISG $ introduced in \ref {introISG}.  Referring to $\ISG $,
instead, will therefore save ourselves of the trouble of introducing a further equivalence relation to our already
long list.

Recall that for each $U$ in $\ISG $,
  $$
  \Weylbisg U,
  $$
  denotes the set of all continuous cross-sections of $\WeylBund $ defined on $U$, and vanishing at infinity.

  Speaking of the cross-section $\xi _{n,a}$ described in \ref {OneSection} in terms of a given $n$ in $N$, and a
given $a$ in $C^b\big (\Dom n\big ) $, we have that $\xi _{n,a}$ is continuous by construction.  Moreover, staring
at \ref {ComputeNorm}, while keeping in mind that $n^*n$ lies in $C_0\big (\Dom n\big )$, and $a$ lies in $C^b\big
(\Dom n\big )$, we see that
  $$
  \xi _{n,a} \in \Weylbisg {\O n}.
  $$

\state Proposition \label DenseGuy
  Given any $n$ in $N$, one has that the set
  $$
  \big \{\xi _{n,a} : a\in A\big \}
  $$
  is a dense linear subspace of $\Weylbisg {\O n}$.

\Proof
  Observing that the correspondence $a\to \xi _{n,a}$ is clearly linear, we see that the above set is a linear
subspace.  Its density then follows from \cite [II.14.1]{FD}.

Incidentally notice that, as a bisection, $\O n$ is homeomorphic to $\Dom n$, and hence a locally compact Hausdorff
space.
  \endProof

As already indicated in the preamble of the present section, besides dealing with Fell bundles over groupoids, as
defined by Kumjian in \cite {KumjianFell}, such as our Fell line bundle $\WeylBund $, we will also be working with
Fell bundles over inverse semigroups \cite {cpisg}, two concepts which are intrinsically linked to each other,
although formally rather distinct.

The terminology ``Fell bundle'' might therefore be a bit confusing, since it refers to these two distinct concepts,
but in this section we will only be dealing with one example of each kind, referring to the first one as \emph {the
Fell line bundle $\WeylBund $ over $\G _N$}, as we have been doing, while we will now introduce the second Fell
bundle of interest, this time over an inverse semigroup.

The inverse semigroup to be used is precisely the semigroup $\ISG $ introduced in \ref {introISG}, and the Fell
bundle over $\ISG $ we have in mind is given by
  \lbldeq IntroFBdlISG
  $$
  \big \{\Weylbisg U \big \}_{U\in \ISG }
  $$
  The structural operations required by \cite [Definition 2.1]{cpisg} may be obtained by viewing each $\Weylbisg U$
as a subset of the twisted groupoid C*-algebra $C^*(\G _N,\WeylBund )$, and borrowing the multiplication operation
  $$
  \Weylbisg U\times \Weylbisg V \to C_0(UV,\WeylBund ),
  $$
  and the adjoint operation
  $$
  * : \Weylbisg U \to C_0(U\inv ,\WeylBund )
  $$
  from $C^*(\G _N,\WeylBund )$.  Regarding the maps $j_{s,t}$ required by \cite [Definition 2.1]{cpisg}, whenever
two given $U$ and $V$ in $\ISG $ satisfy $U\leq V$, namely when $U\subseteq V$, it is clear that
  $$
  \Weylbisg U \subseteq \Weylbisg V,
  $$
  so we take $j_{V,U}$ to be simply the inclusion map.

In order to get a glimpse into the structure of this object, it pays to analyze the behavior of the product and
adjoint operations on the basic local cross-sections of \ref {OneSection}.

\state Proposition \label OperSect
  Let
  $$
  n,m\in N, \quad a\in C^b\big (\Dom n\big ), \and b\in C^b\big (\Dom m\big ).
  $$
  Then, noting that $\pd bna$ (see
  \ref {PartialCompositionUnderlined}) is well defined, bounded and continuous on $\Dom {mn}$, we have
  \izitem
  \zitem $\xi _{m,b} \ \xi _{n,a} = \xi _{mn,\pd bna}$, and
  \zitem $\xi _{n,a}^* = \xi _{n^*, \pd {\bar a}{n^*}}$.

\Proof
  Since $\xi _{m,b}$ is supported in $\O m$ and $\xi _{n,a}$ is supported in $\O n$, their product is supported in
$\O {mn}$.  Picking any $\gamma $ in $\O {mn}$, we may then write $\gamma =[mn,x]$, for some $x$ in $\Dom {mn}$, so
that $x$ lies in $\Dom n$, while $\beta _n(x)$ is in $\Dom m$.  We then have
  \lbldeq Breque
  $$
  (\xi _{m,b} \ \xi _{n,a}) (\gamma ) =
  \sum _{\gamma _1\gamma _2=\gamma }\xi _{m,b}(\gamma _1) \xi _{n,a} (\gamma _2) = \cdots
  $$
  There is only one meaningful way to factor $\gamma $ as $\gamma _1\gamma _2$, i.e., so that $\xi _{m,b}(\gamma
_1)$ and $\xi _{n,a} (\gamma _2)$ are both liable to be nonzero, which is
  $$
  \gamma =[mn,x]= [m,\beta _n(x)] [n, x],
  $$
  so we have that \ref {Breque} equals
  $$
  \cdots = \xi _{m,b}\big ([m,\beta _n(x)]\big ) \xi _{n,a} \big ([n, x]\big ) =
  \big [b(\beta _n(x)), m,\beta _n(x) \big ] \big [a(x), n,x\big ] \quebra =
  \big [(\pd bna)(x), mn,x\big ] = \xi _{mn,\pd bna}(\gamma ),
  $$
  proving (i).  Regarding (ii) observe that both $\xi _{n,a}^*$ and $\xi _{n^*,a^*}$ are supported on $\O {n^*}$, so
let us pick $\gamma $ in $\O {n^*}$, necessarily of the form
  $\gamma =[n^*,x]$, where $x\in \Dom {n^*}$.  Then
  $$
  \xi _{n,a}^* (\gamma ) =
  \xi _{n,a} (\gamma \inv )^* =
  \xi _{n,a} ([n, \beta _{n^*}(x)])^* \quebra =
  [a(\beta _{n^*}(x)), n, \beta _{n^*}(x)]^* =
  [\pd {\bar a}{n^*}(x), n^*,x] =
  \xi _{n^*, \pd {\bar a}{n^*}}[n^*,x].
  \closeProof
  $$
  \endProof

We eventually want to be able to describe our originally given inclusion $\incl {B}{A}$ in terms of the pair $(\G
_N, \WeylBund )$.  A first point of contact between these is given by the following:

\state Proposition \label BasePrepFB
  Given any $n$ in $N$, there is a linear isometry
  $$
  \rho _n: \Weylbisg {\O n} \to B,
  $$
  whose range coincides with $\clsr {nA}$, and such that $\rho _n(\xi _{n,a}) = na$, for all $a$ in $A$.

\Proof
  Denote by $C_n$ the linear subspace of $\Weylbisg {\O n}$ described in \ref {DenseGuy}.  For any given $a$ in $A$,
notice that
  $$
  \big \Vert \xi _{n,a}\big \Vert =
  \sup _{x\in \Dom n}\big \vert \xi _{n,a}([n, x])\big \vert \={ComputeNorm}
  \sup _{x\in \Dom n}|\evx a|\ \evx {n^*n}^{1/2} \quebra =
  \sup _{x\in X}\evx {a^*n^*na}^{1/2} = \Vert a^*n^*na\Vert ^{1/2} = \Vert na\Vert .
  $$
  We therefore get a well defined linear isometric map
  $$
  \rho _n: C_n\to B,
  $$
  by putting $\rho _n(\xi _{n,a}) = na$, for all $a$ in $A$.  Since $C_n$ is dense in $\Weylbisg {\O n}$ by \ref
{DenseGuy}, we may therefore extend $\rho _n$ to $\Weylbisg {\O n}$, obtaining the required map.
  \endProof

We would now like to use the above $\rho _n$ to build a representation of the Fell bundle in \ref {IntroFBdlISG}.
We thus need to show that the correspondence $\O n\to \rho _n$ is well defined.

\state Lemma \label XisRestricted
  Let $n, m\in N$ be such that $\O n\subseteq \O m$, and choose $a,b\in A$ as in \ref {Smallerbisection}, so that
neither $a$ nor $b$ vanish on $\Dom n$, and $na=mb$.  Then, for every $c$ in $C_c\big (\Dom n\big )$, one has that
  $$
  \xi _{n,c}= \xi _{m,c'},
  $$
  where
  $
  c'=\ds {bc\over a}.
  $

\Proof
  Observe that the above expression for $c'$ gives a well defined member of $C_c\big (\Dom n\big )$ because $a$ does
not vanish on $\Dom n$, and $c$ is compactly supported.  Given any element of $\O n$, write it as $[n,x]$, for some
$x$ in $\Dom n$, and observe that $[n,x]=[m,x]$.  Then
  $$
  \xi _{n,c}([n,x]) =
  [\evx {c}, n,x] \explica {(!)}=
  [\evx {c'}, m,x] =
  \xi _{m,c'}([m,x]) =
  \xi _{m,c'}([n,x]),
  $$
  where the equality marked with (!) is justified by the fact that
  $$
  na=mb, \and {\evx a\over \evx a} = {\evx {c'}\over \evx b}.
  $$
  This shows that $\xi _{n,c} = \xi _{m,c'}$ on $\O n$.  For $\gamma $ in $\O m\setminus \O n$, one has that $\xi
_{n,c}(\gamma )$ clearly vanishes, and we claim that also $\xi _{m,c'}(\gamma )=0$.  To see this, write $\gamma
=[m,x]$, for some $x$ in $\Dom m$.  Then $x$ cannot be in $\Dom n$, or otherwise $[m,x]=[n,x]\in \O n$.  So
  $$
  \vert \xi _{m,c'}(\gamma )\vert =
  \vert \xi _{m,c'}([m, x])\vert \={ComputeNorm}
  |c'(x)|\evx {m^*m}|^{1/2} = 0.
  $$
  This concludes the proof.
  \endProof

\state Proposition \label SmallerRhos
  Let $n$ and $m$ be normalizers in $N$ such that
  $
  \O n\subseteq \O m.
  $
  Then the restriction of $\rho _m$ to $C_0\big (\O n,\WeylBund \big )$ coincides with $\rho _n$.

\Proof
  In order to prove that $\rho _n$ and $\rho _m$ agree on $C_0\big (\O n,\WeylBund \big )$, it is enough to prove
that they coincide on the dense subset given by \ref {DenseGuy}, namely to prove that
  \lbldeq SameRhos
  $$
  \rho _n(\xi _{n, c}) = \rho _m(\xi _{n, c}),
  $$
  for every $c$ in $A$.  If $\{u_i\}_i$ is an approximate unit for $C_0\big (\Dom n\big )$ contained in $C_c\big
(\Dom n\big )$, we claim that
  $$
  \li i \xi _{n, cu_i } = \xi _{n, c}.
  $$
  To see this we compute
  $$
  \Vert \xi _{n, cu_i } - \xi _{n, c}\Vert ^2 \={BasePrepFB}
  \Vert ncu_i -nc\Vert ^2 =
  \Vert n(cu_i -c) \Vert ^2 \quebra =
  \Vert (cu_i -c)^*n^* n(cu_i -c) \Vert \leq
  \Vert cu_i -c\Vert \Vert n^* ncu_i -n^* nc) \Vert ,
  $$
  which converges to zero since $n^* nc$ belongs to $C_0\big (\Dom n\big )$.  With this we see that it is enough to
prove \ref {SameRhos} for $c$ in $C_c\big (\Dom n\big )$.

Using that $\O n\subseteq \O m$, we next choose $a$ and $b$ in $C_0\big (\Dom n\big ) $ as in \ref
{Smallerbisection}, so that neither $a$ nor $b$ vanish on $\Dom n$, and $na=mb$.  Defining
  $$
  c' = {bc\over a},
  $$
  we have that $c' $ lies in $C_c\big (\Dom n\big )$, and from \ref {XisRestricted} it follows that
  $\xi _{n,c} = \xi _{m,c'}$, so
  $$
  \rho _m(\xi _{n, c}) =
  \rho _m(\xi _{m, c' }) = mc' =
  m\ {bc\over a} =
  nc= \rho _n(\xi _{n, c}),
  $$
  verifying \ref {SameRhos} and therefore concluding the proof.
  \endProof

As a consequence of \ref {SmallerRhos} we have:

\state Corollary
  If
  $n,m\in N$ are such that
  $
  \O n= \O m,
  $
  then $\rho _m=\rho _n$.

In view of the above result, for every $U$ in $\ISG $, we may unambiguously define
  $$
  \rho _U = \rho _n,
  $$
  where $n$ is any element in $N$ chosen such that $U=\O n$.

\state Proposition
  The collection of maps
  $$
  \Big \{\rho _U: \Weylbisg U\to B\Big \}_{U\in \ISG }
  $$
  forms a representation of the Fell bundle introduced in \ref {IntroFBdlISG} within the C*-algebra $B$.

\Proof
  According to \cite [3.1]{cpisg} we need to show that, for every $U$ and $V$ in $\ISG $, every $\xi $ in
$C_0(U,\WeylBund )$, and every $\eta $ in $C_0(V,\WeylBund )$, one has
  \iaitem
  \aitem $\rho _U(\xi ) \rho _V(\eta ) = \rho _{UV}(\xi \eta )$,
  \aitem $\rho _{U}(\xi )^* = \rho _{U^*}(\xi ^*)$,
  \aitem if $V\subseteq U$, then the restriction of $\rho _{U}$ to $\Weylbisg V$ coincides with $\rho _V$.

\bigskip \noindent Writing $U=\O m$ and $V=\O n$, with $n$ and $m$ in $N$, let us first prove (a).  Using \ref
{DenseGuy} it suffices to take
  $\xi =\xi _{m, b}$ and
  $\eta =\xi _{n, a}$, with $a$ and $b$ in $A$.
  Once these choices are made we have
  $$
  \rho _{UV}(\xi \eta ) = \rho _{mn}(\xi _{m, b}\xi _{n, a}) \={OperSect.i}
  \rho _{mn}(\xi _{mn,\pd bna}) =
  mn\pd bna\={IntertwinerTwo}
  mbna\quebra =
  \rho _m(\xi _{m, b}) \rho _n(\xi _{n, a}) =
  \rho _U(\xi ) \rho _V(\eta ).
  $$
  Leaving (b) as an exercise, we note that (c) has been proved in \ref {SmallerRhos}.
  \endProof

The main result of this section may now be proven.  Even though we have so far only used that $N$ is admissible, it
will now be crucial that $N$ is also generating.

\state Theorem \label MapFromGpg
  Let $\incl {B}{A}$ be a regular inclusion with $A$ abelian, and let $N\subseteq \Norm BA $ be a generating
$*$-subsemigroup.  Then there exists a surjective $*$-homomorphism
  $$
  \Phi : \Cw \to B,
  $$
  restricting to the identity map from
  $$
  C_0 \big (\G _N^{(0)} \big ) = C_0(X) \ \ \longrightarrow \ \ C_0(X) = A,
  $$
  and such that,
  $$
  \Phi (\xi _{n,a}) =na, \for n\in N, \for a\in A,
  $$
  where the local cross-section $\xi _{n,a}$ is given by \ref {OneSection}.

\Proof \def \CstarB {C^*(\Bun \, )}\relax
  Denoting by $\Bun $ the Fell bundle introduced in \ref {IntroFBdlISG}, and by $\CstarB $ its cross-sectional
C*-algebra \cite [3.4]{cpisg}, we may use the universal property \cite [3.5]{cpisg} of $\CstarB $ to obtain a
$*$-homomorphism
  $$
  \Phi : \CstarB \to B,
  $$
  collectively extending all of the $\rho _U$. Since $\ISG $ is wide by \ref {Wide}, we may apply \cite [Corollary
5.6]{BussMeyer} to conclude that $\CstarB $ is naturally isomorphic to $\Cw $, from where the existence of the
required map follows.  Since $N$ is generating, we have that $\NoA $ spans a dense subspace of $B$, so we see that
$\Phi $ is surjective.
  \endProof

\endsection

\startsection Topological freeness of the Weyl groupoid and the main Theorem

\sectiontitle

In this section we continue adopting \ref {StandingChooseN}.  The existence of the map $\Phi $ provided by \ref
{MapFromGpg} is the first step in the search for a groupoid model for the inclusion $\incl {B}{A}$.  However we
should worry that perhaps the null space of $\Phi $ might be to big for us to claim victory.  For example, if
$A={\bf C}$, then any unitary element of $B$ is a normalizer, and if we choose to build the relative Weyl groupoid
out of $N=\Norm BA $, the resulting groupoid will be the whole unitary group of $B$ equipped with the discrete
topology.  Although the map $\Phi $ above is still available, the difference between $\Cw $ and $B$ will be too big,
invalidating any claims that the former could be used to describe the latter.  The key to avoid an excessively large
kernel is to show that the Weyl groupoid is topologically free, and this will be the first medium term goal of this
section.

\definition \label DefineRelTriv
  Given $n$ in $N$, we will say that a given point $x$ in $\Dom n$ is \newConcept {trivial relative to $n$}{trivial
point relative to a normalizer}, if there exists
  $v$ in $A$ such that $\evx v=1$, and $nv\in A$.  The set of all such points will be denoted by $\newsymbol
{T_n}{set of trivial points relative to $n$}$.

For $x$ in $\Dom n$, observe that
  $x$ is trivial relative to $n$
  if and only if the germ
  $[n,x]$ is a unit in the Weyl groupoid.

\medskip The next result is the technical basis for proving the topological freeness of $\G _N$.

\state Lemma \label TrivialPts
  For every $n$ in $\Norm BA$ one has that
  \izitem
  \zitem $T_n$ is open,
  \zitem $T_n\subseteq F_n\cap \Fix n$,
  \zitem If $\incl BA$ is a {light} inclusion, then $T_n=\int F_n\cap \iFix n$,
  \medskip \noindent where $\int F_n$ refers to the interior of $F_n$, and $\iFix n$ is the interior of $\Fix n$.

\Proof (i)\enspace Obvious.

\itmProof (ii) Given $x$ in $T_n$, choose $v$ in $A$ such that $\evx v=1$, and $nv\in A$.  For any two states $\psi
_1$ and $\psi _2$ on $B$ extending $\varphi _x$, we then have that
  $$
  \psi _1(n) \={CondexForState} \Frac {\psi _1(nv)}{\evx v} =
  \Frac {\varphi _x(nv)}{\evx v} =
  \Frac {\psi _2(nv)}{\evx v} \={CondexForState}
  \psi _2(n),
  $$
  so $x\in F_n$.
  Moreover,
  $$
  \beta _n(x) = \beta _n\big (\beta _v(x)\big )=\beta _{nv}(x) = x,
  $$
  so $x\in \Fix n$.  This concludes the proof of (ii).

\itmProof (iii)
  By (i) and (ii) it is clear that
  $$
  T_n\subseteq \int {F_n}\cap \iFix n,
  $$
  so we next set out to prove the reverse inclusion.
  Given
  $x$ in $\int {F_n}\cap \iFix n$, we may choose an open set $W\subseteq X$ with
  $$
  x\in W\subseteq F_n\cap \Fix n.
  $$

  Picking any $v$ in $A$ vanishing off $W$, and such that $\evx v=1$, set $m=nv$. Then
  $$
  \beta _{m} = \beta _{nv} = \beta _{n} \beta _{v} = \beta _{n} |_{\Dom v},
  $$
  so $\beta _m$ is the identity map and we have by \ref {BetaIdentity} that $m$ lies in $A'$.

Regarding the function $\varepsilon _m$ introduced in \ref {DefineRelFree}, observe that it is only defined on
$F_m$, but let us instead consider the complex valued function $g$ defined on the whole of $X$ by setting it to
coincide with $\varepsilon _m$ on $F_m$, and with zero on the complement of $F_m$.  We then claim that $g$ lies in
$C_0(X)$.  To see this observe that, for every $y$ in $F_m$, one has that
  \lbldeq PLessV
  $$
  |\varepsilon _m(y)| =
  |\varepsilon _{nv} (y) | \={UniqueForBA.iii}
  |\varepsilon _n(y) \ev vy| \leq
  \Vert n\Vert \, | \ev vy|.
  $$
  This implies that $\varepsilon _m$ vanishes on $F_m\setminus W$, and hence that $g$ vanishes on $X\setminus W$.

To analyze the behavior of $g$ within $W$, note that
  $$
  W\subseteq F_n\explain {UniqueForBA.ii}\subseteq F_{nv} =F_m,
  $$
  so we see that $g$ is continuous on $W$ by \ref {EbContinuous}.  In addition, notice that the restriction of $g$
to $W$ actually lies in $C_0(W)$ by \ref {PLessV}.  This said, it is now easy to see that $g$ is in $C_0(X)$, as
claimed.  The crucial property of $g$ that will be relevant to us is that
  \lbldeq PropGCricial
  $$
  g(y) = \varepsilon _m(y), \for y\in F_m.
  $$

We next claim that $m-g$ belongs to the {gray} ideal of $B$. To see this, notice that given $y\in F$, we clearly
have that also $y\in F_m$, whence
  $$
  \E yy;(m-g) = \E y;(m-g) \={ThreeValues} \varepsilon _m( y) -g(y) \={PropGCricial} 0.
  $$
  Since both $m$ and $g$ lie in $A'$, we have for every $y\neq z$ in $F$, that
  $$
  \E yz;(m)=0=
  \E yz;(g),
  $$
  by \ref {DiagoComuta}, so in fact
  $$
  \E yz;(m-g)=0, \for y, z\in F,
  $$
  regardless of whether $y=z$ or not.  Consequently $m-g\in \Gamma $, by \ref {IntroGray.iii}, as claimed.  Since we
are assuming that the {gray} ideal $\Gamma $ vanishes, we then conclude that $m=g$, so
  $$
  nv=m=g\in A,
  $$
  thus proving that $x$ is trivial relative to $n$.
  \endProof

The following result is the motivation for the first part of the title of the present section:

\state Proposition \label WeylTopFree
  Under \ref {StandingChooseN} suppose that:
  \izitem
  \zitem $\incl BA$ is a {light} inclusion,
  \zitem $N$ is countable,
  \zitem every element of $N$ is smooth.
  \medskip \noindent Then the relative Weyl groupoid $\G _N$ is topologically free.

\Proof
  For every $n$ in $N$ we have by the definition of smoothness that $\int {F_n}$ is dense in $X$.  Since $N$ is
countable, Baire's Theorem says that
  $$
  F' := \medcap _{n\in N} \int {F_n},
  $$
  is dense in $X$.  The proof will therefore be concluded once we show that the isotropy group $\G _N(x)$ is trivial
for every $x$ in $F'$.  Fixing $x$ in $F'$, let us pick any $\gamma $ in $\G _N(x)$, which we write as $[n,x]$,
where $n$ is in $N$, and $x\in \Dom n$.  We then have that the range of $[n,x]$ coincides with its source, which
translates into saying that $\beta _n(x)=x$, so $x\in \Fix n$.  Moreover
  $$
  x\in
  F_n\explain {TwoExtensions}\subseteq
  X\setminus \bd n,
  $$
  so $x$ is necessarily in $\iFix n$.  Consequently
  $$
  x\in \int {F_n}\cap \iFix n \= {TrivialPts.iii} T_n,
  $$
  so we see that
  $$
  \gamma = [n,x] \in \G _N^{(0)} .
  $$
  This shows that $\G _N(x)$ is trivial for every $x$ in $F'$, so $\G _N$ is topologically free.
  \endProof

The promised characterization of inclusions involving essential groupoid C*-algebras is one of the main results of
this paper.

\state Theorem \label MainTwo
  Let $\incl {B}{A}$ be a weak Cartan inclusion, with $B$ separable and $A$ abelian, the spectrum of the latter
denoted by $X$.  Choose a countable, generating subsemigroup $N\subseteq \Norm BA $ consisting of smooth
normalizers, and let $\G _N$ and $\WeylBund $ be the associated relative Weyl groupoid and line bundle,
respectively.  Then there exists a $*$-isomorphism
  $$
  \Psi :\EssAlg {\G _N}{\WeylBund } \to B
  $$
  carrying $C_0(X)$ onto $A$.

\Proof
  We should first observe that, since $\incl {B}{A}$ is a smooth inclusion, there exists a generating subsemigroup
$M\subseteq \Norm BA $ consisting of smooth normalizers.  Since $B$ is separable, then so is $M$.  Taking any
countable dense subset of $M$, and letting $N$ be the $*$-semigroup generated by this set, we get a semigroup with
the properties listed in the statement.
  Employing \ref {WeylTopFree} we then deduce that $\G _N$ is topologically free, so the inclusion
  $$
  \Incl {\Cw }{C_0(X)}
  $$
  is topologically free by \ref {TopoFreeFull}.
  Considering the $*$-epimomorphism
  $$
  \Phi : \Cw \to B
  $$
  given by \ref {MapFromGpg}, which restricts to the identity on the corresponding copies of $C_0(X)$, we may then
apply \ref {TwoQuotients} which says that $\Phi $ becomes an isomorphism after it is factored through the quotients
by the corresponding {gray} ideals.  In the case of $\Cw $, the quotient is $\EssAlg {\G _N}{\WeylBund }$, by
definition, and in the case of $B$, the quotient is $B$ itself, by hypotheses.  This concludes the proof.
  \endProof

As a consequence of the last result and \ref {MainOne}, we obtain:

\state Corollary
  The class of all weak Cartan inclusions $\incl {B}{A}$, with $B$ separable and $A$ abelian, coincides with the
class of all inclusions
  $$
  \Incl {\Cgf }{C_0(X)},
  $$
  with $\G $ second countable and topologically free.

Besides providing the model for separable weak Cartan inclusions, there are other strong reasons to consider
essential groupoid C*-algebras.  Among them we point out that $\Cgf $ does not suffer from many of the problems
plaguing the classical reduced C*-algebras of non-Hausdorff groupoid, such as those described in \cite {ExelNHaus}.

\state Theorem \label GpdSimplicityCriteria
  Assume that $\G $ is topologically free.  Then,
  \izitem
  \zitem every nonzero ideal of\/ $\Cgf $ has a nonzero intersection with $C_0(X)$,
  \zitem $\Cgf $ is simple if and only if $\G $ is minimal.

\Proof
  The first point is an immediate consequence of \ref {MainOne}, which says that our inclusion is weak Cartan, and
hence also {light}, together with \ref {SummarySimplicity.i}.

Regarding (ii), it is well known that the minimality of $\G $ is equivalent to $C_0(X)$ being $\Cg $-simple, hence
it is also equivalent to $C_0(X)$ being $\Cgf $-simple by \ref {FunctorBSimple.c\&d}.  The conclusion then follows
from \ref {SummarySimplicity.ii}
  \endProof

We remark that a purely algebraic version of the above Theorem has been obtained by Nekrashevych in \cite
[Proposition 4.1]{Nekr}.

\endsection

\startsection {Semi-masa}s

\def \NdotA {N{\cdot }A}

\sectiontitle

The reader might have noticed that the notion of maximal abelian subalgebras, which is absolutely central in the
study of effective Hausdorff groupoids or Cartan subalgebras, has been conspicuously absent until now.  Indeed we
have so far based our development on the notion of topological freeness and this has proved a fruitful approach.
However, should one insist in uncovering the role of masas in the present context, there are many relevant things to
be said.

In this section we will return to \ref {StandingThree}, namely we will assume that $\incl {B}{A}$ is a regular
inclusion, with $A$ abelian.

\definition
  \izitem
  \zitem
  We will say that $A$ is a \emph {maximal abelian} subalgebra of $B$, \newConcept {masa}{} for short, if the \emph
{relative commutant}
  $$
  A' = \{b\in B: ba=ab, \ \forall a\in A\},
  $$
  coincides with $A$.
  \lbldzitem DefNearMasa
  When there exists a generating $*$-subsemigroup $N\subseteq \Norm BA $, such that
  $$
  (\NdotA ) \cap A'\subseteq A,
  $$
  (please see footnote\fn { In all cases so far considered, the product of two subsets of a C*-algebra is meant to
be interpreted as the \emph {closed linear span} of the set of individual products.  We are therefore using the
special notation $\NdotA $, above, to emphasize that this convention is not being used here.  In other words,
$\NdotA = \{na: n\in N, a\in A\}$.})  we will say that $A$ is a \newConcept {{semi-masa} relative to
$N$}{{semi-masa} relative to a $*$-semigroup of normalizers}{}, or simply a \emph {{semi-masa}} if $N$ is
understood.

Observe that if $A$ is a masa, then it is also a {semi-masa} relative to $\Norm BA $.  In other words, the property
of being a {semi-masa} is indeed weaker than being a masa.

Given $n$ in $N$, recall from \ref {DefineRelFree} that $F_n$ denotes the set of all free points relative to $n$,
while $T_n$ refers to the set of trivial points relative to $n$, as defined in \ref {DefineRelTriv}.

\state Lemma \label IntroNu
  Suppose that $A$ is a {semi-masa} relative to a given semigroup $N$ of normalizers.  Given $n$ in $N$ and $a$ in
$A$, put $m=na$.  Then
  \izitem
  \zitem if $\beta _m$ coincides with the identity on $\Dom m$, then $m\in A$,
  \zitem For every $x$ in $\iFix m$, there exists $v$ in $A$, such that $\evx v=1$, and $mv\in A$,
  \zitem $\iFix m = T_m$,
  \lbldzitem FreeAllMinusBd $\bd m= X\setminus F_m$.

\Proof (i)\enspace
  Assuming that $\beta _m$ coincides with the identity on $\Dom m$, we have that $m$ lies in $A'$ by \ref
{BetaIdentity}.  Therefore
  $
  m\in ( \NdotA )\cap A' \subseteq A.
  $

\itmProof (ii) Given $x$ in $\iFix m$, choose $v$ in $A$ vanishing off $\Fix m$, and such that $\evx v=1$.  We then
have that
  $$
  \beta _{mv} =
  \beta _m\circ \beta _v=
  \beta _m\circ \hbox {id}_{\Dom v} =
  \beta _m|_{\Dom v} =
  \hbox {id}_{\Dom v}.
  $$
  Noticing that $mv=nav\in \NdotA $, it then follows from (i) that $mv\in A$.

\itmProof (iii)
  By \ref {TrivialPts.ii} we have that $T_m\subseteq \Fix m$, and hence also that $T_m\subseteq \iFix m$, because
the former is open.  The reverse inclusion then follows from (ii).

\itmProof (iv)
  We have already seen in \ref {TwoExtensions.ii} that
  $
  \bd m\subseteq X\setminus F_m.
  $
  In order to prove the reverse inclusion, pick any $x$ in $X\setminus F_m$.  Using \ref {MandatoryValues.i}, we
have that $x$ lies in $\Fix m$, so it is either in $\iFix m$ or in $\bd m$, and our task is thus to rule out the
first alternative.  But this is easy because
  $$
  \iFix m \explica {(iii)}= T_m \explain {TrivialPts.ii}\subseteq F_m.
  \closeProof
  $$
  \endProof

As a consequence we have:

\state Proposition \label NearMasaSmooth
  If $A$ is a {semi-masa} relative to a given semigroup $N\subseteq \Norm BA $, then every $n$ in $N$ is smooth.

\Proof
  By \ref {FreeAllMinusBd} we have that
  $F_n= X\setminus \bd n$, and since $\bd n$ is a nowhere dense set, the interior of $F_n$ is dense, thus proving
$n$ to be smooth.
  \endProof

From our point of view, the main connection between {semi-masa}s and the theory we have been developing in the
previous sections is as follows:

\state Theorem
  Let $\incl {B}{A}$ be a {light} inclusion such that $A$ is a {semi-masa}.  Then $\incl {B}{A}$ is weak Cartan.

\Proof
  From the previous result it follows that our inclusion is also smooth.
  \endProof

Under the above assumptions, if $B$ is also separable, one may then apply our main Theorem \ref {MainTwo} in order
to describe $B$ as an essential groupoid C*-algebra.  However, there are many inclusions for which our main Theorem
applies, even though the {semi-masa} property fails.  The most striking such examples are those for which $B$ is
commutative, as we shall see in section \ref {PerFunSec}, below.

\medskip For second countable, Hausdorff groupoids, topological freeness \ref {DefTopoFreeGpd} and effectiveness
\ref {EffectiveGpd} are equivalent concepts \cite [Proposition 3.1]{Renault}\fn {The reader should be warned that
the expression \emph {essentially principal} used in \cite {Renault} should be taken to mean \emph {effective},
according to \ref {EffectiveGpd}.}  and in fact they are also equivalent to the fact that $C_0(X)$ is a masa in
$\Cgr $ \cite [II.4.7]{Renault}.  However, for non-Hausdorff groupoids things become a lot more complicated and we
would now like to discuss the relationship between these and several related properties.  We believe condition \ref
{WeakEffective}, below, has first been considered in \cite {BCFS}.

\state Theorem \label CrazyGroupoids
  Under \ref {StandingGPDS}, suppose that $\G $ is second countable, and consider the following statements:
  \izitem
  \lbldzitem Masa
    $C_0(X)$ is a masa in $\Cgr $,
  \lbldzitem Effective
    $\G $ is effective,
  \lbldzitem WeakMasaNg
    $C_0(X)$ is a {semi-masa} relative to $N_\G $,
  \lbldzitem WeakMasa
    $C_0(X)$ is a {semi-masa} (relative to some $*$-semigroup of normalizers),
  \lbldzitem WeakEffective
    $\G '\setminus \Gz $ has empty interior,
  \lbldzitem TopoPrin
    $\G $ is topologically free.
  \smallskip \noindent Then
    \refl {Masa} $\Rightarrow $
    \refl {Effective} $\Leftrightarrow $
    \refl {WeakMasaNg} $\Rightarrow $
    \refl {WeakMasa} $\Rightarrow $
    \refl {WeakEffective} $\Leftrightarrow $
    \refl {TopoPrin}.
  If $\G $ is moreover assumed to be Hausdorff, then all of the above statements are equivalent to each other.

\Proof Throughout this proof we will let $B=\Cgr $, and $A=C_0(X)$, whenever this helps to lighten up notation.

\itemImply {\refl {Masa}}{\refl {Effective}}
  Assuming that $\G $ is not effective, there exists some $\gamma $ in the interior of $\G '$ which is not in $X$.
So we may pick an open bisection $U$ such that $\gamma \in U\subseteq \G '$.  We may then choose $n$ in $\scrCcL U$
such that $n(\gamma )\neq 0$.  Since $U$ is contained in $\G '$, we have that $\beta _n$ is the identity on its
domain, whence $n$ lies in $A'$ by \ref {BetaIdentity}.  But $n\notin A$ because $n(\gamma )\neq 0$.  Therefore,
  $A$ is not a masa.

\itemImply {\refl {Effective}}{\refl {WeakMasaNg}}
  It is clear that $N_\G $ is a generating $*$-subsemigroup of $\Norm BA $, so all we need to do is show that $(N_\G
{\cdot }A) \cap A'\subseteq A$.  Moreover, it is easy to see that $N_\G {\cdot }A=N_\G $, so it actually suffices to
prove that $N_\G \cap A'\subseteq A$.

Picking $n$ in $N_\G \cap A'$, note that we have two descriptions of $\beta _n$, the first one being given by \ref
{Beta}, and under our hypothesis it clearly implies that $\beta _n$ is the identity on $\Dom n$.  The second
description, namely \ref {BetaDois}, then implies that the open support $U_n$ of $n$ is contained in $\G '$.
  Since $U_n$ is open, it is actually contained in the interior of $\G '$, which coincides with $X$ by virtue of $\G
$ being effective.  In other words, we have that $U_n\subseteq X$, so $n$ lies in $C_0(X)=A$, as required.

\itemImply {\refl {WeakMasaNg}}{\refl {Effective}} Choose $\gamma $ in the interior of\ $\G '$ and let $U$ be an
open bisection with $\gamma \in U\subseteq \G '$.  Pick $n$ in $\scrCcL U$, such that $n(\gamma )\neq 0$, and
observe that $\beta _n$ is the identity by \ref {BetaDois}.  Consequently $n$ lies in $A'$ by \ref {BetaIdentity},
and hence $n$ belongs to $A$ by hypothesis.  Since $n(\gamma )\neq 0$, it follows that $\gamma \in X$.

\itemImply {\refl {WeakMasaNg}}{\refl {WeakMasa}} Obvious.

\itemImply {\refl {WeakMasa}}{\refl {TopoPrin}}
  If $C_0(X)$ is a {semi-masa} relative to a $*$-semigroup $N$, then every $n$ in $N$ is smooth by \ref
{NearMasaSmooth}, whence $\incl BA$ is a smooth inclusion.  Since $\G $ is second countable we have that $B$ is
separable so we deduce from \ref {RelTFs.ii} that $\incl BA$ is a topologically free inclusion, and hence that $\G $
is a topologically free groupoid by \ref {TopoFreeRed}.

\subProof {\refl {TopoPrin}$\Leftrightarrow $\refl {WeakEffective}} Follows from \ref {Sexy}.

\itemImply {\refl {WeakEffective}}{\refl {Masa}} Assuming that $\G $ is Hausdorff, this follows from \cite
[Propositions 3.1 \& 4.2]{RenaultCartan}.
  \endProof

Except for the Hausdorff case, the above conditions are not all equivalent to each other.  In particular (\LocalMasa
) is not a consequence of (\LocalEffective ), a counter example being given by \cite [2.3]{ExelNHaus}.  Moreover
\refl {WeakEffective} does not imply \refl {WeakMasa} due to the \emph {two-headed snake}, namely the well known
non-Hausdorff space obtained by adding an extra point $1'$ at the right-end of the unit interval $[0, 1]$.  That
space has the structure of an \'etale groupoid, where all points in $[0,1]$ are units, while the extra point $1'$ is
an involution in the isotropy group of 1.  This groupoid clearly satisfies \refl {WeakEffective} but $C_0(X)$ is not
a {semi-masa} in $C_\red ^*(\G )$ since the latter is abelian.

The only ``\emph {one way implication}'' among the statements of \ref {CrazyGroupoids} not yet discussed is \refl
{WeakMasaNg} $\Rightarrow $ \refl {WeakMasa}.  Unfortunately we have been unable to decide whether \refl
{WeakMasaNg} and \refl {WeakMasa} are equivalent.

\endsection

\startsection Canonical states

\sectiontitle

The groupoid $\G _N$ occurring in Theorem \ref {MainTwo} is not guaranteed to be be Hausdorff, and indeed in many
interesting examples it is not.  However, when $\G _N$ happens to be Hausdorff, e.g.~when $\incl BA$ is a Cartan
pair \cite [Proposition 5.4]{RenaultCartan}, one concludes from \ref {EssGpdHausdorff.iii} that
  $$
  B\simeq \EssAlg {\G _N}{\WeylBund } \simeq \Cwr ,
  $$
  recovering in this way the main result of \cite {RenaultCartan}, namely Theorem (5.6).

The present section is dedicated to providing milder sufficient conditions for the above isomorphism to exist, which
can hold even when $\G _N$ fails to be Hausdorff.  A concrete example illustrating this situation is described in
section \ref {GrayIdealTwisted}, below, near \ref {CondForDiag}.

The standing hypotheses for this section will be \ref {StandingChooseN}.

Given $x$ in $X$, and given any state $\psi $ on $B$ extending $\varphi _x$, let us analyze the possible values of
$\psi (n)$ for a given normalizer $n$ in $N$.  If $x$ is a free point relative to $n$, i.e., if $x$ lies in $F_n$,
then of course $\psi (n)=\varepsilon _n(x)$, so that is the end of the story.  On the other hand, if $x$ is not in
$F_n$, there is no way of guessing the value of $\psi (n)$ without further information about $\psi $, since by
definition there are at least two choices of $\psi $ producing distinct values at $x$.

\definition \label DefineCano
  A state $\psi $ on $B$ will be called an \newConcept {$N$-canonical state}{} relative to a given $x$ in $X$,
provided
  \izitem
  \zitem $\psi $ extends $\varphi _x$, and
  \zitem $\psi (n)=0$, for every normalizer $n$ in $N$, for which $x\not \in T_n$.
  \medskip \noindent In case $N$ is understood, we will simply say that $\psi $ is a canonical state relative to
$x$.

When $N$ is taken to be the entire collection of normalizers $N(A,B)$, every $N$-canonical state $\psi $ is a
compatible state in the sense of \cite [Definition 4.1]{SRI}.  However, the following example shows the collection
of $N$-canonical states may be a proper subset of the compatible states.  Let $S$ be the unilateral shift,
$B=C^*(S)$, and $A=C^*\left (\{I\}\cup \{S^nS^*{^n}: n\in \N \}\right )$. Then any pure state $\psi $ on $B$ which
annihilates the compact operators but not $S$ is a compatible state but is not an $N$-canonical state.

A typical situation where canonical states exist in abundance is as follows:

\state Proposition \label CanoInGpd
  Let $(\G , \Lbpar )$ be as in \ref {StandingGPDS}, and let $\incl {B}{A}$ be the inclusion described in \ref
{TwoAlgebras} (so that $B$ is either the full or the reduced twisted groupoid C*-algebra).  Given $x$ in $X$, define
  $$
  \psi _x(f) = f(x), \for f\in B.
  $$
  Then $\psi _x$ is an $N_\G $-canonical state relative to $x$, where $N_\G $ is the generating $*$-semigroup of
normalizers introduced in \ref {DefineNG}.

\Proof
  It is clear that $\psi _x$ is a state on $B$ extending the evaluation state $\varphi _x$. Given $n$ in $N_\G $,
with
  $x\notin T_n$, we must then prove that $\psi _x(n)=0$.

If $x$ is not in $\Dom n$, then $x$ lies in $BJ_x$, by \ref {NormaJota.i}, so $\psi _x(n)=0$, by \ref {VanishJK}.
So let us now suppose that $x$ is in $\Dom n$.  Let $U_n$ be the open bisection defined in \ref {OpenSuport}, so
that $n\in \scrCzL {U_n}$, and
  $
  \Dom n = s(U_n).
  $
  It then follows that there exists a unique $\gamma $ in $U_n$ such that $x=s(\gamma )$, and we claim that $x\neq
\gamma $.  To see this, assume otherwise, so that $x\in U_n$.  Choosing any element $v$ of $C_0(X)$ vanishing off
$U_n\cap X$, with $\evx v=1$, it is easy to see that $nv$ lies in $C_0(X)$, so $x$ is in $T_n$, contradicting our
assumptions.

This shows that $x\neq \gamma $, whence $x$ is not in $U_n$, as $\gamma $ is already an element of $U_n$ whose
source is $x$.  Therefore $n(x)=0$, which is to say that $\psi _x(n)=0$, as required.
  \endProof

Back to the general situation of \ref {StandingChooseN}, given a state $\psi $ extending $\varphi _x$, let us discus
exactly what it takes to check \ref {DefineCano.ii} for a given normalizer $n$.  Observing that $X=F_n\cup \Fix n$,
by \ref {MandatoryValues}, and partitioning $X$ as
  $$
  \def \;{\kern 10pt}
  X \;=\;
    X \setminus \Fix n \; \sqcup \;
    X \setminus F_n \; \sqcup \;
    \Fix n \cap F_n \cap (X\setminus T_n) \; \sqcup \;
    T_n,
  $$
  there is of course nothing to be done in case $x\in T_n$, but otherwise we need to check that $\psi (n)=0$.
  \bigskip

  \begingroup \noindent \hfill \beginpicture
  \setcoordinatesystem units <0.015truecm, -0.015truecm>
  \setplotarea x from -300 to 300, y from -200 to 200
  \put {\null } at -300 -200
  \put {\null } at -300 200
  \put {\null } at 300 -200
  \put {\null } at 300 200 \plot 50.5 -0.5
 51.35 5.998
 51.9 9.485
 50.153 12.949
 48.108 16.38
 44.769 22.766
 43.14 26.097
 39.225 33.361
 38.027 37.548
 37.554 37.647
 36.81 44.648
 31.803 46.542
 27.54 50.318
 27.03 56.967
 20.281 54.48
 18.303 57.847
 12.105 61.062
 11.698 69.115
 1.093 71.999
 2.302 75.707
 -5.664 78.232
 -8.791 79.568
 -16.068 79.709
 -20.482 83.649
 -25.017 80.384
 -31.661 85.909
 -38.4 84.22
 -42.219 90.315
 -53.104 88.19
 -56.04 87.844
 -64.012 87.275
 -66.005 91.481
 -76.004 89.462
 -82.995 89.218
 -86.961 85.75
 -93.89 89.059
 -99.764 88.146
 -104.571 87.015
 -106.295 87.667
 -114.922 84.106
 -119.438 80.337
 -122.829 79.363
 -134.082 78.189
 -133.183 71.821
 -140.12 73.265
 -147.881 70.527
 -151.453 64.613
 -152.826 60.533
 -157.987 59.292
 -163.928 55.899
 -169.637 53.362
 -169.106 49.692
 -170.327 48.895
 -175.29 38.983
 -177.989 35.964
 -178.416 31.849
 -182.567 30.649
 -186.434 23.373
 -187.015 21.032
 -187.304 14.638
 -190.299 10.201
 -188.997 4.732
 -191.396 4.242
 -188.496 -4.257
 -187.295 -7.754
 -188.796 -11.238
 -191.998 -12.697
 -186.903 -22.121
 -185.516 -21.499
 -181.838 -30.819
 -179.875 -29.07
 -180.631 -36.243
 -181.111 -37.327
 -175.323 -43.311
 -175.272 -46.185
 -168.967 -54.941
 -162.416 -57.567
 -162.628 -62.056
 -157.612 -63.399
 -152.378 -68.587
 -145.937 -69.612
 -141.3 -72.467
 -141.479 -77.145
 -134.485 -75.639
 -126.331 -78.942
 -124.03 -83.049
 -115.596 -82.955
 -110.041 -87.655
 -110.38 -85.144
 -103.627 -84.42
 -98.795 -89.478
 -87.901 -88.316
 -82.958 -92.932
 -75.982 -88.325
 -73.987 -91.493
 -66.988 -92.436
 -60 -89.155
 -56.039 -92.649
 -49.119 -86.921
 -44.254 -90.971
 -36.461 -85.803
 -29.752 -87.419
 -26.143 -81.823
 -19.647 -79.019
 -12.278 -82.011
 -10.05 -76.804
 -0.975 -74.404
 1.933 -71.817
 8.663 -72.049
 11.203 -64.107
 14.54 -63.999
 21.666 -58.732
 19.568 -56.314
 28.238 -54.754
 28.666 -48.062
 33.843 -47.245
 33.762 -39.314
 37.416 -37.279
 43.797 -31.149
 43.901 -28.935
 45.72 -23.647
 46.252 -21.297
 45.492 -13.895
 47.437 -11.451
 50.085 -8.978
 49.434 -2.486
 50.5 -0.5
 / \plot 190 0
 189.4 8.985
 187.608 17.88
 184.64 26.597
 180.527 35.048
 175.31 43.148
 169.04 50.818
 161.781 57.98
 153.605 64.562
 144.593 70.499
 134.836 75.732
 124.432 80.209
 113.483 83.884
 102.1 86.72
 90.396 88.69
 78.488 89.775
 66.496 89.962
 54.539 89.25
 42.736 87.646
 31.205 85.167
 20.062 81.837
 9.418 77.689
 -0.62 72.765
 -9.953 67.113
 -18.487 60.792
 -26.137 53.862
 -32.827 46.395
 -38.489 38.464
 -43.067 30.149
 -46.515 21.532
 -48.799 12.701
 -49.896 3.742
 -49.795 -5.254
 -48.498 -14.197
 -46.016 -22.999
 -42.375 -31.57
 -37.611 -39.827
 -31.772 -47.685
 -24.916 -55.067
 -17.112 -61.899
 -8.437 -68.112
 1.021 -73.645
 11.169 -78.442
 21.904 -82.455
 33.12 -85.644
 44.705 -87.978
 56.542 -89.432
 68.513 -89.993
 80.5 -89.655
 92.381 -88.421
 104.039 -86.303
 115.357 -83.323
 126.222 -79.511
 136.525 -74.904
 146.163 -69.549
 155.04 -63.499
 163.068 -56.814
 170.166 -49.562
 176.262 -41.814
 181.297 -33.649
 185.22 -25.147
 187.992 -16.395
 189.585 -7.478
 190 0
 / \setdots <1.5pt> \plot 30 0
 29.85 3.993
 29.402 7.947
 28.66 11.821
 27.632 15.577
 26.327 19.177
 24.76 22.586
 22.945 25.769
 20.901 28.694
 18.648 31.333
 16.209 33.659
 13.608 35.648
 10.871 37.282
 8.025 38.542
 5.099 39.418
 2.122 39.9
 -0.876 39.983
 -3.865 39.667
 -6.816 38.954
 -9.699 37.852
 -12.484 36.372
 -15.145 34.528
 -17.655 32.34
 -19.988 29.828
 -22.122 27.019
 -24.034 23.939
 -25.707 20.62
 -27.122 17.095
 -28.267 13.4
 -29.129 9.57
 -29.7 5.645
 -29.974 1.663
 -29.949 -2.335
 -29.624 -6.31
 -29.004 -10.222
 -28.094 -14.031
 -26.903 -17.701
 -25.443 -21.193
 -23.729 -24.474
 -21.778 -27.511
 -19.609 -30.272
 -17.245 -32.731
 -14.708 -34.863
 -12.024 -36.647
 -9.22 -38.064
 -6.324 -39.101
 -3.365 -39.748
 -0.372 -39.997
 2.625 -39.847
 5.595 -39.298
 8.51 -38.357
 11.339 -37.033
 14.056 -35.338
 16.631 -33.291
 19.041 -30.911
 21.26 -28.222
 23.267 -25.251
 25.041 -22.027
 26.566 -18.584
 27.824 -14.955
 28.805 -11.177
 29.498 -7.287
 29.896 -3.324
 30 0
 / \setsolid
  \put {$T_n$} at 0 0
  \put {$F_n$} at -120 0
  \put {$\Fix n$} at 120 0
  \put {$X\setminus \Fix n$} at -110 110
  \put {$X\setminus F_ n$} at 130 110
  \put {$\Fix n \cap F_n \cap (X\setminus T_n)$} at 0 -140
  \arrow <0.15cm> [0.35,0.75] from 0 -120 to 0 -55
  \endpicture \hfill \null \endgroup

  \bigskip
  \noindent In case $x$ lies in $X \setminus \Fix n$, then \ref {MandatoryValues.ii} says that any extension $\psi $
necessarily vanishes on $n$, so again there is nothing to be done.  There is not much to be said if $x$ lies in $X
\setminus F_n$, but when $x$ is in the third set in the above partition of $X$, then every extension $\psi $
satisfies
  $$
  \psi (x) = \varepsilon _n(x),
  $$
  but $\varepsilon _n(x)\neq 0$, thanks to \ref {TwoExtensions.i}, so there is no hope of finding an $N$-canonical
state relative to such an $x$.

We therefore have the following:

\state Proposition \label CharacCanonical
  There is no $N$-canonical state relative to any point $x$ lying in
  $$
  \bigcup _{n\in N} \Fix n \cap F_n \cap (X\setminus T_n).
  $$
  For all other points $x$ in $X$, one has that
  \izitem
  \zitem a state $\psi $ is $N$-canonical relative to $x$, if and only if $\psi (n)=0$, for all $n$ in $N$ for which
$x\in X \setminus F_n$,
  \zitem there is at most one $N$-canonical state relative to $x$.

\Proof After what was said in the paragraph before the statement, it suffices to prove (ii).
  Given an $N$-canonical state $\psi _1$, and given $n$ in $N$, observe that by (i)
  $$
  \psi _1(n) = \clauses {
  \cl \varepsilon _n(x) if {x\in F_n},
  \cl 0 if {x\not \in F_n}.
  }
  $$
  If $\psi _2$ is another $N$-canonical state, we then have for all $a$ in $A$, that
  $$
  \psi _1(na)\={CondexForState}
  \psi _1(n) \evx a=
  \psi _2(n) \evx a\={CondexForState}
  \psi _2(na),
  $$
  so $\psi _1$ and $\psi _2$ coincide on $\NoA $, and hence everywhere since $N$ is generating.
  \endProof

For the case of {semi-masa}s, we have the following convenient way to check that states are $N$-canonical:

\state Proposition
  In addition to \ref {StandingChooseN}, assume that $\incl BA$ is a {semi-masa} relative to $N$, and let $x\in X$.
Then
  \izitem
  \zitem given a state $\psi $ on $B$ extending $\varphi _x$, one has that $\psi $ is $N$-canonical relative to $x$,
if and only if $\psi (n)=0$, for every $n$ in $N$, for which $x\in \bd n$,
  \lbldzitem UniqueIsCanonical if $x$ is a free point, then the unique extension of $\varphi _x$ to a state on $B$
is $N$-canonical.

\Proof Focusing on (i), let us first prove that the set $\Fix n \cap F_n \cap (X\setminus T_n)$, appearing in \ref
{CharacCanonical}, is empty for every $n$ in $N$.  Indeed, assuming that $x$ lies in this set, observe that $x$
cannot be in $\bd n$ by \ref {TwoExtensions.ii}, so
  $$
  x\in \iFix n \={IntroNu.iii}T_n,
  $$
  a contradiction.  This shows that $x$ is among the points referred to as ``\emph {all other points}'' in \ref
{CharacCanonical}, and hence $\psi $ is $N$-canonical if and only if $\psi (n)=0$, for all $n$ in $N$, for which
  $$
  x\in X \setminus F_n\={FreeAllMinusBd} \bd n.
  $$
  This concludes the proof of (i).  Point (ii) then follows from (i), since a free point is never in any $\bd n$, by
\ref {TwoExtensions.ii}.
  \endProof

Following Brown and Guentner \cite {BG}, Buss, Echterhoff and Willett \cite {BEW}, \cite {BEWTwo} introduced the
notion of an \emph {exotic crossed product} for a C*-dynamical system $(A,G,\alpha )$, which turns out to be a
C*-algebra obtained by completing $C_c(G,A)$ under a norm sitting in between the full and the reduced norms.  An
exotic crossed product is in fact required to satisfy certain functoriality properties \cite [Definition
3.5.(2)]{BEW}, but let us nevertheless give the following simplified definition which completely ignores any
functoriality aspects.

\def \Cmu #1#2{C^*_\mu \big (#1,#2 \big )} \def \Cwmu {\Cmu {\G _N}{\WeylBund }}

\definition
  Given a twisted \'etale groupoid $(\G , \Lbpar )$, and given any C*-norm $\mu $ on $\scrCcL \G $, with
  $$
  \Vert f\Vert _\red \leq \mu (f) \leq \Vert f\Vert _{\hbox {\sixrm max}}, \for f \in \scrCcL \G ,
  $$
  where $\Vert f\Vert _{\hbox {\sixrm max}}$ is the maximum
  C*-norm\fn {The second inequality above is superfluous since any norm is obviously bounded by the maximum one.},
we will denote by $\Cmu \G \Lbpar $ the completion of $\scrCcL \G $ under $\mu $, and we will call it an \emph
{exotic twisted groupoid C*-algebra} for $(\G , \Lbpar )$.

When $\mu $ is the reduced C*-norm we of course get the reduced C*-algebra, while the same holds for the maximum
norm and the full C*-algebra.  Observe however that the norm $\Vert {\,\cdot \,}\Vert \mnrm $ of \ref {MinNorm} is
technically not included since it does not always sit in between
  $\Vert {\,\cdot \,}\Vert _\red $ and $\Vert {\,\cdot \,}\Vert _{\hbox {\sixrm max}}$.

The following main result is an effective criterion, not requiring any freeness hypothesis, for determining when a
given inclusion admits a groupoid model.  See section \ref {PerFunSec} for an example of its application.

\state Theorem \label MainThree
  Let $\incl {B}{A}$ be an inclusion, with $A$ non-degenerate
  (see \ref {RegAndStuff}) and abelian.  Also let $N\subseteq \Norm BA $ be a generating $*$-subsemigroup and assume
that there exists an $N$-canonical state relative to every $x$ in the spectrum $X$ of $A$. Then:
  \izitem
  \zitem $B$ is isomorphic to some exotic groupoid C*-algebra $\Cwmu $, under an isomorphism
  $$
  \Psi :\Cwmu \to B
  $$
  carrying $C_0(X)$ onto $A$.
  \zitem If the {gray} ideal of $B$ vanishes, the above exotic groupoid C*-algebra can be taken to be the reduced
one.  In particular the {gray} ideal of\/ $\Cwr $ vanishes.

\Proof
  Let
  $$
  \Phi : \Cw \to B
  $$
  be the $*$-epimomorphism given by \ref {MapFromGpg}.

  Given $x$ in $X$, as well as an $N$-canonical state $\psi $ relative to $x$, we claim that
  \lbldeq CrossExpect
  $$
  \psi \big (\Phi (f)\big ) = f(x), \for f\in \Cw .
  $$
  Observe that when we speak of $f(x)$, above, we are of course identifying $x$ with its canonical image in the unit
space of $\G _N$, as in \ref {XAsUnitSpace}.

  By \ref {DenseGuy}, in order to prove our claim we may assume that
  $f=
  \xi _{n,a}
  $, where $n\in N$, and $a\in A$.

Given $n$ and $a$, as above, suppose first that $x$ is not in $T_n$.  Then
  $$
  \psi \big (\Phi (f)\big ) = \psi \big (\Phi (\xi _{n,a})\big ) = \psi (na) \={CondexForState} \psi (n)\evx a= 0,
  $$
  because $\psi $ is $N$-canonical.  We must therefore prove that $f(x)=0$, as well.  Assuming by contradiction that
this is not so, it follows that $x$ lies in the support of $f$, and hence also in $\O n$.  However, this can only
happen if $x\in \Dom n$, and the germ $[n,x]$ is trivial, in which case necessarily $x$ lies in $T_n$, contrary to
our hypothesis.

Let us now assume that $x$ lies in $T_n$, so we may find $v$ in $A$, such that $\evx v=1$, and
  $$
  c:= nv\in A.
  $$
  It follows that $nv^2=cv$, and since $\evx {v^2}\neq 0\neq \evx v$, we deduce that $[n,x] = [c,x]$.  So,
  $$
  \psi \big (\Phi (f)\big ) = \psi \big (\Phi (\xi _{n,a})\big ) = \psi (na) \={CondexForState}
  \Frac {\psi (n{v^2}a)}{\evx {v^2}} =
  \Frac {\psi (cva)}{\evx {v^2}} =
  \Frac {\evx {cva}}{\evx {v^2}} =
  \Frac {\evx {ca}}{\evx {v}}.
  $$
  On the other hand,
  $$
  f(x) = f([c,x]) = \xi _{n,a}([n,x]) = [\evx a, n,x] =
  \left [\Frac {\evx {ca}}{\evx {{v^2} }}, v,x\right ],
  $$
  where the last step is justified by the fact that $n{v^2}=vc$, and
  $$
  \Frac {\evx a}{\evx {v^2}} = \Frac {\Frac {\evx {ca}}{\pilar {5.8pt}\evx {v^2}}}{\evx c}.
  $$
  We then conclude from the above that
  $$
  f(x) =
  \Frac {\evx {ca}}{\evx {v}}[\evx v\inv , v,x] \={FoundUnit}
  \Frac {\evx {ca}}{\evx {v}} {\bf 1},
  $$
  completing the proof of \ref {CrossExpect}.

We next consider the map $P$ from $B$ to the algebra of all bounded, complex functions on $X$, given by
  $$
  P(b)\calcat x= \psi _x(b), \for b\in B, \for x\in X,
  $$
  where $\psi _x$ is the unique $N$-canonical state relative to $x$.  Given $f$ in $\Cw $ we then have for every $x$
in $X$ that
  $$
  P\big (\Phi (f)\big )\calcat x= \psi _x\big (\Phi (f)\big ) \={CrossExpect} f(x),
  $$
  which is to say that
  \lbldeq PGivesRestr
  $$
  P\big (\Phi (f)\big )=f|_X.
  $$

With these preparations we may now prove that
  $$
  \Ker (\Phi ) \subseteq \Ker (\Lambda ),
  $$
  where $\Lambda $ is the left regular representation mentioned in \ref {IntForlLeftReg}, because for every $f$ in
$\Cw $, one has that
  \lbldeq AlmostFaith
  $$
  \Phi (f)=0 \IMPLY
  P\big (\Phi (f)^*\Phi (f)\big )=0 \IFF
  $$
  $$
  \IFF
  P\big (\Phi (f^*f)\big )=0 \explain {PGivesRestr}\IFF (f^*f)|_X=0 \IFF \Lambda (f)=0.
  $$
  It follows that there exists a natural $*$-epimomorphism
  $$
  \def \quad {\ \ }
  {\Cw \over \Ker (\Phi )} \quad \longrightarrow \quad {\Cw \over \Ker (\Lambda )} \quad = \quad \Cwr ,
  $$
  so the algebra in the left-hand-side above, which is evidently isomorphic to $B$, is an exotic groupoid C*-algebra
for $(\G , \Lbpar )$.  This proves (i).

In order to prove (ii), assume that the {gray} ideal of $B$ vanishes.  Recalling that $\Pf $ denotes the free
expectation introduced in \ref {DefineFreeCondExp}, notice that
  $$
  P(b)|_F= \Pf (b), \for b\in B,
  $$
  so, for every $b$ in $B$, we have that
  $$
  P(b^*b)=0 \IMPLY \Pf (b^*b)=0 \explain {FaithGray}\IMPLY b\in \Gamma \IMPLY f=0,
  $$
  so we see that $P$ is faithful.  Therefore the ``$\Rightarrow $'' in \ref {AlmostFaith} may be replaced by a
``$\Leftrightarrow $'', hence showing that actually $\Ker (\Phi ) = \Ker (\Lambda )$.

Once the kernel of $\Phi $ is quotiented out, one is therefore left with an isomorphism
  $$
  \Psi :\Cwr \to B,
  $$
  proving (ii).
  \endProof

We should point out that, unlike \ref {MainTwo}, or any one of the various results in the literature on Cartan
subalgebras, the above Theorem characterizes twisted groupoid C*-algebras, albeit possibly exotic ones, via a
condition (existence of canonical states) which holds for any twisted \'etale groupoid by \ref {CanoInGpd}.
  Summarizing, we therefore have the following main consequence:

\state Corollary
  Let $\incl {B}{A}$ be an inclusion of C*-algebras, with $A$ non-degenerate and abelian.  Then the following are
equivalent:
  \izitem
  \zitem There exists a generating $*$-subsemigroup of normalizers $N\subseteq \Norm BA$ such that, for every $x$ in
the spectrum of $A$, there is an $N$-canonical state relative to $x$,
  \zitem $\incl BA$ is naturally isomorphic to $\INCL {\Cmu {\G } {\Lb }}{\CzGz }$, where
  $\G $ is a (not necessarily Hausdorff) \'etale grou\-poid (with a locally compact Hausdorff unit space), $\Lb $ is
a Fell line bundle over $\G $, and $\Cmu {\G } {\Lb }$ is an exotic groupoid C*-algebra.

\Proof
  The implication (i) $\Rightarrow $ (ii) is nothing but \ref {MainThree}, while the converse may be proved by
observing that \ref {CanoInGpd} applies just as well to exotic groupoid C*-algebras with the exact same proof.
  \endProof

\endsection

\PART {III -- EXAMPLES AND OPEN QUESTIONS}

\startsection Example: non-smooth normalizers

\label NonSmoothNrm

\sectiontitle

As we have seen, the notion of smooth normalizers is a central piece of our strategy for dealing with regular
inclusions.  In this section we would like to exhibit a {light} inclusion of abelian C*-algebras containing a
normalizer $n$, with a big dense set of relative free points which nevertheless has empty interior.  Consequently
$n$ is continuous but not smooth.  This will also illustrate how badly behaved can the set of free points be.

  Let $X$ be a compact space, $D\subseteq X$ be a dense subset, and let us fix a continuous function
  \lbldeq FixH
  $$
  h:D\to \circle ,
  $$
  where $\circle $ denotes the unit circle.  We should remark that we are not requiring $h$ to extend to a
continuous function on $X$, and in fact all of the fun resides on the failure of such an extension to exist.  Define
  $$
  H:x\in D \mapsto \big (x,h(x)\big )\in X\times \circle ,
  $$
  and put $Y=\overline {H(D)}$.  In other words, $Y$ is the closure of the graph of $h$.

Letting $q:X\times \circle \to X$ denote the projection on the first coordinate, observe that $q(Y)=X$, because
$q(Y)$ is compact and contains the dense set $D$.  Therefore the $*$-homomorphism
  $$
  \iota : C(X)\to C(Y),
  $$
  dual to $q$, is injective and hence we may view $C(X)$ as a subalgebra of $C(Y)$.  Precisely, $C(X)$ is identified
with the set of functions in $C(Y)$ which factor through $q$.

  Considering the inclusion
  $$
  \incl {B}{A} = \Incl {C(Y)}{C(X)},
  $$
  observe that, since both algebras are commutative, the only condition for an element $b$ in $B$ to be a normalizer
is that $b^*b$ belong to $A$.
  In particular, every unitary element of $B$ is a normalizer, and since $B$ is linearly spanned by its unitary
elements, the inclusion is evidently regular.

\state Proposition \label InclComH
  \izitem
  \zitem $D$ is contained in the set $F$ of free points.
  \lbldzitem InterprPresent
  For every $b$ in $C(Y)$, and every $x$ in $D$, one has that $\Ex (b) = b\big (x, h(x)\big )$.
  \zitem $\Incl {C(Y)}{C(X)}$ is a {light} inclusion.

\Proof
  By \ref {FreeCommut}, in order to prove (i), it suffices to show that, for $x$ in $D$, one has that $q\inv
(\{x\})$ is a singleton, but since $H(x)$ necessarily lies in $q\inv (\{x\})$, our task is to show that
  \lbldeq Singleton
  $$
  q\inv (\{x\})\subseteq \{ H(x)\}.
  $$
  Given any member of $q\inv (\{x\})$, we may write it as $(x, \lambda )$, for some $\lambda $ in $\circle $.  Since
$Y$ is the closure of $H(D)$, we may also write
  $$
  (x,\lambda )=\lim _i\big (x_i,h(x_i)\big ),
  $$
  with $x_i\in D$.  We then have that $x_i\to x$, and since $h$ is continuous at $x$ we get
  $$
  h(x) = \lim _i h(x_i) = \lambda ,
  $$
  so $(x,\lambda )=(x,h(x))=H(x)$, as desired.  Stressing that $\big (x, h(x)\big )$ is the unique point in $q\inv
(\{x\})$, we see that \refl {InterprPresent} follows from the last sentence in \ref {FreeCommut}.

Turning now to (iii), pick $b$ in the {gray} ideal of $C(Y)$, so that for every free point $x$ of $X$, and in
particular for every $x$ in $D$, one has
  $$
  0 =
  \Pf (b^*b)\calcat x=
  \Ex (b^*b) =
  \big (b^*b) \big (x, h(x)\big ) =
  \big |b\big (x, h(x)\big )\big |^2.
  $$
  If follows that $b$ vanishes on $H(D)$, and hence that $b=0$.  This shows that the {gray} ideal of $C(Y)$
vanishes, so our inclusion is indeed {light}.
  \endProof

The above construction leads to a very natural normalizer, namely the function on $Y$ given by
  \lbldeq IntroBadNormalizer
  $$
  n(x,\lambda )=\lambda , \for (x,\lambda )\in Y.
  $$
  Since $C(Y)$ is commutative, any unitary element of $C(Y)$ is clearly a normalizer of $C(X)$, and it is obvious
that our $n$ above is such an element.
  By \ref {InclComH.i} we have that
  \lbldeq DinExtn
  $$
  D\subseteq F\subseteq F_n,
  $$
  so there is no question that $n$ is a continuous normalizer according to \ref {DefineContSmooth.i}.  Regarding the
function $\varepsilon _n$ defined on $F_n$, we may easily compute it on $D$, observing that, for $x$ in $D$, we have
  \lbldeq enExtendsH
  $$
  \varepsilon _n(x) \={ThreeValues} \Ex (n) \={InterprPresent} n\big (x, h(x)\big ) = h(x).
  $$

Whether or not $n$ is smooth will of course depend on our initial data, namely the function $h$ chosen in \ref
{FixH}.  Our goal here is to produce a non-smooth normalizer by choosing a highly irregular function $h$, and we
will do so based on \cite [Remark 4.31]{Rudin}\fn {We thank Nate Eldredge for pointing out this example to us
through his answer to our question posted in\hfill \break \withfont {cmtt10 scaled
800}{https://mathoverflow.net/questions/317463}.}.
  We in fact need to make a small adaptation to it since we would like $X$ to be compact.

Let $X$ be the interval $[0, 2\pi ]$, and let $D=X\setminus {\bf Q}$.  Choose an enumeration $\{q_n\}_{n\in {\bf
N}}$ of $X\cap \bf Q$ and, for every $x$ in $D$, let us define
  \lbldeq deflogh
  $$
  f(x) = \sum _{n : q_n < x} 2^{-n}.
  $$

We leave it as an exercise for the reader to prove that $f$ is continuous on $D$, and that it cannot be extended
continuously to any proper superset\fn {Should we have chosen the more popular interval $X=[0, 1]$, instead, one
would be able to continuously extend $f$ to $D\cup \{1\}$.  This would in fact be the case for any interval with a
rational right-end point.  For this reason we have decided to choose $[0,2\pi ]$, which is of course almost as
popular, although we have no intention of wrapping our interval around the circle!}  of $D$.

Observing that $f$ is positive and bounded above\fn {In case one thinks of $\N $ as starting at $1$, $f$ would also
be bounded by $1$.}  by $2$, its range is mapped homeomorphically to a subset of $\circle $ via the exponential map
$t\to e^{it}$, so the function
  $$
  h:x\in D \mapsto e^{if(x)}\in \circle
  $$
  is also continuous on $D$ but it cannot be extended continuously to any proper superset.

\state Proposition \label BadNormalizer
  Considering the inclusion $\Incl {C(Y)}{\cclosint 0{2\pi }}$ discussed in \ref {InclComH} and built with the above
function $h$, and considering the normalizer $n$ given by \ref {IntroBadNormalizer}, one has that
  $$
  F = F_n= D = [0, 2\pi ]\setminus {\bf Q},
  $$
  and hence $n$ is a non-smooth normalizer.

\Proof
  We have already seen in \ref {DinExtn} that $D\subseteq F_n$.  Supposing by contradiction that $D$ is a proper
subset of $F_n$, observe that $\varepsilon _n$ is a continuous function on $F_n$ by \ref {EbContinuous}, while
  $\varepsilon _n$ extends $h$ by \ref {enExtendsH}.  However this is a contradiction with the already established
fact that $h$ admits no continuous extension.

In order to prove that $F$ also coincides with $D$, it is enough to notice that
  $$
  D\explain {InclComH.i}\subseteq F\explain {XfAndXb}\subseteq F_n=D.
  \closeProof
  $$
  \endProof

The above result notwithstanding, the inclusion in \ref {BadNormalizer} is not so badly behaved:

\state Proposition
  The inclusion $\Incl {C(Y)}{\cclosint 0{2\pi }}$ mentioned in \ref {BadNormalizer} is a weak Cartan inclusion.

\Proof We first show our inclusion to be essential.  If $x\in D$ and $(x,\lambda )\in Y$ for some $\lambda \in \bf
T$, then $\lambda = h(x)$ by continuity of $h$, so $(x,\lambda )\in H(D)$.  Suppose $R\subseteq Y$ is closed and
$q(R)=X$.  Since $D\subseteq q(R)$, it follows that $H(D)\subseteq R$, whence $R=Y$.  It follows that $\Incl
{C(Y)}{\cclosint 0{2\pi }}$ is essential.

The function $f$ defined in \ref {deflogh} is increasing on $D$ and has the property that for every $q_n\in {\bf
Q}\cap (0,2\pi ]$
  $$
  2^{-n}+ \lim _{t\rightarrow q_n^-,\, t\in D} f(t)=\lim _{t\rightarrow q_n^+,\, t\in D}f(t).
  $$

Let $Z\subseteq [0,2\pi ]\times [0,2]$ be the closure of the graph of $f$, and let $\varpi : Z\rightarrow [0,2\pi ]$
be the projection onto the first component, $(x,t)\mapsto x$.  Then for $x\in [0,2\pi ]$,
  $$
  \varpi ^{-1}(x)=\cases { \{(0,0)\} & if $x=0$\cr \{(x, f(x-)), (x,f(x+))\} & if $x\in {\bf Q} \cap (0,2\pi )$\cr
\{(x, f(x))\} & if $x \in
  (0,2\pi )\setminus \bf Q$.\cr }
  $$
  It follows that $Z$ is a totally disconnected, perfect set (and hence a Cantor set).  As the map $Z\ni
(x,y)\mapsto (x, \exp (iy))\in Y$ is a homeomorphism, $Y$ is also totally disconnected. By \ref {proj=smooth},
$\Incl {C(Y)}{\cclosint 0{2\pi }}$ is smooth, so~\ref {RelTFs} implies $\Incl {C(Y)}{\cclosint 0{2\pi }}$ is a weak
Cartan inclusion.
  \endProof

\endsection

\startsection Example: periodic functions on the interval

\def \1{{\bf 1}}

\label PerFunSec

\sectiontitle

This is an example based on an extremely familiar inclusion, presenting a very unusual and unexpected behavior.

Letting
  $$
  A=C(\circle ), \and B=\cclosint 01,
  $$
  we shall identify $A$ as a subalgebra of $B$ via the familiar point of view according to which a function on the
circle corresponds to a function $f$ on the unit interval such that $f(0)=f(1)$.

We observe that this inclusion may also be obtained via \ref {InclComH}, where $h$ is any continuous branch of the
complex square root function defined on the set $D=\circle \setminus \{1\}$.  Either pursuing this lead, or proving
it directly, one finds that $\Incl {\cclosint 01}{C(\circle )}$ is a {light} inclusion, whose set of free points is
given by $F=\circle \setminus \{\1\}$, where $\1$ is the point in the circle where the two ends of the unit interval
are sent under the restriction of the dual of the inclusion map $C(\circle )\hookrightarrow \cclosint 01$ to the
pure states.

With so many free points, every normalizer is smooth.  Consequently our inclusion is smooth, whence it is also a
weak Cartan inclusion.  Therefore \ref {MainTwo} applies, enabling us to describe it by means of a relative Weyl
groupoid and line bundle, but we must first choose a generating semigroup of normalizers $N$.

As already mentioned the only condition for an element $b$ in $\cclosint 01$ to be a normalizer is that $b^*b$
belong to $A$.  In other words, all we need is
  \lbldeq SemeModulo
  $$
  |b(0)| = |b(1)|.
  $$

Speaking of a generating semigroup of normalizers, there are many to choose from, and an important condition is that
$N$ should include at least one element $b$ with $b(0)\neq b(1)$, or else \ref {DefineAdmissTotal.ii} will not be
satisfied.  Let us therefore assume we have chosen a generating semigroup $N$, and let us briefly describe the
relative Weyl groupoid $\G _N$.

Given $n$ in $N$, and given $x$ in $\circle $ with $\evx {n^*n}\neq 0$, let us analyze the germ $[n,x]$.  Assuming
first that $x\neq \1$, and choosing any $v$ in $B$ such that $\evx v=1$, and $\ev v\1=0$, we see that $nv\in A$, so
  $$
  [n,x] = [nv,x] \in \G _N^{(0)} .
  $$
  In other words, the reduction of $\G _N$ to $\circle \setminus \{\1\}$ is entirely made of units.

The germs of the form
  $
  [n,\1],
  $
  on the other hand, are not necessarily trivial.  In fact, observe that if $n$ and $m$ are two normalizers in $N$,
with
  $$
  \ev {n^*n}\1\neq 0\neq \ev {m^*m}\1,
  $$
  then $[n,\1] = [m,\1]$ if and only if that there are $a$ and $b$ in $A$, such that $\ev {a}\1\neq 0\neq \ev
{b}\1$, and $na=mb$.  Since $a$ and $b$ are periodic, it follows that
  $$
  \Frac {n(1)}{m(1)} = \Frac {\ev {b}1}{\ev {a}1} = \Frac {\ev {b}0}{\ev {a}0} = \Frac {n(0)}{m(0)}
  \quad \Longrightarrow \quad
  \Frac {n(1)}{n(0)} = \Frac {m(1)}{m(0)},
  $$
  so the number
  $$
  \lambda (n):= \Frac {n(1)}{n(0)}
  $$
  is an invariant of the germ $[n,\1]$, and in fact it is a complete invariant.  Observing that \ref {SemeModulo}
forces $\lambda (n)$ to have unit norm, we get an isomorphism
  $$
  \lambda :\G _N(\1) \to \circle ,
  $$
  onto some subgroup of the circle.

Precisely which subgroup of $\circle $ will show up depends on the chosen generating semigroup $N$, but the reader
could certainly imagine many different choices producing quite different isotropy groups at $\1$.  Any countable
subgroup of the circle is in fact possible and, should we temporarily disregard the requirement that $N$ be
countable, any subgroup whatsoever could turn up, although of course $\G _N(\1)$ will always carry the discrete
topology.

As we have seen, the groupoid algebras $\Cw $ vary wildly with each choice of $N$ but, as proven in \ref {MainTwo},
once we quotient by the {gray} ideal, we have
  $$
  \EssAlg {\G _N}{\WeylBund } \simeq \cclosint 01.
  $$

We would now like to discuss a specific choice of $N$ which stands out due to the fact that it leads to a twisted
Weyl groupoid whose reduced C*-algebra has a vanishing {gray} ideal; we will use this specific choice for the
remainder of this section.  Consider the element $n\in \cclosint 01$, given by
  $$
  n(t) = e^{\pi it}, \for t\in [0,1].
  $$
  Observing that $n(0)=1$, and $n(1)=-1$, we see that $n$ satisfies \ref {SemeModulo}, so it is a normalizer,
evidently a unitary one.  Let us therefore pick the $*$-semigroup
  $$
  N=\{n^k:k\in {\bf Z}\},
  $$
  and let us consider the twisted Weyl groupoid built in terms of $N$.

The first important point we'd like to make is that the state $\psi $ on $\cclosint 01$ given by
  $$
  \psi (f) = \Frac {f(0)+f(1)}{2}, \for f\in \cclosint 01
  $$
  is an $N$-canonical state relative to $\1$.  Leaving for the reader to check that $\psi $ extends $\varphi _{\1}$,
observe that for every $p\in {\bf Z}$, one has that
  $$
  n^{2p}(1) = e^{2p\pi i} = 1 = n^{2p}(0),
  $$
  so $n^{2p}$ is periodic and hence it lies in $C(\circle )$.  Consequently the set $T_{n^{2p}}$ of trivial points
is the whole of $\circle $, and in particular $\1\in T_{n^{2p}}$, so \ref {DefineCano.ii} is vacuously satisfied for
$n^{2p}$.  For odd $k$, on the other hand, we have that
  $$
  \psi (n^k) =
  \Frac {1+e^{k\pi i}}{2} =
  \Frac {1-1}{2} =
  0.
  $$
  This shows that $\psi $ is indeed an $N$-canonical state relative to $\1$.

Regarding $N$-canonical states relative to the remaining points, let $x\in \circle \setminus \{\1\}$.  We then claim
that $x$ is a trivial point relative to any normalizer $n^k$ in $N$.  To see this, it suffices to choose $v$ in
$C(\circle )$ with $\evx v=1$ and $\ev v\1=0$, and to observe that $n^kv$ lies in $C(\circle )$.  This said we
conclude that the unique extension of $\varphi _x$ to $C([0,1])$ is an $N$-canonical state, since \ref
{DefineCano.ii} is again vacuously satisfied.  We therefore see that the hypotheses of \ref {MainThree.ii} hold, so
we conclude that
  \lbldeq desirediso
  $$
  \Cwr \simeq \cclosint 01,
  $$
  without the need to mod out the {gray} ideal which, after all, must vanish.

It is interesting to provide an explicit formula for the isomorphism in \ref {desirediso}.  We present an outline,
leaving details for the reader.  We have
  $$
  \G _N=\{[n^0, t]: t\in (0,1)\}\cup \{[n^0, {\bf 1}]\}\cup \{[n, {\bf 1}]\}\simeq (0,1)\cup \{[n^0, {\bf 1}]\}\cup
\{[n, {\bf 1}]\},
  $$
  with open bisections
  $$
  U_0=(0,1)\cup \{[n^0,{\bf 1}]\}\simeq \circle \dstext {and} U_1=(0,1)\cup \{[n, {\bf 1}]\}\simeq \circle .
  $$
  Also, $\Cwr $ is the linear span of $C(U_0,\Lbpar _N)$ and $C(U_1,\Lbpar _N)$.  Furthermore,
  $$
  f_0\in C(U_0,\Lbpar _N)\quebra\Leftrightarrow
  \text { there exists } \tilde f_0\in A \text { such that }
f_0([m,t])=\cases { [\tilde
              f_0(t),n^0,t]& if $[m,t]=[n^0,t]$;\cr
              [0, n, {\bf 1}], & if $[m,t]=[n,{\bf 1}]$}
  $$
  and
  $$
  f_1\in C(U_1,\Lbpar _N)\quebra\Leftrightarrow
  \text { there exists } \tilde f_1\in A \text { such that }
f_1([m,t])=\cases { [\tilde
              f_1(t),n,t]& if $[m,t]=[n,t]$;\cr
              [0, n^0, {\bf 1}], & if $[m,t]=[n^0,{\bf 1}]$.}
  $$
  For $t\neq {\bf 1}$ and $\tilde f_1\in A$, $[\tilde f_1(t),
              n, t]=[\tilde f_1(t)n(t), n^0, t]$ and a computation
              shows that if $f_i\in C(U_i,\Lbpar _N)$, then
  $$
  \Cwr \ni f_0+f_1=0\Leftrightarrow 0=\tilde f_0 +\tilde f_1 n\in
              C[0,1].
  $$
  Therefore the map $\alpha : \Cwr \rightarrow
     C[0,1]$ given by $\alpha (f_0+f_1)=\tilde f_0+ \tilde f_1 n$ is a
     well-defined $*$-monomorphism.  To see it is onto, let
     $\tau :[0,1]\rightarrow [0,1]$ be the
     map $\tau (t)=1-t$.  Then for $h\in C[0,1]$, $\tilde h_0= (h+h\circ
     \tau )/2$ and $\tilde h_1=(h-h\circ \tau )/(2n)$ both belong to $A$
     and $h=\tilde
              h_0 +\tilde h_1 n$.  Thus $\alpha $ is the desired
     isomorphism.

\def \untwist {\Lbpar ^{\hbox {\sixrm \ trivial}}_N}

Before concluding this section, we wish to compare $\EssAlg {\G _N}{\WeylBund }$ (which in the present setting is
$\Cwr $) with its untwisted version, $\EssAlg {\G _N}{\untwist }$, which we denote more simply by $C^*_\ess (\G
_N)$.  Once again, details are left to the reader.

Using the same bisections $U_0$ and $U_1$ as above, we obtain
  $$
  C(U_0)\simeq C(\circle )\simeq C(U_1).
  $$
  A computation gives $C^*_\red (\G _N)=\{g_0\chi _{U_0}+g_1\chi _{U_1}: g_0, g_1\in A\}$.  Let
  $$
  \frak B=A\oplus {\bf C} \dstext {and} \frak A=\{a\oplus a({\bf 1}): a\in A\}\simeq A.
  $$
  For $i=0, 1$ and $g_i\in A$, the map $C^*_\red (\G _N)\ni g_0\chi _{U_0}+g_1\chi _{U_1}\mapsto (g_0+g_1)\oplus
(g_0({\bf 1})-g_1({\bf 1}))$ gives an isomorphism of $C^*_\red (\G _N)$ onto $\frak B$ which carries $C(\G
_N^{(0)})$ onto $\frak A$.  The {gray} ideal for $\incl {\frak B}{\frak A}$ is $\Gamma =0\oplus {\bf C}\subseteq
\frak B$ and the quotient ${\frak B}/\Gamma \simeq A$.  Thus, $\incl {C^*_\ess (\G _n)}{C(\G _N^{(0)})} \simeq \incl
{A}{A}=\incl {C(\circle )}{C(\circle )}$.

{\it This (abelian) example, together with the non-abelian example in the next section, show the absence or presence
of the {gray} ideal can depend upon the line bundle and not just the underlying groupoid, see \ref {diffLB}.}

\endsection

\startsection Example: the {gray} ideal of twisted groupoid C*-algebras

  \def \mat #1#2#3#4{\left [\matrix {#1&#2\cr #3&#4}\right ]}
  \def \tr {\text {tr}}

\sectiontitle

\label GrayIdealTwisted
  One of the most pressing unsolved problems regarding the theory of non Hausdorff grou\-poids is the question of
determining the precise necessary and sufficient conditions the ingredients of \ref {StandingGPDS} should satisfy
for
  $$
  \Incl {\Cgr }{C_0(X)}
  $$
  to be a weak Cartan inclusion.

Topological freeness of $\G $ is certainly a necessary condition but examples such as \cite {ExelNHaus}, show that
it is by no means sufficient.  According to \ref {CharacFreePtsTwo}, topological freeness implies smoothness, so the
missing link is whether or not the {gray} ideal vanishes.

Having spent a significant energy in trying to solve this problem, perhaps the only sensible thing we have to say is
actually something that makes the problem look even more ominous: the vanishing of the {gray} ideal does not depend
only on the structure of the groupoid.  It is influenced by the line bundle as well!

In this section we would therefore like to show an example of a topologically free \'etale groupoid $\G $ such that
the {gray} ideal relative to the inclusion $\Incl {\Cgr }{C_0(X)}$ vanishes for a certain choice of $\Lb $, but not
for all.

The groupoid $\G $ we have in mind is precisely the groupoid of \cite [Section 2]{ExelNHaus}, which we shall now
describe for the convenience of the reader and also because a slightly different point of view is better adapted to
our purposes.

Consider the topological space $X$ consisting of the union of four intervals joined by one end, here denoted
``$0$'', as in the diagram below.

  \begingroup \noindent \hfill 
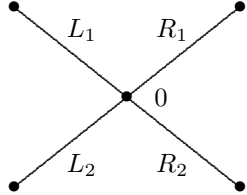
\beginpicture
  \setcoordinatesystem units <0.015truecm, -0.012truecm>
  \setplotarea x from -150 to 150, y from -150 to 150
  \put {\null } at -150 -150
  \put {\null } at 150 150
  \plot 100 100 0 0 -100 100 /
  \plot 100 -100 0 0 -100 -100 /
  \put {$\bullet $} at 0 0
  \put {$\bullet $} at 100 100
  \put {$\bullet $} at -100 100
  \put {$\bullet $} at 100 -100
  \put {$\bullet $} at -100 -100
  \put {$L_1$} at -40 -75
  \put {$R_1$} at 40 -75
  \put {$L_2$} at -40 75
  \put {$R_2$} at 40 75
  \put {$0$} at 30 0
  \endpicture \hfill \null \endgroup

Let us consider an action of the Klein four-group ${\bf Z}_2\times {\bf Z}_2$ on $X$, according to which one of the
two standard generators, say $\sigma =(1,0)$, acts by interchanging the left two arms, and fixing the right ones,
while the second generator $\tau =(0,1)$ interchanges the right two arm and fixes the left ones.

  \begingroup \noindent \hfill 
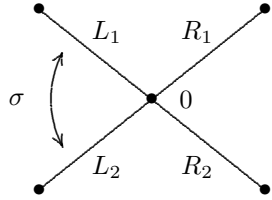
\beginpicture
  \setcoordinatesystem units <0.015truecm, -0.012truecm>
  \setplotarea x from -150 to 150, y from -150 to 150
  \put {\null } at -150 -150
  \put {\null } at 150 150
  \plot 100 100 0 0 -100 100 /
  \plot 100 -100 0 0 -100 -100 /
  \put {$\bullet $} at 0 0
  \put {$\bullet $} at 100 100
  \put {$\bullet $} at -100 100
  \put {$\bullet $} at 100 -100
  \put {$\bullet $} at -100 -100
  \put {$L_1$} at -40 -75
  \put {$R_1$} at 40 -75
  \put {$L_2$} at -40 75
  \put {$R_2$} at 40 75
  \circulararc 55 degrees from -80 -50 center at 0 0
  \arrow <0.15cm> [0.4,1] from -80 -50 to -79.6 -51
  \arrow <0.15cm> [0.4,1] from -80 50 to -79.6 51
  \put {$\sigma $} at -120 0
  \put {$0$} at 30 0
  \endpicture \hfill \null \endgroup

Viewing the Klein group as generating a pseudogroup on $X$, we let $\G $ be the groupoid of germs.  It should
perhaps be stressed that $\G $ is not meant to be the transformation groupoid, but the groupoid built out of the
classical notion of germs.

For $x$ in any one of the right two arms, we therefore have that the germ of $\sigma $ at $x$, denoted by $[\sigma
,x]$, coincides with $[id,x]$, where $id$ denotes the identity map on $X$, simply because there exists a
neighborhood of $x$, e.g.~the arm where it sits, where $\sigma $ coincides with $id$.  In fact there is no group
element fixing such an $x$, other than those having a trivial germ at $x$, so its isotropy group is trivial.  On the
other hand, the isotropy group of the center point $0$ is non trivial as it is easily seen to be isomorphic to the
Klein group itself.

Consequently the set $F$ of free points coincides with $X\setminus \{0\}$, so $\G $ is topologically free.

Choosing the trivial line bundle over $\G $, let us consider the inclusion
  $$
  \Incl {C_\red ^*(\G )}{C_0(X)},
  $$
  and let us show that its {gray} ideal is non trivial.  Given any $g\in {\bf Z}_2\times {\bf Z}_2$, let $U_g$ be
the clopen bisection
  $$
  U_g=\{[g, x]:x\in X\},
  $$
  and denote by $n_g$ the characteristic function of $U_g$.  As discussed in \cite {ExelNHaus}, one has that
  $$
  f:=n_1-n_\sigma -n_\tau +n_{\sigma \tau }
  $$
  vanishes on all points of $\G $, except on those in the isotropy group $\G (0)$.  For every free point $x$ in $X$,
namely for every $x\neq 0$, one therefore has that
  $$
  \Ex (f^*f) = \Ex (4f)\={CharacFreePts.iv} 4f(x)=0,
  $$
  where we have used that $f^*f=4f$, as an easy computation shows.  We then see that $f$ is a nonzero element of the
{gray} ideal $\Gamma $, proving that $\Gamma \neq 0$.

We will next choose a nontrivial line bundle $\Lb $ over the same groupoid $\G $, and we will show that the {gray}
ideal relative to the inclusion
  $$
  \Incl {\Cgr }{C_0(X)},
  $$
  is trivial.  Although not the most direct and geometric way, the easiest description of the line bundle we have in
mind is as the Weyl line bundle $\WeylBund $ for an inclusion whose Weyl groupoid happens to be $\G $.

Consider the C*-algebra
  $$
  B= C\big ([-1, 1],M_2({\bf C})\big )
  $$
  consisting of all continuous functions from the interval $[-1, 1]$ to the algebra of $2\times 2$ complex matrices.
Consider also the subalgebra $A$ of $B$ consisting of all functions $f$ such that
  \lbldeq CondForDiag
  $$
  \def \crr {\cr \pilar {18pt}}
  \matrix {
  f(t)\in D, & \forall t\in [-1, 0], & \hbox {and}\crr
  f(t)\in D',& \forall t\in [0, \hfill 1],\hfill }
  $$
  where $D$ is the set of all diagonal matrices and
  $$
  D' = \left \{ \mat abba : a, b\in {\bf C}\right \}.
  $$

It should be noted that $D'=UDU\inv $, where $U$ is the unitary matrix
  $$
  U={1\over \sqrt 2}\mat 11{-1}1,
  $$
  so what is happening in the left half of the interval $[-1,1]$ is mirrored on the right half by the inner
automorphism determined by $U$.

It is interesting to remark that, since $0$ lies in $[-1,0]$ as well as in $[0,1]$, both conditions in \ref
{CondForDiag} apply, forcing $f(0)$ to be a multiple of the identity.

Since $D$ and $D'$ are commutative algebras, $A$ is a commutative subalgebra of $B$ and one may easily show that its
spectrum is homeomorphic to the space $X$ above.  It is moreover easy to see that the unitary matrices
  $$
  n_\sigma =\mat 0110, \and n_\tau =\mat 100{-1},
  $$
  seen as constant functions on $[-1,1]$, both lie in $\Norm BA $.  Since it is trivial to prove that $A\cup
\{n_\sigma , n_\tau \}$ generates $B$, the inclusion $\incl {B}{A}$ is seen to be regular.

The $*$-subsemigroup $N\subseteq \Norm BA $ generated by $n_\sigma $ and $n_\tau $ is then a group of unitaries, and
it looks very much like the Klein-four group, but not quite.  It is actually a central extension of ${\bf Z}_2\times
{\bf Z}_2$, with a nontrivial cocycle, since
  $$
  n_\sigma n_\tau =- n_\tau n_\sigma .
  $$

As signs do not matter at all in the equivalence relation defining germs in $\G _N$, meaning that
  $$
  [n,x]=[-n,x],
  $$
  whenever $[n,x]$ is a well formed germ, the above cocycle plays no big role in the construction of $\G _N$, and in
fact one may easily see that $\G _N$ is isomorphic to the groupoid $\G $ of our first example above.  The Weyl line
bundle $\WeylBund $, on the other hand, is no longer trivial, in part due to the the above cocycle.

The reader is urged to try to understand the fun structure of $\WeylBund $, but we will instead prove the vanishing
of the {gray} ideal of the inclusion
  $$
  \Incl {C^*_\red (\G , \WeylBund )}{C_0(X)}
  $$
  in an indirect form, using \ref {MainThree.ii}, so we need to verify the hypotheses of this Theorem.

As a first step, we will prove that $A$ masa in $B$.  For this, pick any $f$ in $A'$.  Given $x$ in $[-1, 0)$, let
us prove that $f(x)$ is a diagonal matrix.  In fact, for every $d$ in $D$ one may pick some $a$ in $A$, such that
$a(x)=d$.  Since $f$ commutes with $a$, it follows that $f(x)$ commutes with $a(x)=d$, so $f(x)$ lies in the
commutant of $D$, and hence in $D$, because $D$ is well known to be a masa in $M_2({\bf C})$.  The same argument
applies to show that $f(x)\in D'$, for every $x$ in $(0,1]$, since $D'$ is clearly also a masa in $M_2({\bf C})$.
By continuity we therefore see that $f(0)\in D\cap D'$, so we conclude that $f$ lies in $A$, and hence that $A$ is
indeed a masa in $B$, as claimed.

Evidently it also follows that $A$ is a {semi-masa} relative to the $*$-semigroup $N$, mentioned above.

Let us next prove that $\incl BA$ is a {light} inclusion.  To see this we leave for the reader the task of verifying
that the set $F$ of free points in the spectrum of $A$ coincides with $X\setminus \{0\}$, as before.  Next one
checks that, for each $x$ in $F$, the unique state extension of $\varphi _x$ is given by
  $$
  \psi _x(b) = \clauses {
  \cl b_{11}(x) if x\in L_1,
  \cl b_{22}(x) if x\in L_2,\pilar {20pt}
  \cl \ds \Frac {\ds b_{11}(x) + b_{12}(x) + b_{21}(x) + b_{22}(x)}4 if x\in R_1,\pilar {20pt}
  \cl \ds \Frac {\ds b_{11}(x) - b_{12}(x) - b_{21}(x) + b_{22}(x)}4 if x\in R_2,\pilar {20pt}\stake {10pt}
  }
  $$
   where we view $L_1$ and $L_2$ as two copies of the interval $[-1,0)$, while $R_1$ and $R_2$ are seen as two
copies of $(0, 1]$.

The last two slightly more complicated formulas above have a simple interpretation: since $D'=UDU\inv $, the
expression for $\psi _x(b)$ above, say for $x$ in $R_1$, is obtained by conjugating $b$ by $U$, and then applying
the upper left entry of the resulting matrix of functions to $x$.

Once this description is agreed upon, let us prove that the free expectation $\Pf $ of \ref {DefineFreeCondExp} is
faithful.  For this suppose that $b\in B$ is such that $\Pf (b^*b)=0$.  Then, for every $x$ in $[-1,0)$ we have that
  $$
  \tr \big ((b^*b)(x)\big ) =
  {(b^*b)_{11}(x)+(b^*b)_{22}(x)} =
  {\psi _{x_1}(b^*b)+\psi _{x_2}(b^*b)} = 0,
  $$
  where $x_1$ and $x_2$ are the points of $L_1$ and $L_2$ corresponding to $x$.  Since $\tr $ is a faithful state on
$M_2({\bf C})$, we deduce that $(b^*b)(x)=0$, and hence also that $b(x)=0$.  The same phenomena plays out on the
right-hand-side, although the computations become a bit obscured under conjugation by $U$.  In any case, if $\Pf
(b^*b)=0$, then $b(x)=0$, for all $x$ in $[-1,0)$, as well as for all $x$ in $(0,1]$.  Since $b$ is continuous, we
then deduce that $b(0)=0$, as well, so indeed $b=0$.  This proves that $\Pf $ is faithful, and hence that $\incl
{B}{A}$ is a {light} inclusion.

As we already mentioned, our intention is to apply \ref {MainThree.ii} to our inclusion, so we must provide an
$N$-canonical state relative to each $x$ in $X$.  Besides the case $x=0$, this follows immediately from \ref
{UniqueIsCanonical}, so it is enough to find an $N$-canonical state relative to $0$, but this is easily accomplished
by taking
  $$
  \psi _0(f) = \Frac {\tr \big (f(0)\big )}2,
  $$
  and observing that $\psi _0(n)=0$, for every $n$ in $N$, except for $n=I$ and $n=-I$, these being normalizers in
$N$ for which $0$ lies in $\iFix n$, rather than $\bd n$.

Having verified all of the hypotheses, we may apply \ref {MainThree.ii} to conclude that the inclusions
  $$
  \INCL {C_\red ^*(\G , \WeylBund )}{C_0(X)}, \and \incl BA
  $$
  are isomorphic, and that the {gray} ideal of $C_\red ^*(\G , \WeylBund )$ vanishes.

\state Conclusion \label {diffLB}
  We have exhibited a non-Hausdorff, \'etale groupoid $\G $, equipped with two different Fell line bundles, one
trivial and the other one, called $\WeylBund $, nontrivial.  In the first case, the {gray} ideal of the inclusion
$\Incl {C^*(\G )}{C_0(X)}$ is nonzero, while in the second case, the {gray} ideal of $\Incl {\Cgr }{C_0(X)}$
vanishes.

We remark that Nekrashevych has shown that the algebraic version of the {gray} ideal of a certain groupoid algebra
associated to the Grigorchuk group \cite [Example 4.5]{Nekr} is nontrivial.

\endsection

\startsection Some open questions

  \def \Lbf {\Lbsys \free }
  \def \Sgf {\mathscr {C}\kern 2pt ^*\kern -1pt(\G \free ,\Lbf )}
  \label AltDefs

\sectiontitle

Let us begin this section assuming \ref {StandingGPDS}.  The description of the {gray} ideal of $\Cg $ given in \ref
{DescribeGray.ii}, contrasted with the definition of the essential groupoid C*-algebra as the quotient of $\Cg $ by
its {gray} ideal, gives a somewhat intriguing picture of $\Cgf $ which we have not pursued with sufficient vigor,
but which perhaps deserves some further thought.  So let us begin by briefly outlining a few consequences of our
results, emphasizing this point of view, as well as some unresolved questions.

  Viewing the free component $\G \free $ as a discrete groupoid, it is clearly principal, i.e., each element $\gamma
$ in $\G \free $ is determined by the pair $(r(\gamma ), s(\gamma ))$, so in fact $\G \free $ is isomorphic to the
principal groupoid obtained by restricting the orbit equivalence relation ``$\sim $'' of \ref {OrbEquiv} to the set
of free points $F$, so we are allowed to make the identification
  $$
  \G \free = \big \{(x, y)\in F:x\sim y\big \}.
  $$
  Let us denote by $\Lbf $ the restriction of $\Lb $ to $\G \free $, so that, once its topology is ignored, $\Lbf $
becomes a bona fide Fell line bundle over $\G \free $.

So far so good but, from now on, things start to look a bit strange: let
  $$
  \Sgf
  $$
  be the set of all cross-sections of $\Lbf $ obtained by restricting the elements of $\Cg $ to $\G \free $.  The
correspondence
  $$
  f\in \Cg \to f|_{\G \free }\in \Sgf
  $$
  clearly defines a linear map onto $\Sgf $, whose null space is precisely the {gray} ideal by \ref
{DescribeGray.ii}, so it factors through the quotient, providing a linear bijection
  $$
  \Cgf \to \Sgf .
  $$

  One could use this bijection to export the $*$-algebra structure from $\Cgf $ to $\Sgf $, but in view of \ref
{MatrixMult}, the operations may be described intrinsically by the formulas
  $$
  (f{*}g)(x,y) = \sum _{z\in \Orb x} f(x,z)g(z,y),
  $$
  and
  $$
  f^*(x,y) = f(y,x)^*,
  $$
  which hold for every $f$ and $g$ in $\Sgf $, and for all $(x,y)\in \G \free $.

The fact that $\Cgf $ and $\Sgf $ are isomorphic therefore allows for an equivalent definition of the former, but if
one is to give a free standing alternative definition, one should be able to spell out precisely which sections of
$\Lbf $ belong to $\Sgf $.

One possibility is of course to do as above, namely to say that the sections we are interested in are the ones
obtained by restricting the elements of $\Cg $, but it would be desirable to give a more intrinsic definition using
only the algebraic and topological information that can be read out of $\G \free $ and $\Lbf $.

Given the dire topological properties of $F$ illustrated by the example on section \ref {NonSmoothNrm}, we feel this
might be a difficult task, but let us nevertheless pose the following:

\state Question
  \rm
  Is it possible to construct $\Sgf $, and hence also $\Cgf $, using only algebraic and topological information that can
be read out of $\G \free $ and $\Lbf $?

\bigskip
  A similar question regards the Fourier coefficients of section \ref {FourierSection}.  In the context of \ref
{StandingThree}, consider the Fell bundle $\Bun $ over the discrete groupoid $\Xt $ given by \ref {BxyFellBun}.
Every element $b$ in $B$ gives rise to a cross-section of $\Bun $ by means of the correspondence
  $$
  (y,x)\in \Xt \mapsto \ft yxb\in \Byx .
  $$

  When $b$ is in the {opaque} ideal, the above cross-section is zero by \ref {BlackNoFourrier}, so one should
actually see the $\ft yxb$, above, as depending on the class of $b$ modulo $\Delta $.  The question is then:

\state Question \label FellBundleOverEqulvRel
  \rm
  Under \ref {StandingThree}, can the quotient C*-algebra $B/\Delta $ be described as some sort of cross
sectional C*-algebra of $\Bun \, $?

\bigskip In the context of \ref {StandingThree}, recall from \ref {ManyTopFree.ii} that an inclusion $\incl BA$ is
called smooth when there exists a generating $*$-semigroup consisting of smooth normalizers.  This implies that the
set of all smooth normalizers spans a dense subspace of $B$, so it makes sense to ask:

\state Question \label SmoothSgp
  \rm
  Let $\incl {B}{A}$ be a regular inclusion with $A$ is abelian, and suppose that the set of all smooth normalizers
spans a dense subspace of $B$.  Is $\incl {B}{A}$ a smooth inclusion?

The difficulty here is that the set of all smooth normalizers is not necessarily closed under multiplication, so it
may not form a semigroup.  In order to illustrate this phenomena, consider the action $\theta $ of ${\bf Z}\times
{\bf Z}$ on $\circle $ such that $\theta _{(1,0)}$ is a given irrational rotation while $\theta _{(0,1)}$ is the
identity map.  Regarding the inclusion
  $$
  \Incl {\ C(\circle )\rtimes ({\bf Z}\times {\bf Z})}{C(\circle )}
  $$
  it is clear that each implementing unitary $u_{p,q}$ is a normalizer, and that $\beta _{u_{p,q}}$ coincides with
$\theta _{p,q}$.  Therefore, whenever $p$ is nonzero, the set $\Fix {u_{p,q}}$ of fixed points for $\beta
_{u_{p,q}}$ is empty, whence the set $F_{u_{p,q}}$ of relative free points coincides with $\circle $ by \ref
{MandatoryValues.i}, and consequently $u_{p,q}$ is a smooth normalizer.  In particular both $u_{1, 0}$ and $u_{-1,
1}$ are smooth, but the reader will have no difficulty in checking that their product, namely
  $$
  u_{1, 0} \, u_{-1, 1} = u_{0, 1}
  $$
  is not smooth.  This shows that the product of two smooth normalizers may indeed fail to be smooth, but it does
not seem to provide a counter example for \ref {SmoothSgp} since the smooth normalizers in this case do not appear
to generate $C(\circle )\rtimes ({\bf Z}\times {\bf Z})$ as a right $C(\circle )$ module, as required by \ref
{ManyTopFree.ii} (although
  $
  \{u_{p,q} : p\neq 0\}\cup C(\circle )
  $
  certainly generates $C(\circle )\rtimes ({\bf Z}\times {\bf Z})$ as a C*-algebra).

Smoothness seems to be connected with the ideal intersection property, and the relationship between these properties
is unclear, even in the commutative case.  We therefore pose the following test question.

\state Question
  If the inclusion $\incl {C(Y)}{C(X)}$ has the ideal intersection property, must it be smooth?

\noindent We note that \ref {smoothext} and \ref {Mlambda} give some evidence that this may be the case.

Finally, we have seen that when $\incl {C(Y)}{C(X)}$ has the ideal intersection property and is generated by
projections, it is smooth \ref {proj=smooth}.  Recall from~\cite {PZ} that when an inclusion has a unique
pseudo-expectation which is also faithful, the inclusion has the {\it faithful unique pseudo-expectation property.}  The
inclusion $\incl {C(Y)}{C(X)}$ has the ideal intersection property if and only if it has the faithful unique
pseudo-expectation property \cite [Corollary~3.22]{PZ}.  Because of this and \cite [Proposition~5.5]{SRITwo}, the unique
faithful pseudoexpectation property seems likely to be the correct property instead of the ideal intersection property
for the non-abelian context.  An affirmative answer to the following would be a non-commutative analog of~\ref
{proj=smooth}.

\state Question
  Let $\incl BA$ be a regular inclusion with $A$ abelian and having the faithful unique pseudo-expectation property.
If the linear span of the normalizing partial isometries is dense in $B$, must $\incl BA$ be smooth?

\endsection

\PART {IV -- APPENDIX}

\startsection Isotropy projection

\sectiontitle

Our goal in this section is to prove \ref {IsotropyProjection} below, which is a basic result in the theory of
twisted groupoid C*-algebras.  We have been unable to locate this result in the literature in the generality we
require.

Throughout this section we shall fix an \'etale groupoid $\G $ equipped with a Fell line bundle $\Lb $, as described
in \ref {StandingGPDS}.

\fix We shall also choose a point $x$ in $\Gz $, to be considered fixed for this entire section.

\medskip \noindent The restriction of $\Lb $ to the isotropy group $\G (x)$ will be denoted by $\scrM $, so that
$\scrM $ is a Fell line bundle over $\G (x)$ and we can then consider the twisted group C*-algebra $\CstarGx $.

Our aim is to prove that the map
  \lbldeq TheRestriction
  $$
  \rho : \scrCcL \G \to \CcGx ,
  $$
  sending any $f$ in $\scrCcL \G $ to its restriction to $\G (x)$, extends to a bounded, positive linear map from
$\Cg $ to $\CstarGx $, which we shall call the \emph {isotropy projection}.

In the following we present two general results which will be instrumental in attaining this goal.

\state Lemma
  \izitem
  \lbldzitem LemaBanach
  Let $E$ and $F$ be normed spaces, and let $T:E\to F$ be a linear map such that, for every $\psi $ in the
topological dual $F'$, the composition $\psi \circ T$ is continuous on $E$.  Then $T$ is continuous.
  \lbldzitem BanachToAlgebra
  Let $A$ and $B$ be C*-algebras, and let $E\subseteq A$ be a dense linear subspace. Also let $T:E\to B$ be a linear
map such that, for every state $\psi $ on $B$, the composition $\psi \circ T$ extends to a positive linear
functional on $A$.  Then $T$ extends to a bounded, positive linear map from $A$ to $B$.

\Proof
  Regarding (i), we claim that the transpose map $T':F'\to E'$ is continuous.  Even though $E$ and $F$ might not be
complete, $E'$ and $F'$ are always complete, so we may employ the Closed Graph Theorem.  We therefore let $\{\psi
_n\}_n$ be a sequence in $F'$, converging to zero, and such that the sequence $\{T'(\psi _n)\}_n$ converges to some
$\varphi $ in $E'$.
  Given any $x$ in $E$, we then have
  $$
  \varphi (x)= \lim _{n\to \infty } T'(\psi _n)(x) = \lim _{n\to \infty } \psi _n(T(x)) = 0,
  $$
  so $\varphi =0$.  This proves the claim that $T'$ is continuous.

In order to prove the continuity of $T$, pick any $x$ in $E$, and take $\psi $ in $F'$ with unit norm, such that
  $\big |\psi \big (T(x)\big )\big | = \|T(x\|$.  Then
  $$
  \|T(x)\| = \big |\psi \big (T(x)\big )\big | = |T'(\psi )(x)| \leq \|T'(\psi )\|\|x\| \leq \|T'\|\|x\|.
  $$

  In order to prove (ii), consider the map
  $$
  T^\star : \psi \in B'\ \mapsto \ T\circ \psi \in E^\star ,
  $$
  where $E^\star $ stands for the algebraic dual of $E$.  Given a state $\psi $ of $B$, observe that the hypothesis
implies that $T^\star (\psi )$ actually lies in $E'$.  Since the states of $B$ span $B'$, we see that the range of
$T^\star $ is a subset of $E'$.

This verifies the conditions of (i), so we conclude that $T$ is bounded, and hence it extends to a bounded linear
map $\tilde T:A\to B$.  To see that $\tilde T$ is positive, choose any $a\in A_+$.  Then, given any state $\psi $ on
$B$, let $\varphi $ be a positive linear functional on $A$ extending $\psi \circ T$.  Then
  $$
  \varphi (x)=\psi (\tilde T(x)),
  $$
  for all $x$ in $E$, and this in turn implies that the same holds for every $x$ in $A$.  It then follows that
  $$
  \psi (\tilde T(a)) = \varphi (a) \geq 0.
  $$
  Since $\psi $ is arbitrary, we conclude that $\tilde T(a)\geq 0$, so $\tilde T$ is positive.
  \endProof

Returning to our stated goal of extending \ref {TheRestriction} to the respective completions, and in order to
prepare for an application of \ref {BanachToAlgebra}, let us fix a state $\psi $ on $\CstarGx $.  This said,
consider for each $g$ in $\G (x)$, the linear functional $p_g$ defined on $\scrM _g$ (the fiber of $\scrM $ over
$g$, as in~\ref {StandingGPDS}) by
  $$
  p_g:u\in \scrM _g\mapsto \psi (u{\delta _g})\in {\bf C},
  $$
  where we denote by $u\delta _g$ the element of $\CcGx $ given by
  \lbldeq DeltaNotation
  $$
  (u\delta _g)(h) = \left \{\matrix {
  u, & \text {if } h=g, \hfill \cr
  0, & \text {otherwise}.
  } \right .
  $$
  So, for any $a$ in $\CcGx $, we have that
  \lbldeq FormulaPsi
  $$
  \psi (a) = \psi \Big (\sum _{g\in \G (x)}a(g)\delta _g\Big ) = \sum _{g\in \G (x)}p_g\big (a(g)\big ).
  $$

\state Proposition \label PhiPositive
  Suppose that, for each $g$ in $\G (x)$, we are given $u_g$ in $\scrM _g$, such that $u_g=0$ for all but finitely
many $g$.  Then
  $$
  \sum _{g, h\in \G (x)}p_g(u_h^*u_{hg})\geq 0.
  $$

\Proof Letting $a=\sum _{g\in \G (x)}u_g\delta _g$, we have by the positivity of $\psi $ that
  $$
  0\leq \psi (a^*a) =
  \psi \left ( \Big (\sum _{g\in \G (x)}u_g\delta _g\Big )^*\Big (\sum _{g\in \G (x)}u_g\delta _g\Big )\right )
\quebra =
  \sum _{g, h\in \G (x)}\psi (u_h^*u_g\delta _{h\inv g}) =
  \sum _{g, h\in \G (x)}p_{h\inv g}(u_h^*u_g) =
  \sum _{g, h\in \G (x)}p_g(u_h^*u_{hg}).
  \closeProof
  $$
  \endProof

Let us now define a linear functional $\varphi $ on $\scrCcL \G $ by
  \lbldeq DefinePhi
  $$
  \varphi (b) =
  \sum _{g\in \G (x)} p_g\big (b(g)\big ), \for b\in \scrCcL \G .
  $$
  It is worth noticing that $\varphi $ is defined via the exact same expression as in \ref {FormulaPsi}, although
the variable $a$, used there, is a function defined only on $\G (x)$, while $b$ is defined on the whole of $\G $.
Nevertheless only the restriction of $b$ to $\G (x)$ affects the outcome of \ref {DefinePhi}.

\state Lemma
  Given any state $\psi $ on $\CstarGx $, the linear functional $\varphi $ defined in terms of $\psi $ by \ref
{DefinePhi} is positive in the sense that $\varphi (b^*b)\geq 0$, for all $b$ in $\scrCcL \G $.

\Proof Given $b$, one has that
  $$
  \varphi (b^*b) =
  \sum _{g\in \G (x)} p_g\big ((b^*b)(g)\big ) =
  $$
  \lbldeq ContaPositivo
  $$
  =
  \sum _{g\in \G (x)} p_g\Big (\sum _{\gamma \in \G _x}b^*(\gamma \inv ) b(\gamma g)\Big ) =
  \sum _{(g, \gamma )\in \G (x)\times \G _x} p_g\big (b(\gamma )^*b(\gamma g)\big ).
  $$
  The set of indices for the above sum, namely $\G (x)\times \G _x$, can be split as the disjoint union
  $$
  \G (x)\times \G _x =\bigcup _{y\in \Orb x } \G (x)\times \G _x^y,
  $$
  where $\Orb x $ is the orbit of $x$.
  Moreover, if for each $y$ in $\Orb x $ we choose a representative $\gamma _y\in \G _x^y$, we may parametrize $\G
_x^y$ as
  $$
  \G _x^y = \{\gamma _yh: h\in \G (x)\},
  $$
  and hence
  \lbldeq ParametrizeGs
  $$
  \G (x)\times \G _x =\{(g, \gamma _yh): y\in \Orb x, \ g, h\in \G (x)\}.
  $$
  So, using \ref {ContaPositivo}, we have that
  \lbldeq ContaPositivoTwo
  $$
  \varphi (b^*b) =
  \sum _{y\in \Orb x }\ \ \sum _{(g, h)\in \G (x)\times \G (x)}p_g\big (b(\gamma _yh)^*b(\gamma _yhg)\big ).
  $$
  In order to prove that this is non-negative, it is enough to show that, for each $y$ in $\Orb x $, the
corresponding summand above, namely
  \lbldeq DefineS
  $$
  S_y := \sum _{g, h\in \G (x)}p_g\big (b(\gamma _yh)^*b(\gamma _yhg)\big ),
  $$
  is non-negative.  So, fixing such a $y$, choose $e_y\in \Lb _{\gamma _y}$, with $e_ye_y^*=1$, and observe that
  $$
  S_y =
  \sum _{g, h\in \G (x)} p_g\big (b(\gamma _yh)^*e_ye_y^*b(\gamma _yhg)\big ) =
  \sum _{g, h\in \G (x)} p_g\big (u_h^*u_{hg}\big ),
  $$
  where
  $$
  u_h:=e_y^*b(\gamma _yh)\in \Lb _{\gamma _y\inv }\Lb _{\gamma _yh} = \Lb _h.
  $$
  That $S_y$ is non-negative then follows from \ref {PhiPositive}.
  \endProof

The GNS representation associated to a given positive linear functional can be proven to exist in a variety of
situations, but quite often one has to assume completeness of the algebra involved.  In our case of $\scrCcL \G $,
this condition is not granted so we must prove it by hand.

\state Proposition \label AdHocGNS
  Given any state $\psi $ on $\CstarGx $, there exists a triple $(\pi , \Hilb , \xi )$, where $\Hilb $ is a Hilbert
space, $\xi $ is a unit vector in $\Hilb $, and $\pi $ is a *-representation of\/ $\scrCcL \G $ in $\BH $, such that
  $$
  \langle \pi (b)\xi , \xi \rangle =\varphi (b), \for b\in \scrCcL \G ,
  $$
  where $\varphi $ is the linear functional defined in terms of $\psi $ by \ref {DefinePhi}.

\Proof We begin by claiming that for every $b\in \scrCcL \G $, and every pointwise nonnegative $a\in \CurlyCcGz $,
one has that
  \lbldeq HidenNormalizer
  $$
  \varphi (b^*ab) \leq \|a\| \varphi (b^*b),
  $$
  where $\|a\|$ refers to the supremum norm.
  In order to verify this, we shall retain most of the ingredients introduced in the above proof.  More specifically
we will retain the choice of $\gamma _y\in \G _x^y$, for each $y$ on $\Orb x $, and we will again use the
parametrization of
  $\G (x)\times \G _x$ given by \ref {ParametrizeGs}.  Moreover we shall use the fact that $S_y$, as defined in \ref
{DefineS}, is non-negative, and that
  $$
  \varphi (b^*b) = \sum _{y\in \Orb x } S_y.
  $$
  This said, we have
  $$
  \varphi (b^*ab) =
  \sum _{g\in \G (x)} p_g\big ((b^*ab)(g)\big ) =
  \sum _{g\in \G (x)} p_g\Big (\sum _{\gamma \in \G _x}b^*(\gamma \inv )\ \langle a, r(\gamma )\rangle \ b(\gamma
g)\Big ) \quebra =
  \sum _{(g, \gamma )\in \G (x)\times \G _x} \langle a, r(\gamma )\rangle \ p_g\big (b(\gamma )^*b(\gamma g)\big ) =
  \sum _{y\in \Orb x }\langle a, y\rangle \sum _{g, h\in \G (x)} p_g\big (b(\gamma _yh)^*b(\gamma _yhg)\big )
\quebra \leq
  \|a\|\sum _{y\in \Orb x } S_y =
  \|a\| \varphi (b^*b),
  $$
  hence proving our claim.

  Following the first few steps of the classical construction of the GNS representation, let
  $$
  I=\{x\in B:\varphi (x^*x)=0\},
  $$
  and let $\Hilb _0:=\scrCcL \G /I$ be equipped with the inner product
  $$
  \langle x+I, y+I\rangle = \varphi (y^*x), \for x, y\in \scrCcL \G .
  $$
  For each $b$ in $\scrCcL \G $, consider the linear operator $\pi _0$ on $\Hilb _0$, given by
  $$
  \pi _0(b)(x+I) = bx+I, \for x\in \scrCcL \G .
  $$
  We next claim that $\pi _0(b)$ is a bounded operator.  To see this, suppose first that $b$ lies in $\scrCcL U$,
for some open bisection of $\G $.  Then it is easy to see that $b^*b\in \CurlyCcGz $, so we may use \ref
{HidenNormalizer} to obtain for every $x$ in $\scrCcL \G $, that
  $$
  \|\pi _0(b)(x+I)\|^2 = \varphi (x^*b^*bx) \leq \|b^*b\| \varphi (x^*x) = \|b^*b\|\, \|x+I\|^2,
  $$
  thus proving that $\pi _0(b)$ is bounded.  Since $\scrCcL \G $ is spanned by the set of elements $b$ considered
above, we conclude that $\pi _0(b)$ is bounded for every $b$ in $\scrCcL \G $.

We may then extend each $\pi _0(b)$ to a bounded operator on the completion $\Hilb $ of $\Hilb _0$, thus obtaining a
*-representation of $\scrCcL \G $ on $\Hilb $.

Finally, choose any $a$ in $\CurlyCcGz $, such that $\langle a, x\rangle =1$, and let $\xi =a+I$.  Then
  $$
  \|\xi \|^2 = \langle a+I, a+I\rangle = \varphi (a^*a) \={DefinePhi}
  \sum _{g\in \G (x)} p_g\big ((a^*a)(g)\big ) = p_x\big (\langle a^*a, x\rangle \big ) = p_x(1) = 1,
  $$
  and if $b\in \scrCcL \G $, we have that
  $$
  \langle \pi (b)\xi , \xi \rangle =
  \langle ba+I, a+I\rangle =
  \varphi (a^*ba) \={DefinePhi}
  \sum _{g\in \G (x)} p_g\big ((a^*ba)(g)\big ) \quebra =
  \sum _{g\in \G (x)} p_g\big (\overline {\langle a, x\rangle }b(g) \langle a, x\rangle \big ) =
  \sum _{g\in \G (x)} p_g\big (b(g)\big ) =
  \varphi (b).
  \closeProof
  $$
  \endProof

We are now ready to prove the main result of this appendix.

\state Theorem \label IsotropyProjection
  Given any twisted \'etale groupoid $(\G ,\Lbpar )$, as in \ref {StandingGPDS}, and given any $x$ in $\Gz $, the
restriction map
  $$
  \rho : \scrCcL \G \to \CcGx
  $$
  extends to a bounded, positive linear map from the full groupoid C*-algebra $\Cg $, to the full group C*-algebra
$\CstarGx $.  We call $\rho $ the \newConcept {isotropy projection}{isotropy projection}.

\Proof
  As already indicated we shall prove the statement via an application of \ref {BanachToAlgebra}, where we set
  $$
  A=\Cg , \quad E=\scrCcL \G , \and
  B=\CstarGx ,
  $$
  while $T$ will be taken to be $\rho $, seen as a map into $\CstarGx $.  In order to verify the hypotheses of \ref
{BanachToAlgebra}, let us therefore fix a state $\psi $ on $\CstarGx $, and let us consider the linear functional
$\varphi $ introduced in \ref {DefinePhi}, as well as the representation $\pi $ of \ref {AdHocGNS}.

By the universal property of $\Cg $, we have that $\pi $ extends to a *-representation $\Pi $ of $\Cg $ on $\Hilb $,
so we may define a state $\Phi $ on $\Cg $ by
  $$
  \Phi (y) = \big \langle \Pi (b)\xi , \xi \big \rangle , \for y\in \Cg ,
  $$
  which evidently extends $\varphi $.  Since it is obvious that $\varphi =\psi \circ \rho $, the hypotheses of \ref
{BanachToAlgebra} are verified, so its conclusion becomes available, thus finishing the present proof.
  \endProof

\endsection

\SymbolIndex

\ConceptIndex

\vskip 1.5cm \goodbreak

\references

\AddToTableOfContents {}{References}

\Article BCFS
  J. H. Brown, L. O. Clark, C. Farthing, and A. Sims;
  Simplicity of algebras associated to \'etale groupoids;
  Semigroup Forum, 88 (2014), no. 2, 433-452

\Article BG
  N. P. Brown and E. Guentner;
  New C*-completions of discrete groups and related spaces;
  Bull. Lond. Math. Soc., 45 (2013), 1181-1193.

\Bibitem BEW
  A. Buss, S. Echterhoff and R. Willett;
  Exotic crossed products;
  Operator algebras and applications -- the Abel Symposium 2015, 67--114, Abel Symp., 12, 2017

\Article BEWTwo
  A. Buss, S. Echterhoff and R. Willett;
  Exotic crossed products and the Baum-Connes conjecture;
  J. Reine Angew. Math., 740 (2018), 111-159

\Article BussExel
  A. Buss and R. Exel;
  Fell bundles over inverse semigroups and twisted \'etale groupoids;
  J. Operator Theory, 67 (2012), 153-205

\Article BussMeyer
  A. Buss and R. Meyer;
  Inverse semigroup actions on groupoids;
  Rocky Mountain J. of Math., 47 (2017), 53-159

\Bibitem CEP
  L. O. Clark, R. Exel, E. Pardo, A. Sims and C. Starling;
  Simplicity of algebras associated to non-Hausdorff groupoids;
  Trans.\ Amer.\ Math.\ Soc., 372 (2019), no.\ 5, 3669--3712.

\Article DKR
  V. Deaconu, A. Kumjian and B. Ramazan;
  Fell bundles associated to groupoid morphisms;
  Math. Scand., 102 (2008), 305-319

\Article Dixmier
  J. Dixmier;
  Sur certains espaces consideres par M. H. Stone;
  Summa Bras. Math., 2 (1951), 151-181

\Article DonsigPitts
  A. Donsig and D. Pitts;
  Coordinate systems and bounded isomorphisms;
  J. Operator Theory, 59 (2008), 359-416

\Article DFPvn
  A. Donsig, A. Fuller and D. Pitts;
  Von Neumann algebras and extensions of inverse semigroups;
  Proc. Edinb. Math. Soc. (2), 60 (2016), 57-97

\Bibitem DFP
  A. Donsig, A. Fuller and D. Pitts;
  Cartan Triples;
  International Mathematics Research Notices, (2019), \hfill \break
  https://doi.org/10.1093/imrn/rnz340

\Bibitem EffRuan
  E. G. Effros and Z-J. Ruan;
  Operator spaces;
  London Mathematical Society Monographs, 23.
  Oxford University Press, 2000

\Article cpisg
  R. Exel;
  Noncommutative Cartan sub-algebras of C*-algebras;
  New York J. Math., 17 (2011), 331-382

\Article actions
  R. Exel;
  Inverse semigroups and combinatorial C*-algebras;
  Bull. Braz. Math. Soc. (N.S.), 39 (2008), 191-313

\Article ExelNHaus
  R. Exel;
  Non-Hausdorff \'etale groupoids;
  Proc. Amer. Math. Soc., 139 (2011), 897-907

\Bibitem FD
  J. M. G. Fell and R. S.  Doran;
  Representations of $*$-algebras, locally compact groups, and Banach $*$-algebraic bundles, vols. 1 \& 2;
  Pure and Applied Mathematics, 126, Academic Press, 1988

\Article Gleason A. Gleason; Projective topological spaces; Illinois J. Math., 2 (1958), 482-489

\Article FM  J. Feldman and C. Moore;
Ergodic equivalence relations, cohomology, and von {N}eumann
algebras. {II};
Trans.\ Amer.\ Math.\ Soc., 234 (1977), no. 2, 325--359.

\Article HadPaul
  D. Hadwin and V. Paulsen;
  Injectivity and projectivity in analysis and topology;
  Sci. China Math., 54 (2011), no. 11, 2347-2359

\Bibitem Halmos
  P. Halmos;
  Lectures on {B}oolean algebras;
  Van Nostrand Mathematical Studies, No. 1, 1963

\Article Hahn
  P. Hahn;
  Haar measure for measure groupoids;
  Trans. Amer. Math. Soc., 242 (1978), 1-33

\Bibitem CH
  E. Hewitt and K. A. Ross;
  Abstract harmonic analysis II;
  Academic Press, 1970

\Bibitem Hofmann
  K. H. Hofmann;
  Bundles and sheaves are equivalent in the category of Banach spaces;
  pp. 53-69, Lecture Notes in Mathematics, 575, Springer-Verlag, 1977

\Article Kumjian
  A. Kumjian;
  On C*-diagonals;
  Canad. J. Math., 38 (1986), 969-1008

\Article KumjianFell
  A. Kumjian;
  Fell bundles over groupoids;
  Proc. Amer. Math. Soc., 126 (1998), 1115-1125

\Article KM
  B. K. Kwa\'sniewski and R. Meyer;
Essential crossed products for inverse semigroup actions:
              simplicity and pure infiniteness;
              Doc. Math., 26 (2021), 271--335

\Article Nekr
  V. Nekrashevych;
  Growth of \'etale groupoids and simple algebras;
  Int. J. of Alg. and Comp., 26 (2016), 375-397

\Bibitem Mrcun
  I. Moerdijk and J. Mrcun;
  Introduction to Foliations and Lie Groupoids;
  Cambridge University Press, 2003

\Bibitem Paterson
  A. L. T. Paterson;
  Groupoids, inverse semigroups, and their operator algebras;
  Birkhauser, 1999

\Bibitem Paulsen
  V. Paulsen;
  Completely bounded maps and operator algebras;
  Cambridge Studies in Advanced Mathematics, 78. Cambridge University Press, Cambridge

\Bibitem PittsSR
  D. R. Pitts;
  Structure for regular inclusions;
  2012. arXiv:1202.6413v2

\Article SRI
  D. R. Pitts;
  Structure for regular inclusions. I;
  J. Operator Theory, 78 (2017), 357-416

\Bibitem SRITwo
  D. R. Pitts;
  Structure for regular inclusions. II: Cartan Envelopes, Pseudo-Expectations and Twists;
  J. Funct. Anal., 281 (2021), DOI 10.1016/j.jfa.2021.108993

\Article PZ
  D. R. Pitts and V. Zarikian;
  Unique Pseudo-Expectations for C*-Inclusions;
  Illinois J. Math., 59 (2015), 449-483

\Bibitem Renault
  J. Renault;
  A groupoid approach to C*-algebras;
  Lecture Notes in Mathematics vol.~793, Springer, 1980

\Article RenaultCartan
  J. Renault;
  Cartan subalgebras in C*-algebras;
  Irish Math. Soc. Bulletin, 61 (2008), 29-63

\Article Rieffel
  M. A. Rieffel;
  Induced representations of C*-algebras;
  Adv. Math., 13 (1974), 176-257

\Bibitem Rudin
  W. Rudin;
  Principles of mathematical analysis;
  McGraw-Hill, 1976

\Bibitem Sims
 A. Sims;
  Hausdorff \'etale groupoids and their C*-algebras;
  in A. Sims, G. Szab\'{o}, and D. Williams,
  {\it Operator algebras and dynamics: groupoids, crossed products,
              and {R}okhlin dimension} (F. Perera, ed), pp. 58--120, Advanced
 Courses in Mathematics - CRM Barcelona, Springer, (2020)

\Article Willett
  R. Willett;
  A non-amenable groupoid whose maximal and reduced C*-algebras are the same;
  M\umlaut unster J. of Math., 8 (2015), 241-252

\Bibitem Yg
  S. Yamagami;
  On primitive ideal spaces of C*-algebras over certain locally compact groupoids;
  Mappings of operator algebras (H. Araki and R. Kadison, eds.), Progress in Math. Vol. 84,  Birkh\umlaut au\-ser, Boston (1991),
199-204

\Article Zettl
  H. Zettl;
  A characterization of ternary rings of operators;
  Adv. Math., 48 (1983), 117-143

  \endgroup

\close